\documentclass[english,parskip=half]{scrartcl} %Deutschsprachiges Dokument mit der passenden Absatzkennzeichnung in der KOMA-Klasse article

\usepackage[utf8]{inputenc} %Unicode Eingabekodierung für Umlaute und Sonderzeichen

\usepackage{babel} %Sprachunterstützung mit babel
\usepackage{lmodern} %Schriftart Latin Modern
\usepackage[T1]{fontenc} %Fontencoding T1 = skalierbare Vektorschrift
\usepackage{microtype} %Verbesserte Trennregeln mithilfe von Mikrotypografie
\usepackage{graphicx} %Standard-Paket zum einbinden von Grafiken
\usepackage{csquotes} %Paket für kontextsensitive Zitate
\usepackage{verbatim} %fuer comment
\usepackage[output-decimal-marker={,},exponent-product=\cdot]{siunitx} %Paket für Einheiten mitsamt der deutschen Anpassungen
\usepackage[sorting=nyt, giveninits=true, backref=false, isbn=false, doi=false, url=false]{biblatex} %Lade biblatex für die Erstellung der Bibliografie mit biber
\renewbibmacro{in:}{}
\usepackage{tabularx} %Paket für Tabellen mit variabler Breite um diese gleich der Seitenbreite zu setzen
\usepackage{url}
\usepackage{mathtools,amssymb} %Standard-Mathepakete
\usepackage{amsmath}   % enthaelt nuetzliche Makros fuer Mathematik
\usepackage{amsthm}    % Sätze, Beweise, etc.
\usepackage{bbm}
\usepackage{braket}
\usepackage{enumitem}
\usepackage{relsize}
\usepackage[section]{placeins} %Verhindert, dass Abbildungen außerhalb der aktuellen section gesetzt werden
\usepackage{listings}
\usepackage{gnuplot-lua-tikz}
\usepackage{caption}
\usepackage{subcaption}
\usepackage{tikz}
\usepackage{pgfplots}
\usepackage{pgf}
\usepackage{pgfplotstable, booktabs} % Tabellen aus csv-Datei setzen
\usepackage{hhline}
\usepackage{xcolor,colortbl}
\usepackage{longtable}
\usetikzlibrary{hobby}
\usepackage{enumitem}
\usepackage{soul} % niceres underline \ul
\usepackage{algorithm}
\usepackage{algorithmic}
\usepackage{xfrac}

\usepackage[normalem]{ulem} % für \sout, kann danach gelöscht werden

\usepackage{afterpage}

\usepackage{accents}
\newcommand{\ubar}[1]{\underaccent{\bar}{#1}}

\def\XXint#1#2#3{{\setbox0=\hbox{$#1{#2#3}{\int}$ }
		\vcenter{\hbox{$#2#3$ }}\kern-.6\wd0}}

\DeclarePairedDelimiter{\abs}{\lvert}{\rvert}
\DeclarePairedDelimiter{\norm}{\lVert}{\rVert}
\NewDocumentCommand{\normL}{ s O{} m }{%
	\IfBooleanTF{#1}{\norm*{#3}}{\norm[#2]{#3}}_{L^2(\Omega)}%
}

\DeclareFontFamily{U}{mathc}{}
\DeclareFontShape{U}{mathc}{m}{it}%
{<->s*[1.03] mathc10}{}

\DeclareMathAlphabet{\mathscr}{U}{mathc}{m}{it}
\usepackage{hyperref}
\usepackage{cleveref}

\newtheorem{theorem}{Theorem}[section]
\newtheorem{corollary}[theorem]{Corollary}
\newtheorem{remark}[theorem]{Remark}
\newtheorem{assumptions}[theorem]{Assumptions}

\newtheorem{lemma}[theorem]{Lemma}

\crefname{remark}{Remark}{Remarks}
\crefname{assumptions}{Assumptions}{Assumptions}

\bibliography{paper-bib.bib} %Einbinden der Bibliotheksdatenbank

\AtEveryBibitem{\clearfield{month}}
\AtEveryBibitem{\clearfield{day}}
\AtEveryBibitem{\clearfield{number}}

\addto\extrasenglish{%  
} 

\makeatletter	% Zentrierte Bilder
\g@addto@macro\@floatboxreset\centering
\makeatother

\definecolor{whitesmoke}{rgb}{0.96, 0.96, 0.96}
\definecolor{timberwolf}{rgb}{0.86, 0.84, 0.82}

\def\luise#1{\textbf{\\ \colorbox{green}{\vbox{#1}}}}
\def\johannes#1{\textbf{\color{blue} #1}}
\def\change#1{#1}
\def\changemy#1{#1}
\def\changetwo#1{#1}
\def\kommentar#1{}

\def\anm#1{}

\def\del#1{}

\def\y{\text{\ttfamily y}}

\begin{document}

% FUER LUISE
%%%%%%%%%%%%%%%%%%%%
%\boldmath \LARGE
%%%%%%%%%%%%%%%%%%%%

%Titelseite-----------------------------------------------------------------------------------------
%\begin{titlepage}
%\centering
%
%\Large\bfseries Parabolic optimal control with strongly monotone quasilinearity and its time discretization
%%\Large\bfseries Optimal control of parabolic equations with strongly monotone quasilinearity and their time discretization
%
%%\vspace{2cm}
%%%\includegraphics[width=.5\linewidth]{urlogo}
%%
%%\vspace{2cm}
%
%\normalsize\normalfont
%Luise Blank\\
%Johannes Meisinger\\
%Department of Mathematics\\
%University of Regensburg
%
%\vspace{1.8cm}
%\bfseries
%...\\
%[\baselineskip]
%\today
%\end{titlepage}
%Ende der Titelseite--------------------------------------------------------------------------------
%\afterpage{\blankpage}
%
%\begin{keywords}
%	quasilinear parabolic equation, Allen-Cahn equation, anisotropy, optimal control, implicit discretization
%	%	quasilinear parabolic partial differential equation, phase field approach, double obstacle potential, nonsmooth optimization
%\end{keywords}

% REQUIRED
%\begin{AMS}
%	35K59, 49J20, 49M41, 65M12, 65M60
%\end{AMS}

\begin{center}
{\Large\bfseries Optimal control of anisotropic Allen-Cahn equations}
	
	Luise Blank\footnote{
		Department of Mathematics,
		University of Regensburg, D-93040 Regensburg, Germany (luise.blank@ur.de, johannes.meisinger@ur.de)}, 
	Johannes Meisinger\footnotemark[1]
\end{center}

\newif\ifabstract
\abstracttrue

\newif\ifgraphics
\graphicstrue

\bigskip
% !TeX spellcheck = en_US

\begin{abstract}
\ifabstract
\textbf{Abstract.}
\fi
This paper aims at solving an optimal control problem governed by an anisotropic Allen-Cahn equation numerically.
Therefore we first prove the Fr\'echet differentiability of an in time discretized parabolic control problem
under certain assumptions on the involved quasilinearity {and formulate the first order necessary conditions}.
{As a next step, since the anisotropies are in general not smooth enough, the convergence behavior of the optimal controls are studied for a sequence of (smooth) approximations of the former quasilinear term.
%\change{In addition driving simultaneously the approximation and the time step size to zero is considered.}
\change{In addition the simultaneous limit in the approximation and the time step size is considered.}
  % and also  with respect to the time discretization}.
}%
For a class covering a large variety of anisotropies we introduce a certain regularization and show the previously formulated requirements.
%for differentiability and convergence with respect to regularization.
Finally, a trust region Newton solver is applied
to various anisotropies and configurations, and numerical evidence for mesh independent behavior
and convergence with respect to regularization is presented.
%the required equations for the second derivative is given formally
%\johannes{126 words (ohne meins)}

\kommentar{
In this paper we build on the results of our recent publication (L. Blank and J. Meisinger, "Optimal control of a quasilinear parabolic equation and its time discretization", $\ldots$), where we introduced the optimal control of a general quasilinear parabolic state equation. We now restrict the quasilinearity and the nonlinearity to a specific form, more concretely to an anisotropic Allen-Cahn equation that is used as a model for crystal growth.
Furthermore we introduce a certain regularization of the quasilinear term and show that the assumptions from our previous publication are fulfilled.
This allows us to show Fr{\'e}chet differentiability of the implicit in time discretized state equation and further to derive the first order necessary condition.
At this point we have to restrict the space dimension to $d\leq 3$ since some strong boundedness properties of the state are needed for the proof.
The equations underlying the definition of the second derivative of the control-to-state operator will be given formally.
They are solved in our trust region algorithm that computes a solution to the optimization problem.
A discussion of the thereby obtained numerical findings concludes the paper.
\johannes{190 Wörter}
}
\end{abstract}
\bigskip



\bigskip

%\tableofcontents
%\clearpage
%\newpage
% !TeX spellcheck = en_US

\section{Introduction}%The optimization problem}
\label{sec:intro}
The goal of this paper is to study the optimal control of anisotropic phase field models describing interface evolution. These models are successfully applied e.g. for anisotropic solidification processes like crystal growth (see also \cite{DDE} and references therein).
The defining equations are given by a gradient flow of a Ginzburg-Landau energy. Here several ansatzes exist to incorporate anisotropy.
%Incorporating interfacial energy anisotropy in phase field models already starts at the level of defining an appropriate energy functional.
In the pioneering paper \cite{Kobayashi} the author
considers convex anisotropies in order to obtain a well-posed problem.
% replaces the laplacian by a certain anisotropy function.
% In order for this leading to a well posed problem one has to demand some restrictions that keep the anisotropy convex.
In \cite{Eggleston,Miranville,RV,TLVW,Wise} various approaches are taken to enlarge this also to non-convex anisotropies by adding regularization terms
or changing the structure of the energy functional. 
%To enlarge this also to non-convex anisotropies, the authors of \cite{Wise} (nachschauen, war nicht korrekt angegeben in proposal wo ich hab) or \cite{RV} added regularization terms to the energy functional.
%In \cite{TLVW} and \cite{Miranville}, these solutions are enhanced to overcome a problem of directional dependence of the interface thickness.
%These approaches lead to higher dimensional terms which make the analysis complicated.
%The authors of \cite{Eggleston} try to tackle the problem by convexifying the anisotropies, but also there mathematical difficulties appear.
%Another way---that also we pursue---is
We pursue 
to define the Ginzburg-Landau functional resembling \cite{Kobayashi}, i.e.
\begin{equation}
\label{GinzburgLandauEnergyAniso}
\mathcal{E}(y) \coloneqq \int_{\Omega} {\varepsilon}
A(\nabla y) + \varepsilon^{-1}
\psi(y) \, \mathrm dx,
\end{equation}
but using a different class of anisotropy functions
$A:\mathbb{R}^d \to \mathbb{R}$
\change{where
  $A(p) = \tfrac{1}{2}|\gamma(p)|^2 $
with
a so-called density function $\gamma : \mathbb{R}^d \to \mathbb{R}_{\geq 0}$ (see e.g. \cite{Elliott1996}).}
%is chosen such that it avoids directional dependence of the interface thickness and higher dimensional terms.
The first part of the functional represents the surface energy
while the potential $\psi $ drives the order parameter $y$ to the pure phases given by the local minimizers
 $\pm1$. The interface thickness is proportional to the variable $\varepsilon>0$.	
The scaled $L^2$-gradient flow of \eqref{GinzburgLandauEnergyAniso}
yields the anisotropic Allen-Cahn equation  which defines the state equation for our control problem. This reads as \anm{$u\in L^2(Q)$ entfernt}
\begin{equation}
\label{eq1a}
\min_{}  J(y,u) \coloneqq \dfrac{1}{2} \norm{y(T) - y_\Omega}^2_{L^2(\Omega)} + \dfrac{\lambda}{2\varepsilon} \norm{u}^2_{L^2(Q)}
\end{equation}
subject to the quasilinear parabolic state equation with potentially nonsmooth $A'$
\begin{equation}\label{weakAC}
\begin{aligned}
  \int_Q \varepsilon\partial_t y \eta + \varepsilon A'(\nabla y)^T \nabla \eta +
  \dfrac{1}{\varepsilon} \psi'(y) \eta &=  \int_Q u\eta 
     \quad  \forall \eta \in L^2(0, T; H^1(\Omega))\\
  y(0)&=y_0     \quad \mbox{ in } \Omega
\end{aligned}
\end{equation}
where
$\Omega \subset \mathbb{R}^d$ is a bounded Lipschitz domain,
$Q \coloneqq [0,T]\times \Omega$ denotes the space-time cylinder,  $\Sigma \coloneqq [0,T] \times \partial \Omega$ its boundary and the target function
$y_\Omega \in L^2(\Omega)$
as well as the initial state
$y_0 \in H^1(\Omega)$ are given.
%$u\in L^2(Q)\cong L^2(0,T;L^2(\Omega))$ and $y \in L^2(0,T;H^1(\Omega))\cap H^1(0,T;H^1(\Omega)')$.
% 
{\color{blue} %ALTE FORMULIERUNGEN
\kommentar{The variable $\varepsilon >0$ is related to the interfacial thickness as the interface typically shows a behaviour $\sim \tanh(\tfrac{x}{\sqrt{2}\varepsilon})$ \cite{DDE}.
\\
For further introduction to phase field models we refer to {\cite%[Chapter 7.9]
  {eck2017mathematical}} and references therein.
\\
$J$ is well defined due to the embedding $L^2(0,T;H^1(\Omega))\cap H^1(0,T;H^1(\Omega)') \hookrightarrow C([0,T]; L^2(\Omega))$.
}
}% ALTES
\changetwo{The weak formulation implies the boundary condition $A'(\nabla y)^T\nu=0$ on $\Sigma$ where $\nu$ is the outer normal. Furthermore note}
that the weight of the control cost is divided by $\varepsilon$ since from practical observations one expects contributions only in vicinity of the interface.

In this paper we focus on the control problem discretized in time and aim at differentiability of the reduced cost functional such that efficient control solvers can be applied.
Regarding differentiability considerations
and first order conditions we point also to
\cite{Casas1988, Casas1991, Casas199320} where the authors consider quasilinear elliptic problems
related to ours as well as
to 
\cite{WachsmuthII, WachsmuthIII} for a problem with similar time discretization that is also regularized.
For \crefrange{eq1a}{weakAC} the time discretization is given as in our paper \cite{BlMe21} where existence of optimal controls and the convergence with respect to the time discretization is shown.
For convenience we shortly repeat the dG(0) discretization here.
Let $I_j:=(t_{j-1},t_j]$ and \anm{$n$ auf $j$ geändert}
\begin{equation}
\begin{aligned}
Y_{\tau}&:=\{ y_\tau :Q\rightarrow\mathbb{R} \mid
y_\tau(t,.)\in H^1(\Omega), y_\tau(.,x) \text { constant in } I_j   \text { for } j=1,\ldots,N \} ,\\
U_{\tau}&:=\{ u_\tau :Q\rightarrow\mathbb{R} \mid
u_\tau(t,.)\in L^2(\Omega), u_\tau(.,x) \text { constant in } I_j  \text { for } j=1,\ldots,N \} ,
\end{aligned}
\end{equation}
and 
for each interval we label the constant by a subscript, e.g. $y_j:=y_\tau \vert_{I_j}$.
% The vector containing these constants will be denoted in boldface, e.g. $\boldsymbol{y} \coloneqq (y_n)_{n = 1,\ldots, N} \in H^1(\Omega)^N$. 
The time discretized control problem is then given by 
\begin{equation}
\label{discrete_problem1}
%\min_{H^1(\Omega)^N \times L^2(\Omega)^N} \mathscr{J}(\boldsymbol{y},\boldsymbol{u})
\min_{Y_\tau\times U_\tau} {J}({y_\tau},{u_\tau})
= \dfrac{1}{2}\|y_N-y_\Omega\|^2 + \dfrac{\lambda}{2\varepsilon}\sum_{j=1}^N \tau_j \|u_{j}\|^2
\end{equation}
subject to %the time discretized state equation
\begin{equation}
\begin{aligned}
\label{scheme_state1}
\tfrac\varepsilon{\tau_{j}}(y_{j}, \varphi_{}) + {\varepsilon}  (A'(\nabla y_{j}), \nabla \varphi_{}) +
\tfrac{1}{\varepsilon} (\psi'(y_{j}), \varphi_{}) =
(u_{j}, \varphi_{}) + \tfrac\varepsilon{\tau_{j}}(y_{j-1},\varphi_{}) \quad \forall \varphi &\in H^1(\Omega),\\
j &= 1,\ldots,N
\end{aligned}
\end{equation}
where  %$y_\tau(0,.):=
$y_0 \in H^1(\Omega)$ is the given initial value.

Usually the density function $\gamma $
in $A(p)= \tfrac12 |\gamma(p)|^2$ is assumed 
to be a positive $1$-homogeneous function in
$C^2(\mathbb{R}^d \textbackslash \{0\})\cap C(\mathbb{R}^d)$ (see\change{,} e.g.\change{,} \cite{Barrett2008, DDE, Elliott1996, Graser2013})
providing absolutely 2-homogeneity of $A$.
%Consequently $A$ is in general not twice differentiable at the origin and hence the control-to-state operator is not differentiable.
\change{Consequently $A''$ is absolutely 0-homogeneous and therefore it does not exist at the origin unless $\gamma $ is an energy norm. Hence the control-to-state operator may not be differentiable.}
\change{Since numerical methods for nonsmooth optimal control problems are still in its infancy this is  problematic for efficient solvers.
For a nonsmooth quasilinear elliptic
control problem a semismooth Newton method  is
applied to a relaxed optimality system in \cite{clason2021}.
%Here directional derivatives are required.
To the best of our knowledge globally convergent methods for parabolic equations without extra regularity
requirements do not exist.} 
%\change{This is problematic for efficient numerical methods which involve derivatives.}
%This is problematic for numerical methods
%that require the existence of a gradient of the reduced cost functional like gradient or Newton type methods.
%  Alternatively one might consider active set methods that can be applied to the variational inequalities that arise as first order conditions in the nonsmooth case.}
To circumvent this problem%
\del{ for the numerics}, the present approach is to consider a regularized $A$ by modifying the function $\gamma$.
We give the details later, when we have introduced the specific form of $A$ which was first proposed in \cite{BGN2007, Barrett2008}.
\kommentar{die kamen ungefähr zur gleichen Zeit raus(recieved, published, etc. überschneiden und sind nicht genau rauszukriegen). In keinem wird auf ein vorheriges verwiesen. In Barrett 2008 wird die Wahl motiviert, s. auch meinen Text aus dem entsprechenden Kapitel}%
Furthermore, while for
the numerical experiments we use the smooth double-well potential $\psi(y) = \tfrac{1}{4}(y^2-1)^2$,
the analysis holds for more general $\psi$.
This is also the case for $A$. The functions $A$ and $\psi$
shall fulfill at least the assumptions from \cite{BlMe21} listed below under a) to guarantee the existence of the optimal control and the existence and convergence of the time discretized optimization problem.
The assumptions are further restricted to obtain differentiability of the reduced cost functional. These conditions are met by the regularized $A$ as can be seen in \cref{sec:regularization}.
\begin{assumptions}\label{assu}
\hfill
  \begin{itemize}
\item[a.]
Assume $A\in C^1(\mathbb{R}^d)$
with $A'$ \change{being} strongly monotone and \change{fulfilling} the growth condition $ |A'(p)| \leq C|p| $.\\
Let $\psi \in C^1(\mathbb{R})$
be bounded from below and \change{such that it} can be approximated by 
$f_n $ satisfying
%\begin{equation}
%		\label{f_approximation_property}
$
f_n \in C^2(\mathbb{R}), %\quad
f_n \to \psi %\quad
\text{ in } C^1_{\text{loc}}, %\quad
 -c \leq f_n \leq c (\psi+1) , %\quad 
 f_n'' \geq -C_\psi,  %\quad
 |f_n''| \leq C_n ,
$%		\end{equation}
		with $c,C_\psi, C_n \geq 0$
		and  $\psi(y_0)\in L^1(\Omega)$ for the given initial data $y_0\in H^1(\Omega)$.
\\
 Furthermore  for the time discretization the restriction $\max_j \tau_j < \varepsilon^2/C_\psi$ on the time steps $\tau_j:=t_j-t_{j-1}$ holds.
%		\johannes{Furthermore recall that $\max_j \tau_j < 1/C_\tau$ is required to obtain a unique solution to the discrete problem with Lipschitz-continuous dependence on the data (for details see \cite{BlMe21})}
              \item[b.]
Assume in addition $A\in C^{2}(\mathbb{R}^d)$
with bounded $A''$ and \change{let}
$\psi \in C^{2}(\mathbb{R})$ \change{where the Nemytskii operator given by
$\psi''$ is continuous from $H^1(\Omega)$ to $L^q(\Omega)$ for some \changetwo{$q>\max\{d/2,1\}$}.}

% For the space dimension it holds $d\leq 3$.
\kommentar{
		Assume $A\in C^{2}(\mathbb{R})$
		with $A'$ is strongly monotone ($\Leftrightarrow$ $A''$ is uniformly positive definite)%
		\kommentar{nach https://xingyuzhou.org/blog/notes/strong-convexity Proposition iii) $\Leftrightarrow$ $A$ strongly convex (die Def. dort entspricht uniformly convex in Optimierung 1, s. (3.7) dort) $\Leftrightarrow$ $A''$ glm. pos. def. wenn Gebiet (hier $\mathbb{R})$ offen, mach Opt. 1 Theorem 3.4.c}%
		, and fulfills the growth condition $ |A'(p)| \leq C|p| $. Also let $A''$ be bounded.\\
	Furthermore let $\psi \in C^{2}(\mathbb{R})$
        be bounded from below, $\psi'' \geq -C_\psi$, $\lim_{t\to \pm \infty} \psi''(t)=\infty$
        %$\psi'(0) = 0$ \johannes{?} für $L^\infty$ aussage von Casas
		and  $\psi(y_0)\in L^1(\Omega)$ for the given initial data $y_0\in H^1(\Omega)$.
                Additionaly let
\kommentar{$\psi'': L^6(\Omega)\to L^{q}(\Omega)$ with $q>3/2$}
                $\psi'': L^p(\Omega)\to L^{q}(\Omega)$ with $p=2d/(d-2)$ and $q>d/2$
                be a Nemytskii operator.}
\end{itemize}
\end{assumptions}
\changetwo{Let us mention that one can find $p>2$ with $H^1(\Omega)\hookrightarrow L^p(\Omega)$ and $\tfrac{1}{q} + \tfrac{2}{p}<1$, e.g., when $d\in\{1,2\}$ choose some $p\in(\tfrac{2q}{q-1},\infty)$ and for $d\geq3$ choose $p=\tfrac{2d}{d-2}$. Such a $p$ will be used in the following.}
%\change{ Let us mention that $q$ is chosen to satisfy $\tfrac{1}{q} + \tfrac{2}{p}<1$ where $p$ is defined in this paper as $p=2d/(d-2)$ for $d\neq 2$ and $p=\infty$ for $d=2$
%such that $H^1(\Omega)\hookrightarrow L^p(\Omega)$ holds.}
Note that the assumptions imply that $A''$ is uniformly positive definite and $\psi'' \geq -C_\psi$ holds.
%the assumptions on $\psi''$ from b. guarantee that it can be approximated as stated in a. (see \cite[Remark 2.8]{BlMe21}).
\change{Furthermore, the double-well potential fulfills the condition if $d\leq 3$
  since $\psi''$ induces a continuous Nemytskii operator from $L^{2q}(\Omega)$ to $L^q(\Omega)$
  % by a known result on the continuity of superposition operators
 (see, e.g., \cite[Proposition 26.6]{zeidler2013nonlinear}) and the imbedding $H^1(\Omega) \hookrightarrow L^{2q}(\Omega)$  is only valid for $d\leq 3$.}
  %and a known result on the continuity of superposition operators (see, e.g., \cite[Proposition 26.6]{zeidler2013nonlinear}).}

The outline of this paper is as follows.
In \cref{sec2} we study under above  assumptions the Fr\'echet differentiability of the  reduced cost functional for \crefrange{discrete_problem1}{scheme_state1}. As a first step we analyze the differentiability of the state equation in one time step. Then, due to the implicit discretization one can {successively} prove differentiability of the control-to-state operator and of the reduced cost functional. The corresponding time discrete adjoint equation is deduced \changemy{rigorously}.
Subsequently, in \cref{sec3} we give sufficient conditions on the regularization $A_\delta$ of $A$ such that the corresponding states converge to the solution of the originally given state equation. Furthermore, also the convergence of a subsequence of \change{global minimizers $u^\delta_\tau$ with respect to the regularization parameter $\delta$ and the time discretization coarseness $\tau$ is addressed}.
While these results hold under above assumptions on $A$, in the \change{subsequent} section we study the class of anisotropies given in \cite{BGN2007,bgn13}. We introduce a regularization  for this class by adjusting $\gamma$.
Furthermore we show that $A$ fulfills \cref{assu}.a, and that in addition \labelcref{assu}.b
as well as the conditions in \cref{sec3} hold for the regularizations.
In the final section we first set up formally the linearized equations needed for a trust
region Newton solver
applied to the only in time discretized control problem.
We provide numerical evidence for convergence with respect to the regularization and for iteration numbers independent of the discretization level.
Finally numerical results for various facets of the anisotropy and different configurations are presented.

%End kommentar
%\newpage
% !TeX spellcheck = en_US

\section{Fr\'echet differentiability of the reduced cost functional for the time\changemy{-}discretized problem}\label{sec2}
\anm{$\phi$ zu $\varphi$ geändert\\}%
{In this section we investigate the Fr{\'e}chet differentiability of the cost functional
\change{$j_\tau(u_\tau):= J(y_\tau(u_\tau), u_\tau)
$}
for the time discretized optimal control problem reduced to the control $u_\tau$ when the \cref{assu} hold.
Hence the first order optimality system
can be shown rigorously. Furthermore its derivative is needed for the numerical optimization solver.
}

As a first step the Fréchet differentiability of the discrete control-to-state operator
${S}_\tau: U_\tau \to Y_\tau$ of \cref{scheme_state1}
is shown. Here, the idea is to prove it for a single time step and then to apply the chain rule.
Let us recall, that
the solution operator ${S}_\tau: U_\tau \to Y_\tau$ of  \cref{scheme_state1} is given by mapping
$u_{\tau}$, correspondingly $(u_1,\ldots,u_N)$,  to
$ y_\tau={S}_\tau(u_{\tau})$ determined by
$ (y_1,\ldots,y_N)$  with
\begin{equation}
\label{eq:define_y_j_via_S}
y_j=S(\tfrac{1}{\varepsilon}u_j + \tfrac{1}{\tau_j}y_{j-1}) \quad \forall j=1,\ldots N .
\end{equation}

Here  $S:L^2(\Omega)\to H^1(\Omega), g\mapsto \y$ is defined as the solution operator of the quasilinear elliptic problem
\begin{equation}
\label{modelprob}
(A'(\nabla \y), \nabla \varphi) +  (\zeta(\y), \varphi) = (g, \varphi)	\qquad \forall \varphi \in H^1(\Omega)
\end{equation}
%with $y\in H^1(\Omega)$ and $g \in L^2(\Omega)$.
with % $\zeta$ is given by
\begin{equation}
\label{eq:def_zeta}
\zeta(s) := \tfrac{1}{\varepsilon^2}\psi'(s) + \tfrac{1}{\tau_j}s
% \stackrel{\text{here}}{=}  \tfrac1{\varepsilon^2} s^3 +	(\tfrac1{\tau_j}-\tfrac1{\varepsilon^2} ) s
.
\end{equation}
%since $\psi'(\y)=\y^3-\y$,
%Each $y_j$ in \cref{scheme_state1} is then obtain by $S(g_j)$ where $g_j:=\tfrac1{\varepsilon}u_j+\tfrac1{\tau_j}y_{j-1}$.
Note under \cref{assu}.a
\kommentar{for $\tau_j < 1/C_\psi \stackrel{\text{here}}{=} \varepsilon^2$}
the left-hand side defines a strongly monotone operator. In \cite{BlMe21} we have shown the unique existence of the solution.
Let us mention that with 
%For the analysis we often use
the restriction on the space dimension $d\leq 3$ %.
%In this case
a result from \cite{Casas1995} for quasilinear elliptic equations with controls on the Neumann boundary
provides solutions in $L^\infty(\Omega)$  if in addition $A'(0)=0$. Here the  restriction on $d$ is due to the use of Stampacchia's method. %This leads to the following regularity result.
%From \cite{Casas1995} we even get an $L^\infty$-bound if the right hand side takes on more regularity:
\kommentar{\\
	* $\zeta(\y) \in C^1(\mathbb{R})$\\
	* $\zeta(0) = 0$ (soweit ich herausgefunden habe braucht man das um z.z. dass der Operator Koerziv ($\zeta(x)x = (\zeta(x)-\zeta(0)(x-0)\geq 0$). Ich denke man kann die Vor. wie üblich umgehen, da die Neumann rbd nur von $A'$ abhängt (d.h. man braucht $A'(0)=0$)\\
	* $c \leq \zeta'(\y) \leq h(|\y|)$ 	where $h$ shall be a positive increasing function (das ist immer erfüllt; der Casas hat im Pendant zu $\zeta$ noch eine $x$-Abhängigkeit)\\
	* $A'(0)=0$ !\\
}
% \begin{theorem}
\kommentar{
	\label{th:y_in_L_infty}
	Let \cref{assu} (toDO: stimmen die für das Casas-Theorem) hold and $d \leq 3$.
	Then $\y:=S(g) \in H^1(\Omega)\cap L^\infty(\Omega)$ for $g \in L^2(\Omega)$.
}%end Kommentar

%
%The next auxiliary lemma we will use several times in what follows.
The next auxiliary lemma is obtained by subtracting the defining equations, testing with
\change{$\y-\tilde{\y}$} and using {strong monotonicity} of $A'$ and $\zeta$.
\begin{lemma}
	\label{lemconv}
	Let \cref{assu}.a hold. 
	Then the solution operator $S$ for \cref{modelprob} is Lipschitz continuous \change{with a constant independent of $\tau$ for small enough $\tau $}, i.e. to be more precise
	for \change{$g, \tilde{g}\in L^{2}(\Omega)$} and
	\change{$\y=S(g)$, $\tilde{\y}=S(\tilde{g})$} it holds
	\begin{equation}
	\label{eq:modelprob_Lipschitz}
	\|\y-\tilde{\y}\|_{H^1(\Omega)} \leq C \|g-\tilde{g}\|_{H^{1}(\Omega)'}.
	\end{equation}
\end{lemma}
%\begin{proof}
\kommentar{
	Subtracting the defining equations, testing with
	$\y_1-\y_2$ and using {strong monotonicity}, we obtain
	\begin{align*}
	\begin{aligned}
	%		C\|y_m-y\|_{H^1(\Omega)}^2 &\leq \underbrace{C\|\nabla y_m-\nabla y\|_{L^2(\Omega)}^2}_{\leq (A'(\nabla y_m) - A'(\nabla y), \nabla y_m - \nabla y)} + \underbrace{\|y_m-y\|_{L^2(\Omega)}^2}_{\leq C (\zeta(y_m)-\zeta(y), y_m-y)}\\&\leq C\braket{u_m-u,y_m -y}_{H^{1\ast}, H^1} \leq  C\|u_m-u\|_{H^{-1}(\Omega)} \|y_m-y\|_{H^1(\Omega)},
	\|\y_m-\y\|_{H^1(\Omega)}^2 &= \norm{\nabla \y_m-\nabla \y}_{L^2(\Omega)}^2 + \norm{\y_m-\y}_{L^2(\Omega)}^2\\
	&\leq C(A'(\nabla \y_m) - A'(\nabla \y), \nabla \y_m - \nabla \y) + C (\zeta(\y_m)-\zeta(\y), \y_m-\y)\\
	&= C\braket{g_m-g,\y_m -\y}_{{H^1}', H^1} \leq  C\|g_m-g\|_{H^{1}(\Omega)'} \|\y_m-\y\|_{H^1(\Omega)}
	\end{aligned}
	\end{align*}
	and therefore the assertion after dividing by $\|\y_m-\y\|_{H^1(\Omega)}$.
	%	and therefore the convergence $\y_m \to \y$ in $H^{1}(\Omega)$.
}%
\change{Let us mention that 
  for $y_{\tau}= S_\tau(u_\tau) $ and $\tilde y_{\tau}= S_\tau(\tilde u_\tau) $ it holds also
  (see \cite[Theorem 2.4]{BlMe21})
\begin{equation}
	%	\begin{aligned}
\label{eq:Lipschitz_discrete_inequ2}
	\norm{y_\tau-\tilde y_\tau}_{L^\infty(0,T;L^2(\Omega))} + \norm{\nabla y_\tau-\nabla \tilde y_\tau}_{L^2(0,T;L^2(\Omega))} 
        \leq {C}_{A,\psi,T}
\norm{u_\tau - \tilde u_\tau}_{L^2(0,T;H^1(\Omega)')}  .
	%	\end{aligned}
\end{equation}
}%
\kommentar{
   We note that the constant does not depend on $\tau$ if it is small enough.
  However, on successive application of this theorem there enters a $\tfrac{1}{\tau}$ dependence through the appearance of the previous time step on the right-hand side. This leads to issues in taking the limit $\tau\to0$ for the whole system.
	If one were to consider the latter convergence, one should resort to \cite[Theorem 2.4]{BlMe21}, which shows Lipschitz continuous dependence of the whole system with a constant independent from $\tau$.
The latter convergence was shown in \cite{BlMe21} for a further (minor) restriction on the time step, where also a similar result for the whole discretized equation is given by Theorem 2.4 there.
This result also implies \Cref{lemconv}, we decided however to formulate it in the given form to better fit the needs of this chapter, where the Fr\'echet differentiability will be shown time step wise.\\
}%
Due to difficulties related to a required norm-gap for the differentiability of the $A'$-term, the implicit function theorem is not applicable directly (cf. \cite{Wachsmuth2014}). Since it is more illuminating, we follow the approach in \cite{Casas199320,Casas1988,Casas1995} to show Gâteaux differentiability. We have to add some work afterwards in order to upgrade to Fr{\'e}chet differentiability.

\begin{theorem}
	\label{th_disc_frechet1}
	Let \cref{assu} hold.
	Then the solution operator $S: L^2(\Omega) \to H^1(\Omega)$ of \cref{modelprob}
	is G\^ateaux differentiable and the directional derivative $S'(g)v = z$ is given with $\y = S(g)$ by \change{$z\in H^1(\Omega)$ such that}
	\begin{equation}
	\label{ellGdiff1}
	(A''(\nabla \y) \nabla z, \nabla \varphi) + (\zeta'(\y)z, \varphi) = (v, \varphi)	\qquad \forall \varphi \in H^1(\Omega). 
	\end{equation}
	Furthermore there exists a $C$ independent of
	$g \in L^2(\Omega) $ \change{and $\tau$} with
	\begin{equation}
	\label{zestimate}
	\| z\|_{H^1(\Omega)} \leq C  \| v \|_{L^2(\Omega)} .
	\end{equation}  
\end{theorem}
\anm{$t$ zu $\rho$ geändert im folgenden Beweis}%
\begin{proof}
	\kommentar{Damit geht es auch: Recall that $\y \in H^1(\Omega)\cap L^\infty(\Omega)$,}
	% if $d<4$ (cf. \cref{th:y_in_L_infty}).
	%	\johannes{VORSICHT! Für $d=4$ muss $u$ regulärer sein! Dann aber nicht $S:L^2 \to H^1$ möglich, oder?}
	% alter Kommentar der ursprünglich in dem Remark drüber stand: \johannes{s. auch Trölzsch Seite 283}
	Due to the \cref{assu} the bilinear form defined by the \changemy{left-}hand side of \eqref{ellGdiff1}
	is elliptic
	with an ellipticity constant independent of
	$\changemy{\y}=S(g)$ and is continuous
	given that
	$|(\zeta'(\y)z, \varphi)|  \leq \change{C}\|\zeta'(\y)\|_{L^q(\Omega)}
	\|z\|_{L^p(\Omega)} \|\varphi\|_{L^p(\Omega)}
	\leq C \|\zeta'(\y)\|_{L^q(\Omega)}
	\|z\|_{H^1(\Omega)} \|\varphi\|_{H^1(\Omega)}$.
	%(due to \cref{proof_regA}) and $\zeta'(\y)\geq0$.
	Hence the Lax-Milgram theorem  provides existence and uniqueness of the solution of \cref{ellGdiff1} for $v\in L^2(\Omega)$ and\change{---using $\zeta'(s)\geq - \tfrac 1{\varepsilon^2} C_\psi +\tfrac 1\tau \geq c>0$ for small enough $\tau$---}the estimate \eqref{zestimate} holds for the solutions independently of $g$ \change{ and $\tau$.}\\
	% , as $\y\in H^1(\Omega)\cap L^\infty(\Omega)$ (cf. \cref{th:y_in_L_infty}) is considered fixed, 
	For $v \in L^2(\Omega)$ \change{and $\rho>0$} let us
	consider
	\begin{equation}
	(A'(\nabla \y_\rho), \nabla \varphi) + (\zeta(\y_\rho), \varphi) = (g +\rho v, \varphi) \;.
	\end{equation}	
	Subtracting the equation with $\rho=0$ and dividing by $\rho$, we obtain
	\begin{equation}
	\label{ellGdiff2}
	\left(\dfrac{A'(\nabla \y_\rho) - A'(\nabla \y)}{\rho}, \nabla \varphi \right) + \left(\dfrac{\zeta(\y_\rho) - \zeta(\y)}{\rho}, \varphi\right) = (v, \varphi) \qquad \forall \varphi \in H^1(\Omega).
	\end{equation}
	
	\Cref{lemconv} yields for $z_\rho := \tfrac{1}{\rho}(\y_\rho-\y)$% \in H^1(\Omega)$
	\kommentar{
		\begin{equation}
		C\|z_t\|_{H^1(\Omega)}^2 \leq C\|\nabla z_t\|_{L^2(\Omega)}^2 + \|z_t\|_{L^2(\Omega)}^2\leq (v,z_t)_{L^2(\Omega)} \leq  \|v\|_{L^2(\Omega)} \|z_t\|_{L^2(\Omega)} \leq  \|v\|_{L^2(\Omega)} \|z_t\|_{H^1(\Omega)},
		\end{equation}
	}
	\begin{equation}\label{ztbound}
	\|z_\rho\|_{H^1(\Omega)} \leq C \|v\|_{L^2(\Omega)} \qquad \forall \rho.
	\end{equation}
	% meaning that $z_t$ is bounded since $v$ is considered fixed.
	Therefore there exists a subsequence with $z_{\rho_n} \rightharpoonup z$ in $H^1(\Omega)$. % for a $z \in H^1(\Omega)$.
%	\change{ANM: soll hier noch ein Kommentar bzgl. der Ambiguität rein?}
	We now show that this solves \cref{ellGdiff1} by taking the limit in \cref{ellGdiff2}. \change{This implies that $z$ is in fact the desired G\^ateaux derivative.}
	%	\begin{equation}
	%		\label{ellGdiff3}
	%		\int_{\Omega} \tfrac{1}{t}\left(A'(\nabla y_{t_n}) - A'(\nabla y)\right)\nabla \varphi \, dx +\int_{\Omega}\tfrac{1}{t}(\zeta(y_{t_n})-\zeta(y),\varphi)_{L^2(\Omega)} = (v,\varphi)_{L^2(\Omega)} \quad \forall \varphi \in H^1(\Omega)% C_0^\infty(\Omega).
	%	\end{equation}
	For the first term we have
	\begin{equation}
	\label{limitpass}
	\int_{\Omega} \tfrac{1}{\text{\change{$\rho_n$}}}\left(A'(\nabla \y_{\rho_n}) - A'(\nabla \y)\right)\nabla \varphi \, dx = \int_{\Omega} \nabla z_{\rho_n}^T A''(w_{\rho_n})\nabla \varphi\, dx \stackrel{n \to \infty}{\longrightarrow} \int_{\Omega} \nabla z^T A''(\nabla \y) \nabla \varphi \, dx,
	\end{equation}
	where
	%in the first relation we use the mean value theorem
	$w_{\rho_n}(x) = \nabla \y(x) + s(x)(\nabla \y_{\rho_n}(x)-\nabla \y(x))$) with $s(x)\in[0,1]$ is some intermediate point. Since $y_\rho\to y$ in $H^1(\Omega)$ as $\rho\to0$ (see \cref{lemconv}) it holds $w_{\rho_n}\to \nabla \y$ in $L^2(\Omega)$. The convergence follows since $\nabla z_{\rho_n}$ converges weakly, and 
	\kommentar{and hence $|A''(w_{t_n})\nabla \varphi |\leq C |\nabla \varphi|$}%
	$A''(\cdot) \nabla \varphi$ is a \change{continuous} Nemytskii operator from $(L^2(\Omega))^d$ to $(L^2(\Omega))^d$
	given $A'':\mathbb{R}^d\to \mathbb{R}^{d\times d}$ is continuous and bounded.  Therefore it holds $A''(w_{\rho_n})\nabla \varphi \to A''(\nabla \y) \nabla \varphi$ in $(L^2(\Omega))^d$.
	We proceed analogously with the second term
	\kommentar{
		\begin{equation}
		\label{limitpass2}
		\int_{\Omega} \tfrac{1}{t}\left(\zeta(\y_{t_n}) - \zeta(\y)\right) \varphi \, dx = \int_{\Omega}  z_{t_n} \zeta'(s_{t_n}) \varphi\, dx \stackrel{n \to \infty}{\longrightarrow} \int_{\Omega} z \zeta'(\y) \varphi \, dx,
		\end{equation}	
		where now $s_{t_n}$ is a series of intermediate points from the mean value theorem converging to $\y$ in $H^1(\Omega)$.}%
              with intermediate values  $s_{\rho_n}$ between $\y_{\rho_n}$ and $\y$ using
\del{the imbedding $H^1(\Omega) \hookrightarrow L^p(\Omega)$ \kommentar{nicht kompakt}  and }that
  \change{$\zeta': H^1(\Omega)\to L^{q}(\Omega)$} is a \change{continuous} operator.
	\kommentar{for $d\leq 3$ we got that $z_{t_n} \zeta'(s_{t_n}) \varphi \to z \zeta'(\y) \varphi$ in $L^1(\Omega)$.%
	}
	\kommentar{$\zeta'(s_{t_n})\varphi\to \zeta'(\y)\varphi$ in $L^{6/5}$ (da $2/3 + 1/6 = 5/6$) und $z_{t_n} \rightharpoonup z$ in $L^6(\Omega)$. Achtung in $d=3$ ist die Einbettung $H^1(\Omega) \hookrightarrow L^6(\Omega)$ nicht kompakt.}%
	
	%	So we obtain the weak formulation of \cref{ellGdiff1}.
	%
	Hence $z$ fulfills \cref{ellGdiff1}.
	Since %the subsequence above was arbitrary and
	the limit is given uniquely by the latter equation, the whole sequence $(z_\rho)_{\rho\geq 0}$ converges weakly to $z$ in $H^1(\Omega)$.
	
	It remains to show the strong convergence in $H^1(\Omega)$. Due to the compact imbedding into $L^2(\Omega)$ only the part $\nabla z_\rho \to \nabla z$ in $L^2(\Omega)^\changetwo{d}$ is left. For this we consider the sequence $\{L^T_\rho \nabla z_\rho\}_{\rho>0}$ where $L_\rho$ is the Cholesky-decomposition of the s.p.d.-matrix $A''(w_\rho)$.
	%	$L_t$ is uniformly positive definite \johannes{nicht spd, evtl andere formulierung, z.b. es ex. konst. s.d. eigenwerte groesser als konst.} and bounded (w.r.t. $t$ and $x$), as $A''(w_t)$ is.
	%	From Lemma \ref{proof_regA}.c we obtain \johannes{hier allgemeines A''}
	From boundedness and uniformly positive definiteness of $A''$, we obtain $c \leq \norm{L_\rho(x)} \leq C$ with constants independent from $\rho$ and $x$. Since $A''(w_\rho) \to A''(\nabla \y)$ in $L^2(\Omega)^{d\times d}$, %(keyword Nemytskii)
	from the resulting almost everywhere convergence and just stated boundedness one can verify by dominated convergence that
	% Cholesky decomposition is a continuous function from A'', since 
	% see also https://en.wikipedia.org/wiki/Cholesky_decomposition#The_Cholesky_algorithm
	% all appearing operations are continuous functions from entries in A'' (a_i,i>0 from positive definiteness implying also no pivoting)
	% so we get that also L_t -> L pointwisely
	$$L_\rho \to L \text{ and } L^{-1}_\rho \to L^{-1} \quad \text{in } L^2(\Omega)^{d\times d},$$
	where $L$ is the Cholesky-decomposition of $A''(\nabla \y)$. Furthermore we have
	\begin{align*}
	\|L^T_\rho \nabla z_\rho\|^2_{L^2(\Omega)} &= \int_{\Omega} \nabla z_\rho^T A''(w_\rho) \nabla z_\rho \, dx \leq\int_{\Omega} \nabla z_\rho^T A''(w_\rho) \nabla z_\rho \, dx + (\underbrace{\zeta'(s_{\rho})}_{>0}z_\rho, z_\rho)
	\\ &= (v,z_\rho)_{L^2(\Omega)} \leq C\|v\|^2_{L^2(\Omega)}
	\end{align*}
	using \eqref{ellGdiff2} in the intermediate value formulation and   \eqref{ztbound}.
	So we can extract from the sequence $\{L^T_\rho \nabla z_\rho\}_{\rho>0}$ %is bounded in $L^2(\Omega)^d$, and we can extract
	a weakly convergent subsequence whose limit  is $L^T \nabla z$,
	due to the strong convergence of $L_\rho$.
	\kommentar{(using that the limits in $L^1(\Omega)$ and $L^2(\Omega)$ have to coincide).}%
	% here we used f_n -> f strong g_n -> g weak and f_n, g_n bounded => f_n g_n -> fg weak
	% see https://math.stackexchange.com/questions/1629822/convergence-of-product-of-weakly-converging-sequences-in-lp?noredirect=1&lq=1
	% or my notes
	% the short statement done in the text actually means
	% $L_t \to L$ in $L^2(\Omega)$ + $\nabla z_t \rightharpoonup \nabla z$ in $L^2(\Omega)$ $\implies$ \int_\Omega L_t \nabla z_t \to \int_\Omega L \nabla z
	Due to the uniqueness of the limit also the whole sequence converges weakly  in $L^2(\Omega)^d$.
%
%\change{Given $\zeta': H^1(\Omega)\to L^{q}(\Omega)$}  is a \change{continuous} operator.
\change{Furthermore, there exists a $p'<p$  with
  $\tfrac{1}{q} + \tfrac{2}{p'} \leq 1$. 
Hence the compact imbedding $H^1(\Omega) \hookrightarrow L^{p'}(\Omega)$
provides $z_\rho \to z$ in $L^{p'}(\Omega)$. Then
using $s_\rho\to \y $ in $H^1(\Omega)$ and given $\zeta': H^1(\Omega)\to L^{q}(\Omega)$ is a continuous operator}
%\change{but still $\tfrac{1}{q} + \tfrac{2}{p} \leq 1$ (recall $q>\tfrac{d}{2})$,}
we have
\kommentar{ $s_t \to \y$ not in $L^\infty$ daher nicht Nemytskii mit $L^\infty$ nutzbar!}
	% hier könnte man auch wieder starke L^2 und weak* L^\infty Konv. ausnutzen
	\begin{eqnarray*}
		\|L^T \nabla z\|^2_{L^2(\Omega)} &\leq& \liminf_{\rho\to0} \|L^T_\rho \nabla z_\rho\|^2_{L^2(\Omega)} = \lim_{\rho\to0} \left[ (v,z_\rho) - (\zeta'(s_\rho)z_\rho, z_\rho)\right] \\
		&=& (v, z) - (\zeta'(\y)z, z) = \|L^T \nabla z\|^2_{L^2(\Omega)}
	\end{eqnarray*}
	and with that we can even deduce $L_\rho^T \nabla z_\rho \to L^T \nabla z$ in $L^2(\Omega)^d$. Furthermore there exists some dominating function $m\in L^2(\Omega)$ with $ |L_\rho^\changetwo{T} \nabla z_\rho| \leq m $. Finally, from the pointwise relations
	$$\nabla z_\rho = (L_\rho^\changetwo{{-T}})(L_\rho^\changetwo{{T}}\nabla z_\rho)
	\to (L^\changetwo{{-T}})(L^\changetwo{T}\nabla z)= \nabla z,
	\qquad | \nabla z_\rho | = |L_\rho^\changetwo{{-T}}L_\rho^\changetwo{T}\nabla z_\rho| \leq C |L_\rho^\changetwo{T} \nabla z_\rho| \leq Cm$$
	% left: we have pointwise convergence of L_t and L_t\nabla z_t due to the shown strong L^2 convergences
	% right: use that L_t is bounded independent from x,t
	we get by dominated convergence that $\nabla z_\rho \to \nabla z$ in $L^2(\Omega)$.
	
	%	 Adding the strong convergence of $L_t^{-1}$, from this we deduce the $L^1$-convergence (and therefore a.e.-convergence) of $\{\nabla z_t\}_{t>0}$ and together with the boundedness of the sequence in $L^2(\Omega)$ we finally have $\nabla z_t \to \nabla z$ in $L^2(\Omega)$. \johannes{toDo: ich glaub hier stimmt etwas nicht wie ich es formuliert habe}
	
	\changemy{Ultimately},
	\kommentar{from
		%	\begin{equation}
		$	\|z\|_{H^1(\Omega)} = \lim_t \|z_t\|_{H^1(\Omega)} \leq C \|v\|_{L^2(\Omega)}
		$%	\end{equation}
		we obtain }%
	the desired continuity of $S'(g): L^2(\Omega)\to H^1(\Omega)$ needed for the G\^ateaux differentiability
	is given by \eqref{zestimate}.
\end{proof}

The following theorem upgrades the last result to Fr\'echet differentiability.
\begin{theorem}
	Let \cref{assu} hold. Then the mapping	$S: L^2(\Omega)\to H^1(\Omega)$ is Fr\'echet differentiable.
\end{theorem}
\begin{proof}
	For $g_n \to g$ in $L^2(\Omega)$  \cref{lemconv} provides $\y_n \coloneqq S(g_n)\to \y\coloneqq S(g) $ in $H^1(\Omega)$ and
	given $v \in L^2(\Omega)$  we set $z_n := S'(g_n)v$, $z := S'(g)v$.
	Subtracting the defining equations for $z_n$ and $z$, testing with $(z_n - z) \in H^1(\Omega)$ and inserting 0 terms yields   
	\begin{align*}
	0=  &\int_{\Omega} \left(A''(\nabla \y_n)\nabla z_n - A''(\nabla \y_n) \nabla z\right)^T(\nabla z_n - \nabla z) \, dx
	\nonumber \\ &+ \int_{\Omega} \nabla z^T\left(A''(\nabla \y_n) - A''(\nabla \y)\right)(\nabla z_n - \nabla z) \, dx
	\\ &+ (\zeta'(\y_n)z_n-\zeta'(\y_n)z,z_n-z) + (\zeta'(\y_n)z - \zeta'(\y)z, z_n- z)
	\nonumber\\
	\geq & C\|\nabla z_n - \nabla z\|^2_{L^2(\Omega)}
	-  \|A''(\nabla \y_n) - A''(\nabla \y) \|
	\|\nabla z \|_{L^2(\Omega)} \|\nabla z_n - \nabla z\|_{L^2(\Omega)}
	\nonumber \\ &
	+ C \|z_n-z\|^2_{L^2(\Omega)}                                          - \change{C} \|\zeta'(\y_n) - \zeta'(\y)  \|_{L^q(\Omega)}
	\|z\|_{L^p(\Omega)} \|z_n - z\|_{L^p(\Omega)}
	\end{align*}
	using the ellipticity of $A''(s)$ and $\zeta'(s)$ with constants independent of $s$.
	Given the estimate \eqref{zestimate} $ \|z\|_{H^1(\Omega)}, \|z_n - z\|_{H^1(\Omega)} \leq c \| v\|_{L^2(\Omega)}$  and  it follows
	\begin{align*}
	\|A''(\nabla \y_n) - A''(\nabla \y)\|_{}
	+ \|\zeta'(\y_n) - \zeta'(\y)  \|_{\text{\change{$L^q(\Omega)$}}}
	\geq C
	\frac{\|z_n - z\|_{H^1(\Omega)}^2}{\|v\|_{L^2(\Omega)} ^2} \quad \forall v\in L^2(\Omega).
	\end{align*}
	Since
	$A''(.)$ is a \change{continuous} operator 
	from $(L^2(\Omega))^d \to (L^2(\Omega))^{d\times d}$ and
	$\zeta'(.)$ is a \change{continuous} operator 
	from $L^p(\Omega) \to L^q(\Omega)$,
%	their resulting continuity together with
	$y_n\to y$ in $H^1(\Omega)$ provides that the
	left\changemy{-}hand side \change{goes} to 0 as $n\to 0$.
	Hence $S$ is Fr\'echet differentiable.
\end{proof}

We now consider the solution operator ${S}_\tau: U_\tau \to Y_\tau$
% , (u_1,\ldots u_N) \mapsto (y_1, \ldots, y_N)$
of \cref{scheme_state1}.
\kommentar{
	By using the identifications $U_\tau \cong L^2(\Omega)^N$ and $Y_\tau \cong H^1(\Omega)^N$, we note that using ${S}: L^2(\Omega) \to H^1(\Omega)$ one can write% for the model problem \cref{modelprob} ${S}_\tau$ can be written as
	\begin{equation}
	\left(\begin{matrix}
	y_1\\
	y_2\\
	\vdots\\
	y_N
	\end{matrix}\right)
	=S_\tau
	\mathcal{\left(\begin{matrix}
		u_1\\
		u_2\\
		\vdots\\
		u_N
		\end{matrix}\right)}
	=
	\left(\begin{matrix}
	S(u_1 + \tfrac{1}{\tau_1}y_0)\\
	S(u_2 + \tfrac{1}{\tau_2}y_1)\\
	\vdots\\
	S(u_N + \tfrac{1}{\tau_N}y_{N-1})
	\end{matrix}\right)
	=
	\left(\begin{matrix}
	S(u_1 + \tfrac{1}{\tau_1}y_0)\\
	S(u_2 + \tfrac{1}{\tau_2}S(u_1 + \tfrac{1}{\tau_1}y_0))\\
	\vdots\\
	%	S(S(\ldots (S(u_1 + \tfrac{1}{\tau_1}y_0) + \tfrac{1}{\tau_2}y_1)+\ldots) + y_{J-1}) (falsch, steht nur da falls ich noch was davon brauche (wahrscheinlich nicht))
	S(u_N + \tfrac{1}{\tau_N}S(u_{N-1} + \tfrac{1}{\tau_{N-2}}S(\ldots)\ldots))
	\end{matrix}\right).
	\end{equation}
}% end Kommentar
On each time interval $I_j$ we have
$y_\tau(u_\tau)| _{I_j}
=y_j(u_1,\ldots u_j)=S(\tfrac1{\varepsilon} u_j+\tfrac1{\tau_j} y_{j-1}(u_1,\ldots u_{j-1}))$.
Using the previous shown result for $S$ as well as the chain rule we obtain for $j=1,\ldots , N$
\begin{align}
\label{eq:zj}
z_j&:= {\frac{dy_\tau(u_\tau)| _{I_j}}{du_\tau}}v_\tau
=S'(\tfrac1{\varepsilon}u_j+\tfrac1{\tau_j} y_{j-1})(\tfrac1{\varepsilon} v_j
+ \tfrac1{\tau_j} \frac{dy_\tau(u_\tau)|_{I_{j-1}}}{du_\tau}v_\tau)
\nonumber   \\   &= S'(\tfrac1{\varepsilon}u_j+\tfrac1{\tau_j} y_{j-1})(\tfrac1{\varepsilon}v_j
+ \tfrac1{\tau_j} z_{j-1})
\end{align}
where we used $z_0:=0$
and by induction we  can state the following theorem.
\anm{$j \to j_\tau$ for later clarity}
\begin{theorem}
	\label{th:Fdiff_total_S}
	Let \cref{assu} hold. Then the operator ${S}_\tau: U_\tau \to Y_\tau$ is Fr\'echet differentiable and
	consequently also the reduced cost functional
	\change{$j_\tau: U_\tau \to \mathbb{R}$}
	with
        \change{$j_\tau(u_\tau):=J({S}_\tau(u_\tau),u_\tau)$} is  Fr\'echet differentiable with
	\change{$j_\tau'(u_\tau)v_\tau= (y_N-y_\Omega, z_N) + \tfrac\lambda\varepsilon (u_\tau ,v_\tau)$},
	where $z_N$ is given by the solution of the following sequence with $z_0:=0$
\change{
\begin{equation}
\label{diff}
(\varepsilon A''(\nabla y_j) \nabla z_j, \nabla \varphi)
+ (\tfrac1\varepsilon\psi''(y_j)z_j +
\tfrac\varepsilon{\tau_j} z_j, \varphi) = (v_j + \tfrac\varepsilon{\tau_j} z_{j-1} , \varphi)	\qquad \forall j=1,\ldots N, \; \; \varphi \in H^1(\Omega). 
      \end{equation}
}
\end{theorem}
\kommentar{
	We note that the derivative can be written out as (FEHLER: DAS SOLLTE EINE MATRIX SEIN, SO IST NUR NACH $u_1$ ABGELEITET)
	\begin{equation}
	{S}_\tau'\left(\begin{matrix}
	u_1\\
	u_2\\
	\vdots\\
	u_N
	\end{matrix}\right)
	=
	\left(\begin{matrix}
	S'(u_1 + \tfrac{1}{\tau_1}y_0)\\
	S'(u_2 + \tfrac{1}{\tau_2}y_1)\\
	\vdots\\
	S'(u_N + \tfrac{1}{\tau_N}y_{N-1})
	\end{matrix}\right)
	=
	\left(\begin{matrix}
	S'(u_1 + \tfrac{1}{\tau_1}y_0)\\
	S'(u_2 + \tfrac{1}{\tau_2}S(u_1 + \tfrac{1}{\tau_1}y_0))\tfrac{1}{\tau_1}S'(u_1 + \tfrac{1}{\tau_1}y_0)\\
	\vdots\\
	S'(u_N + \tfrac{1}{\tau_N}S(\ldots))\tfrac{1}{\tau_N}S'(u_{N-1} + \tfrac{1}{\tau_{N-1}}S(\ldots))\tfrac{1}{\tau_{N-1}}S'(\ldots)\ldots\tfrac{1}{\tau_1}S'(u_1 + y_0)
	\end{matrix}\right).
	\end{equation}
}

We note that the $z_j$ %defined in \cref{eq:zj}
satisfy a linearized state equation given
by $z_\tau=S_\tau'(y_\tau)v_\tau\in Y_\tau$.
% Using \cref{diff} this may be viewed as
Furthermore, \cref{diff} is the dG(0) discretization (as used for the state equation) of
% the time discrete variant
the in time continuous, linearized state %following parabolic
equation
\begin{equation}
\label{eq:lin_equation_cont_with_z}
\begin{aligned}
\varepsilon(\partial_t z, \eta) + \varepsilon (A''(\nabla y) \nabla z, \nabla \eta) + \dfrac{1}{\varepsilon}(\psi''(y)z, \eta) &= (v, \eta) \quad &\forall \eta \in L^2(0, T;H^1(\Omega))\\
z(0) &= 0 \qquad &\text{in } \Omega.
\end{aligned}
\end{equation}
Given \cref{assu} the unique solvability of \cref{eq:lin_equation_cont_with_z} is guaranteed by standard results on parabolic equations.

Let us mention that, while for the forward problem it may be more efficient to use the semi-implicit scheme of \cite{bgn13} where  $A'(\nabla y_i)$ is approximated by $M(\nabla y_{j-1})\nabla y_j$, to show Fr\'echet differentiability with the above technique
of applying the solution operator $S$ recursively,
higher regularity properties are required. In particular, to our best knowledge, the gradient of the previous time step solution has to be bounded in $L^\infty(\Omega)$ {as it appears in the ellipticity coefficient} \cite{Wolff2018}.

Let us now define for given $y$ the adjoint equation  in the time continuous setting:
\begin{equation}
\label{eq:cont_adjoint}
\begin{aligned}
- \varepsilon (\eta,\partial_t p) + \varepsilon (A''(\nabla y)\nabla \eta, \nabla p) + \dfrac{1}{\varepsilon}(\psi''(y)\eta,p) &= 0 \quad &\forall \eta \in L^2(0, T;H^1(\Omega))\\
p(T) &= y(T)-y_\Omega \quad &\text{in } \Omega
.
\end{aligned}
\end{equation}
After the substitution $t \mapsto -t$ as for the linearized equation \cref{eq:lin_equation_cont_with_z} the existence of a unique adjoint $p$ as a solution of \cref{eq:cont_adjoint} follows. In analogy to the discretization of the state equation, but taking into account the backward-in-time nature,
we use the piecewise constant time discrete
$p_\tau \in P_{\tau}$,
where
$$
P_{\tau}:=\{ p_\tau :Q\rightarrow\mathbb{R} \mid
p_\tau(t,.)\in H^1(\Omega), p_\tau(.,x) \text { constant in } \hat{I}_j   \text { for } j=1,\ldots,N \} 
$$
with  $\hat{I}_j \coloneqq [t_{j-1}, t_j)$ and use the notation
$p_{N+1}	:=p(T)$.
The Galerkin scheme
\change{\begin{equation}	
\begin{aligned}
\label{eq:disc_adjoint}
(\varepsilon A''(\nabla y_{j})\nabla \varphi, \nabla p_{j}) +
(\tfrac1\varepsilon\psi''(y_j) \varphi +
\tfrac\varepsilon{\tau_j} \varphi ,p_{j}) &
=\tfrac\varepsilon{\tau_j} (\varphi, p_{j+1}) 
\qquad  \text{ for } j=N,\ldots, 1
\end{aligned}
\end{equation}
}
starting with $p_{N+1}=y_N-y_\Omega$
then determines the approximation $p_\tau$ of $p$.
Given \eqref{ellGdiff1} and the symmetry of
$S'(g_j)$ with
$g_j= \tfrac1{\varepsilon}u_j+\tfrac1\tau_j y_{j-1}$
and $y_j=S(g_j)$ we have
$$
p_j=\tfrac1{\tau_j} S'(g_j)p_{j+1}  \qquad  \text{ for } j=N,\ldots, 1.
$$
With 
\eqref{eq:zj}
this leads to
\begin{equation}
\begin{aligned}
(z_j,p_{j+1})=(S'(g_j)(\tfrac1{\varepsilon}v_j+ \tfrac1{\tau_j}z_{j-1}),p_{j+1})
=\tfrac{\tau_j}{\varepsilon}(v_j,p_j)+(z_{j-1},p_j)
\end{aligned}
\end{equation}
for $j=N,\ldots, 1$ and consequently we have
$$
(y_N-y_\Omega,z_N)= (p_{N+1},z_N)=
\tfrac1{\varepsilon} \sum_{j=1}^N {\tau_j} (p_j,v_j)
= \tfrac1{\varepsilon} (p_\tau,v_\tau).
$$
Altogether, we have shown 
\begin{corollary}
	Under \cref{assu} the reduced cost functional \change{${j_\tau}: U_\tau \to \mathbb{R}$} is Fr{\'e}chet differentiable
	and the derivative can be represented as
\change{	\begin{equation*}
	\nabla {j_\tau}({u_\tau}) = \dfrac{1}{\varepsilon}({\lambda} {u_\tau} + p_\tau)
      \end{equation*}
    }
	where $p_\tau$ is the solution of the discrete adjoint equation \cref{eq:disc_adjoint}.
      \end{corollary}
%\newpage
% !TeX spellcheck = en_US

\section{Convergence with respect to a regularization of $A$}\label{sec3}
\anm{Die Einleitungssätze hier wurden etwas verbessert.\\
  ANM: Die Theoreme wurden neu formuliert; der zweite Teil von dem letzten Theorem ist komplett neu\\
}%
{In the previous section $A$ had to fulfill \cref{assu}.a and b. However, as mentioned in the beginning, anisotropy functions $A$ typically fulfill only \cref{assu}.a.
In order to guarantee Fr\'echet differentiability for the numerical approach we regularize such an $A$ to $A_\delta$ so that in addition \cref{assu}.b hold.}
%In order to guarantee Fr\'echet differentiability for the numerical approach we regularize $A$ which fulfills only \cref{assu}.a to $A_\delta$ for which in addition \cref{assu}.b have to hold.
{%For an example consider \cref{our_choice_of_gamma} and \cref{derAreg}, which we will discuss later.
  An example of regularization is given and discussed in \cref{derAreg}.\\
}%
In this section we consider the dependence on $\delta$ of the solutions of the in time discretized optimization problem \cref{discrete_problem1}.
To consider convergence with $\delta\to0$ we need that $A'_\delta \to A'$.
However, the results of this section do not require Fr\'echet differentiability yet, such that \cref{assu}.a on $A_\delta$ {are} sufficient.
We denote by  $y_\tau \in Y_\tau$ the solution of \cref{scheme_state1} with $A$, while $y^\delta_\tau \in Y_\tau$ shall be given as the solution of the regularized equation
%\johannes{oder \cite{CASAS199320} Th. 4.3 statt folgendes Lemma???}
\begin{equation}
\label{regstate1_v2}
\begin{aligned}
\varepsilon\dfrac{(y^\delta_{j}, \varphi) - (y^\delta_{j-1},\varphi)}{\tau_{j}}+ \varepsilon(A_\delta'(\nabla y^\delta_{j}), \nabla \varphi) + \dfrac{1}{\varepsilon} (\psi'(y^\delta_{j}), \varphi)&=(u_{j}, \varphi) \quad &j &= 1,\ldots,N\\
%(y_0, \varphi) &= (\hat{y}_0, \varphi) \quad &j&=0.
\end{aligned}
\end{equation}
and $y^\delta_\tau(0,\cdot) = y_0$.
As before we define the reduced cost functional by
{\begin{equation}
  \label{jdelta}
j_{\tau,\delta}(u_\tau) = \dfrac{1}{2} \norm{y^\delta_N(u_\tau) - y_\Omega}_{L^2(\Omega)}^2 + \dfrac{\lambda}{2\varepsilon}\norm{u_\tau}_{L^2(Q)}^2.
\end{equation}}%
{We note that to not overload the notation, $j_\delta$ is used in place of $j_{\tau,\delta}$ as long as it is clear from the context that $\tau$ is considered fixed.}
The goal of this section is to derive a convergence result {for} minimizers of a sequence of
{$j_{\tau,\delta}$ to minimizers of $j_{\tau,0}$}
in the limit $\delta \to0$
{and  to minimizers of $j$
when additionally $\tau\to0$ holds}.
Therefore some convergence behavior of the $\delta$-dependent solution $y^\delta_\tau$ is needed {that then is combined with results concerning $\tau\to0$ from \cite{BlMe21}}. This will be covered by the following two auxiliary results.

\begin{theorem}
	\label{th:conv_state_disc_wrt_delta_full_system}
	Let the \cref{assu}.a hold  and in addition
	let  $A'_\delta$
	be strongly monotone with a constant $C_A$
	independent of $\delta $
	and let
	$|A'_\delta(p) - A'(p)| \leq   \eta(\delta) $ for all $p\in \mathbb{R}^d$.
        %\\
	\kommentar{
		i)	
		$A_\delta^{'}(\cdot) : \mathbb{R}^d \to \mathbb{R}^d$ is strongly monotone, where the constant does not depend on $\delta$, i.e.
		\begin{equation*}
		(A_\delta^{'}(p)-A_\delta^{'}(q), p-q) \geq C_A |p-q|^2 \qquad \forall \delta \leq R, p,q\in \mathbb{R}^d
		\end{equation*}
		ii)
		$A_{(\cdot)}^{'}(q)$ is Lipschitz continuous, i.e.
		\begin{equation*}
		|A_\delta^{'}(\nabla \y) - A^{'}(\nabla \y)| \leq
		% C_q |\delta|.
		g(\delta) f(\| \nabla \y \|)
		\end{equation*}
	}%
        Then, for fixed $y_0$ and $u_\tau$
        and $\max_j\tau_j$ %such that
%$ \tfrac{C_\psi}{\varepsilon}\max_j\tau_j\leq c <1$
%holds,
the solutions $y_\tau(u_\tau)$ and $y_\tau^\delta(u_\tau)$ of \cref{scheme_state1} and \cref{regstate1_v2}, respectively, satisfy the following estimate
%if $ \tfrac{C_\psi}{\varepsilon}\max_j\tau_j\leq c <1$
	%	or even
\begin{equation}
\label{eq:dependence_of_y_on_tau_whole_system}
\norm{y_\tau(u_\tau)-y^{\delta}_\tau(u_\tau)}_{L^\infty(0,T;L^2(\Omega))} + \norm{\nabla y_\tau(u_\tau)-\nabla y^{\delta}_\tau(u_\tau)}_{L^2(0,T;L^2(\Omega))}
\leq {C}_{A,\psi,T}\eta(\delta).
\end{equation}
\end{theorem}
\begin{proof}
	% Let $u_{\tau}^{(i)} \in U_\tau$, $y_{0}^{(i)}\in H^1(\Omega)^N$  for $i=1,2$ and denote the corresponding solutions $y_{\tau}^{(i)}\in Y_\tau$.
We note down the differences by a prescript $\Delta$, e.g. $\Delta y_\tau \coloneqq y_{\tau} - y_{\tau}^{\delta}$.
	With $\tfrac{1}{2}(a^2-b^2) \leq (a-b)a$, %$(a-b)a = \tfrac{a^2}{2} - \tfrac{b^2}{2} + \tfrac{1}{2}(a-b)^2$
	testing the defining equations \cref{scheme_state1} and \cref{regstate1_v2} with $\Delta y_j$ and using that $A'_\delta$ is strongly monotone as well as {$(\psi'(x)-\psi'(y),x-y) \geq -C_\psi|x-y|^2$}, we obtain
	\begin{align*}
		\tfrac{1}{2}&\left(\norm{\Delta y_j}^2 - \norm{\Delta y_{j-1}}^2\right) + \tau_j C_A \norm{\nabla \Delta y_j}^2\\
		&\leq \left(\Delta y_j - \Delta y_{j-1}, \Delta y_j\right) + \tau_j\left(A'_\delta(\nabla y_{j})-A'_\delta(\nabla y^\delta _j), \nabla \Delta y_j \right)\\
		&\leq \left(\Delta y_j - \Delta y_{j-1}, \Delta y_j\right) + \tau_j\left(A'(\nabla y_{j})-A'_\delta(\nabla y^\delta _j), \nabla \Delta y_j \right) + \tau_j\left(A'_\delta(\nabla y_{j})-A'(\nabla y _j), \nabla \Delta y_j \right)\\
		&=  - \tfrac{\tau_j}{\varepsilon^2} \left(\psi'(y_{j}) -\psi'(y_j^{\delta}), \Delta y_j\right) + \tau_j\underbrace{(A'_\delta(\nabla y_j)-A'(\nabla y_j), \nabla \Delta y_j)}_{\leq \|A'_\delta(\nabla y_j)-A'(\nabla y_j)\|\|\nabla \Delta y_j\|
			\leq |\Omega|\eta( \delta)  \|\nabla \Delta y_j\|}\\
                      % &\leq \tfrac{\tau_j}{\varepsilon^2}C_\psi\norm{\Delta y_j}^2 + \tfrac{\tau_j\epsilon}{2}\norm{\nabla \Delta y_j}^2 + \tfrac{|\Omega|^2\tau_j}{2\epsilon} \eta(\delta)^2.
 &\leq \tfrac{\tau_j}{\varepsilon^2}C_\psi\norm{\Delta y_j}^2 + \tfrac{\tau_jC_A}{2}\norm{\nabla \Delta y_j}^2 + \tfrac{|\Omega|^2\tau_j}{2C_A} \eta(\delta)^2.
	\end{align*}
In the last step we used scaled Young's inequality with the scaling $C_A$.
        % $0<\epsilon <2C_A$ (not to be confused with the interface width~$\varepsilon$).%
	\kommentar{$\epsilon<1$ braucht man für $C_{\psi, \tau}>0$ später, damit nur $\tau \leq \ldots$ und nicht $\tau<\ldots$ gefordert ist. Das wiederum wird benötigt, um ohne Probleme $C_{\psi,\tau}$ unabhängig von $\tau $ abzuschätzen.}%
	%	We note that at this step we could equally well have used Cauchy-Schwarz inequality directly to the scalar product instead of the duality product to obtain \cref{eq:Lipschitz_discrete_inequ}. In the following we will only show the proof for \cref{eq:Lipschitz_discrete_inequ2}. For \cref{eq:Lipschitz_discrete_inequ} just replace the ${H^1}'$ and $H^1$ norm with the $L^2$ norm and set $\epsilon=1$ and $\tilde{C}_A = C_A$.
We now sum over $j = 1,\ldots,J$ and get
	\begin{align}
\label{eq:discr_cont_Ausgangsgleichung}
\tfrac{1}{2}\norm{\Delta y_J}^2 +
 \tfrac{C_A}{2}      %  \tilde{C}_A
\sum_{j=1}^{J} {\tau_j} \norm{\nabla \Delta y_j}^2
\leq \tfrac{1}{2}
\tfrac{|\Omega|^2}{C_A} \eta(\delta)^2\sum_{j=1}^{J} \tau_j
% \left(\sum_{j=1}^{J} \tfrac{|\Omega|^2\tau_j}{\epsilon} \eta(\delta)^2\right)
+ \tfrac{1}{2}\tilde{C}_\psi \sum_{j=1}^{J} \tau_j \norm{\Delta y_j}^2
	\end{align}
        for all $1 \leq J \leq N_\tau$. Here we defined
        % $\tilde{C}_A\coloneqq C_A-\tfrac{\epsilon}{2}>0$ and
$\tilde{C}_\psi\coloneqq\tfrac{{C}_\psi}{\varepsilon^2}$.
Omitting the gradient term on the left, absorbing the $J$-th term from the right
and using
$\tfrac{1}{
  (1-\tilde{C}_\psi\tau_J)}\leq c $,
we obtain
\begin{align*}
  \norm{\Delta y_J}^2 &\leq
 \dfrac{c |\Omega|^2}{C_A} T \eta(\delta)^2
+ c \tilde{C}_\psi \sum_{j=1}^{J-1} \tau_j \norm{\Delta y_j}^2 .
  % \dfrac{1}{	(1-\tilde{C}_\psi\tau_J)} \left(\sum_{j=1}^{J} \tfrac{|\Omega|^2\tau_j}{\epsilon} \eta(\delta)^2\right)
% + \dfrac{\tilde{C}_\psi}{1-\tilde{C}_\psi\tau_J}\sum_{j=1}^{J-1} \tau_j \norm{\Delta y_j}^2\\
%		& \leq C_{\psi,\tau} \left(\sum_{j=1}^{N_\tau} \tfrac{|\Omega|^2\tau_j}{\epsilon} \eta(\delta)^2\right) +
%		C_{\psi,\tau}{\tilde{C}_\psi}\sum_{j=1}^{J-1} \tau_j \norm{\Delta y_j}^2,
	\end{align*}
%	where we defined $C_{\psi,\tau}\coloneqq \tfrac{1}{1-\tilde{C}_\psi\max_j\tau_j}>0$.
	To this we apply the discrete Gronwall Lemma
	% (see, e.g. \cite[Lemma A.3]{Kruse2014}),
	which yields
	\begin{equation}
		\label{eq:th_disc_Gronwall_est_1}
		\begin{aligned}
                  \norm{\Delta y_J}^2 &\leq
 %                 \left( \sum_{j=1}^{N_\tau} \tfrac{|\Omega|^2\tau_j}{\epsilon} \eta(\delta)^2 \right) C_{\psi,\tau} \exp\left(C_{\psi,\tau}{\tilde{C}_\psi} \sum_{j=1}^{J-1} \tau_j\right)\\
%			&\leq \left( \sum_{j=1}^{N_\tau} \tfrac{|\Omega|^2\tau_j}{\epsilon} \eta(\delta)^2 \right)
%C_{\psi,\tau} \exp\left(C_{\psi,\tau}{\tilde{C}_\psi T}\right).
\dfrac{c |\Omega|^2}{C_A} T \eta(\delta)^2
\exp\left(c \tilde{C}_\psi T \right).
\end{aligned}
	\end{equation}
	Inserting this into \cref{eq:discr_cont_Ausgangsgleichung} we finally get for all $J=1,\ldots,N_\tau$
\begin{equation}
\label{eq:th_disc_Gronwall_est_2}
{C_A}\sum_{j=1}^{J} {\tau_j} \norm{\nabla \Delta y_j}^2 \leq
\dfrac{ |\Omega|^2}{C_A} T \eta(\delta)^2
\left(1 + c{\tilde{C}_\psi T}
  \exp\left(c \tilde{C}_\psi T \right) \right)
%\left(\sum_{j=1}^{N_\tau} \tfrac{|\Omega|^2\tau_j}{\epsilon} \eta(\delta)^2 \right)
%\left(1 + C_{\psi,\tau}{\tilde{C}_\psi T}
%\exp\left(C_{\psi,\tau}{\tilde{C}_\psi T}\right)\right),
	\end{equation}
	%	\begin{align*}
	%		\tfrac{1}{2}\norm{\Delta y_J}^2 + C_A \sum_{j=1}^{J} {\tau_j} \norm{\nabla \Delta y_j}^2 \leq \left(\norm{\Delta y_0}^2 + \sum_{j=1}^{N_\tau} \tfrac{\tau_j}{\epsilon} \norm{\Delta u_j}_{{H^1}'}^2\right)\left(\tfrac{1}{2}+ \tfrac{(\tfrac{1}{2}+C_\psi)T}{(1-\tilde\tilde{C}_\psi\tau)} \exp\left(\dfrac{\tilde\tilde{C}_\psiT}{1-\tilde\tilde{C}_\psi\tau}\right)\right),
	%	\end{align*}
	which together with \cref{eq:th_disc_Gronwall_est_1}
%	and the boundedness of $C_{\psi,\tau}$ independently of $\tau$
	yields the inequality \cref{eq:dependence_of_y_on_tau_whole_system}. % after taking the maximum.
\end{proof}

%
%\begin{thrm}
%	\label{th:conv_state_disc_wrt_delta}
%Let the \cref{assu}.a hold  and in addition
% let  $A'_\delta$
% be strongly monotone with a constant $C_A$
% independent of $\delta $
% and let
% $|A'_\delta(p) - A'(p)| \leq   \eta(\delta) $ for all $p\in \mathbb{R}^d$.
%
%
%Then it holds my{ANM: ich hab da bei meinem noch eine Anmerkung drin die ich nicht mehr versteh}
%	\begin{equation}
%	\label{ass1}
%	\| y^\delta_\tau (u_\tau) - y_\tau(u_\tau)\|_{H^1(\Omega)^N}
%        \leq C_\tau \eta( \delta) 
%	\end{equation}
%where $y^\delta_\tau(u_\tau)$ are the solutions of \cref{regstate1_v2}.
%\end{thrm}
%

We note that here the
constant $ {C}_{A,\psi,T} $ depends exponentially on the interface thickness $\varepsilon$ as can be seen in~\cref{eq:th_disc_Gronwall_est_2}.
When studying the dependence on $\varepsilon$%
---which is not subject of this paper---%
a more careful analysis in terms of $\varepsilon$
possibly not based on the Gronwall Lemma
is necessary.

Combining the estimate
\eqref{eq:dependence_of_y_on_tau_whole_system}
with results in \cite{BlMe21}
%, respectively \cref{eq:Lipschitz_discrete_inequ2} and
% \cref{th:conv_state_disc_wrt_delta_full_system}
we obtain:
  \kommentar{Entspricht \cite[Th. 4.3]{CASAS199320}, wo es aber ganz anders bewiesen wird}%
 \begin{corollary}
  \label{th:conv_y_delta_if_conv_u}
Let the assumptions of \cref{th:conv_state_disc_wrt_delta_full_system} be fulfilled and $u_\tau,\tilde{u}_\tau\in U_\tau$ be given.
Then the estimate
\begin{equation}
\label{eq:dependence_of_y_on_tau_and_u_whole_system}
%\begin{aligned}
	\norm{y_\tau(u_\tau)-y^{\delta}_\tau(\tilde{u}_\tau)}_{L^\infty(0,T;L^2(\Omega))} + \norm{\nabla y_\tau(u_\tau)-\nabla y^{\delta}_\tau(\tilde{u}_\tau)}_{L^2(0,T;L^2(\Omega))}
	\leq {C}_{A,\psi,T}\big(\eta(\delta) + \norm{u_\tau-\tilde{u}_\tau}_{L^2(0,T;H^1(\Omega)')} \big)
%\end{aligned}
\end{equation}
holds.
Hence, given a sequence $(u_\tau)_\tau$ with
$u_\tau \rightharpoonup u$ in $L^2(0,T;L^2(\Omega))$ for $\tau \to 0$, there exists $\sigma(\tau)$ with $\lim_{\tau\to 0} \sigma(\tau) =0 $ such that
\begin{equation}
\label{eq:dep_tau_delta_L2_L2}
%\norm{y(u) - y_\tau^\delta(u_\tau)}_{L^\infty (0,T;L^2(\Omega))}
 \max_{t\in [0,T]} \norm{y(u)_{|t}  - y_\tau^\delta(u_\tau)_{|t} }_{L^2(\Omega)} 
 \leq C\left(\eta(\delta)+\sigma(\tau)\right).
 % \rho(\tau, \delta),
\end{equation}
%Furthermore, we have for all $t\in[0,T]$ that
%$y_\tau^\delta(u_\tau)_{|t} \to y(u)_{|t}$ in $L^2(\Omega)$
% and $y_\tau^\delta(u_\tau)\to y(u)$ in $L^2(0,T;L^2(\Omega))$
%for $\tau,\delta \to0$.
\end{corollary}

\begin{proof}
  Estimate \cref{eq:dependence_of_y_on_tau_and_u_whole_system} follows by zero completion with $y_\tau(\tilde{u}_\tau)$, triangle inequality and estimating the resulting terms by \cref{th:conv_state_disc_wrt_delta_full_system} and
  \cite[Theorem 2.4]{BlMe21}, respectively \cref{eq:Lipschitz_discrete_inequ2}.\\
  \kommentar{$\norm{y_\tau(u_\tau)-y^{\delta}_\tau(u_\tau^\delta)} \leq \norm{y_\tau(u_\tau)-y_\tau(u_\tau^\delta)} + \norm{y_\tau(u_\tau^\delta)-y^{\delta}_\tau(u_\tau^\delta)}$}%
  For the second estimate we recall from \cite[Theorem 2.6]{BlMe21} that
  % we have $y_\tau \to y$ in $L^2(0,T;L^2(\Omega))$
  if $u_\tau \rightharpoonup u$ in $L^2(0,T;L^2(\Omega))$ there exists $\sigma(\tau)$  with $\lim_{\tau\to 0} \sigma(\tau) =0 $ such that
%  $\sigma(\tau) \to 0$ for $\tau\to0$ and $\norm{y_\tau(u_\tau)-y(u)}_{L^2(0,T;L^2(\Omega))} \leq \sigma(\tau)$.
  $\max_{t\in [0,T]} \norm{y(u)_{|t}  - y_\tau(u_\tau)_{|t} }_{L^2(\Omega)} 
 \leq C \sigma(\tau)
  $.
 By inserting $y_\tau(u_\tau)$ and using the triangle inequality together with the first estimate one obtains %note that the right hand side vanishes since initial value and control are the same
%  \begin{equation*}
%  	\norm{y(u)-y_\tau^\delta(u_\tau)}_{L^2(0,T;L^2(\Omega))} \leq C\left(\eta(\delta)+\sigma(\tau)\right),
%  \end{equation*}
%defining $\rho(\tau,\delta)$.
\cref{eq:dep_tau_delta_L2_L2}.%\\
%Finally, the last convergences follow from
%\begin{equation*}
%	\norm{y_\tau^\delta(u_\tau)_{|t} - y(u)_{|t}}_{L^2(\Omega)} \leq \underbrace{\norm{y_\tau^\delta(u_\tau)_{|t} - y_\tau(u_\tau)_{|t}}_{L^2(\Omega)}}_{\leq C_{A,\psi,T}\eta(\delta) \text{ due to \cref{eq:dependence_of_y_on_tau_and_u_whole_system}}} + \underbrace{\norm{y_\tau(u_\tau)_{|t} - y(u)_{|t}}_{L^2(\Omega)}}_{\to 0 \text{ due to \cite[Theorem 2.6]{BlMe21}}} \to 0,
%\end{equation*}
%where we used that the constant in \cref{eq:dependence_of_y_on_tau_and_u_whole_system} is independent from $\tau$ and analogously for the convergence in $L^2(0,T;L^2(\Omega))$.
\end{proof}
%
%The just shown result in particular demonstrates that if $\eta(\delta)\to0$, $y_0^\delta\to y_0$ and $u_\tau^\delta \rightharpoonup \ubar{u}_\tau$ in $U_\tau$ for $\delta \to 0$, then it holds $y^\delta_\tau(u^\delta_\tau) \to y_\tau(\ubar{u}_\tau)$ in $Y_\tau$.
%
{Finally, we finish this section with the following convergence result of global minimizers.}
\kommentar{The proof of the first part of the theorem uses ideas from \cite[Theorem 4.4]{Casas199320}.
The second part is a generalization of \cite[Theorem 3.3]{BlMe21}.}%
\begin{theorem}
  \label{th:conv_op_delta}
Let the assumptions of  \cref{th:conv_state_disc_wrt_delta_full_system} be fulfilled and $\lim_{\delta \to 0}\eta(\delta)=0$. Denote by $u^\delta_\tau$ a global minimizer of $j_{\tau,\delta}$.
Then it holds:
   % Furthermore let $u^\delta_\tau \in L^2(\Omega)^N$ be a minimizer of $j_\delta$.
   \begin{enumerate}
   	\item 
   	Considering $\delta\to0$ for fixed $\tau>0$, there exists a subsequence such that
   	 it holds $u^\delta_\tau \to \ubar{u}_\tau$ in $U_\tau$,
   	$y^\delta_\tau (u^\delta_\tau) \to y_\tau (\ubar{u}_\tau)$ in $Y_\tau$ and $j_{\tau,\delta}(u^\delta_\tau) \to j_\tau(\ubar{u}_\tau)$ for $\delta \to 0$. Furthermore $\ubar{u}_\tau$ is a global minimizer of $j_{\tau}$.
   	\item
   	Considering $\tau, \delta\to0$, there exists a subsequence such that it holds $u^{\delta}_{\tau}  \to \ubar{u}$ in $L^2(0,T;L^2(\Omega))$, $y^{\delta}_{\tau} (u^{\delta}_{\tau}) \to y (\ubar{u})$ in $L^2(0,T;L^2(\Omega))$ and $j_{{\tau,\delta}}(u^{\delta}_{\tau}) \to j(\ubar{u})$. Furthermore $\ubar{u}$ is a global minimizer of $j$.
   \end{enumerate}
\end{theorem}
\begin{proof}
\begin{enumerate}
\item	  Take $\bar{u}_\tau\in U_\tau$ fixed. From \cref{th:conv_state_disc_wrt_delta_full_system}
  we obtain $y_\tau^\delta(\bar{u}_\tau) \to y_\tau(\bar{u}_\tau) $ {in $Y_\tau$} for $\delta \to 0$ and therefore from the boundedness of this sequence we obtain
{  $$\dfrac{\lambda}{2\varepsilon}
  \norm{u^\delta_\tau}_{L^2(Q)}^2 \leq j_{\tau,\delta}(u^\delta_\tau) \leq j_{\tau,\delta}(\bar{u}_\tau) \leq C.$$}
Hence $u^\delta_\tau \rightharpoonup \ubar{u}_\tau \in U_\tau$ for a subsequence, which is considered in the following, and consequently
          $y^\delta_N(u^\delta_\tau) \to y_N(\ubar{u}_\tau)$
          in $L^2(\Omega ) $my{,} see \cref{th:conv_y_delta_if_conv_u}.
          Using the definition of $j_{\tau,\delta}$ in \cref{jdelta} leads to
	\begin{equation}
	\label{eq:inequ_chain_j_delta}
		j_{\tau}(\ubar{u}_\tau) \leq \liminf_{\delta\to0} j_{\tau,\delta}(u^\delta_\tau) \leq \limsup_{\delta\to0} j_{\tau,\delta}(u^\delta_\tau) \leq \lim_{\delta\to0} j_{\tau,\delta}(u_\tau) = j_\tau(u_\tau) \qquad \forall u_\tau \in U_\tau .
	\end{equation}
Hence %$j(\ubar{u}) \leq j(u)\ \forall u\in L^2(\Omega)^N$ and
	$\ubar{u}_\tau$ is a minimizer.
Since we can also choose $\ubar{u}_\tau$ on the righter part of \cref{eq:inequ_chain_j_delta}, in addition we obtain $j_{\tau,\delta}(u^\delta_\tau) \to j_\tau(\ubar{u}_\tau)$.
\\
{Since we already have $u^\delta_\tau \rightharpoonup \ubar{u}_\tau$,} to obtain the strong convergence $u^\delta_\tau \to \ubar{u}_\tau$ {in $U_\tau$} {it remains to} check that the norms converge. This follows from {
\begin{equation*}
  \tfrac{\lambda}{2\varepsilon}\norm{u^\delta_\tau}^2_{L^2(Q)}
  = j_{\tau,\delta}(u^\delta_\tau) - \tfrac{1}{2}\norm{y^\delta_N(u^\delta_\tau) - y_\Omega}_{L^2(\Omega)}^2 \to j_\tau(\ubar{u}_\tau) - \tfrac{1}{2}\norm{y_N(\ubar{u}_\tau) - y_\Omega}_{L^2(\Omega)}^2 = \tfrac{\lambda}{2\varepsilon}\norm{\ubar{u}_\tau}^2_{L^2(Q)}.
 % \qedhere
\end{equation*}
}
%\end{proof}
\kommentar{um $u_\delta \to \ubar{u}$ strong gegen eine bestimmtes Optimum zu erhalten muss man den $\norm{u - \ubar{u}}^2$ Term wie in Casas in dem Kostenfunktional haben. S. z.B. Kapitel 4 in \cite{CASAS199320} oder andere Casas-papers.}%
\kommentar{aside: $(A'(\nabla \cdot, \nabla v))$ erfüllt die $S^+$-property und ist damit auch vom Typ-M (umgekehrt denke ich folgt aber nicht aus Typ-M die $S^+$-property, weil man aus den Vor. von $S^+$ nicht die boundedness und so etwas wie $Ay_n \rightharpoonup z$ folgern kann was man aber brauch um Typ-M anwenden zu können). Ich hab das hier notiert, weil Casas Typ-M in seinem analog zu Theorem 5.2. braucht.}%
%
%As a generalization of \cref{th:conv_op_delta} and \cite[Theorem 3.3]{BlMe21}, in conclusion we show the following result.
%\begin{thrm}
%	\label{th:conv_op_delta_and_tau}
%	Let the assumptions of  \cref{th:conv_state_disc_wrt_delta_full_system} be fulfilled and $\lim_{\delta \to 0}\eta(\delta)=0$.
%	% Furthermore let $u^\delta_\tau \in L^2(\Omega)^N$ be a minimizer of $j_\delta$.
%	Then there exists a sequence $(\tau_l,\delta_l) \to (0,0)$ such that for the corresponding minimizers $u^{\delta_l}_{\tau_l}$ of $j_{\delta_l}$ it holds $u^{\delta_l}_{\tau_l}  \to \ubar{u}$ in $L^2(0,T;L^2(\Omega))$, $y^{\delta_l}_{\tau_l} (u^{\delta_l}_{\tau_l}) \to y (\ubar{u})$ in $L^2(0,T;L^2(\Omega))$ and $j_{{\delta_l}}(u^{\delta_l}_{\tau_l}) \to j(\ubar{u})$ for $l \to \infty$. Furthermore $\ubar{u}$ is a minimizer of $j$.
%\end{thrm}
%\begin{proof}
\item
%	For the sake of presentation we omit the index $l$ and denote subsequences the same way.
First we choose an arbitrary but fixed $\bar{u} \in L^2(0,T;L^2(\Omega))$ and a sequence $\bar{u}_\tau\in U_\tau$ with {$\lim_{\tau\to0}\bar{u}_\tau = \bar{u}$} in $L^2(0,T;L^2(\Omega))$.
	%	From the Lipschitz-continuity (Theorem \ref{th:lipschitz_disc}) we obtain
 Hence {$ \lim_{\delta,\tau\to0}y^\delta_{\tau}(\bar{u}_\tau)_{|T} = y(\bar{u})_{|T}$ in	$ L^2(\Omega)$}
due to \cref{eq:dep_tau_delta_L2_L2}.
	%	\begin{equation*}
	%		\norm{y^\ast_{\tau, j}}_{L^2(\Omega)} \leq C (\norm{y_0}_{L^2(\Omega)} + \underbrace{\norm{u^\ast_\tau}_{L^2(0,T;L^2(\Omega))}}_{\leq C}) \leq C
	%	\end{equation*}
	%	for all $j=1,\ldots,N$. So in particular $y^\ast_{\tau,N} \leq C$ independently from $\tau$.
\kommentar{	Now let $(u^\delta_\tau)_{\tau,\delta}$ be a sequence of global minimizers to \cref{jdelta} with respect to \cref{regstate1_v2}.}%
	%	It holds
	%		\begin{equation*}
	%		J_\tau(y_\tau, u_\tau) \leq J_\tau(y^\ast_\tau, u^\ast_\tau) = \dfrac{1}{2}\sum_{j=1}^N \tau_j \|y^\ast_{\tau,j}-y_{Q,j}\|^2 + \dfrac{\lambda}{2 \varepsilon} {\norm{u^\ast_\tau}^2} \leq CT+ C,
	%		\end{equation*}
	%		where the first inequality follows from global optimality and in the last inequality we use the just shown boundedness and that $u^\ast_\tau$ is bounded due to convergence.
As above it holds
	\begin{equation*}
\dfrac{\lambda}{2\varepsilon}\norm{u^\delta_\tau}^2_{L^2(Q)} \leq\dfrac{1}{2}\|y^\delta_{\tau}(\bar{u}_\tau)_{|T}-y_{\Omega}\|^2_{L^2(\Omega)} + \dfrac{\lambda}{2\varepsilon} {\norm{\bar{u}_\tau}^2_{L^2(Q)}} \leq  C.
\end{equation*}
%{which is bounded due to the convergences.}
%which is now bounded independent of $\tau,\delta$ for small enough $\tau$ and $\delta$.
\kommentar{
  It holds
{
	\begin{equation*}
\dfrac{\lambda}{2\varepsilon}\norm{u^\delta_\tau}^2_{L^2(Q)} \leq J(y^\delta_\tau(u^\delta_\tau), u^\delta_\tau) \leq J(y^\delta_\tau(\bar{u}_\tau), \bar{u}_\tau) = \dfrac{1}{2}\|y^\delta_{\tau}(\bar{u}_\tau)_{|T}-y_{\Omega}\|^2_{L^2(\Omega)} + \dfrac{\lambda}{2\varepsilon} {\norm{\bar{u}_\tau}^2_{L^2(Q)}}\leq c .
              \end{equation*}
            }
            }%
\kommentar{where the second inequality follows from global optimality and in the last inequality we use the boundedness of the appearing quantities due to the just stated convergences.
  This implies that $(u^\delta_\tau)_{\tau,\delta}$ is bounded in $L^2(0,T;L^2(\Omega))$ and}%
Hence we can deduce a subsequence denoted in the same way with $u^\delta_\tau \rightharpoonup \ubar{u}$ in $L^2(0,T;L^2(\Omega))$.
	%Since $u_\tau \in U_\tau$, from the proof of Theorem \ref{existence_state}
Then \cref{th:conv_y_delta_if_conv_u} yields that {$y^\delta_\tau(u^\delta_\tau)_{|T}\to y(\ubar{u})_{|T}$} in $L^2(\Omega)$ and hence
$j(\ubar{u})\leq \liminf_{\tau,\delta\to 0}
j_{\tau,\delta}(u^\delta_\tau)$.
\kommentar{DENKE DRAN, dass Konvergenz  $\delta,\tau \to 0$ meint:
$u^\delta_\tau \rightharpoonup \ubar{u}$ bedeutet dass es eine Folge $(\delta_n,\tau_n)_n $ gibt mit $u^{\delta_n}_{\tau_n} \rightharpoonup \ubar{u}$ und somit gilt  $y^{\delta_n}_{\tau_n}(u^{\delta_n}_{\tau_n})_{|T}\to y(\ubar{u})_{|T}$ für $n\to \infty$
}%
	%as specified in \cref{eq:various_convergences_of_y_tau}.
	%Note that the limit also is $y(u)$ since for the analysis done there the weak convergence of $u_\tau$ is actually sufficient.
Respectively, given some arbitrary
{$ \tilde{u} \in L^2(0,T;L^2(\Omega))$ and a sequence $\tilde{u}_\tau$ with $\tilde{u}_\tau \to \tilde{u}$
  % in $L^2(0,T;L^2(\Omega))$
  }
  we obtain
  $\lim_{\tau,\delta\to 0} j_{\tau,\delta}(\tilde{u}_\tau)
=j(\tilde u) $.
\kommentar{$\lim_{\delta,\tau \to 0} {y}^\delta_\tau(\tilde{u}_\tau)|_T = {y}(\tilde{u})|_T$.}%
  % the latter also for ${y}^\delta_\tau(\tilde{u}_\tau)$ and ${y}(\tilde{u})$.
  \kommentar{This yields
	\begin{equation}
	\label{eq:conv_disc_and_delta_to_cont_inequality_chain}
	J(y(\ubar{u}),\ubar{u}) \leq \liminf_{\tau,\delta\to 0} J(y^\delta_\tau(u^\delta_\tau), u^\delta_\tau) \leq \liminf_{\tau,\delta\to 0} J({y}^\delta_\tau(\tilde{u}_\tau), \tilde{u}_\tau) = J({y}(\tilde{u}), \tilde{u}).
      \end{equation}
      }%
{Then the assertions follows as in (1).\qedhere}
      \kommentar{
	where in the first inequality we used the weak lower-semicontinuity of the norm and strong convergence of $y^\delta_{\tau}(u^\delta_\tau)_{|T}$, then the global optimality of $(y^\delta_\tau(u^\delta_\tau), u^\delta_\tau)$ and finally the strong convergence of $y^\delta_\tau(\tilde{u}_\tau)_{|T}$ in $L^2(\Omega)$ and $\tilde{u}_\tau$ in $L^2(0,T;L^2(\Omega))$, respectively.
	Since $\tilde{u}$ was arbitrary this yields the global optimality of~$u$. Plugging in $\tilde{u} = \ubar{u}$ yields the convergences $j_{\tau,\delta}(u^\delta_\tau)\to j(\ubar{u})$, $\|u^\delta_\tau\| \to \|\ubar{u}\| $
	and therefore
	with the weak convergence also
	the strong convergence $u^\delta_\tau \to \ubar{u}$ in $L^2(0,T;L^2(\Omega))$.
	The stated convergence of $y^\delta_\tau(u^\delta_\tau)$ follows from \cref{th:conv_y_delta_if_conv_u}.}%
\end{enumerate}
\end{proof}
%\newpage
% !TeX spellcheck = en_US

\section{The regularization of a %special
class of anisotropies}% and its properties.}
\label{sec:regularization}
%In the chapters that follow, we will consider the practical realization of a solver for the optimization problem \crefrange{eq1a}{weakAC}. Solvers typically search for a solution by employing gradients or considering the first order conditions. The arising adjoint equation is to be expected to contain the second derivative $A''$. Since for 2-homogeneous $A$ this is a 0-homogeneous function it is in general not continuous at the origin unless $A$ is induced by a scalar product. This discontinuity does not allow the standard formulation of the first order condition and it is an open problem how they are given.
%\johannes{Approaches could be oriented to \cite{Casas1991} where the author considers another kind of anisotropy and a time-independent state equation.}
%Furthermore for the numerical solution one has to use a nonsmooth optimization solver or to regularize e.g. the function $A$. Here we choose the latter approach. Let us mention, that although $A''$ does not exist in the origin only, one has to keep in mind that $\nabla y$ is (near) zero except on the interface between two phases.

Before we continue with {simulations for optimal control of anisotropic phase field models}
we have to specify the {anisotropy} function $A$. As mentioned in the introduction, this function %related to anistropic phase field models
typically is $2$-homogeneous. This in general however 	conflicts with the requirement of $A$ being twice continuously differentiable.
\kommentar{Our approach would be to apply a regularization, but the process is unclear without further information on $A$. Therefore this section's goal is to restrict it to a special form where we can apply a certain regularization
}%end kommentar
Therefore this section's goal is to specify the employed $A$,
to introduce an appropriate regularization $A_\delta$ 
and to show that
$A$ satisfies \cref{assu}a. and $A_\delta$  fulfills in addition  \cref{assu}b.
%both---the regularized and the original anisotropy function---satisfy \cref{assu}a. and b.
This guarantees that the results from \cite{BlMe21} and the preceding chapters can be applied.
%As mentioned in the introduction we restrict it to be of the form (cf. \cref{eq:form_of_A})
First, recall from the introduction that $A$ can be written as
\begin{equation}
  \label{eq:form_of_A}
  A(p) = \tfrac{1}{2}|\gamma(p)|^2 \qquad \forall p \in \mathbb{R}^d,
  \end{equation}
where the so-called density function $\gamma: \mathbb{R}^d \to \mathbb{R}_{\geq0}$ with $\gamma \in C^2(\mathbb{R}^d\setminus\{0\}) \cap C(\mathbb{R}^d)$ shall be {positive} 1-homogeneous.
{The terminology `density function' goes back to the study of sharp interface models, where the surface energy of the interface between a solid and liquid phase, say, is given by
%\begin{equation}
%	\label{eq:surface_energy_aniso}
$	|\Gamma|_\gamma = \int_{\Gamma} \gamma (\nu) \ ds$.
%\end{equation}
In the isotropic case $\gamma(p) = |p|$ this would reduce to the area of the interface $|\Gamma|_\gamma = |\Gamma|$. The authors of \cite{alfaro2009motion, Elliott1996} show for the Allen-Cahn equation \cref{weakAC}---with $A$ defined as in \cref{eq:form_of_A}---that in the limit $\varepsilon\to 0$ the zero level sets converge to a sharp interface $\Gamma$ moving with $V = \gamma(\nu)\kappa_\gamma$ if $u=0$.}
%--which as explained in the introduction arises from a gradient flow of an energy containing a term of the form \cref{eq:form_of_A}--
{While there exist several approaches to define $\gamma$, like e.g. in \cite{Kobayashi}
\kommentar{Kobayashi anisotropy}%
or in \cite{DDE}},
\kommentar{regularized $l1$-norm laut BGN; auf Seite 53 sind aber mehrere Ansätze}%
%\cite{Graser2013} enthält beide als Beispiel
we constrain ourselves to a class of anisotropies for which the density function $\gamma $ is introduced in \cite{BGN2007}. The corresponding phase field ansatz is studied e.g. in \cite{bgn13}.
%\johannes{and an even further generalization we do not discuss here is given in \cite{Barrett2008}. ODER IST DANN DER NÄCHSTE SATZ VERWIRREND}
In the following they are referred to as BGN-anisotropies. They allow for the modelling and approximation of a large class of common anisotropies. Also they are well suited to model crystal growth, since crystals build characteristical faces.
The basic observation is that for the metric $(\cdot,\cdot)_{\tilde{G}}$ defined by symmetric positive definite $\tilde{G}$, the surface area element can be expressed as $\gamma(\nu) = \sqrt{\nu^T G \nu}$ with $G = \det(\tilde{G})^{\tfrac{1}{d-1}} \tilde{G}^{-1}$ (see \cite{Barrett2008}). This motivates the choice of the class of density functions $\gamma$ given by
\begin{equation}
\label{our_choice_of_gamma}
	\quad \gamma(p) = \sum_{l=1}^{L} \gamma_l(p),
        \quad \mbox{ where }
        \gamma_l(p) = \sqrt{p^T G_l p}
      \end{equation}
      and $G_l \in \mathbb{R}^{d\times d}$ are symmetric and positive definite.
Note that for $p\neq 0$ the derivative of $A$ can then be computed as
\begin{equation}
\label{derA}
A'(p) = \gamma(p){\gamma}'(p)=\sum_{l,m} \dfrac{\gamma_m(p)}{\gamma_l(p)} G_l p
\end{equation}
and $A'$ is continuous also at $p=0$ with $A'(0)=0$.

The second derivative exists for $p\neq 0$  and is given by
\begin{equation}
\label{derAprime0__}
%A''_\delta(p) = \gamma_\delta(p) \gamma''_\delta(p) + \gamma'_\delta(p) \otimes \gamma'_\delta(p),
A''(p) = \gamma(p) \gamma''(p) + \gamma'(p) \gamma'(p)^T,
\end{equation}
where
%$$\gamma''_\delta(p) = \sum_{l} \left(\dfrac{G_l}{\gamma^\delta_l(p)} - \dfrac{G_l p \otimes G_l p}{\gamma^\delta_l(p)^3}\right).$$
$$\gamma''(p) = \sum_{l} \left(\dfrac{G_l}{\gamma_l(p)} - \dfrac{G_l p (G_l p)^T}{\gamma_l(p)^3}\right).$$

We note that $A'': \mathbb{R}^d\setminus\{0\} \to \mathbb{R}^{d\times d}$ is continuous  and
$A''(p)$ is positive definite with constants independent of $p$. %and $A''$ even induces a uniformly equivalent norm on $\mathbb{R}^d$.
Moreover we have Lipschitz-continuity and strong monotonicity of $A'$.
%, and consequently, \johannes{given $A'(0)=0$, (kann das weg?)} \cref{assu} are fulfilled.
These properties follow from results in \cite{Graser2013} where the authors need in addition to the given properties  of $\gamma$, namely continuity on $\mathbb{R}^d$, twice continuously differentiability on
$\mathbb{R}^d \setminus \{0\}$, positive homogeneity of degree one, $\gamma(p)>0$ for $p\neq 0$ and
the following relation
\begin{equation*}
	q^T\gamma''(p)q \geq C |q|^2 \qquad \forall p,q \in \mathbb{R}^d \text{ with } p^T q = 0 \text{ and } |p| = 1.
\end{equation*}
The latter can be shown by an application of the Cauchy-Schwarz inequality
\begin{equation*}
	q^T \gamma''(p) q = \sum_{l} \left(\dfrac{q^TG_lq}{\gamma_l(p)} - \dfrac{(q^TG_l p)^2}{\gamma_l(p)^3}\right) \geq \sum_{l} \left(\dfrac{q^TG_lq}{\gamma_l(p)} - \dfrac{\gamma_l(p)^2\gamma_l(q)^2}{\gamma_l(p)^3}\right) = 0,
\end{equation*}
where equality does not hold for $p^T q = 0$ and the compactness of the set given by $p^T q = 0$, $\| p\|=\|q\|=1$.
%%%%%%%%%%%%%%%%%%%%%%%%%

Our goal for regularizing $A$ is that $A_\delta\in C^2(\mathbb{R}^d) $  shall fulfill the requirements for the existence of an optimal control, that the derivative shall be simple to evaluate and that the influence on the interfacial region (i.e. $\nabla y \not \approx 0$) shall be little.
Our approach is to modify the $\gamma_l$, but one could also think of regularizing e.g. the quotient appearing in the sum in \cref{derA}.
%One obvious regularization is given by $\gamma_l \to (\gamma_l)^{1+\delta}$, which yields the smooth $A_\delta'(p) = (1+\delta) \sum_{l,m} \tfrac{\gamma_m(p)^{1+\delta}}{\gamma_l(p)^{1-\delta}} G_l p$. However,
%this changes the asymptotical as well as the behaviour of $A_\delta'$ near zero rendering \cref{assu} no longer to be satisfied.
%Furthermore we observed numerically that the modified power properties cause the problematic almost-zero values away from the interface to become even smaller and at the interface where the gradient of the state admits large norms the function will deviate strongly from the unregularized one. \johannes{was ist gelöst worden. Neu formulieren. das ist so, warum muss das numerisch beobachtet werden?}
%
%In the following we will regularize the anisotropy and the density function by substituting $\gamma_l $ with
Among various choices we considered the most promising was to alter the functions $\gamma_l$ by a small shift of $\delta$, i.e.	
\begin{equation}
\label{reg_choice}
 \gamma_l^\delta \coloneqq \sqrt{\gamma_l^2+\delta}
      \end{equation}
      where $\delta>0$.
This we use in the following and denote the resulting regularizations by $A_\delta$ and $ \gamma_\delta $.  Both are now in $C^\infty(\mathbb{R}^d)$.
%Here the asymptotical behavior stays the same and the small addition under the square-root should be irrelevant at the interface.
A very  convenient property  for this choice is that $\gamma_\delta(p) = \tilde{\gamma}((p, \sqrt{\delta})^T)$ where $\tilde \gamma$ is defined employing the matrices $\tilde{G}_l \coloneqq \left(\begin{smallmatrix}
G_l &\\&1
\end{smallmatrix}\right)$.
Hence one can also view the regularized anisotropy $A_\delta$ on $\mathbb{R}^d$ as an unregularized BGN-anisotropy $\tilde A$ on $\mathbb{R}^{d+1}$ for which above properties hold.
%     Note that one can also view the regularized anisotropy $A_\delta$ as an unregularized one of the same form for $\gamma$ as in \cref{our_choice_of_gamma} just on a higher dimensional space. Therefore define $\tilde{G}_l \coloneqq \left(\begin{smallmatrix}
% G_l &\\&1
% \end{smallmatrix}\right)$ and  $\tilde{A}$ the same way as before just with the new matrices. Both formulations are connected by $\gamma_\delta(p) = \tilde{\gamma}((p, \sqrt{\delta})^T)$.

The derivatives still have the same structures as in  \cref{derA} and \cref{derAprime0__} namely
\begin{align}
\label{derAreg}
A'_\delta(p) &= \sum_{l,m} \dfrac{\gamma^\delta_m(p)}{\gamma^\delta_l(p)} G_l p ,
\\
% \end{equation}
%
%For completeness and since needed later, we also gi%ve the second derivative
%\begin{equation}
\label{derAprime0}
%A''_\delta(p) = \gamma_\delta(p) \gamma''_\delta(p) + \gamma'_\delta(p) \otimes \gamma'_\delta(p),
A''_\delta(p) &= \gamma_\delta(p) \gamma''_\delta(p) + \gamma'_\delta(p) \gamma'_\delta(p)^T,
\\
\mbox{with } \quad
%$$\gamma''_\delta(p) = \sum_{l} \left(\dfrac{G_l}{\gamma^\delta_l(p)} - \dfrac{G_l p \otimes G_l p}{\gamma^\delta_l(p)^3}\right).$$
\label{der2gamma}
\gamma''_\delta(p) &= \sum_{l} \left(\dfrac{G_l}{\gamma^\delta_l(p)} - \dfrac{G_l p (G_l p)^T}{\gamma^\delta_l(p)^3}\right),
\end{align}
though these hold in the regularized version for all $p\in \mathbb{R}^d$.
\change{Note that for $L=1$, i.e. for $A(p)= \tfrac 12 p^TG p$ which includes the isotropic case, it holds $A'_\delta=A'$.
This is a particularly convenient property
as in this case $A$ is already smooth by itself and hence there is no need for regularization anyway.}
%This is especially convenient, since there is no need for regularization in this case anyway.
Due to 
\begin{equation}\label{AAtilde}
  A_\delta'(p) =\left(\tilde{A}'\left(\begin{smallmatrix}
	p\\\sqrt{\delta}
      \end{smallmatrix}\right)\right)_{1,\ldots,d},
  \quad
    A_\delta''(p) =\left(\tilde{A}''\left(\begin{smallmatrix}
	p\\\sqrt{\delta}
\end{smallmatrix}\right)\right)_{{1,\ldots,d}\atop{1,\ldots,d}}
\end{equation}
for $\delta > 0$ the Lipschitz-continuity and strong monotonicity of $\tilde A'$ provide these properties for $A_\delta'$  with constants independent of $\delta$.
%Note also that \johannes{$A'_\delta(0) = 0$ and therefore (kann weg?)} $|A'_\delta(p)|<C|p| $ holds.
Moreover, since $\tilde A''$ induces uniformly
equivalent norms on $\mathbb{R}^{d+1}$ the same holds for $A_\delta''$.
Hence $A_\delta''(p)$ is bounded independent of $p$.
Using $A_\delta(0) \geq A(0)=0$ %
\kommentar{from Taylors expansion}%
we obtain $c|p|^2 \leq A_\delta(p)  $ with $c>0$.
This inequality holds also for $A$ due to the 2-homogeneity.
Moreover, all constants can be chosen independently of $\delta$. 
The only exception is the upper bound in the growth condition
$A_\delta (p) \leq A_\delta(0) +C|p|^2$ due to
$A_\delta(0)$.
%gilt das auch für die wo man nicht direkt aus Lipschitz und strong monotonicity ableiten kann?
%ja z.b.
%\begin{eqnarray}
%	|A'_\delta(p)| = |\tilde{A}'(\begin{smallmatrix}
%	p\\\sqrt{\delta}
%	\end{smallmatrix}) \geq C |\begin{smallmatrix}
%	p\\\sqrt{\delta}
%	\end{smallmatrix}| \stackrel{\text{normequiv.}}{\geq} C(|p| + \sqrt{\delta})) \geq C|p|\\
%	\sqrt{p^TG_lp + \delta} \geq \sqrt{p^TG_lp} \geq C|p| \implies A^2 \geq C|p|
%\end{eqnarray}
%bis auf A<=C(1+p^2), der geht nicht
Finally,  Hölder-continuity of $A_\delta'(p)$ with respect to $\delta$ follows also with the formulation \eqref{AAtilde} and the Lipschitz continuity of $\tilde A'$.
Summarized we can state

\begin{lemma}
  \label{proof_regA}
 The mappings $A_\delta: \mathbb{R}^d \to \mathbb{R}$ for $\delta \geq 0$ (with $A_{\delta = 0} \coloneqq A$ as shorthand notation) have the following properties:
  \begin{itemize}
  \item[a)]  $A_\delta$ fulfill the growth condition
    $ c|p|^2 \leq A_\delta(p) \leq C_\delta+C|p|^2$ for all $p$ with positive constants $c,C,C_\delta$, where only $C_\delta$ may depend on $\delta$. 
  \item[b)] $A_\delta'$ are Lipschitz-continuous and strongly monotone on  $\mathbb{R}^{d}$
 with constants independent of $\delta$
    and $A_\delta'(0)=0$.
 \item[c)] $A_\delta''$ induce uniformly equivalent norms on $\mathbb{R}^{d}$ for $\delta>0$, i.e. there exist constants $c_0, C$ such that
   $$
c_0 \| q\|^2\leq	q^TA_\delta''(p)q \leq C \| q\|^2 \qquad \forall p,q \in \mathbb{R}^d, \delta> 0
$$
and $A_\delta''(0)= L\sum_{l=1}^LG_l$.
Furthermore, if $\delta=0$ the same holds true 
%and this inequality holds also
for all $p\neq 0$.
% \item[c)] The assumptions in Theorem \ref{existencecontrol} hold, where the upper bound  in the growth condition on $A$ may depend on $\delta$.
\item[d)]  
$A'_{(.)}(p)$ is Hölder-continuous with exponent $1/2$ and with a constant independent of $p$. Especially it holds
	\begin{align}\label{Akonv}
          |A_\delta{'}(p) - A{'}(p)|  \leq C \sqrt{\delta}
\qquad \forall p \in \mathbb{R}^d, \delta> 0
          .
	\end{align}
      \end{itemize}
In particular the \cref{assu} are fulfilled if $\delta >0$, \cref{assu}a. hold for $A$
and the convergence assumption with respect to $\delta \to 0$ in \cref{th:conv_op_delta} hold.

%In particular the assumptions on $A'$ in \cite[Assumptions 1.1]{BlMe21} are fulfilled.
%In particular \cref{assu}
%for the existence of an optimal control are fulfilled for the original nonsmooth and for the regularized control problems.
\end{lemma}

\section{Numerical results}
\label{sec:numerical_resultsPart1}
\anm{trust-region auf trust region vereinheitlicht; $A$ zu $A_\delta$, etc., geändert\\}%
In the last part of this paper we report some numerical findings.
\change{When we consider fixed delta in this chapter, we drop the corresponding index in the relevant quantities to keep the presentation lucid.}
We also consistently use the smooth double-well potential $\psi(s) = \tfrac{1}{4}(1-s^2)^2$ which defines $C_\psi=1$.
Our numerical approach for solving the regularized optimization problem is to first discretize in time and
then apply an optimization algorithm on this semi-discretized problem. {The arising equations and first order condition have rigorously been analyzed in the previous chapters.} Finally, each step in the algorithm is discretized also in space where we use global continuous, piecewise linear finite element approximations.

Preliminary numerical results have been obtained by a line search method based on the gradient $\nabla j_\tau$. However,
since we have not seen any relevant differences in the computed controls and states
we only present here results using second order informations---which are formally derived in the following---to gain efficiency in the solver.
\change{As algorithm to solve for a local minimizer we apply the trust region Newton method \cite{conn2000trust} to the reduced cost functional $j_\tau:U_\tau\to \mathbb{R}$.
  It is a common globalization of Newton's method which is needed
  due to
  % the nonlinearity of the state equation and hence
  the non-convexity of the problem.
The main idea is to determine at each iterate $u_\tau$ an approximate solution $\delta \bar u_\tau $ of the quadratic subproblem
\begin{equation}
\label{eq:tr_subprob_paper2}
\min_{\norm{\delta u_\tau} \leq \sigma} (\nabla j_\tau(u_\tau), \delta u_\tau) + \tfrac{1}{2}(\nabla^2 j_\tau(u_\tau)\delta u_\tau, \delta u_\tau),
\end{equation}
where $\sigma >0$ parameterizes the size of a trust region in $U_\tau$ where this model is considered to be sufficiently valid.
The proper choice of $\sigma$ is controlled by the trust region method.
The solution to \cref{eq:tr_subprob_paper2} is determined by the Steihaug-CG method \cite{Steihaug83}, which iteratively applies the CG-steps to the the first order condition of the unconstrained version of the subproblem
% solves the first order condition
$\nabla^2 j_\tau(u_\tau)\delta u_\tau= - \nabla j_\tau(u_\tau)$, and additionally handles the cases when a CG-iterate $\delta \tilde u_\tau$ exceeds the trust region boundary or
%$\delta \tilde u_\tau$ generates a nonpositive curvature
\changetwo{$j$ has at $u_\tau$ nonpositive curvature in direction $\delta u_\tau$},
i.e. $(\nabla^2 j_\tau(u_\tau)\delta \tilde u_\tau, \delta \tilde u_\tau)\leq 0$.
%by the same steps as the CG method. In addition it suitably handles the cases when an iterate exceeds the trust region boundary or has a negative energy `norm' with respect to the Hessian.
Finally the new trust region iterate is set to $u_\tau + \delta \bar u_\tau $.}\\
%
%As algorithm
%we apply the trust region Newton method \cite{conn2000trust} to the reduced cost functional $j:U_\tau\to \mathbb{R}$.
%The need for a globalization arises due to the nonlinearity of the state equation.
%The approximate solution $\delta \bar u_\tau $
%of the quadratic subproblem at the current iterate $u_\tau$ 
%$$ \min_{\norm{\delta u_\tau} \leq \sigma} (\nabla j(u_\tau), \delta u_\tau) + \tfrac{1}{2}(\nabla^2 j(u_\tau)\delta u_\tau, \delta u_\tau)$$
%is determined by the Steihaug-CG method \cite{Steihaug83} and the new iterate is set to $u_\tau + \delta \bar u_\tau $.
%
% A pseudocode sketch of both algorithms is given in \cref{algTRN} and \cref{algSCG}.
% In the last section we have shown differentiability of $j$. Let us summarize:
Let us summarize for convenience the \change{formulas} of the last sections:
\change{\begin{align}
	{j_\tau}({u_\tau}) &= \dfrac{1}{2}\|y_\tau(T)-y_\Omega\|^2 + \dfrac{\lambda}{2\varepsilon}\|u_{\tau}\|^2\qquad \\
  \nabla {j_\tau}(u_\tau) &= \tfrac{1}{\varepsilon}({\lambda} {u_\tau} + {p_\tau}),
\end{align}}%
where $y_\tau=S_\tau(u_\tau)$ is given by the state equation \cref{scheme_state1} and  $p_\tau=(S'_\tau(u_\tau))^*(y_\tau(T)-y_\Omega)$ is given by the time discrete adjoint equation \cref{eq:disc_adjoint}.\\
The derivative of $y_\tau$ in direction of $\delta u_\tau$ we denote by $\delta y_\tau$, i.e. $\delta y_\tau$ solves the linearized state equation \cref{diff}
with $v=\delta u_\tau$. This derivative is employed to compute the Hessian as can be seen below.\\
The second derivative \change{required only here for the Newton approach (see subproblem \eqref{eq:tr_subprob_paper2})}, in particular its action on an $L^2$-function, we deduce formally.
%Altogether the second derivative is given by
We obtain
\change{\begin{align}
  \nabla^2 {j_\tau}({u_\tau}) {\delta u_\tau} &=
                                           \tfrac{1}{\varepsilon}(\lambda {\delta u_\tau} + \delta {p_\tau}).
\end{align}}
Here, for the given  solution $p_\tau(u_\tau)$ of   \cref{eq:disc_adjoint}
%Given $p_\tau(u_\tau)$ is the solution of   \cref{eq:disc_adjoint}
the derivative
\begin{align*}
  \delta p_\tau := \frac{d p_\tau(u_\tau)}{du_\tau} \delta u_\tau
\end{align*}
fulfills
the so called additional adjoint equation.
In time discrete form it is given by
$\delta p_{N+1}:= \delta y_{N}$ and
\begin{equation}
\label{addadjdisc}
\begin{aligned}
&\varepsilon(A_\delta''(\nabla y_{j})\nabla \varphi,\nabla \delta p_{j}) + (\tfrac{1}{\varepsilon}\psi''(y_{j})\varphi + \tfrac{\varepsilon}{\tau_j}\varphi,\delta p_{j}) =
(\varphi, \tfrac{\varepsilon}{\tau_j} \delta p_{j+1})
\\
& \qquad\qquad -\varepsilon(A_\delta'''(\nabla y_{j})[\nabla \varphi, \nabla \delta y_{j}], \nabla p_j)
- \tfrac{1}{\varepsilon}(\psi'''(y_{j})\varphi \delta y_{j}, p_{j})  \quad
\text{ for } j=1,\ldots,N.
\end{aligned}
    \end{equation}
As for  \cref{eq:disc_adjoint} the unique existence of the solution is guaranteed.
The time-continuous counterpart is %respectively
given by
\begin{equation}
\label{eq:cont_add_adjoint}
\begin{aligned}
&  - \varepsilon (\eta,\partial_t \delta p) + \varepsilon (A_\delta''(\nabla y)\nabla \eta, \nabla \delta p) + \dfrac{1}{\varepsilon}(\psi''(y)\eta,\delta p) = \\
& \qquad\qquad - \varepsilon (A_\delta'''(\nabla y)[\nabla \eta,\nabla \delta y], \nabla p) - 
  \dfrac{1}{\varepsilon}(\psi'''(y) \eta \delta y , p)
   \quad &\forall \eta \in L^2(0, T;H^1(\Omega))\\
&  \quad	\delta p(T) = \delta y(T) \quad \text{in } \Omega   
\end{aligned}
    \end{equation}
and they are related by the discontinuous Galerkin time discretization as for the adjoint equation.
Note that each of the discretized equations has a time continuous counterpart (cf. \cref{weakAC,eq:cont_adjoint,eq:lin_equation_cont_with_z,eq:cont_add_adjoint}).
In fact, \cref{eq:cont_adjoint,eq:lin_equation_cont_with_z,eq:cont_add_adjoint} are the equations you would expect to get as the adjoint and corresponding linearized equations
%\luise{auslassen:? from doing a formal treatment of the analysis}
for \cref{weakAC}. Thus at least on a formal level it holds that the approaches \textit{first discretize then optimize} and \textit{first optimize then discretize} commute for the implicit time discretization \change{in the sense of \cite[chapters 3.2.2 and 3.2.3]{hinze2008optimization}}. Also discretization and optimization are interchangeable for the spatial discretization if one chooses the same ansatz spaces for $y$ and $p$.
Consequently, one expects to obtain for the optimization solver iteration numbers independent of the discretization level. This is strengthened by numerical observations in \cref{sec:nummesh}.
\kommentar{From a more practical perspective (sagt man das so? oder besser another)---and also strengthened by numerical observations---this should yield the problem to be independent of the discretization coarseness to a certain level.}%
Note that in general the \mbox{3-tensor} $A_\delta'''$
appearing in \cref{eq:cont_add_adjoint}
is not \change{symmetric}% although $A''$ is
, so we have to keep care of order in the corresponding term.
In the implementation $A_\delta'''$ is determined by automatic differentiation---stating an explicit formula does give no new insight.
Moreover, since $\nabla y=0$ in the pure phases one is particularly interested in the behavior of $A_\delta'''(p)$ when $p\to 0$. However,
$\tilde A'''(\sqrt\delta ({p \atop 1})) $ behaves like $1/\sqrt\delta$ due to the 2-homogeneity 
of $\tilde A$ (see \cref{AAtilde}).
Hence $\lim_{p\to 0} A'''_\delta(p)$ cannot be bounded independently of $\delta$.
Numerically we see this  problem for values 
$\delta < 10^{-13}$.%1.e-12: konv. nach 5 Schritten, 1.e-13 konv. nach 11 Schritten, darunter hab ich immer bei 10 oder so abgebrochen, da man dort erst bei Residuum 1.e-6 oder so war und es nicht superlinear aussah
% bei isotrop kann man bis ~ 200 runtergehen und es geht one Probleme durch (Fehler erst bei 210)
\kommentar{when
the optimization algorithm does not converge due to inaccuracies in the adjoint and addititional adjoint equation.}%

\kommentar{Note the appearance of the third derivative of $A$ in \cref{eq:cont_add_adjoint}. Together with the second derivative already arising in the adjoint and linearized state equation this leads to some numerical restriction on the so far arbitrarily small parameter $\delta$.
If we would not regularize, then while $A''(p)$ was bounded independently of $p$ {for $p\neq 0$}, this would not be the case for $A'''$. This follows from the fact that $A'''(\lambda p) $ behaves like $1/\lambda$ due to the 2-homogeneity of the unregularized $A$.
Hence the third derivative clearly diverges in the limit $|p|\to 0$, but already the second derivative possesses issues since for small $p$ it contains fractions of small numbers.
Hence one has to choose $\delta$ appropriately to circumvent these problems.
We further note that in general the 3-tensor $A'''$ is not symmetric% although $A''$ is
, so we have to keep care of order in the corresponding term.
In the implementation $A'''$ is determined by automatic differentiation???---stating an explicit formula does give no new insight.
%We will not give the expression for it explicitly, since it gives no new insights. In the implementation the derivative was automatically determined by the tools of FEniCS.
}%

%If one uses a Newton type solver as we do this leads to problems when considering the limit $\delta \to 0$ numerically as the solver fails in evaluating a linearization of the adjoint equation if $\delta$ is to small.
%However one typically takes one small but fixed value for $\delta$ there the solver still succeeds.
\bigskip

To keep the computational cost moderate we set $d=2$ in all experiments.
Furthermore,
throughout this section we use
as a spatial domain the square $\Omega = (-1,1)^2$
and as the time horizon $[0,T]$ with $T=1.625\cdot 10^{-2}$,
we set the parameters $\varepsilon = \tfrac{1}{14\pi}$ and $\lambda = 0.01$ and---if not mentioned otherwise---the regularization parameter $\delta=10^{-7}$.
We choose the constant time step size $\tau = 1.625\cdot 10^{-4}$, which fulfills the condition
$\tau \leq \varepsilon^2/C_\psi $ %= 5.16 \cdot 10^{-4}$
and $\Omega$ is 
uniformly discretized with $129 \times 129$ grid points.
  % $2.4\cdot10^{-4}\approx (2/129)^2 = \Delta x^2
\kommentar{It holds $h^2 \geq \tau$ for the spatial mesh size $h$. \johannes{ist das wirklich die genaue Bedingung?, weil z.B. bei heat equation hing es ja noch von dem Koeffizienten von Laplace ab\\}}%
Moreover by numerical evidence we know that the interface is resolved sufficiently with
6-14 mesh points accross the interface.
Each computation started with $u^{(0)} \equiv 0$.

We have looked at various set-ups that mainly vary by $y_0$, $y_\Omega$ and
the final time $T$.
For the anisotropies determined by $G_l$ we used three different choices, that are listed in the following.
\kommentar{
	Eigentlich habe ich
	$$
	A'(p) = \tfrac{1}{L}\gamma(p){\gamma}'(p)=\dfrac{1}{L}\sum_{l,m} \dfrac{\gamma_m(p)}{\gamma_l(p)} G_l p
	$$
	in diesem Kapitel benutzt. Wie soll man das einführen? Schon in anderen Kapiteln darauf hinweisen? Oder `rescaling of the $G_l$'?\\
	Mögliche Begründung $B_L(q)\coloneqq \dfrac{1}{L}\sum_{l,m} \dfrac{\gamma_m(p)}{\gamma_l(p)} G_l \sim C\sum_{l} G_l$ for $|p|\to0$. ($=id$ für reg$l1$norm und $diag(1.5,1.5)$ für hexa; $C$ variiert je nach $p$) \\
}

\begin{figure}[htbp]
	\includegraphics[scale=.3]{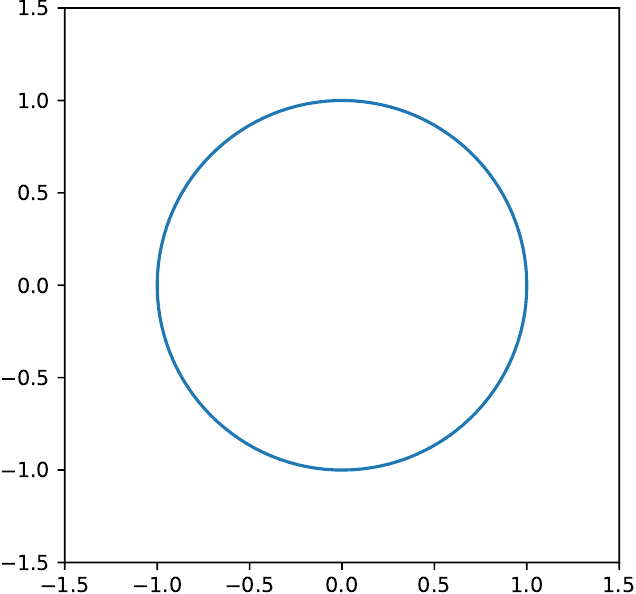}~~~~
	\includegraphics[scale=.3]{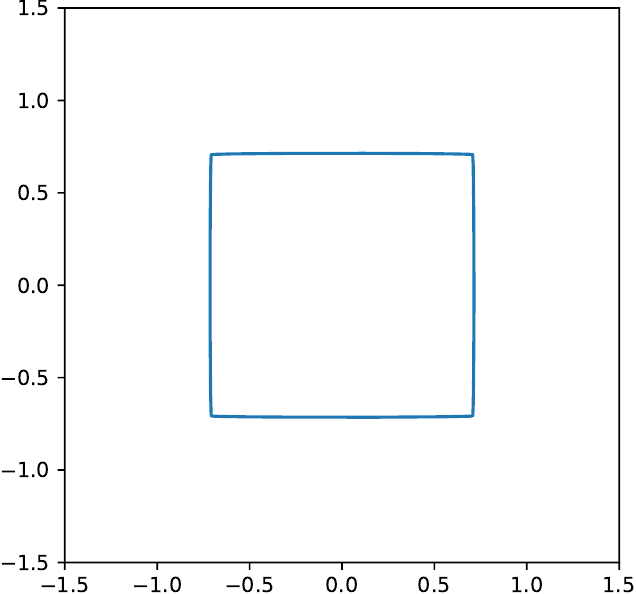}~~~~
	\includegraphics[scale=.3]{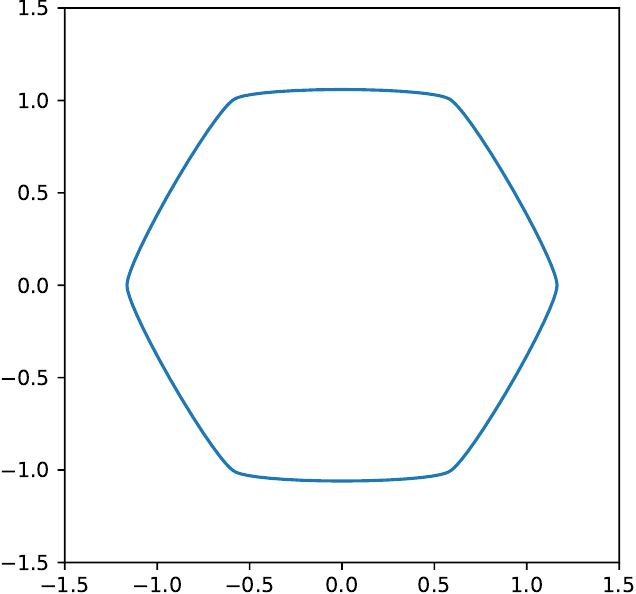}
	\caption{Wulff shapes of isotropic, $l_1$- and hexagon-anisotropy.
	}
	\label{fig:Wulff}
\end{figure}

\begin{enumerate}
	\item \uline{isotropic case:} $\gamma(p) = \norm{p}_2$
          this would belong to the choice
\kommentar{WEGLASSEN; die zweite überhaupt angeben?, weil die ist nur zu testzwecken ob alles richtig implementiert ist, wird aber nicht mehr erwähnt}
	\begin{equation}
	\label{eq:choice_isotropic}
	G = \left(\begin{matrix}
	1&0\\0&1
	\end{matrix}\right).
%	\qquad
%	\text{or}
%	\qquad 
%	%\begin{matrix}
%	G_1 = G_2 = 
%	\left(\begin{matrix}
%	0.25&0\\0&0.25
%	\end{matrix}\right),
%	%&&&&&& G_2 = 
%	%\left(\begin{matrix}
%	%0.25&0\\0&0.25
%	%\end{matrix}\right),
%	%\end{matrix}
	\end{equation}
%	but in the implementation we used a special version where terms were sufficiently simplified, such that we do not need to regularize.
\change{Note that in this case regularization is not necessary.}
%	Note that in this case one does not need to apply regularization.
\change{In addition, the regularization would cancel out as can be seen in \cref{derAreg} and as it is also discussed in the corresponding text thereunder.}	

\item \uline{regularized $l_1$-norm:}
	\begin{equation}
	\label{eq:choice_aniso1}
	\begin{matrix} G_1 = \dfrac{1}{2}
	\left(\begin{matrix}
	1&0\\0&\epsilon
	\end{matrix}\right),
	&&&&&& G_2 = \dfrac{1}{2}
	\left(\begin{matrix}
	\epsilon&0\\0&1
	\end{matrix}\right),
	\end{matrix}
	\end{equation}
	with some small parameter $\epsilon$ that we set to $\epsilon =0.01$ (not to be confused with the interface parameter $\varepsilon$). For $\epsilon = 0$ this reduces to \changetwo{$\gamma(p)=\tfrac{1}{\sqrt2}\|p\|_{\ell_1}$}.% Its Wulff shape is depicted in ???.	

      \item \uline{form of a smoothed hexagon:}
		\begin{equation}
		\label{eq:choice_aniso2}
	\begin{matrix} G_l = \dfrac{1}{3}
	\left(\begin{matrix}
	\cos(\alpha_l)&-\sin(\alpha_l)\\\sin(\alpha_l)&\cos(\alpha_l)
	\end{matrix}\right)
	\left(\begin{matrix}
	1&0\\0&\epsilon
	\end{matrix}\right)
	\left(\begin{matrix}
	\cos(\alpha_l)&\sin(\alpha_l)\\-\sin(\alpha_l)&\cos(\alpha_l)
	\end{matrix}\right),
	\end{matrix}
%	=								% this only holds if $\epsilon=0$, so just for checking
%	\left(\begin{matrix}
%	\cos(\alpha_l)^2& \cos(\alpha_l)\sin(\alpha_l)\\
%	\sin(\alpha_l)\cos(\alpha_l) & \sin(\alpha_l)^2
%	\end{matrix}\right)
	\end{equation}
	where $\alpha_l = \tfrac{\pi}{3}l$, $l=1,2,3$ and $\epsilon$ as before.% The Wulff shape is depicted in ???.
\end{enumerate}
{Note that in contrast to the choices in \cite{Barrett2008} we divide the matrices by their total number $L$. By this scaling the costs between the different anisotropies becomes more comparable, since by numerical observation the velocity of the shrinkage is approximately equal. This can also be seen on the Wulff shapes, which are defined as
  % $z_W(z) \coloneqq \gamma'(z)$ for $z\in S^{d-1}$, see \cite{Gurtin, M.E.:Thermodynamics ...}
  $\mathcal{W} \coloneqq \left\{q \in \mathbb{R}^d: \sup_{p\in\mathbb{R}^d\setminus\{0\}} \tfrac{p^Tq}{\gamma(p)}\leq1 \right\}$, see \cite{WulffXXVZF},
  and which are visualized for above choices of $\gamma$
  %a certain visualization of the anisotropies given
  in \cref{fig:Wulff}.
  Without the rescaling the Wulff-shape
of the hexagon anisotropy would extend approximately to the label $2.0$ on each axis.}

%The following sections each treat one of the above choices for the anisotropy respectively.
%We will always show one example that takes benefit of the considered anisotropy by complicated scenarios that still can be handled with it. We will also show two examples of topology change, one time by expanding a phase to the whole domain and one time by merging two \johannes{phase clusters???}.
As computing framework we used FEniCS \cite{AlnaesBlechta2015a} or rather its C++ interface DOLFIN~\cite{LoggWells2010a}.
The simulations were carried out on an HP EliteDesk 800 G4 workstation containing an Intel Core i7-8700 CPU with 12 cores à 3.20GHz and 16 GB of RAM.

In the following subsections we support the convergence result of \cref{th:conv_state_disc_wrt_delta_full_system} concerning the regularization, we give numerical evidence for mesh independent behavior in the solution process and we present optimal control results for different  anisotropies and different desired states, including star like objects and necessary topology changes.

\kommentar{
\luise{ab hier noch neu, bzw. zu überarbeiten}
In the following \luise{subsections} we look at several examples and compare the different anisotropies among them.
First we investigate how the growing of finger like structures compare between the isotropic and anisotropic cases.
Then we present examples with  topology changes, namely
splitting a mass accumulation into two or merging them the other way round.
This is a situation for which the phase field method is better suited than other methods for moving boundary problems.
%Finally we will look at an example where we split a mass accumulation into two or merge them the other way round.
Finally we give numerical evidence of mesh independent iteration numbers and considerations concerning the regularization parameter $\delta$.
\luise{$\uparrow$DAS ist noch zu verbessern, einzufügen...$\uparrow$}
}%ENDE Kommentar
      
%\newpage
% !TeX spellcheck = en_US

\subsection{\changemy{Dependence on} the regularization parameter $\delta$}%\change{\sout{Behavior of the state equation in }}
\begin{figure}[h]
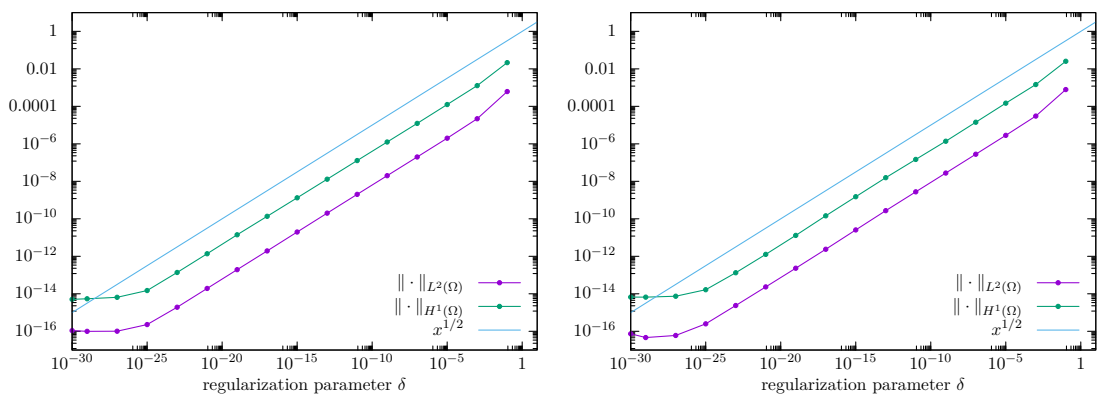
%\hspace*{-1.5cm}
  \begin{subfigure}{0.5\textwidth}
    \begin{center}
      \input{delta_dep_l1}
      %\begin{tikzpicture}[gnuplot, scale=0.6, every node/.style={scale=0.6}]
	\label{fig:delta_dep_l1}
	\caption{Results for the regularized $l_1$-norm.}
      \end{center}
    \end{subfigure}%\hspace{0.1\textwidth}
\begin{subfigure}{0.5\textwidth}
    \begin{center}
	\input{delta_dep_hexa}
	%\begin{tikzpicture}[gnuplot, scale=0.6, every node/.style={scale=0.6}]
	\label{fig:delta_dep_hexa}
	\caption{Results for the hexagon anisotropy.}
\end{center}
\end{subfigure}
\caption{Comparison of \change{$\|y^\delta_\tau(T)-y_\tau(T)\|$} in the $L^2(\Omega)$ and $H^1(\Omega)$-norms for different values of $\delta$.
\label{fig:delta_dep}}
\end{figure}

%\change{In the beginning, we discuss how the dependence on the regularization parameter $\delta$ enters in the simulations. Therefore, we first analyze the solution of the state equation.}
% Here we \change{analyze} numerically the dependence of the solution of the state equation on the parameter $\delta$.
\change{First  we analyze numerically
  the dependence of the solution of the state equation on the parameter $\delta$.}
As a setting we start from a circle of radius $0.5$ and look at the evolution of the state only
using $u=0$.
% until $T=1.625\cdot 10^{-2}$. 
\kommentar{Here $T$ and the remaining parameters remain the same as above.}%
A plot containing the results for the choices of both anisotropies respectively is given in
Figure \ref{fig:delta_dep}.
\kommentar{
  Here we look at the difference only at the end time point $T$. Since errors accumulate during the time evolution this should be a good metric for comparison.
  We have plotted both the difference in $L^2(\Omega)$ as well as $H^1(\Omega)$, where the latter is of course always bigger.
 }%
We have plotted both the difference of the states
in $L^2(\Omega)$ as well as $H^1(\Omega)$ at the end point $T$.
Since errors accumulate during the time evolution the errors at $T$ should be a good metric for comparison.
With the additionally plotted function $f(x)=x^{1/2}$
the figure clearly exhibits the convergence order $1/2$ which is expected according to \cref{eq:dependence_of_y_on_tau_whole_system,Akonv},
\change{i.e. it equals the approximation order of $A'_\delta$ to $A'$.}
\kommentar{
For $\delta$ larger than about $10^{-25}$ the errors are given by parallel lines however.
As a comparison we show the line the function $f(x)=x^{1/2}$ which on the log-log-plot is given by a line with slope $\tfrac{1}{2}$.
This is what one would expect from \cref{ass1,Akonv}.
}
% Fitting these values for the $H^1(\Omega)$-norm to a function of the form $f(x) = a+bx^c$ gives a value of $c\approx 0.63$ which lies close to the $0.5$ one would expect from \cref{ass1,Akonv}. % Fehler ist bloß 0.1355% also nicht erwähnenswert
\kommentar{One can observe deviations from the straight line for values below $10^{-25}$.
Here errors arising from computational inaccuracies dominate.
%For very low values of $\delta$ (like $10^{-100}$) one can reach errors of $4.0\times10^{-14}$ and $1.4\times10^{-12}$ for the $L^2(\Omega)$ and $H^1(\Omega)$ cases, respectively.
However the optimization algorithm does not converge starting from about $\delta=10^{-16}$ \johannes{genauen Wert nochmal überprüfen}, since then inaccuracies in the adjoint and addititional adjoint equation dominate.}%

\change{
When considering the numerical solution of the optimization problem,
% Furthermore, also with regard to the next subsection,let us mention that
%our experience shows at most
according to our experience there is only
weak dependence of the number of the trust region as well as of the Steihaug-CG iterations when varying $\delta$.
They stay nearly the same as in
\cref{dep1,tab:dep2} and are therefore not listed here.
If $\delta$ is such small that rounding errors accumulate for $A''_\delta$ %(and hence for the solution of the adjoint equation)
and even more for $A'''_\delta$ (and consequently for the solutions of the respective equations) %(consequently for the solution of the linearized state and additional adjoint equation)
the algorithm may not converge. %break down and the results may deviate due to that.
However, for $\delta\geq 10^{-10}$ the algorithm was always robust.}

%\newpage
% !TeX spellcheck = en_US
\subsection{Mesh independent behavior}
\label{sec:nummesh}
In this section we numerically investigate the mesh dependence of the problem solver.
More concretely we look at the number
of trust region iterations, called {\em TR steps} in the following tables, as well as
at the number of Steihaug-CG steps that are needed to solve the quadratic subproblems.
Since Steihaug algorithm consists of early stopping criteria given by the trust region algorithm their amount might change drastically during the progress of the algorithm.
%Also, the total amount %for one Steihaug-CG step also
%depends strongly on the initial value.
Therefore %only looking at the average amount is no good metric to measure performance of the method---
we rather look at the average amount of steps, called {\em mean CG} in the following tables, {that are needed} to decrease the residuum by $6$ orders of magnitude for trust region steps where this kind of measurement is possible. In addition we also take the maximum amount, called {\em max CG}  as an indicator. These numbers of CG iterations reflect on the conditioning of the linear systems corresponding to the quadratic subproblems.
%Since we have to vary the mesh size for both the spatial and time discretization there are some deviating choices for the parameters.
As final time we choose in these experiments $T=2\times10^{-3}$ \kommentar{erwähnen, dass $T$ leicht anders bei z.b. $\tau=1.e-4.5$ da sonst nicht teilt? The deviating choice for $T$ from before comes from the fact we would like to vary it in the following.
  }
The remaining parameters are left unchanged.
%AUSSER:
%The absolute tolerance of the Steihaug-CG method was set to $10^{-11}$.
We inspect the dependence on the space discretization by fixing $\tau=10^{-4}$ and varying $h=2/N$.
For analyzing the dependence on the step size $\tau$ we fix the spatial mesh size using $N=128$.
%Since we look at a procedure that is nested in the outer trust region algorithm its behavior might vary due to several reasons.
%The Steihaug algorithm consists of other stopping criteria than the tolerance and the starting value will improve over time. As a measure we give the maximum number of iterations among the whole algorithm as well as the mean number of iterations it took for improving the residuum by $6$ orders of magnitude for trust region steps where this kind of measurement is possible.

As model problem for the isotropic case, we consider the control of a circle from radius $r=0.5$ to $r=0.55$, where the results can be found in
\cref{dep1}.
%in time $T=2\times10^{-3}$ \johannes{erwähnen, dass $T$ leicht anders bei z.b. $\tau=1.e-4.5$ da sonst nicht teilt?}.
{One cannot observe a clear tendency that would suggest dependency
of the maximal or mean number of CG iterations and of the trust region steps on the granularity,
as is expected by the discussion in the introduction of this section.}
Only the amount of total computing time increases with the number of unknowns. %which is owed to the growing computational costs.
Let us mention that in case of $\tau=10^{-4}$ and $N=516$ the reduced optimization problem has around $5.4$ million unknowns given by the amount of discretization points of $u$. If $\tau=10^{-6}$ and $N=128$ the number of unknowns is roughly
$33.3$ million. %33282000
Due to the parallelization of the algorithm and the non-commutativity of floating point operations
the results might vary slightly among runs sharing the same configuration.
\begin{table}[h]
\centerline{\begin{minipage}{1.05\textwidth}
		\begin{tabular}{lrrrr}
		\toprule
		{N} &    64 &   128 &   256 &   512 \\
		\midrule
		max CG     &    38 &    48 &    38 &    39 \\
		mean CG    &  21.2 &  22.7 &  18.8 &  20.4 \\
		TR steps &    12 &    16 &    11 &    12 \\
		time (s) &    17 &    72 &   235 &  1196 \\
		\bottomrule
	\end{tabular}~~~
%	\begin{tabular}{llll}
%		\toprule
%		{$\tau$} & $10^{-4}$ & $10^{-5}$ & $10^{-6}$ \\
%		\midrule
%		max      &    48 &    34 &    34 \\
%		mean     &  22.7 &  18.4 &  17.2 \\
%		TR steps &    16 &     9 &     5 \\
%		time (s) &    72 &   213 &  1275 \\
%		\bottomrule
%	\end{tabular}
	\begin{tabular}{lrrrrr}
		\toprule
		{$\tau$} & $10^{-4}$ & $10^{-4.5}$ & $10^{-5}$ & $10^{-5.5}$ & $10^{-6}$ \\
		\midrule
		max CG     &      48 &        60 &      34 &        34 &      34 \\
		mean  CG   &    22.7 &      22.2 &    18.4 &      18.0 &    18.5 \\
		TR steps &      16 &        11 &       9 &         8 &       8 \\
		time (s) &      72 &       105 &     213 &       706 &    2032 \\
		\bottomrule
	\end{tabular}
\caption{Dependence on $N$ and $\tau$ for the isotropic case. 	\label{dep1}}
\end{minipage}}
\end{table}

Next we do the same analysis for the anisotropic Allen-Cahn equation with the regularized $l_1$-norm. Here we choose $y_0$ and $y_\Omega$ the same, i.e. we try to keep a square constant.
 The outcomes are listed in \cref{tab:dep2}.
 Again, almost no dependence on the discretization parameters is observed. Numbers for the average Steihaug steps as well as for the total trust region steps rather seem to ameliorate for more accurate computations.
 %The values are comparable to the isotropic case---the table comparing different $\tau$ may show less tendencies. Only the amount of trust region steps is generally lower. But since we chose configurations adapted to the anisotropies, it is no surprise that they differ a bit.
\begin{table}[h]
\centerline{\begin{minipage}{1.05\textwidth}
	%	\begin{tabular}{lllll}
%		\toprule
%		{N} &    64 &   128 &   256 &   512 \\
%		\midrule
%		max      &    43 &    41 &    42 &    43 \\
%		mean     &  19.8 &  18.8 &  18.8 &  19.0 \\
%		TR steps &     6 &     5 &     5 &     3 \\
%		time (s) &     8 &    24 &  5295 &  5953 \\
%		\bottomrule
%	\end{tabular}
	\begin{tabular}{lrrrr}
		\toprule
		{N} &    64 &   128 &   256 &   512 \\
		\midrule
		max CG     &    60 &    40 &    40 &    39 \\
		mean  CG   &  30.0 &  22.0 &  21.3 &  21.0 \\
		TR steps &    10 &     6 &     6 &     6 \\
          time (s) &    24 &   66 % 38, habe die Zeit aus der zweiten Tabelle übernommen
                                    &   194 &  1193 \\
		\bottomrule
	\end{tabular}~~~
%	\begin{tabular}{llll}
%		\toprule
%		{$\tau$} & 1.e-4 & 1.e-5 & 1.e-6 \\
%		\midrule
%		max      &    41 &    53 &    54 \\
%		mean     &  18.8 &  20.5 &  20.5 \\
%		TR steps &     5 &     6 &     6 \\
%		time (s) &    24 &   130 &  5528 \\
%		\bottomrule
%	\end{tabular}
	\begin{tabular}{lrrrrr}
		\toprule
		{$\tau$} & $10^{-4}$ & $10^{-4.5}$ & $10^{-5}$ & $10^{-5.5}$ & $10^{-6}$ \\
		\midrule
		max CG     &    40 &      40 &    39 &      35 &    35 \\
		mean CG    &  22.0 &    21.8 &  24.0 &    20.0 &  20.7 \\
		TR steps &     6 &       6 &     7 &       5 &     5 \\
		time (s) &    66 &     161 &   537 &    1127 &  3306 \\
		\bottomrule
	\end{tabular}
%	\caption{dependence on $\tau$          for the anisotropic case 2.        \label{tab:dtdep2}}
\caption{Dependence on $N$ and $\tau$ for the regularized $l_1$-norm where $y_\Omega=y_0$.\label{tab:dep2}}
\end{minipage}}
\end{table}
Finally, we list in \cref{tab:star}
the results for the control of a circle to a star with four fingers as can be seen in \cref{fig:l1iso_star4}. Here more control is necessary and the solution process takes significantly more trust region steps. Also the number of CG iterations are increased {compared to \cref{tab:dep2}}. While this indicates the dependency on the control configuration the results concerning the CG-iterations still show a behavior independent of the discretization level. A slight increase in the number of trust region steps is present.

\begin{table}[h]
\centerline{\begin{minipage}{1.15\textwidth}
\begin{tabular}{lrrrr}
\toprule
{$N$} &     64 &    128 &    256 &    512 \\
\midrule
max CG     &    185 &    207 &    243 &    167 \\
mean  CG   &  129.0 &  127.7 &  139.7 &  117.0 \\
TR steps &     71 &     83 &    106 &    125 \\
time (s) &     95 &    641 &   4068 &  23657 \\
\bottomrule
\end{tabular}~~~
%   diese Kopie ist zur Erinnerung wie die Spalte 128 früher aussah (hab die gleiche genommen wie im tau plot nachträglich)
%\begin{tabular}{lrrrr}
%	\toprule
%	{$N$} &     64 &    128 &    256 &    512 \\
%	\midrule
%	max CG     &    185 &    186 &    243 &    167 \\
%	mean  CG   &  129.0 &  106.0 &  139.7 &  117.0 \\
%	TR steps &     71 &     87 &    106 &    125 \\
%	time (s) &     95 &    467 &   4068 &  23657 \\
%	\bottomrule
%\end{tabular}
%
%  	python3 analyze.py output_1.e-4/out.txt 1.e-4 output_1.e-4.5/out.txt 1.e-4.5 output_1.e-5_other/out.txt 1.e-5 output_1.e-5.5/out.txt 1.e-5.5 output_1.e-6/outs2/out_all.txt 1.e-6
	\begin{tabular}{lrrrrr}
		\toprule
		{$\tau$} & $10^{-4}$ & $10^{-4.5}$ & $10^{-5}$ & $10^{-5.5}$ & $10^{-6}$ \\
		\midrule
		max CG   &    207 &     180 &    175 &     181 &    175 \\
		mean CG  &  127.7 &   153.0 &  148.0 &   130.5 &  137.3 \\
		TR steps &     83 &      95 &    146 &     126 &    108 \\
		time (s) &    641 &    1462 &   5476 &   13753 &  46297 \\
		\bottomrule
	\end{tabular}
\caption{Dependence on $N$ and $\tau$ for the simulation circle to 4-star for the regularized $l_1$-norm.
  % in the regularized $l_1$ anisotropic setting.
\label{tab:star}}
\end{minipage}}
\end{table}

\kommentar{ die folgende Tabelle wurde für eingrößeres  $T$ berechnet
\begin{table}[h]
%\begin{adjustwidth}{-0cm}{}
%  \begin{center}
  \hspace*{-0.05\textwidth}
    \begin{minipage}{0.3\textwidth}
 % \begin{center}
	\begin{tabular}{llll}
		\toprule
		{N} &     64 &    128 &     256 \\
		\midrule
		max CG     &    607 &    783 &     801 \\
		mean CG    &  308.0 &  479.0 &   480.4 \\
		TR steps &    225 &    221 &     623 \\
		time (s) &   2371 &   8298 &  143147 \\
		\bottomrule
	\end{tabular}
		\caption{dependence on $h$
			for the simulation circle to star in the anisotropic setting.
		\label{tab:dhstar}}
\end{table}
}%END KOMMENTAR    

%\newpage
% !TeX spellcheck = en_US
\subsection{Numerical examples for different desired states and anisotropies}
Finally we present solution for three different objectives: the evolution to star-like structures, the splitting of geometries and the merging of geometries.
For the presented figures which show the evolution of the control $u$ we employed a scaling of the color that was adjusted to the values at $t\approx T/2$.
%corresponding to the obtained values at $t\approx T/2$. 
Hence this allows to see where in $\Omega$ the control is present although its values may be clipped on some images.
To see how much the system is controlled at which time, we in addition include figures showing the $L^2(\Omega)$-norm of the control over time. 
\subsubsection{Evolution to star-like structures}
In the first experiment we start from a circle of radius $0.5$ and try to steer it to a star-like structure with $4$ %(regularized $l1$-norm)
or $6$ %(hexagon)
fingers respectively.
The images of the time evolution of the corresponding states and controls 
can be found in \cref{fig:l1iso_star4,fig:isohex_star6}.
The state is given by the Allen-Cahn equation with
the regularized $l_1$-norm anisotropy for the `4-star' target, and with  the `hexagon' anisotropy for the `6-star' target. In addition we present the results in both cases for the isotropic evolution equation.
In practice the choice of Allen-Cahn equation is given by the model equation and not by the desired state.
%
%We compare the applied control costs to the isotropic case with same start and end configurations.
%The images of the time evolution can be found in
%Figure \ref{fig:l1_star4}-\ref{fig:iso_star6}.

Qualitatively the main observation is that
in all cases
the control takes place in a neighborhood of the interface.
Moreover,
the evolution is controlled essentially in the second half of the time interval.
\kommentar{Grund (Hinweis für mich): ich hab die Bilder so skaliert, dass sie für den Schritt 3 oder 4 die beste Farbe haben; auch unten ergänzen wenn so passt\\}%
This can be particularly seen in \cref{fig:star_iso_compare} where the $L^2(\Omega)$-norms
%for several time steps of the controls corresponding to the various anisotropies are plotted.
of the control are plotted over the time $t$.
On the $t$-axis we indicated the times at which the states and controls were sampled for \cref{fig:l1iso_star4,fig:isohex_star6}.
%At the beginning the mass evolves almost according to the Allen-Cahn equation with $u=0$.
%\luise{das ist in den Bildern nicht zu sehen}
%One can observe that the control cost gets bigger during the time evolution.
While for the isotropic case the middle part seems to grow nearly linearly, for the cases of the regularized $l_1$-norm and the hexagon anisotropy one observes bigger jumps intersected by approximately constant parts. The first plateau comes from the fact that at the first part the evolution follows the nearly uncontrolled Allen-Cahn flow to get a square-like, respectively a hexagon-like shape. Only then the control truly enters to initiate the development of the fingers
with the strongly non-convex parts.
From that point onwards more control is needed for the anisotropic cases than 
% the curve lies above the one
for the isotropic case but towards the end they approximately overlap.
The last peak arises from the fact that much of the control is spent to form the details of the fingers in the last few time steps.
% the main part of the control does not take part before half of the evolution is over.
% At the beginning the mass is left to take on a form that it would get if one would evolve the Allen-Cahn equation with $u=0$.
%Much of the control is also spent to form the details of the fingers in the last few time steps.
%Concludingly we can say: while the evolution to a (smoothed) square does not need much control with the $l_1$-norm it is not clear that the isotropic case is inferior when fingers shall build due to the strongly nonconvex parts. However, the choice of anisotropy is given by the model and not by the desired state. \johannes{??}

\ifgraphics
% für die Skalierung des control-plots sollte ich eigentlich immer den 4. Zeitschritt benutzt haben
\begin{figure}[htbp]\begin{center}
\begin{subfigure}{1.\textwidth}
\begin{center}    
	\includegraphics[trim={540px 112px 540px 112px}, clip, scale=.08]{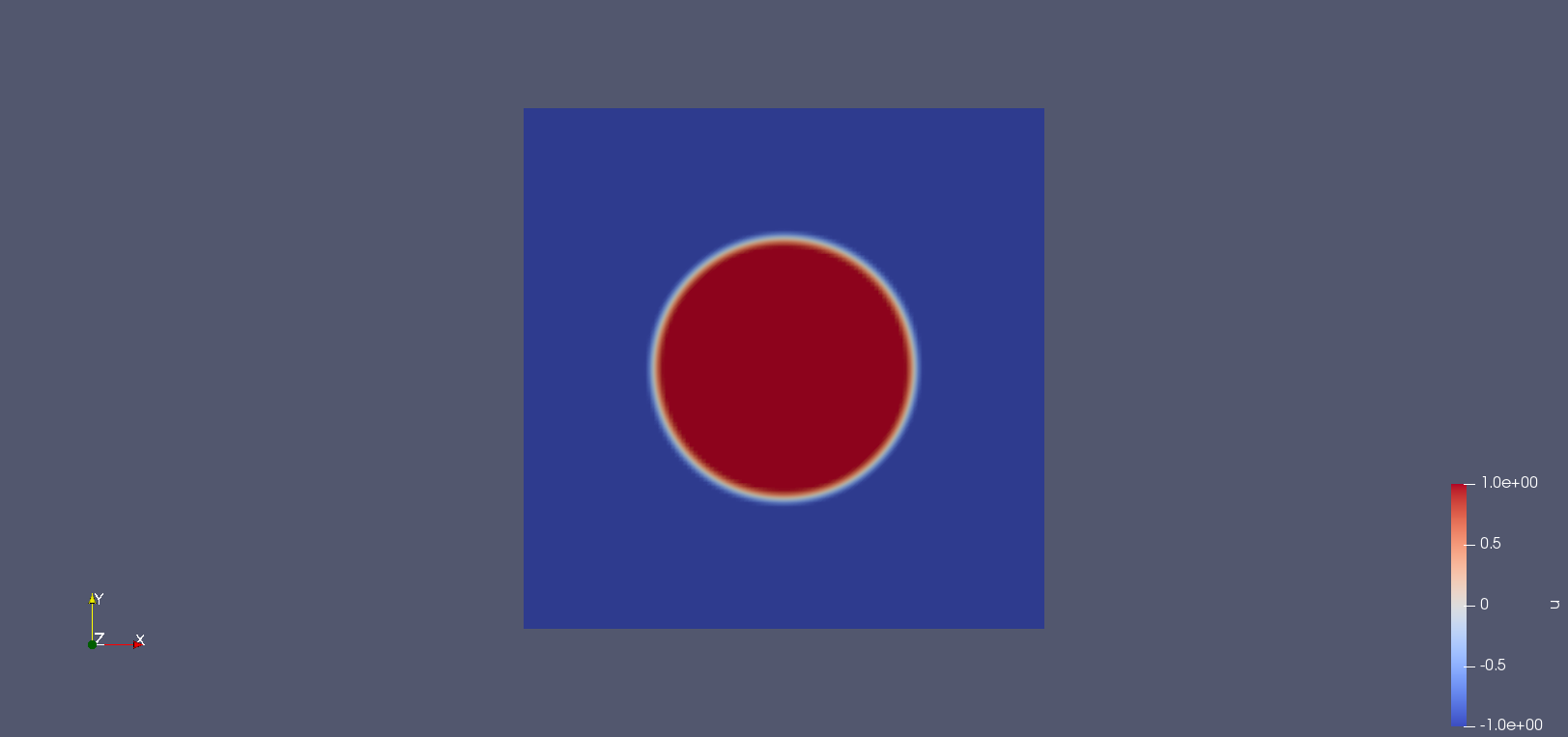}
	\includegraphics[trim={540px 112px 540px 112px}, clip, scale=.08]{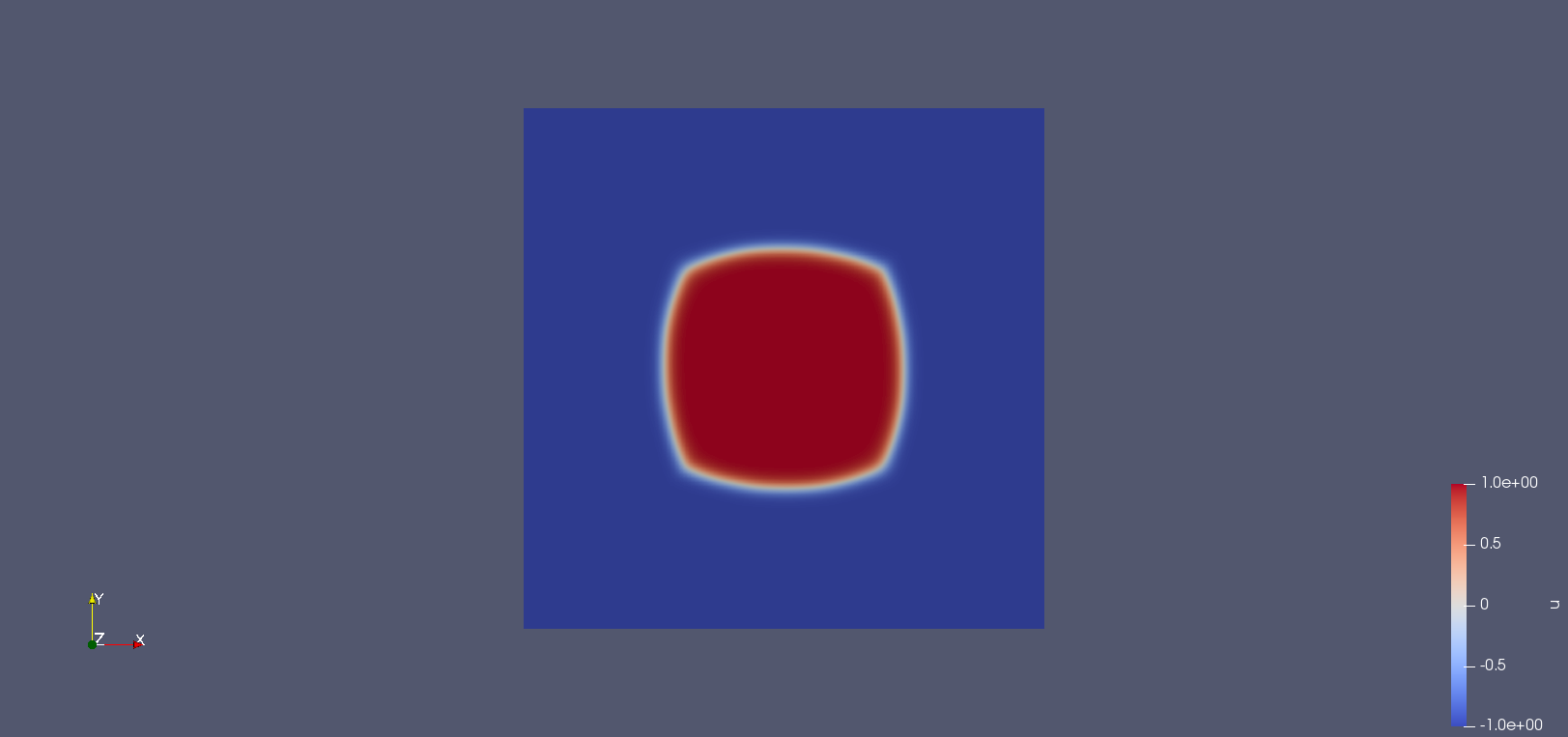}
	\includegraphics[trim={540px 112px 540px 112px}, clip, scale=.08]{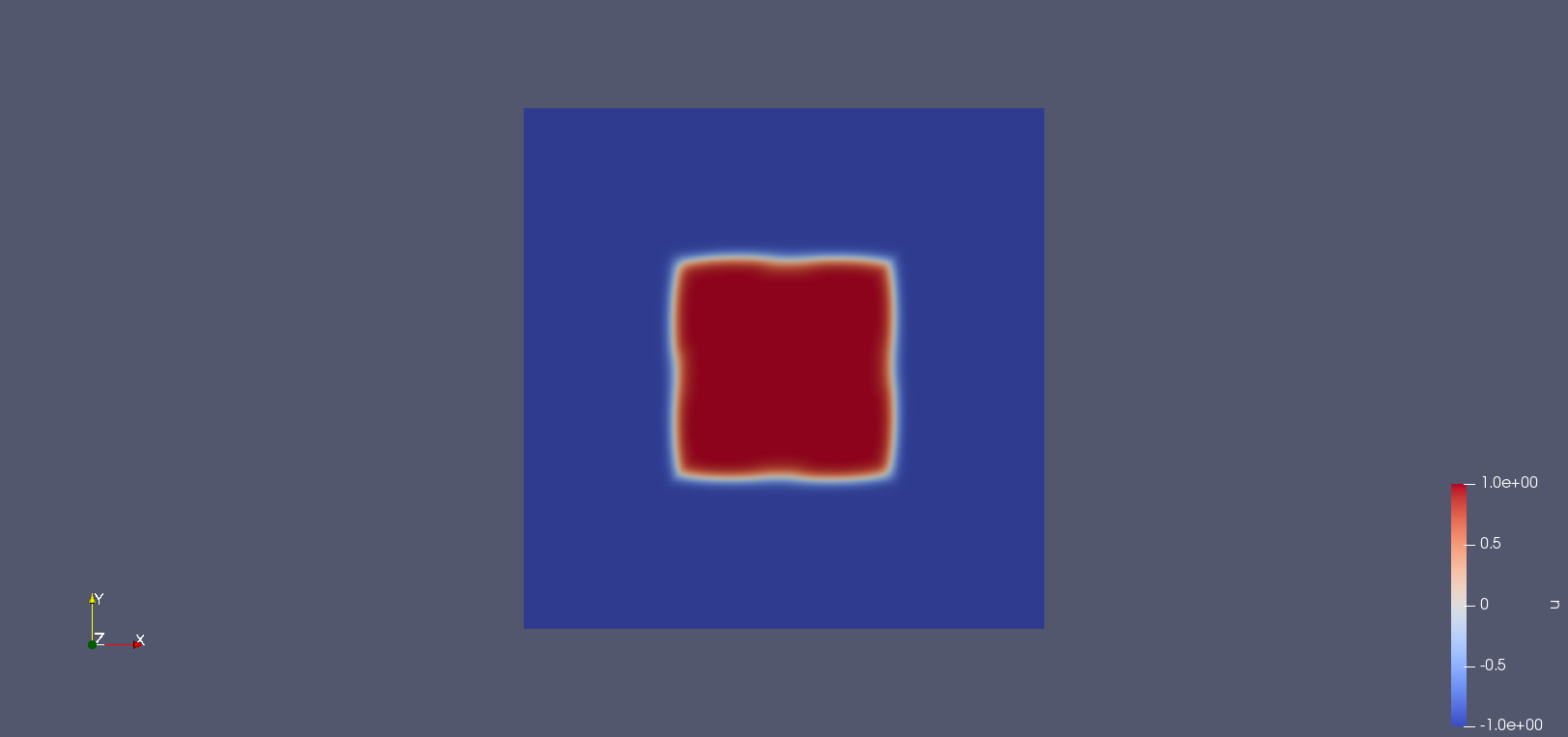}
	\includegraphics[trim={540px 112px 540px 112px}, clip, scale=.08]{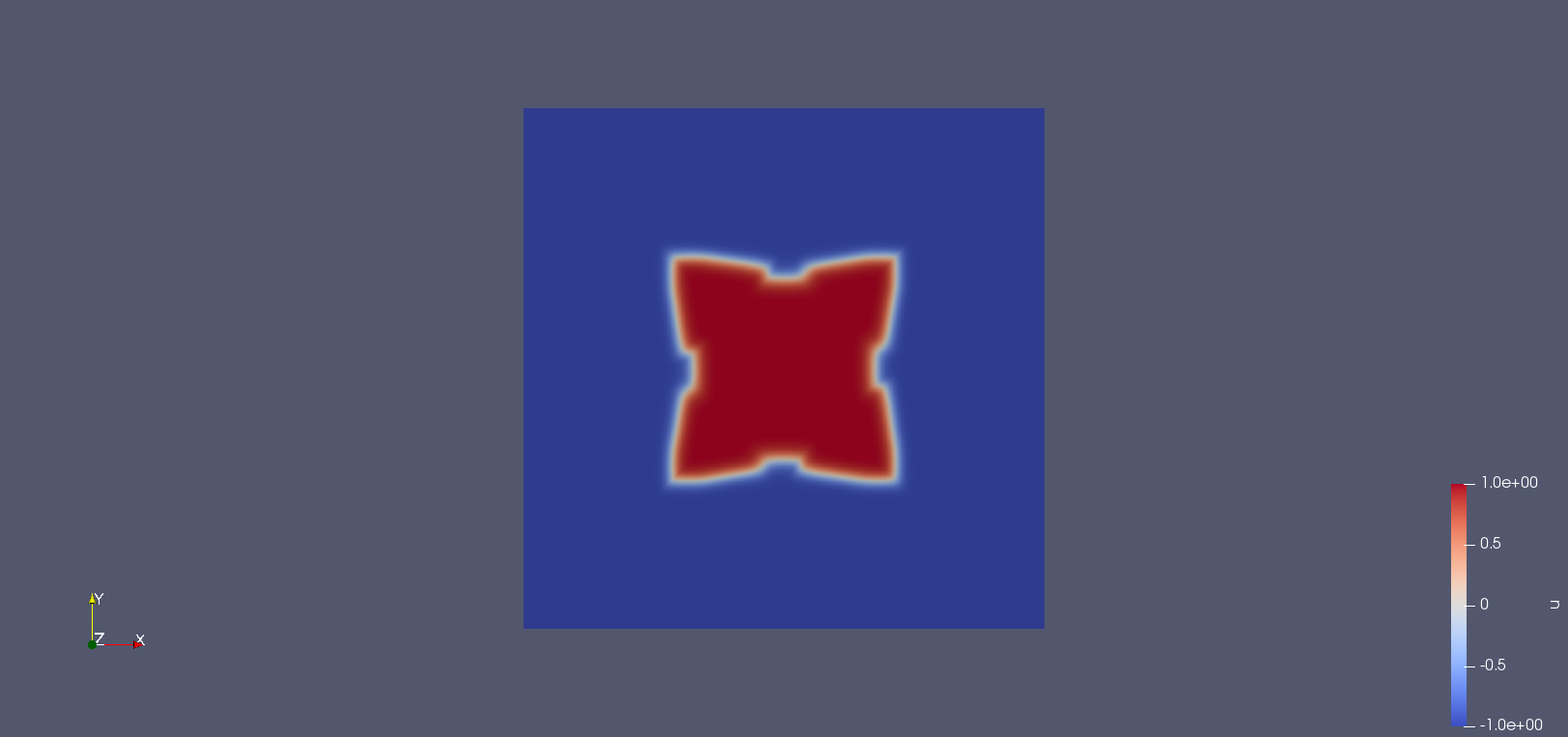}
	\includegraphics[trim={540px 112px 540px 112px}, clip, scale=.08]{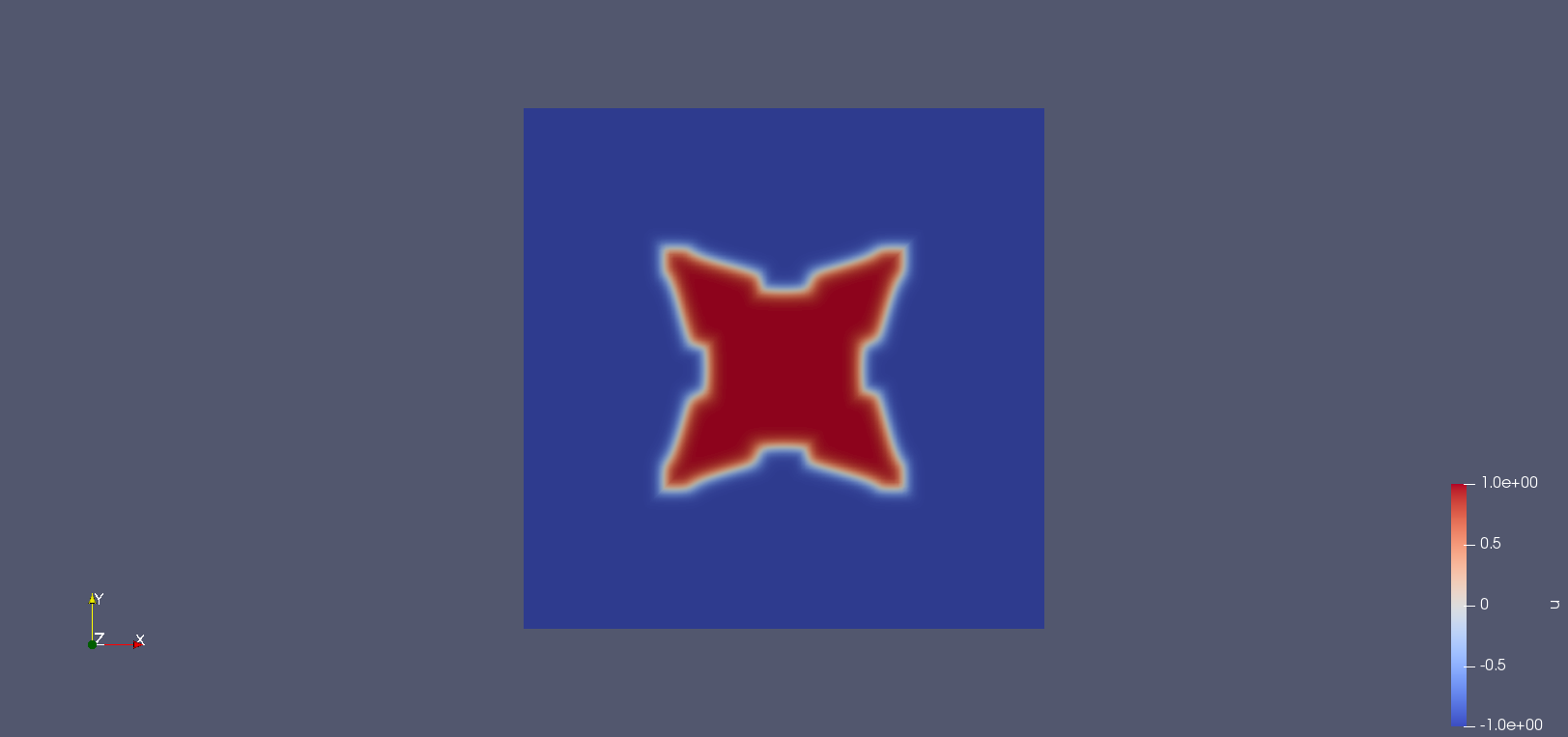}
	\includegraphics[trim={540px 112px 540px 112px}, clip, scale=.08]{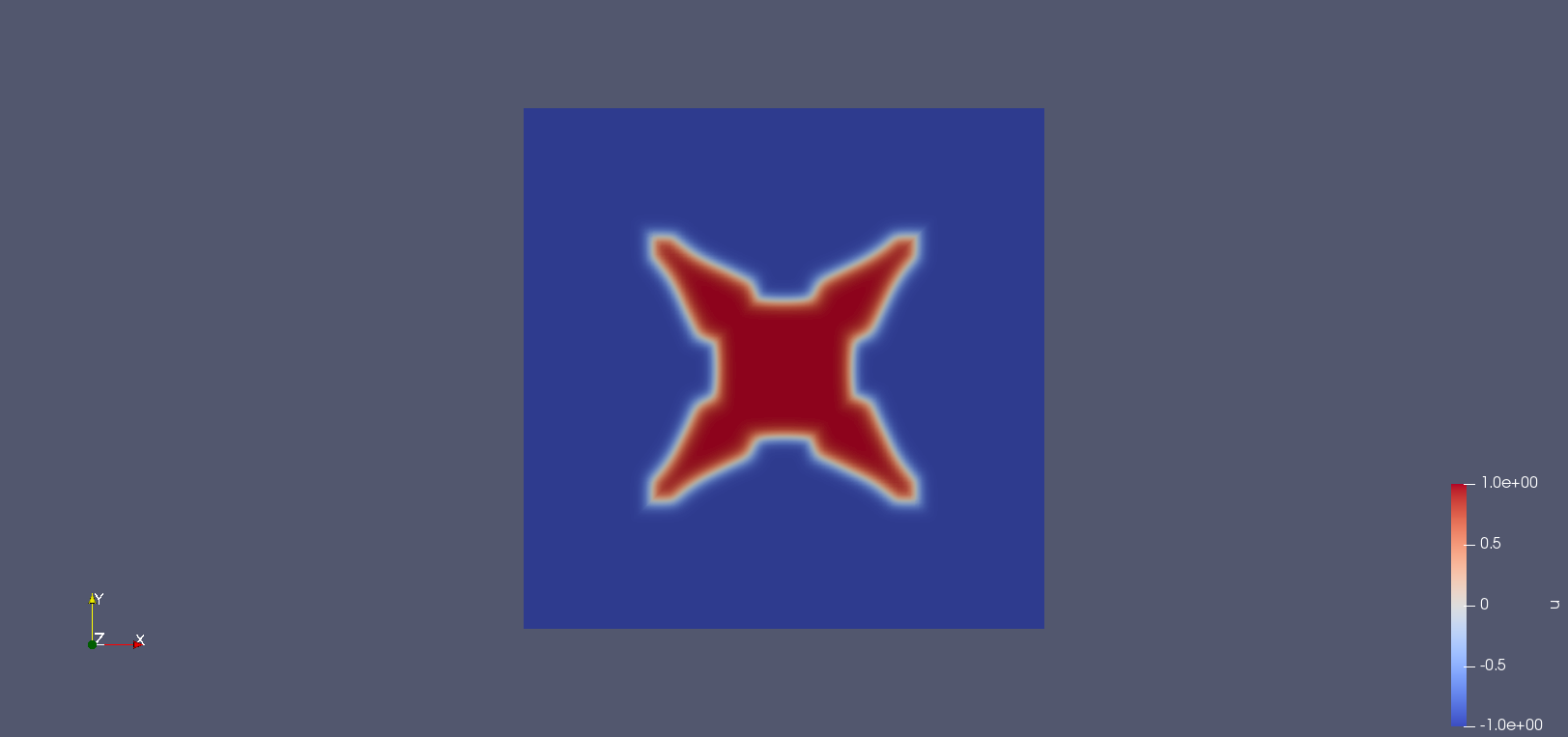}
	\includegraphics[trim={540px 112px 540px 112px}, clip, scale=.08]{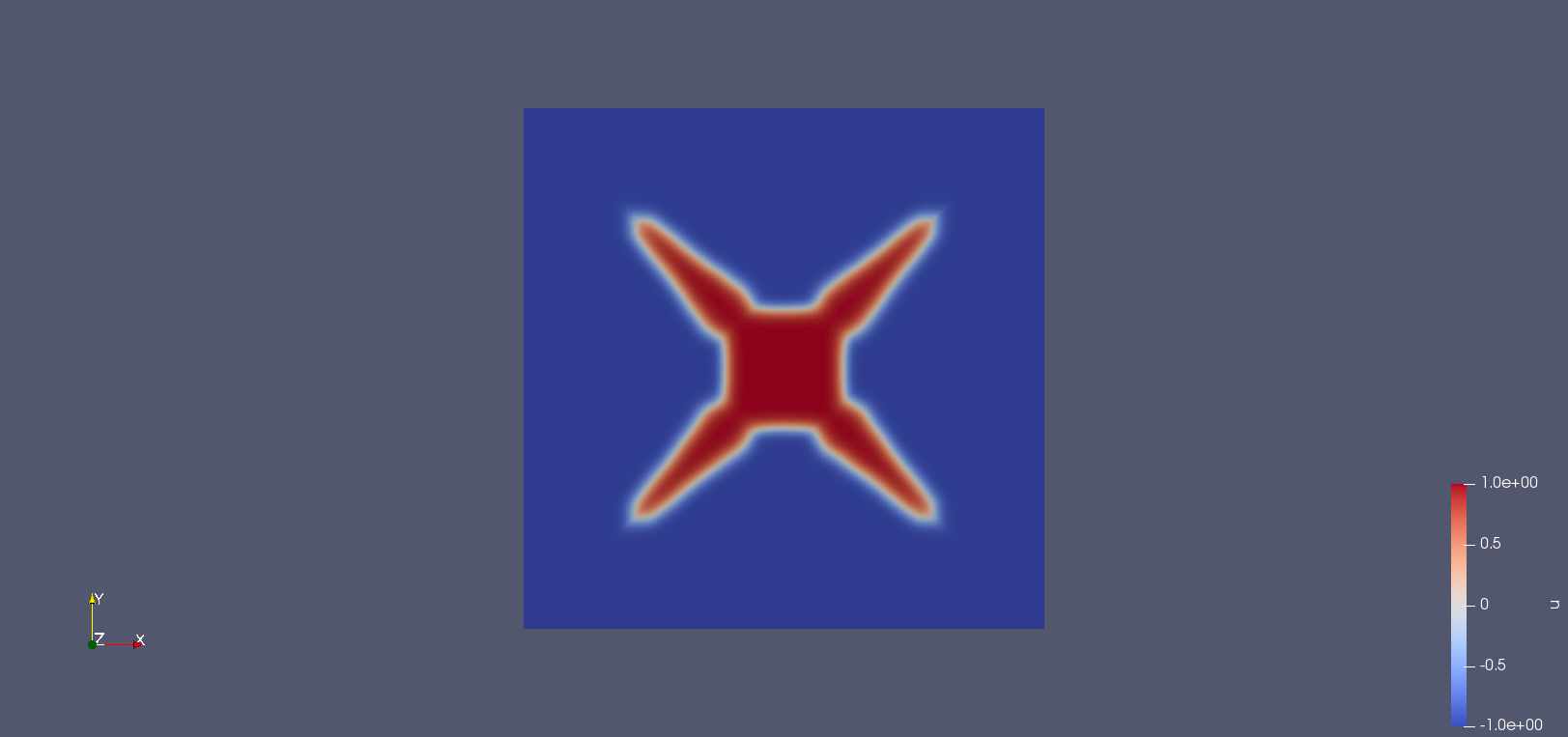}
	\includegraphics[trim={540px 112px 540px 112px}, clip, scale=.08]{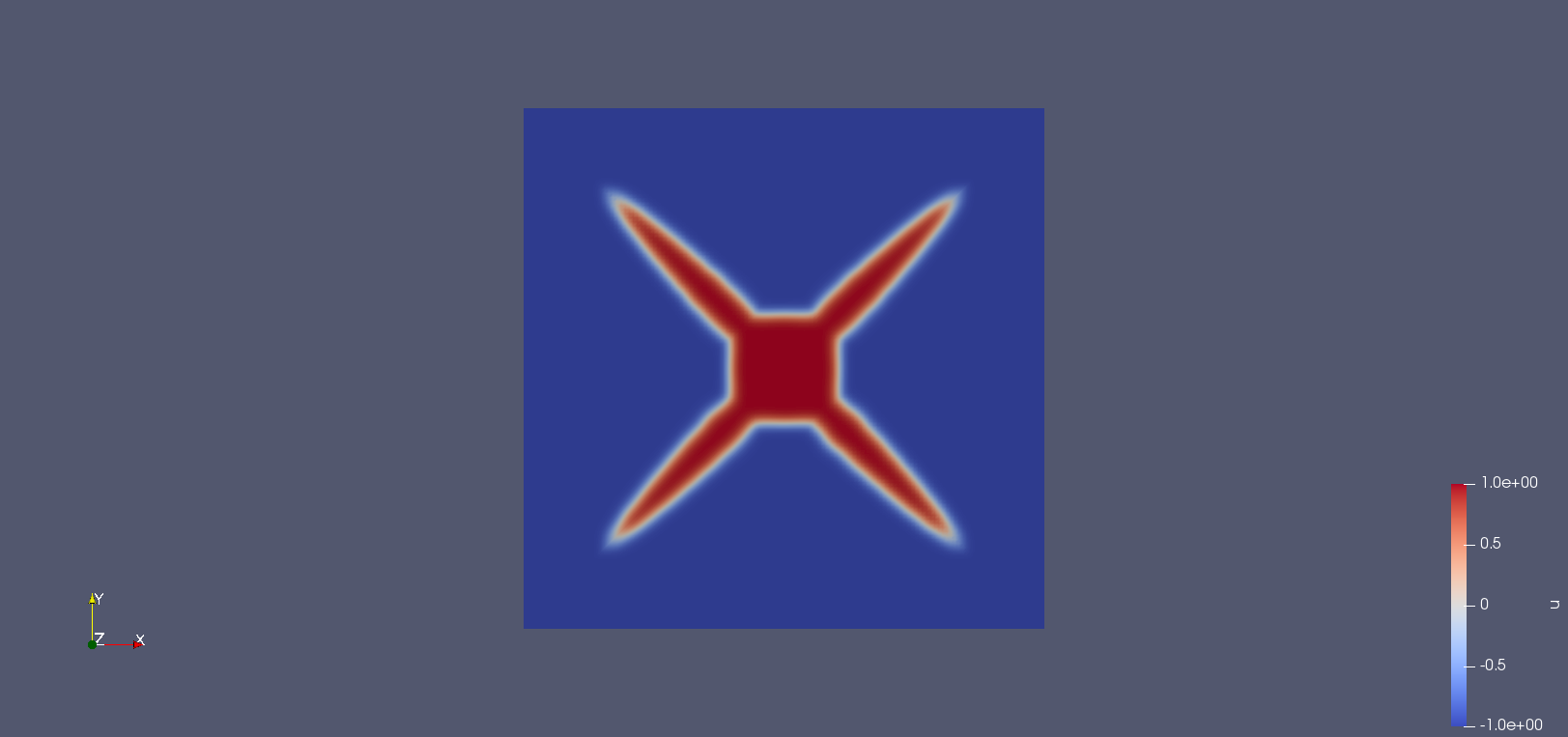}~~~
	\includegraphics[trim={1490px 0px 20px 490px}, clip, scale=.16]{{l1_star_state.0000}.png}~\\~\\
	\includegraphics[trim={540px 112px 540px 112px}, clip, scale=.08]{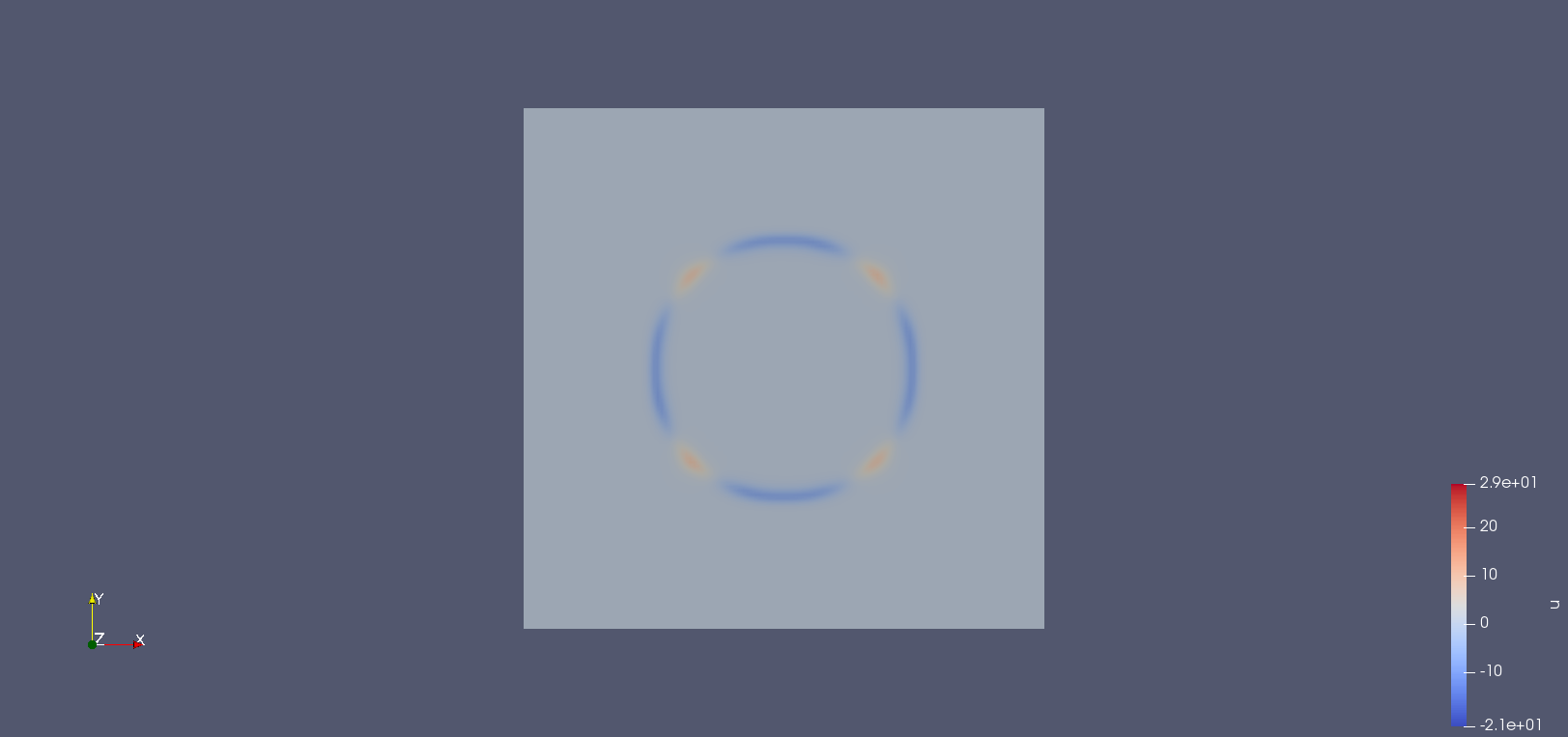}
	\includegraphics[trim={540px 112px 540px 112px}, clip, scale=.08]{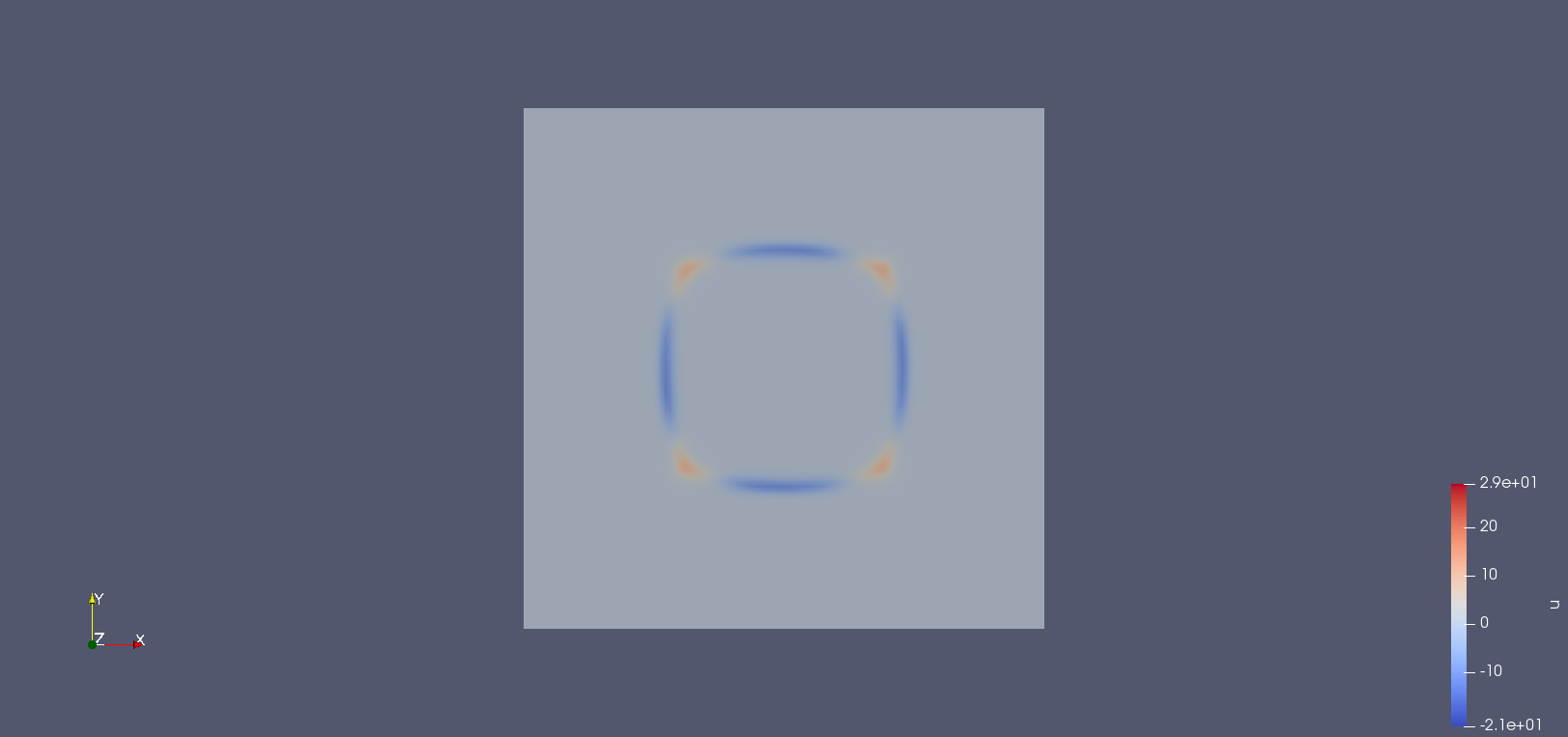}
	\includegraphics[trim={540px 112px 540px 112px}, clip, scale=.08]{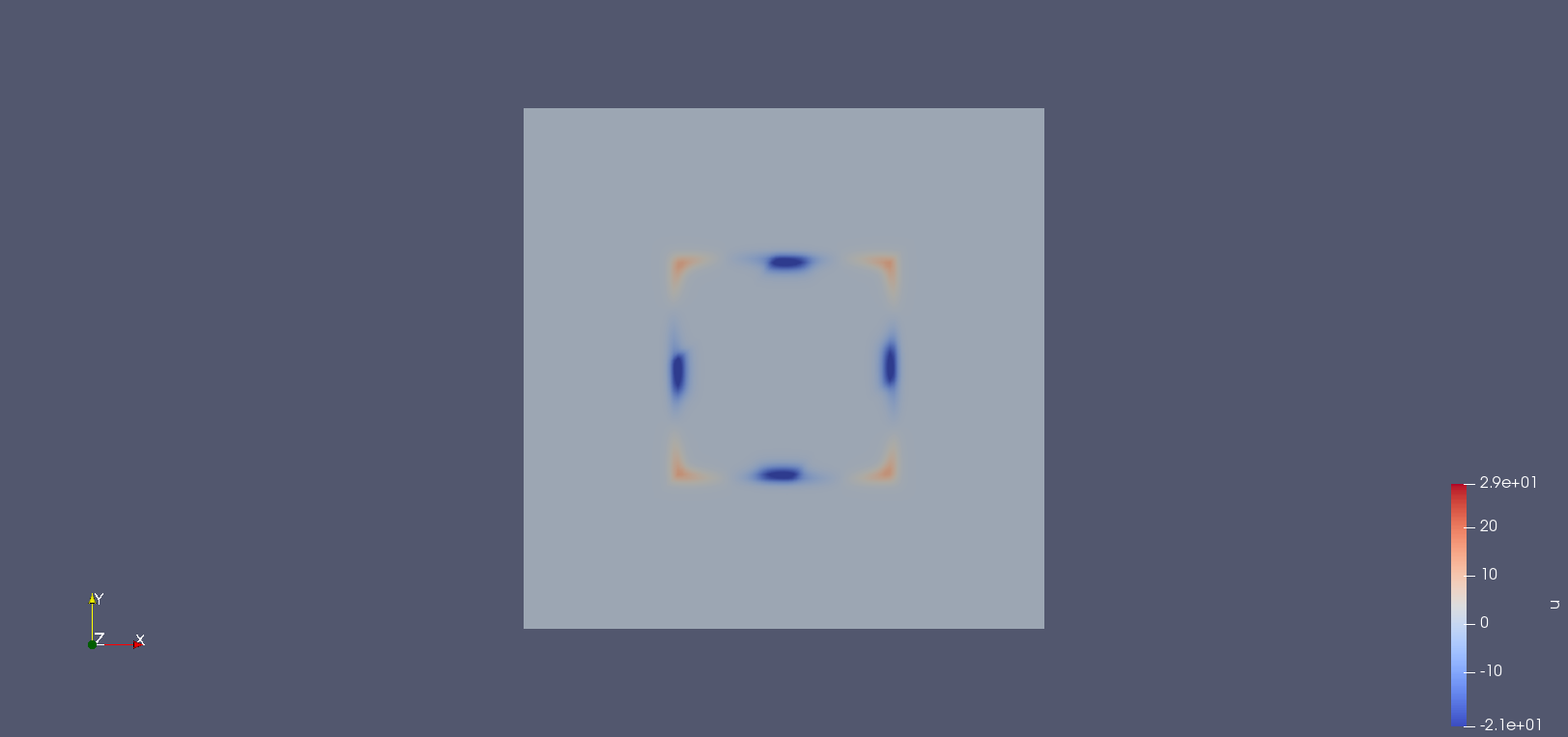}
	\includegraphics[trim={540px 112px 540px 112px}, clip, scale=.08]{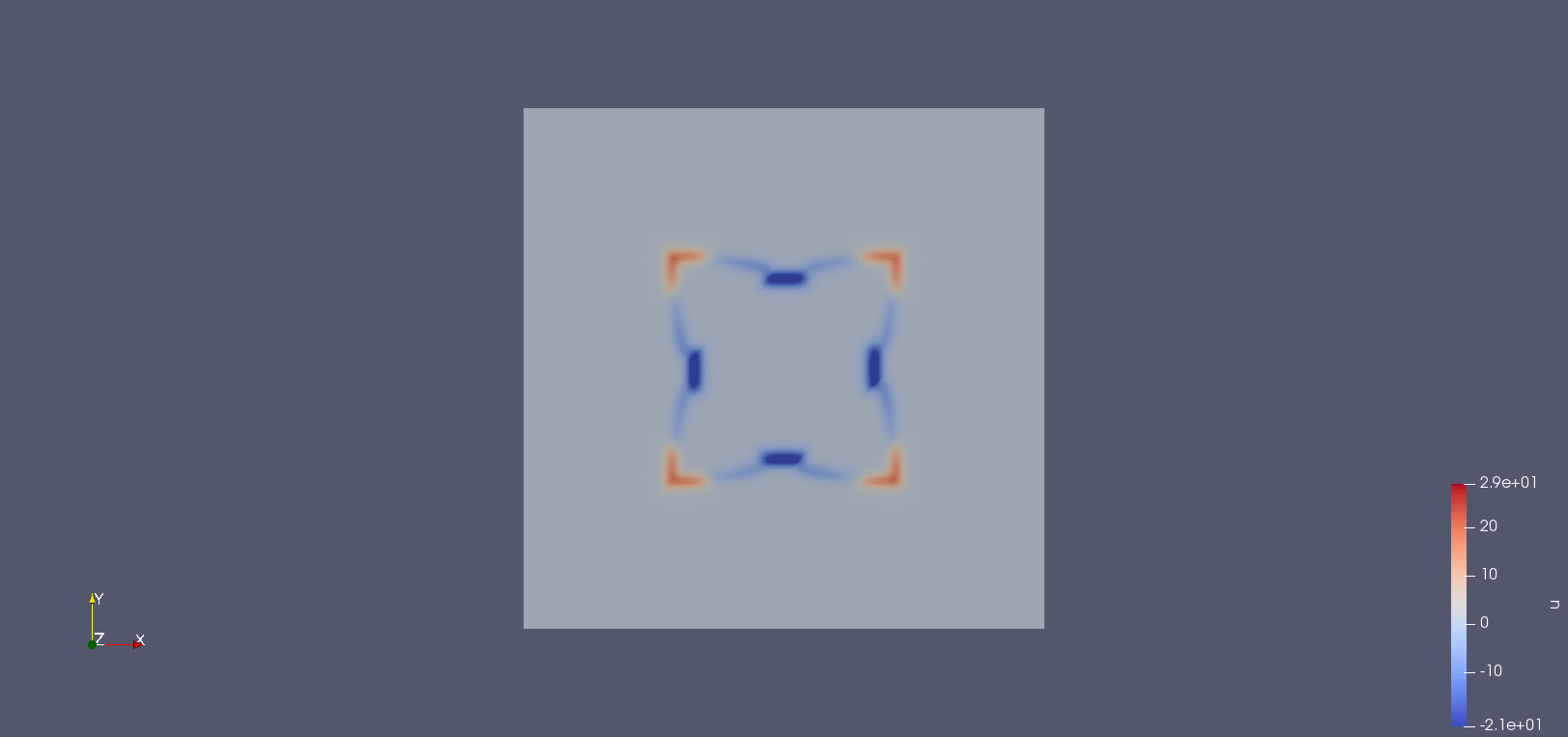}
	\includegraphics[trim={540px 112px 540px 112px}, clip, scale=.08]{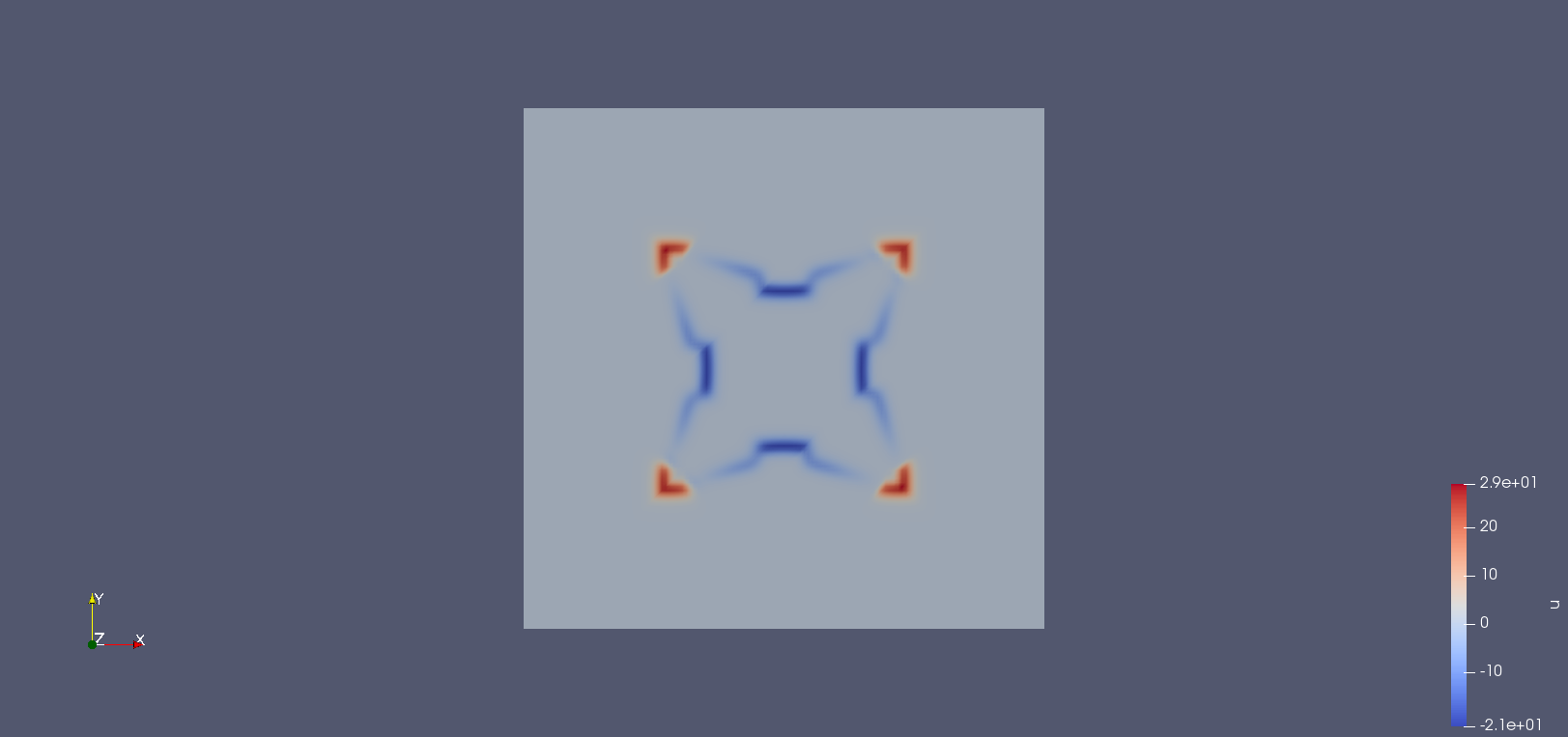}
	\includegraphics[trim={540px 112px 540px 112px}, clip, scale=.08]{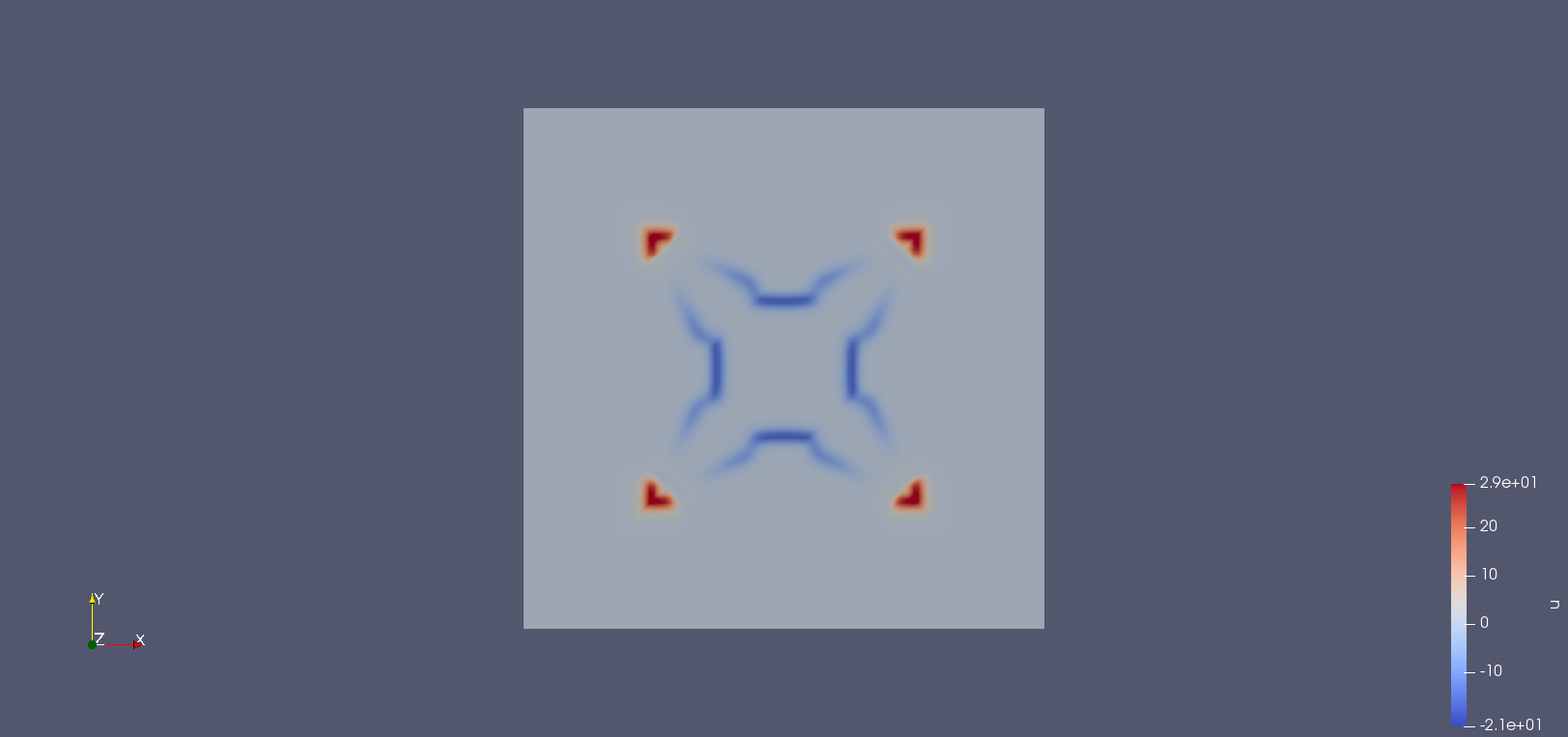}
	\includegraphics[trim={540px 112px 540px 112px}, clip, scale=.08]{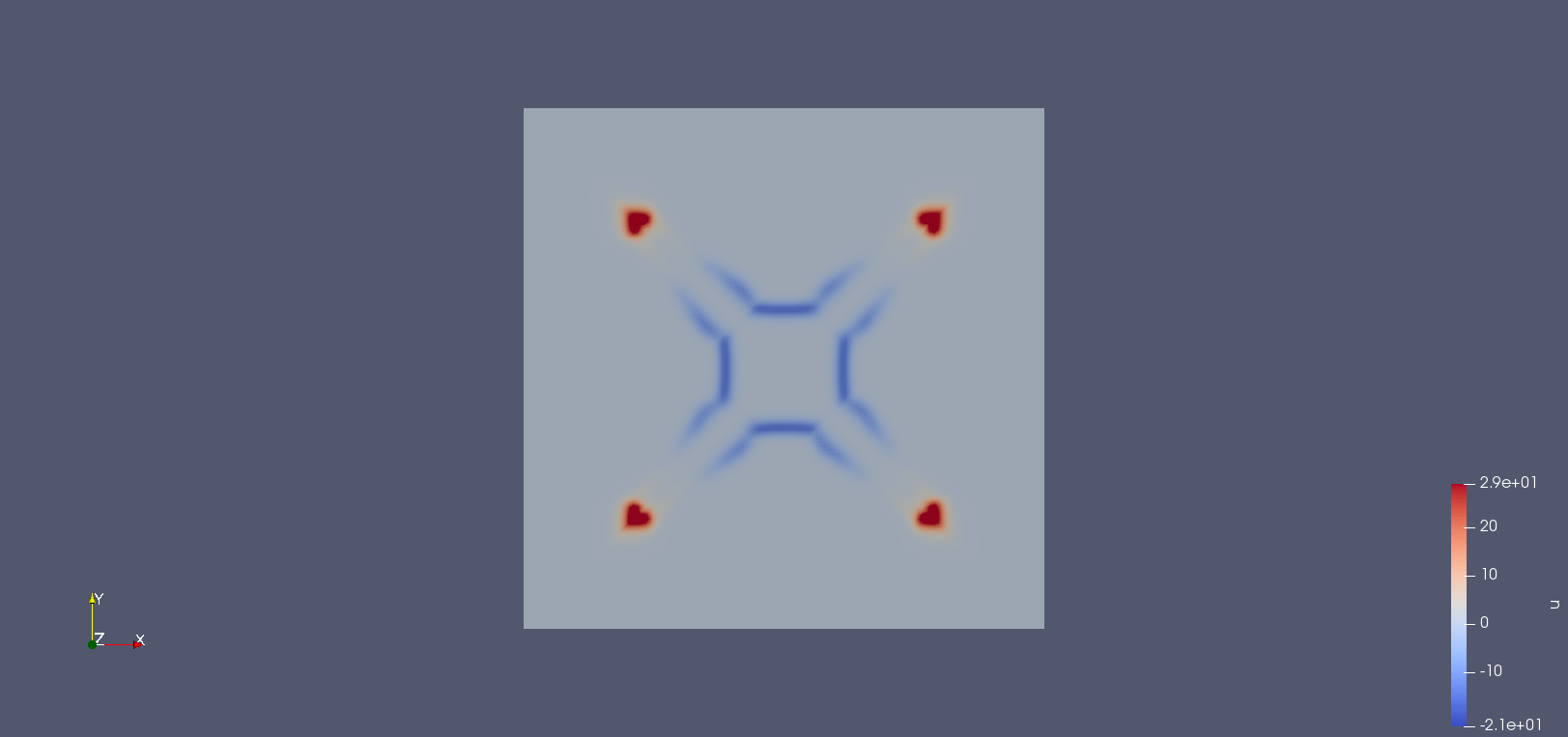}
	\includegraphics[trim={540px 112px 540px 112px}, clip, scale=.08]{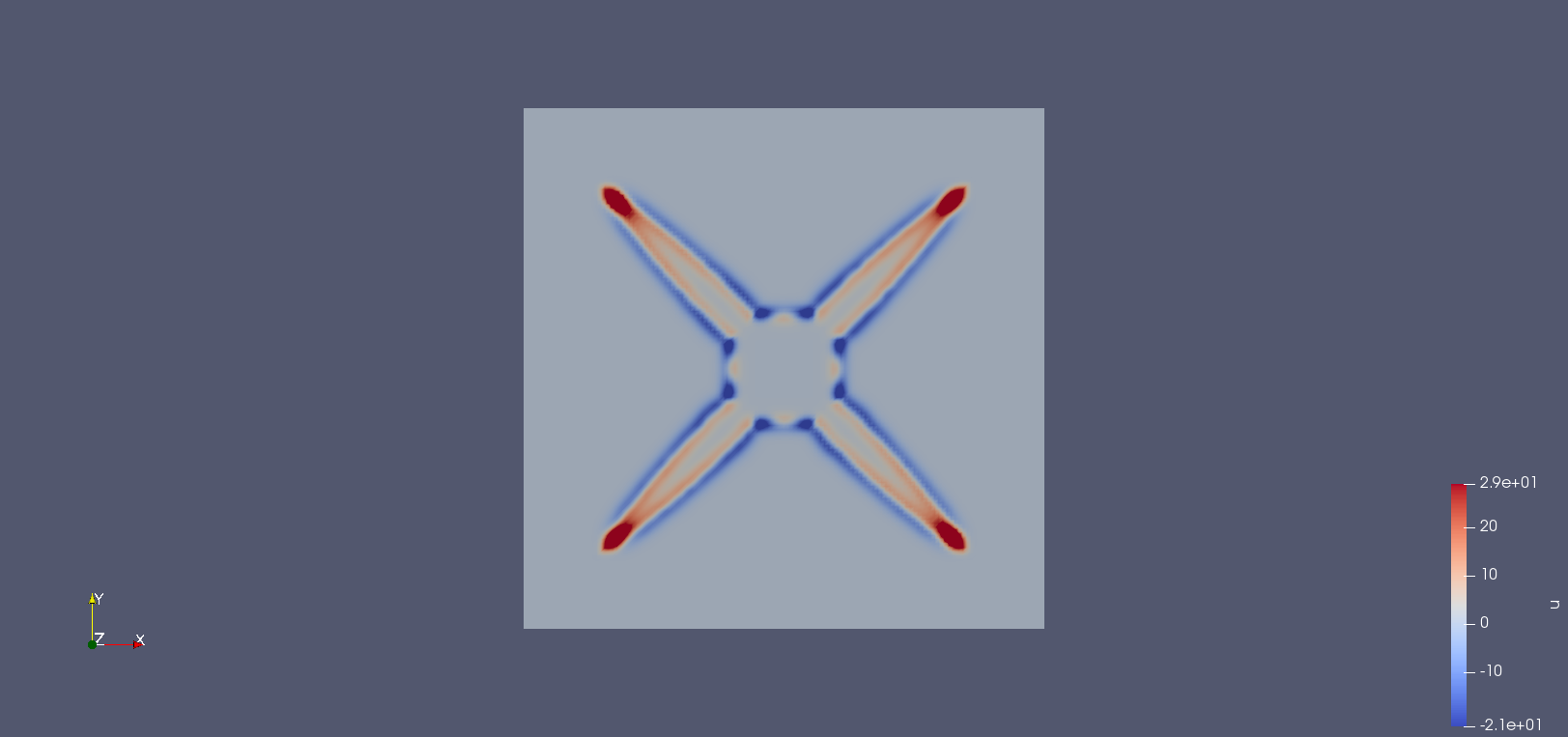}~~~
	\includegraphics[trim={1490px 0px 20px 490px}, clip, scale=.16]{{l1_star_control.0000}.png}
\caption{Results for the regularized $l_1$-norm.} 
      \label{fig:l1_star4}
\end{center}    
\end{subfigure}
\begin{subfigure}{1.0\textwidth}
\begin{center}    
	\includegraphics[trim={540px 112px 540px 112px}, clip, scale=.08]{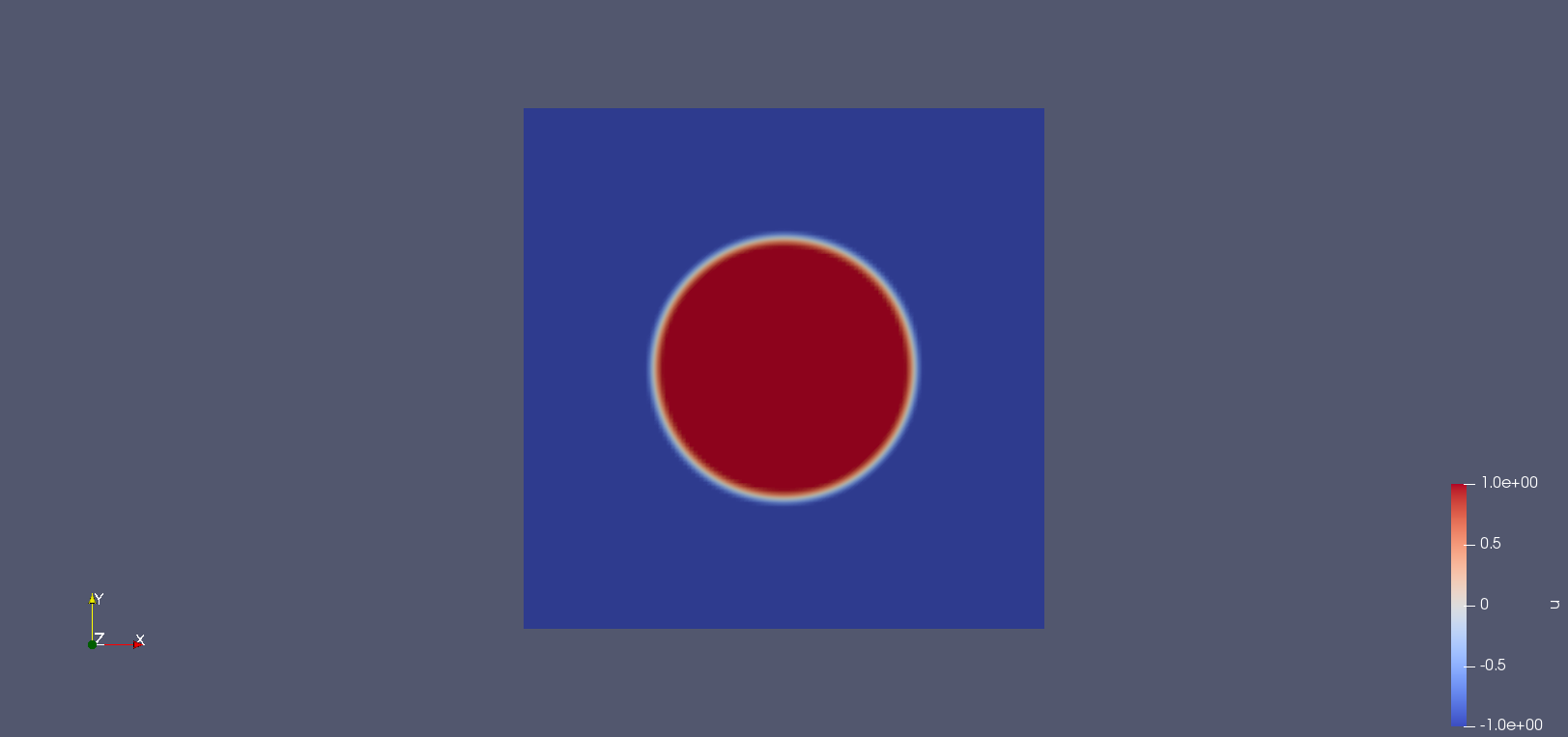}
	\includegraphics[trim={540px 112px 540px 112px}, clip, scale=.08]{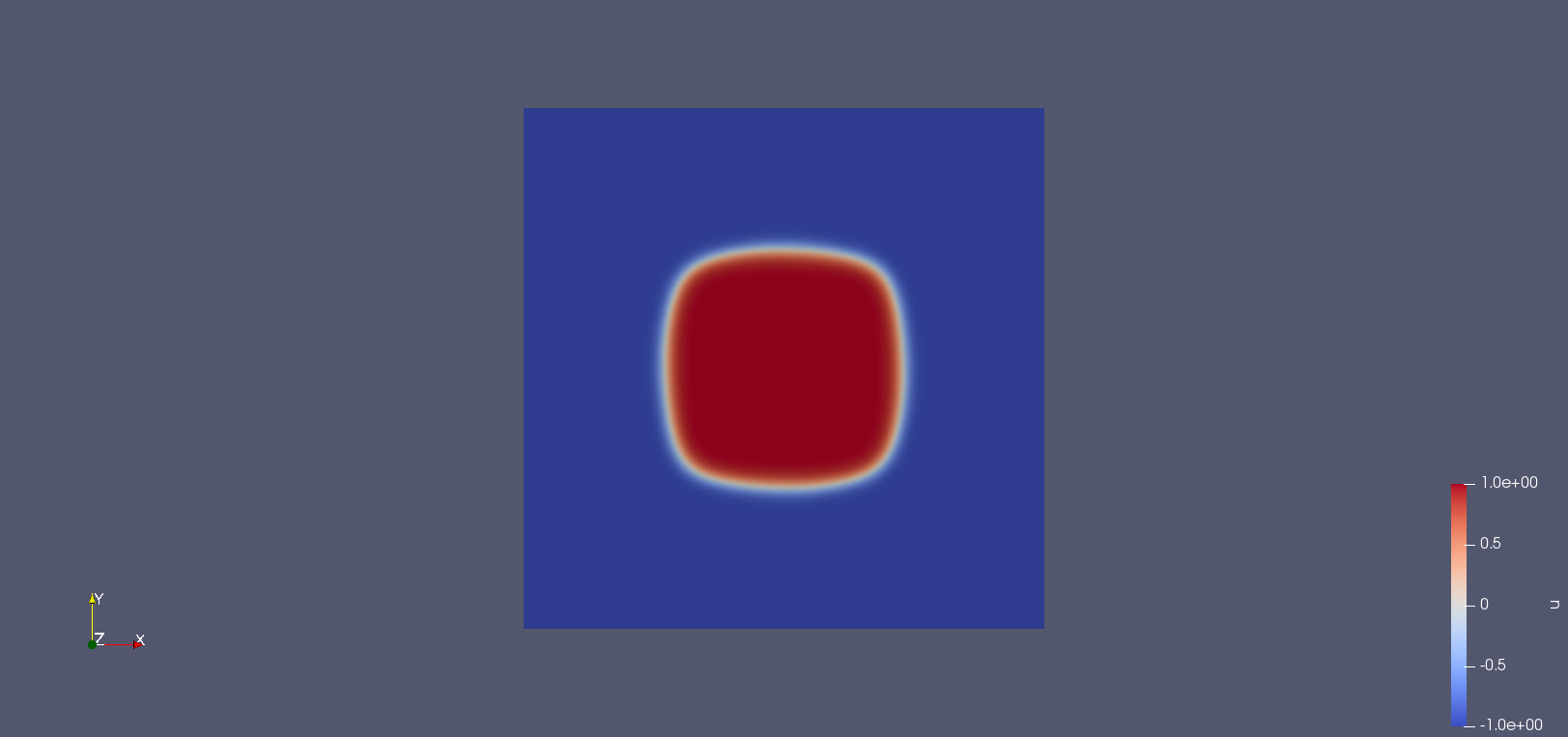}
	\includegraphics[trim={540px 112px 540px 112px}, clip, scale=.08]{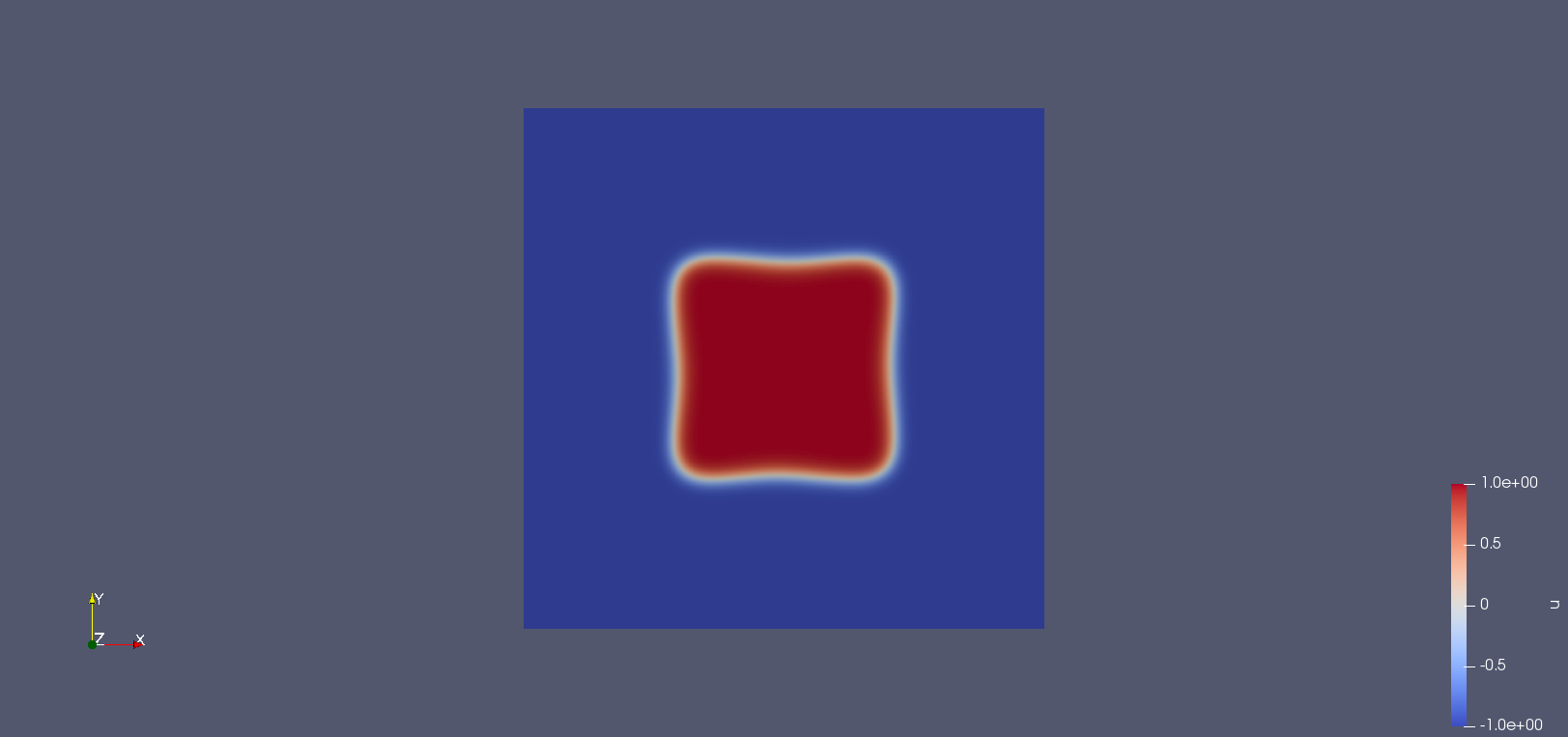}
	\includegraphics[trim={540px 112px 540px 112px}, clip, scale=.08]{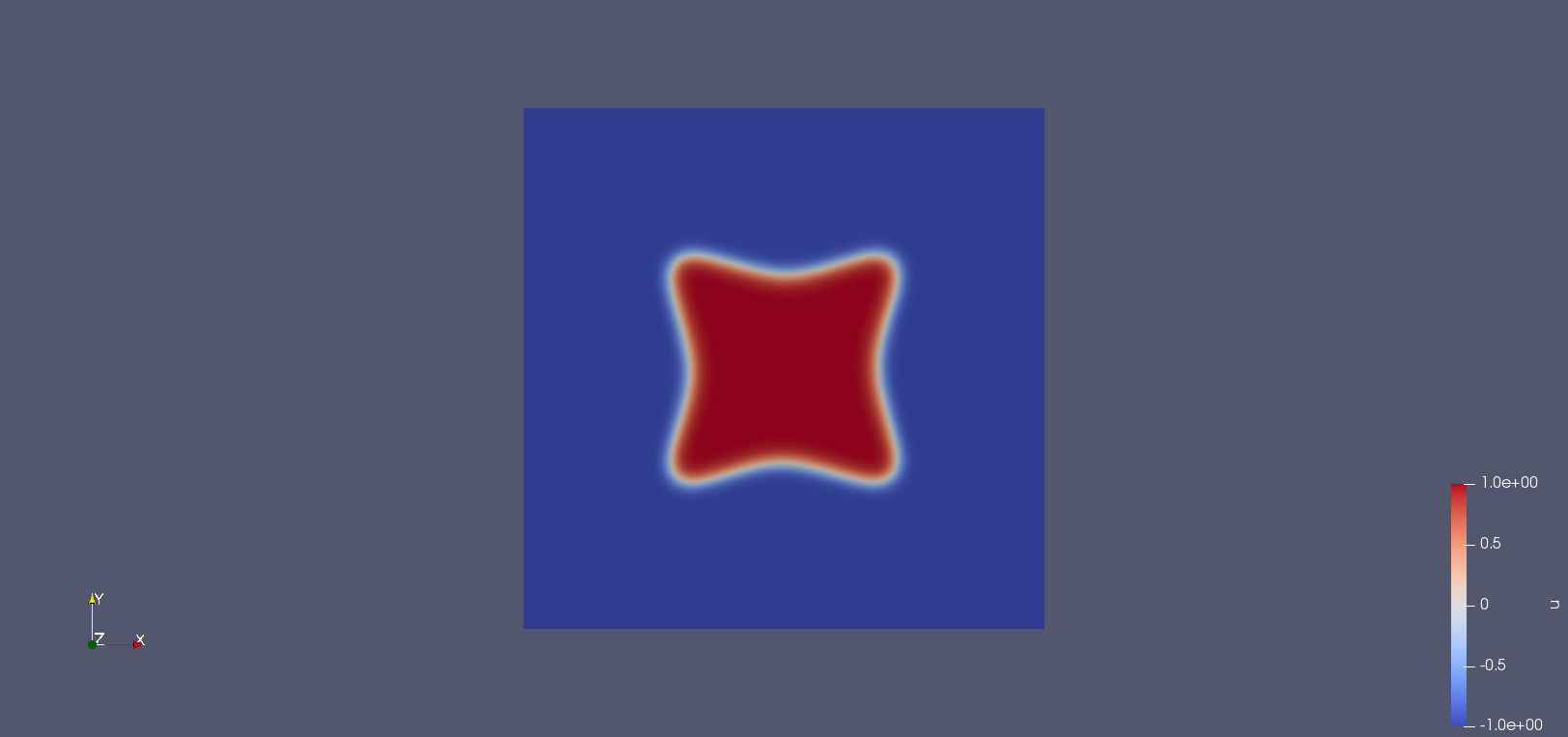}
	\includegraphics[trim={540px 112px 540px 112px}, clip, scale=.08]{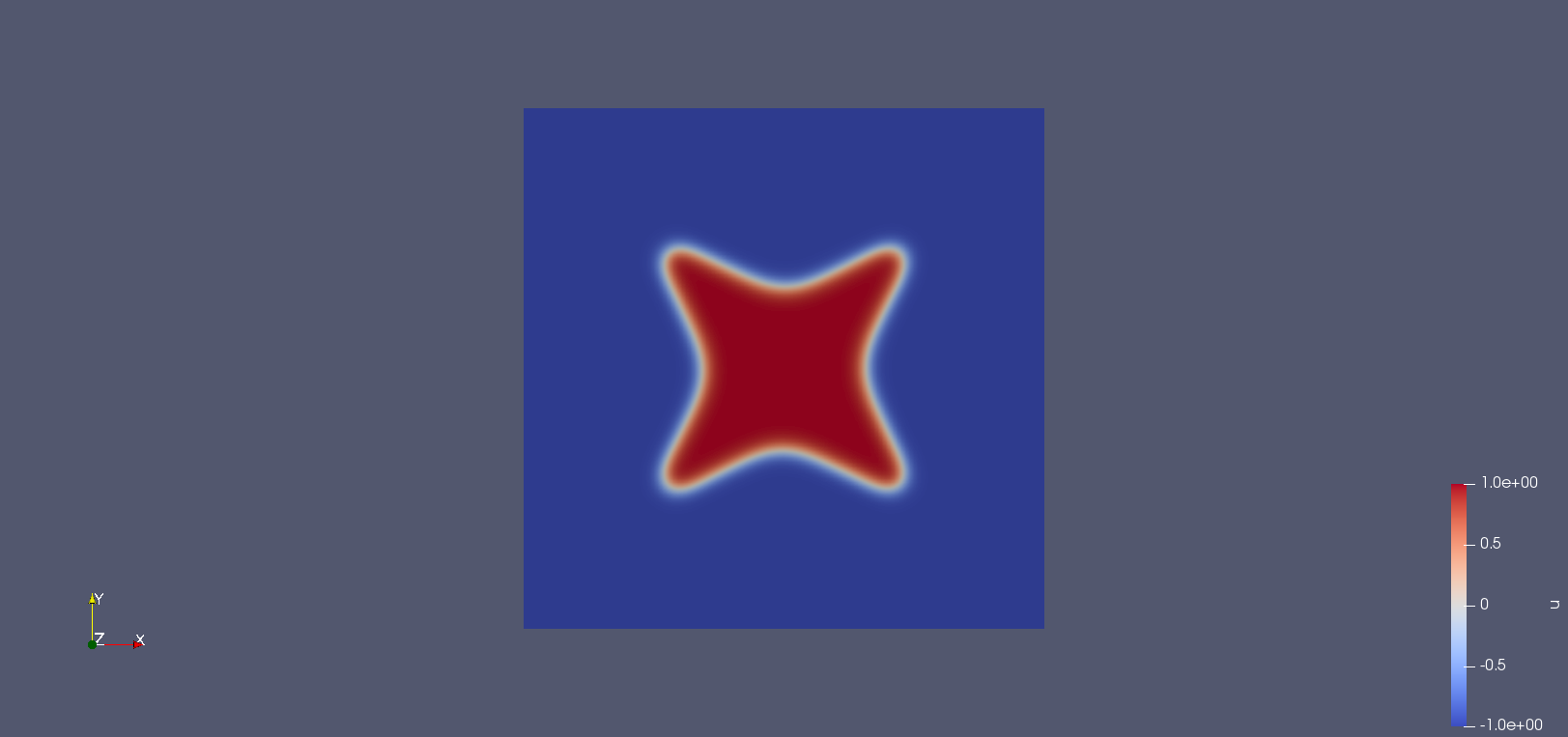}
	\includegraphics[trim={540px 112px 540px 112px}, clip, scale=.08]{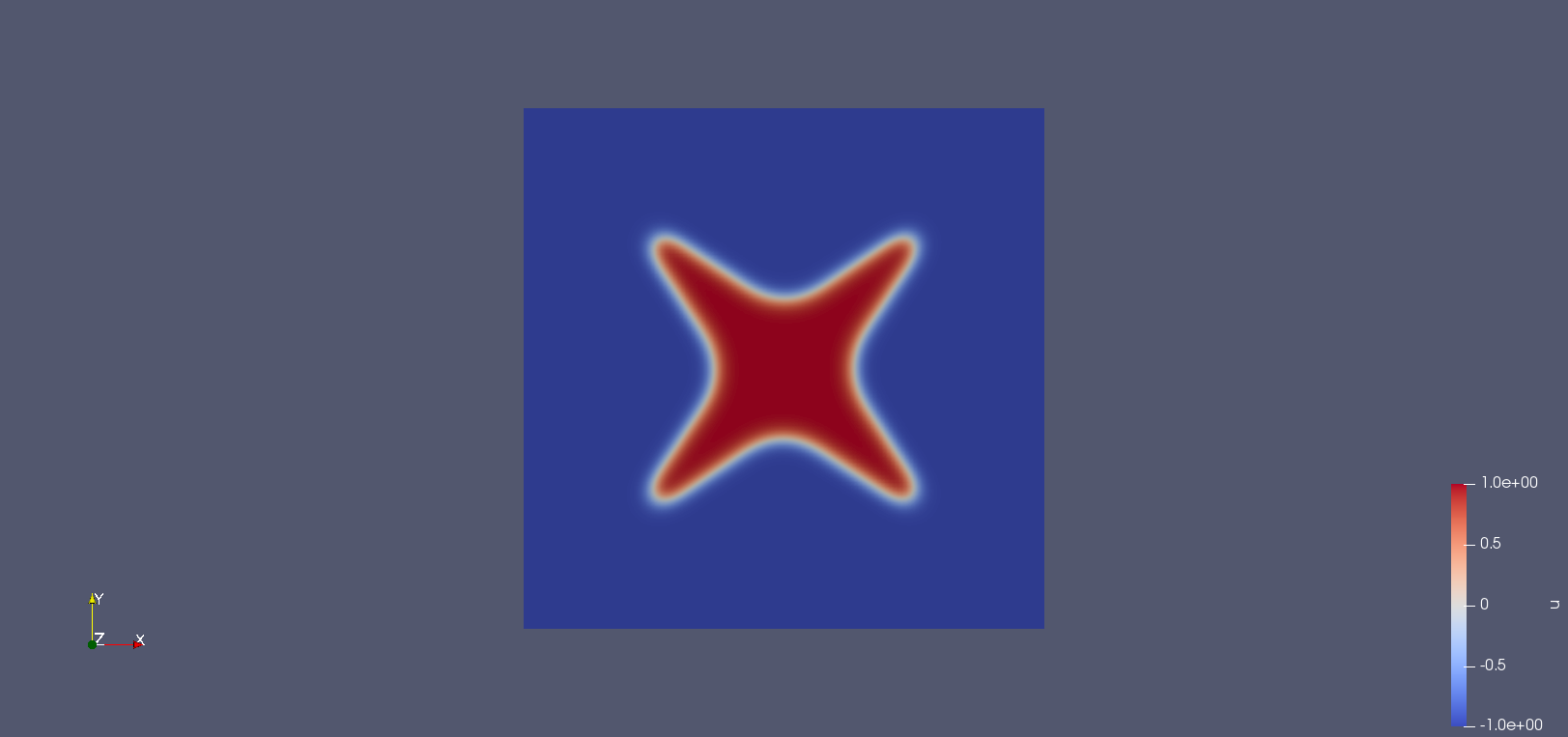}
	\includegraphics[trim={540px 112px 540px 112px}, clip, scale=.08]{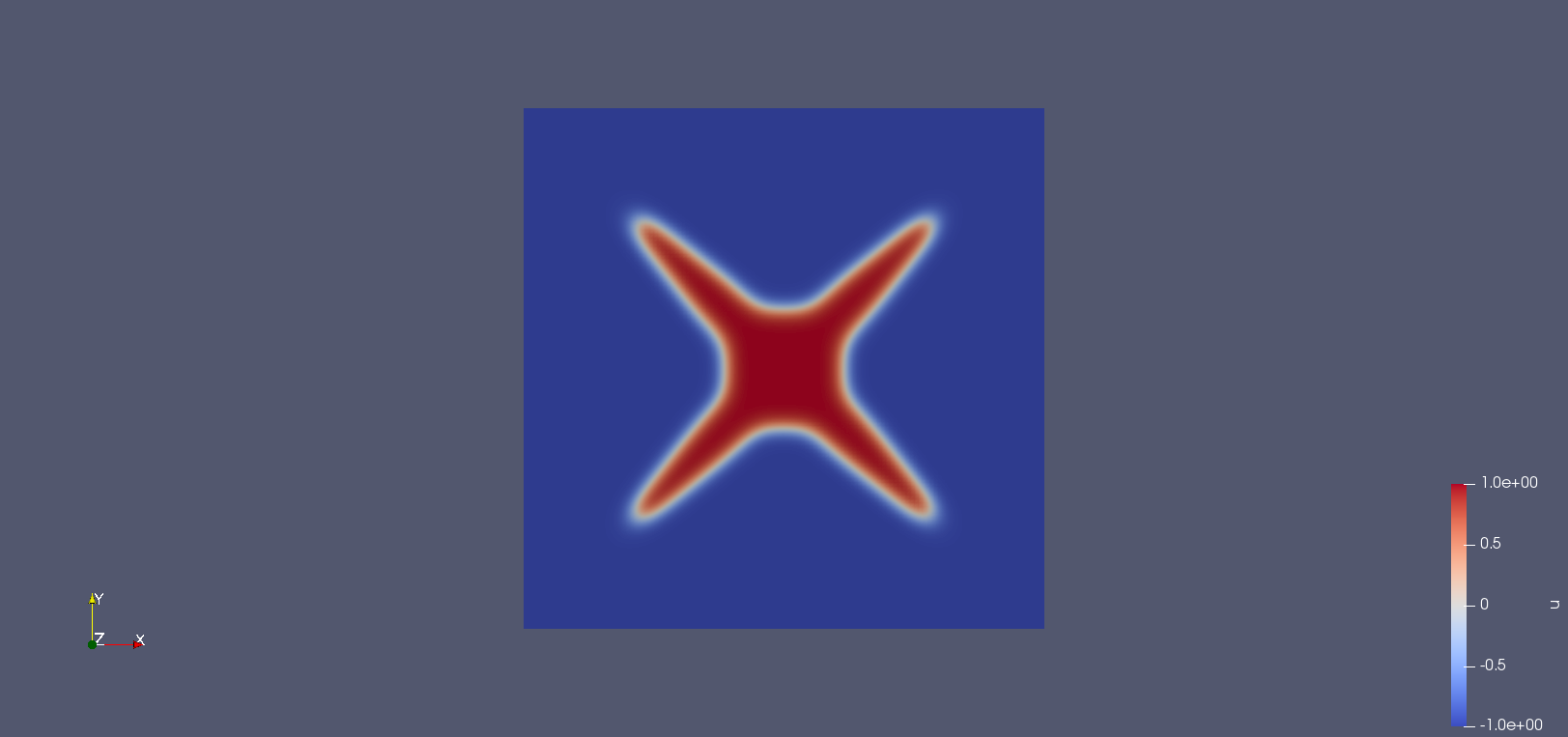}
	\includegraphics[trim={540px 112px 540px 112px}, clip, scale=.08]{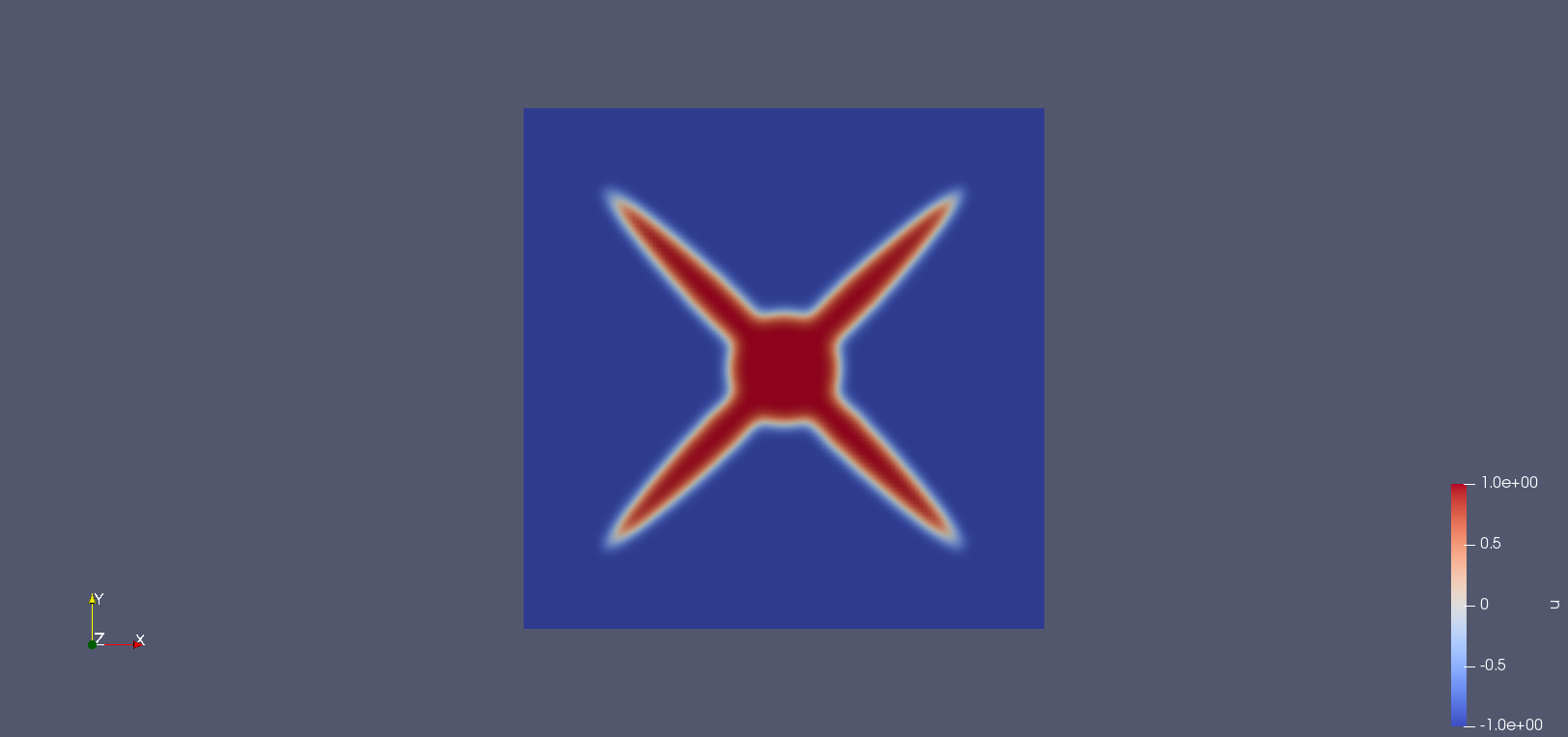}~~~
	\includegraphics[trim={1490px 0px 20px 490px}, clip, scale=.16]{{iso_starl1_state.0000}.png}~\\~\\
	\includegraphics[trim={540px 112px 540px 112px}, clip, scale=.08]{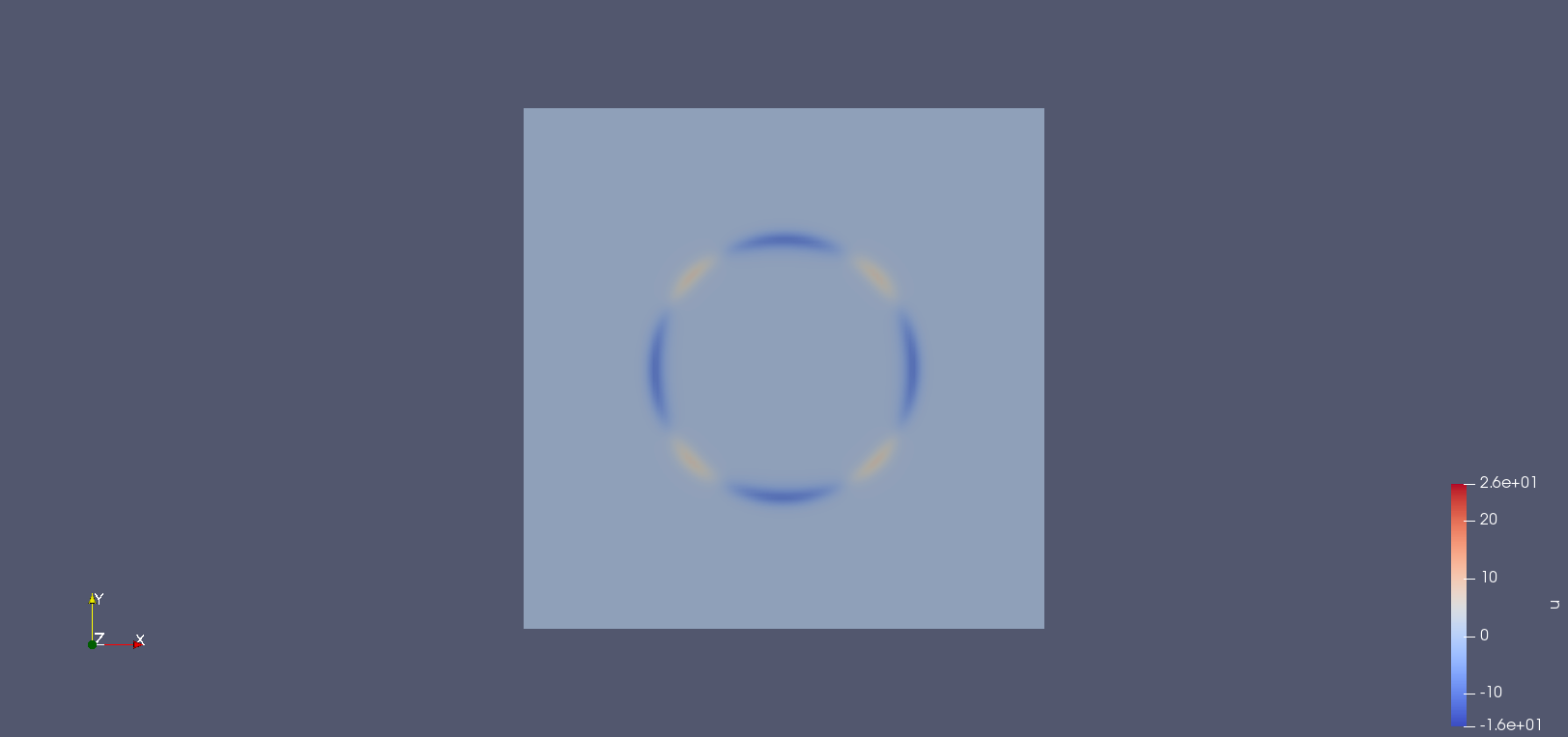}
	\includegraphics[trim={540px 112px 540px 112px}, clip, scale=.08]{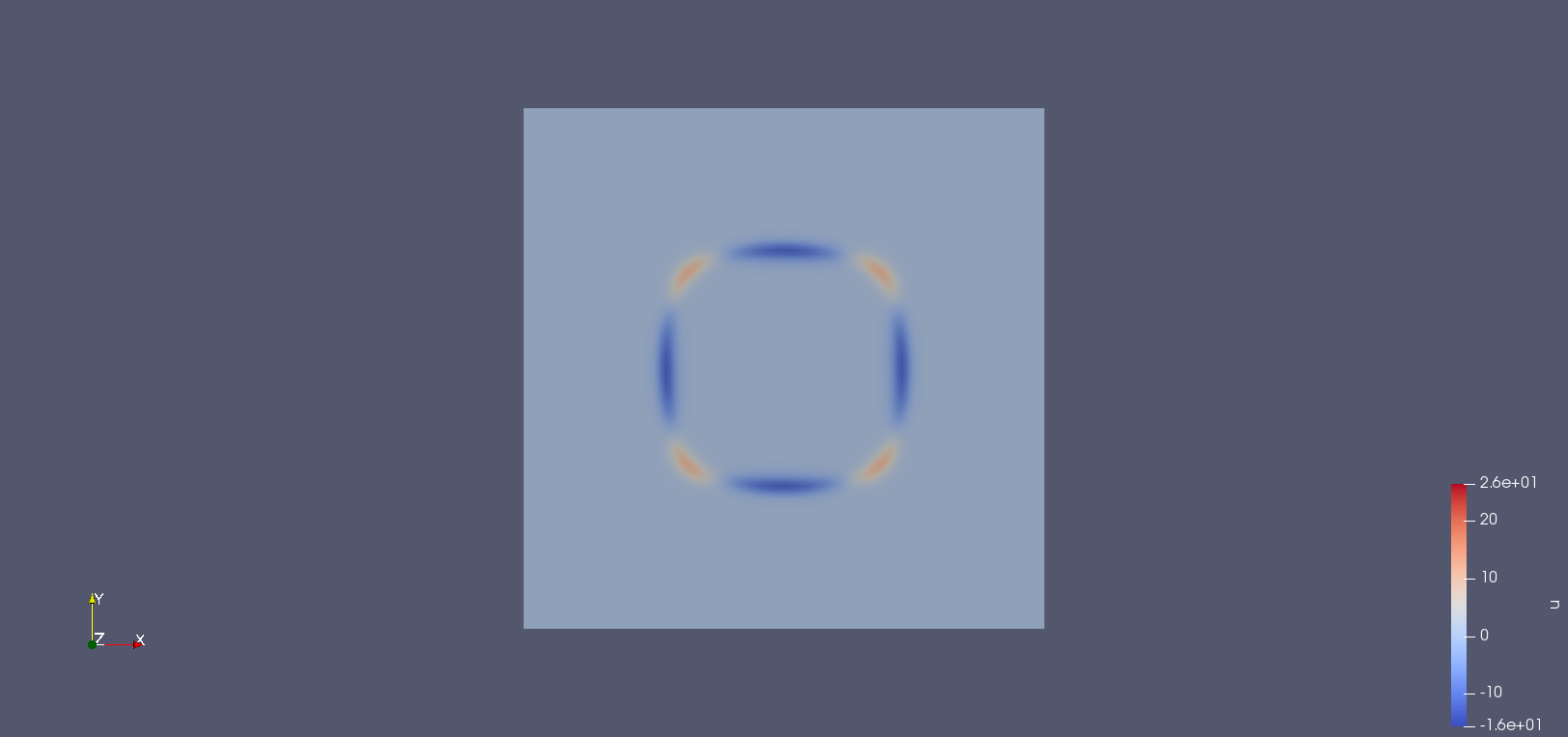}
	\includegraphics[trim={540px 112px 540px 112px}, clip, scale=.08]{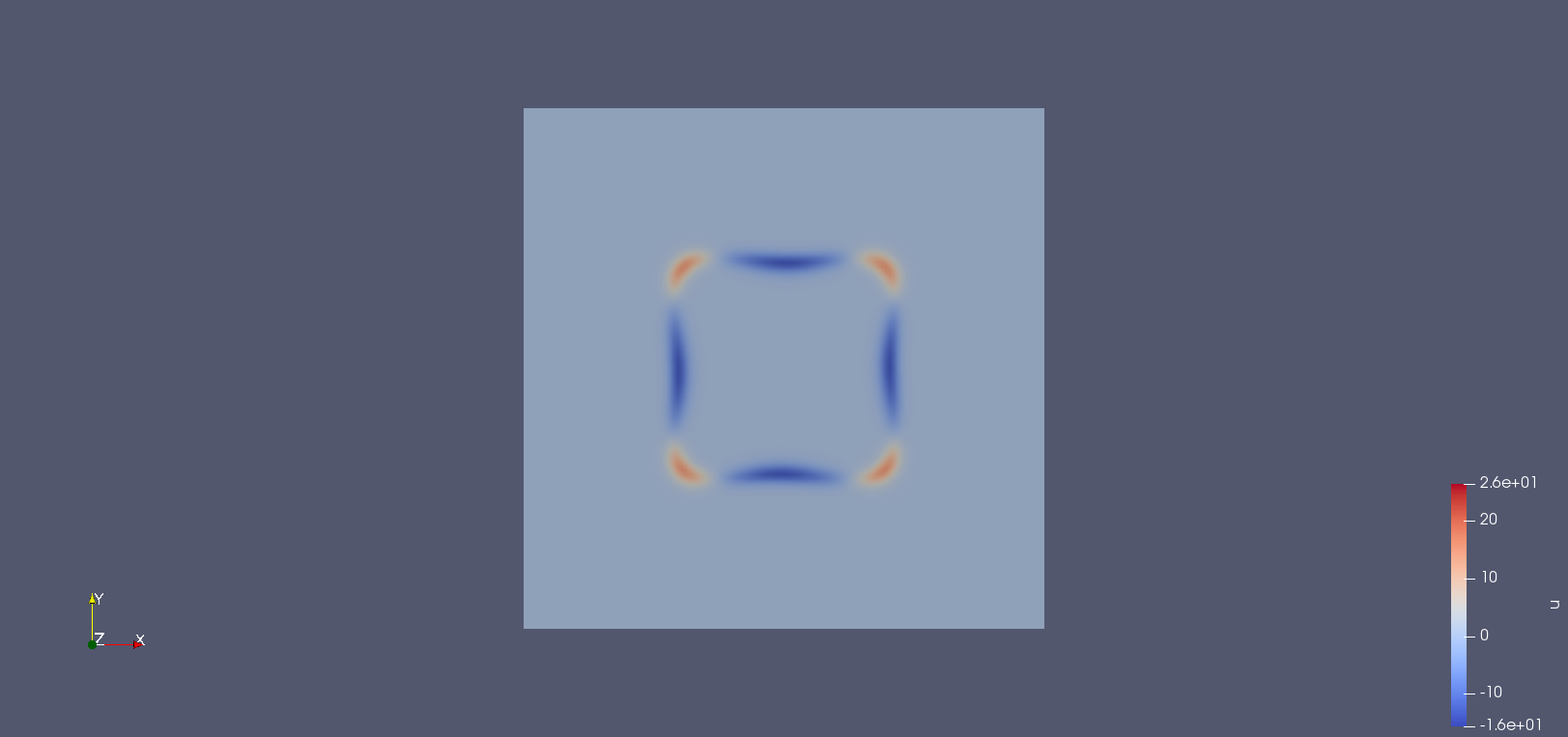}
	\includegraphics[trim={540px 112px 540px 112px}, clip, scale=.08]{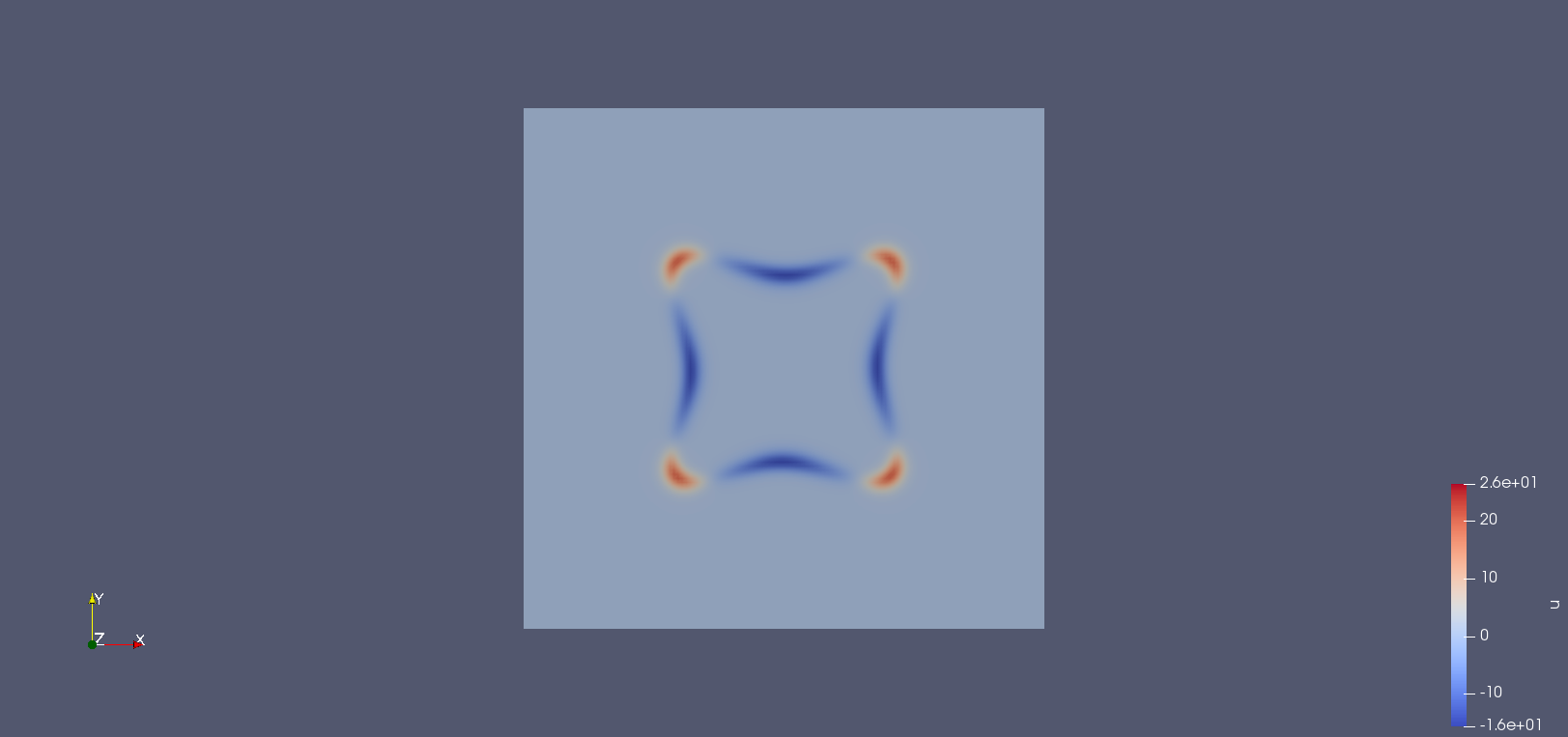}
	\includegraphics[trim={540px 112px 540px 112px}, clip, scale=.08]{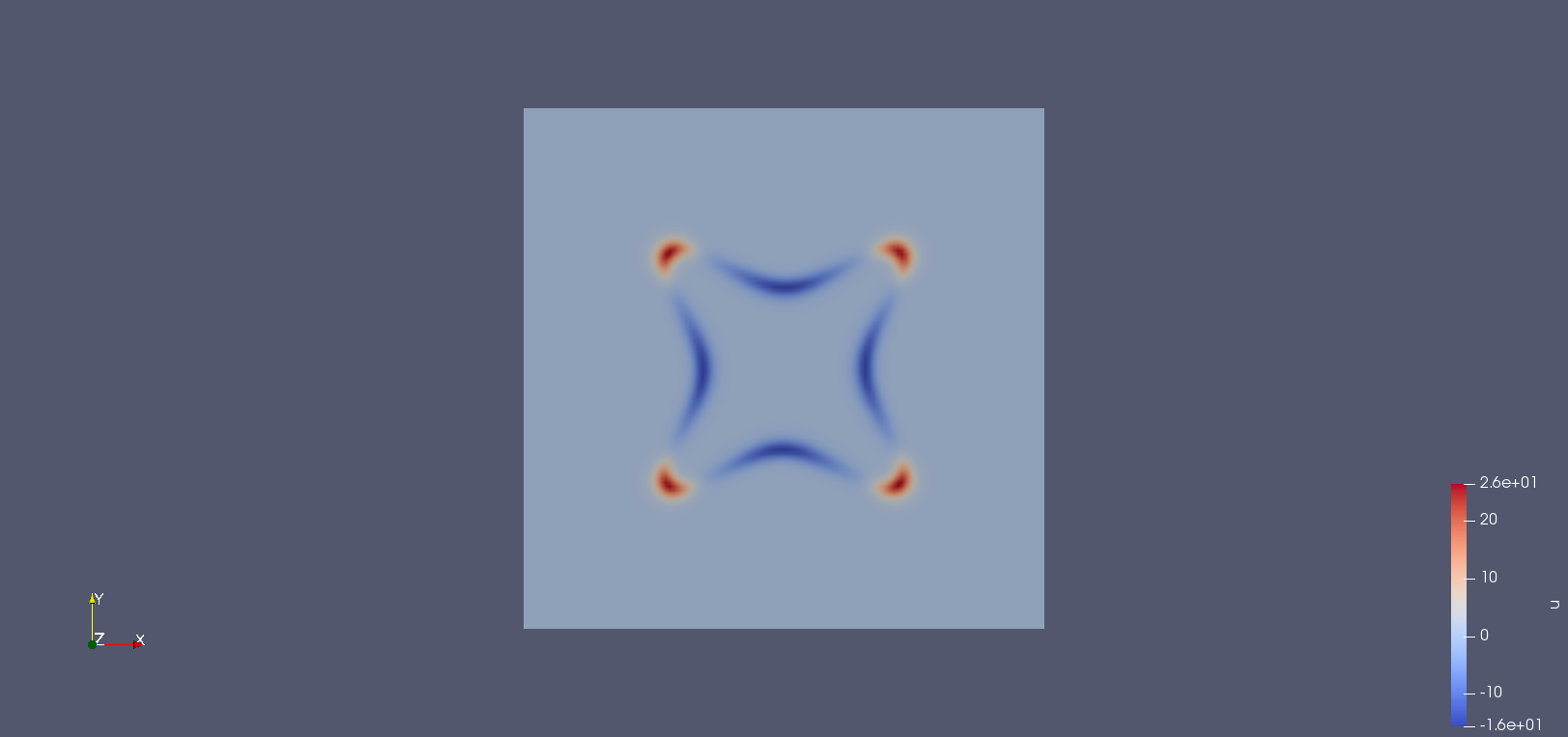}
	\includegraphics[trim={540px 112px 540px 112px}, clip, scale=.08]{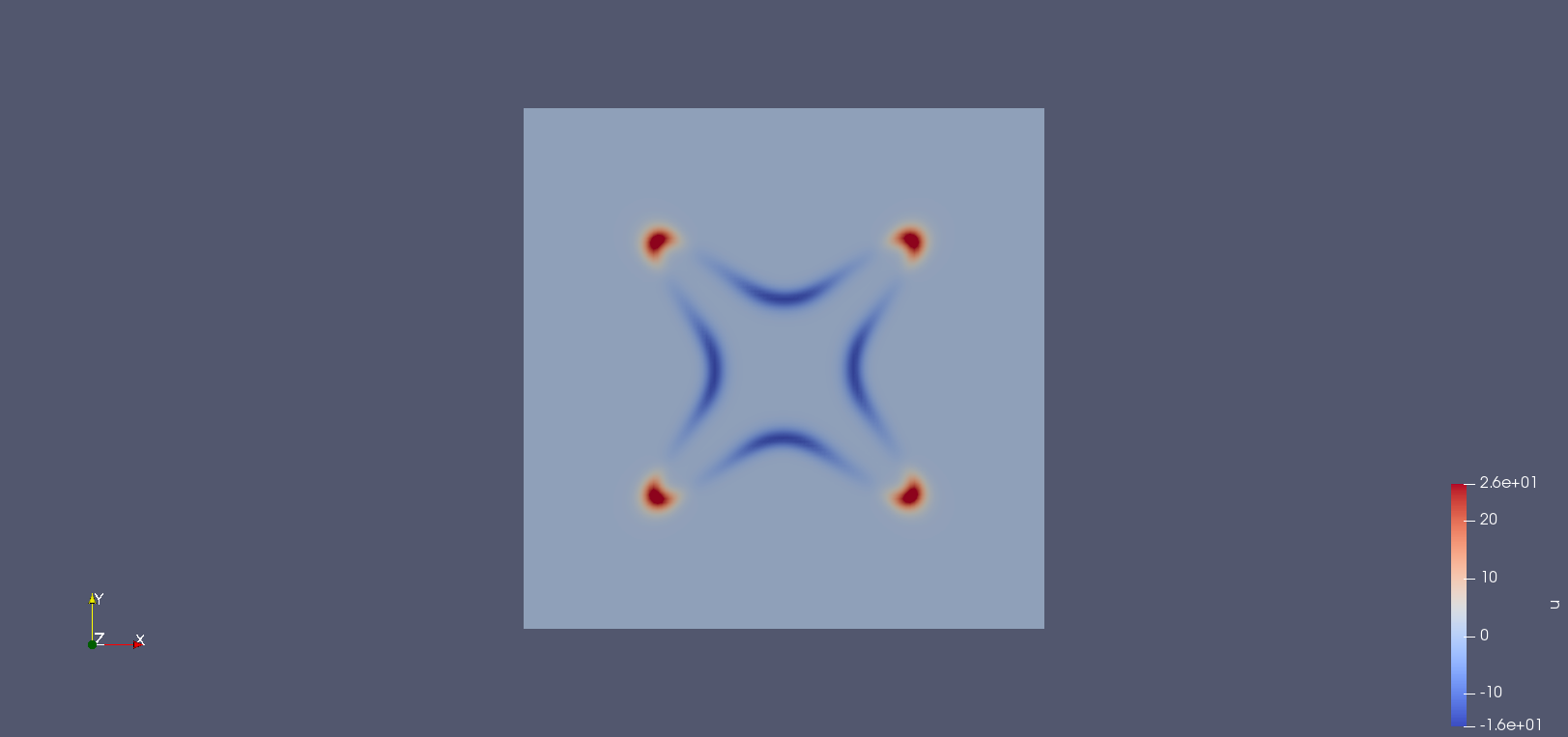}
	\includegraphics[trim={540px 112px 540px 112px}, clip, scale=.08]{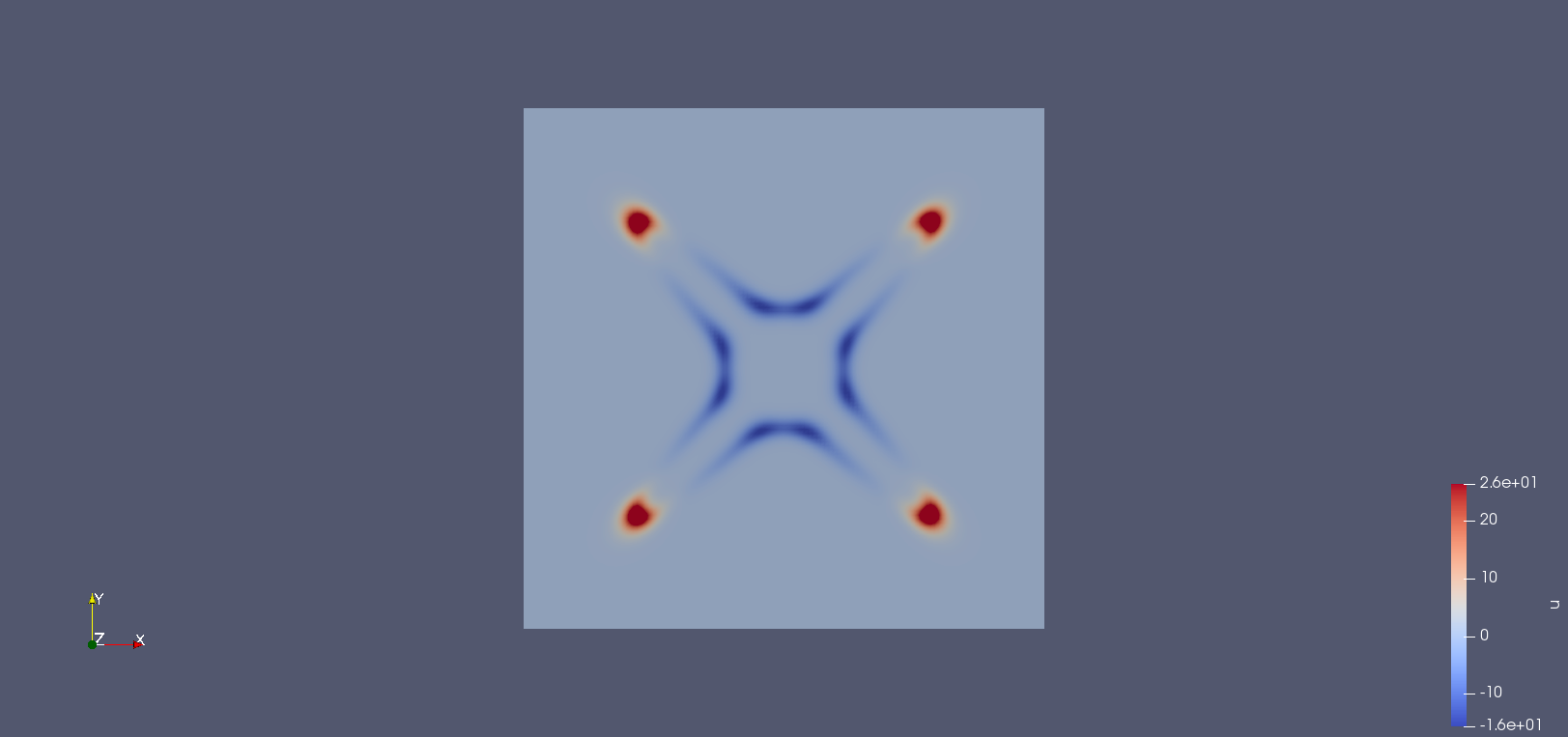}
	\includegraphics[trim={540px 112px 540px 112px}, clip, scale=.08]{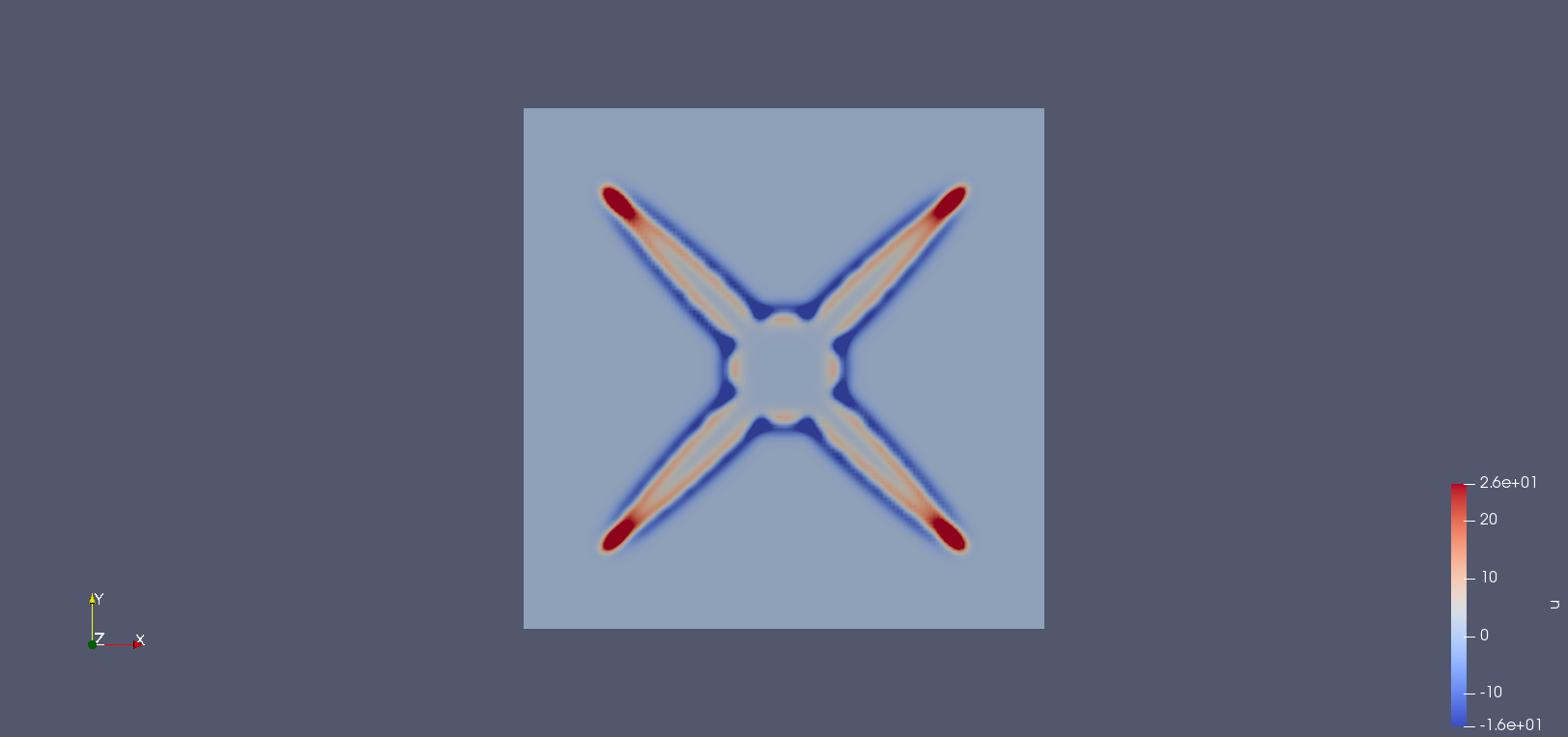}~~~
	\includegraphics[trim={1490px 0px 20px 490px}, clip, scale=.16]{{iso_starl1_control.0000}.png}
\caption{Results for the isotropic case.} 
\label{fig:iso_star4}
\end{center}
\end{subfigure}
\caption{`circle to 4-star' solutions: states in the first and third row and controls in the second and fourth row \label{fig:l1iso_star4}.}
%  \johannes{    $t = timesnotlookedupsofar$.}
% Hinweis: die Zahlen entsprechen ungefähr diesen Zeitschritten, ich hab die Bilder von 0-7 darauf umgerechnet und auf "schöne" Zahlen gerundet
\end{center}
\end{figure}
\fi

\ifgraphics
\begin{figure}[htbp]
  \begin{center}
\begin{subfigure}{1.\textwidth}
\begin{center}    
	\includegraphics[trim={540px 112px 540px 112px}, clip, scale=.08]{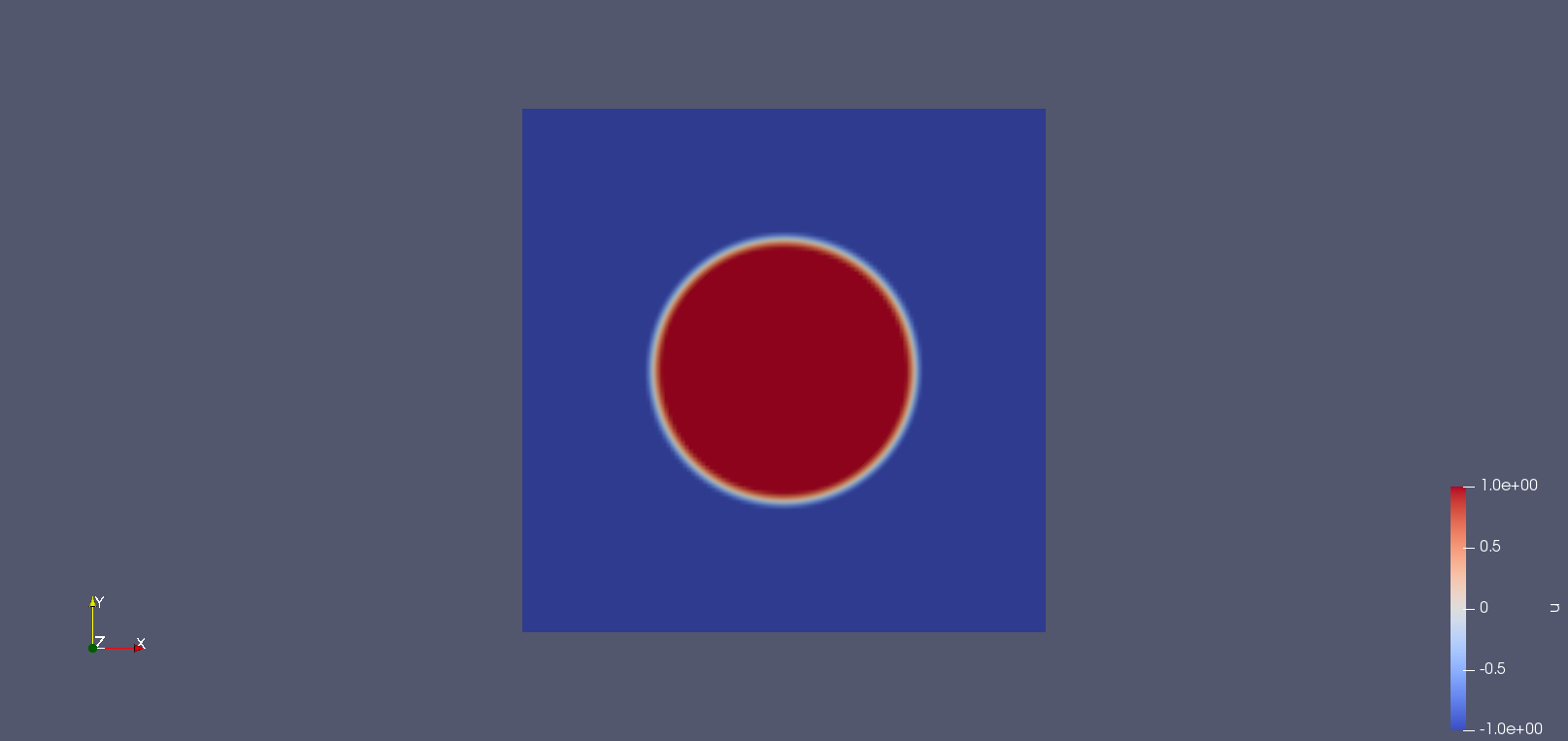}
	\includegraphics[trim={540px 112px 540px 112px}, clip, scale=.08]{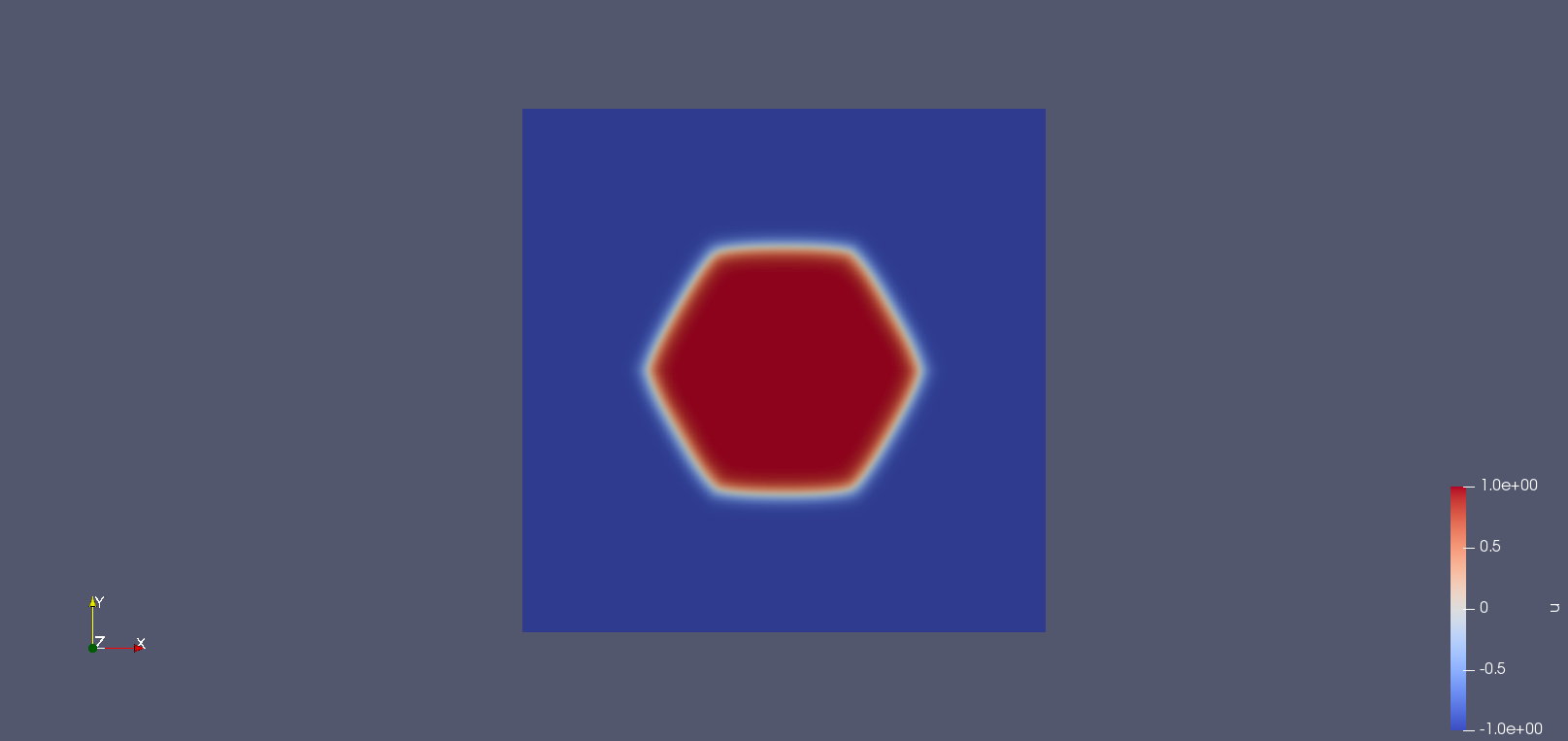}
	\includegraphics[trim={540px 112px 540px 112px}, clip, scale=.08]{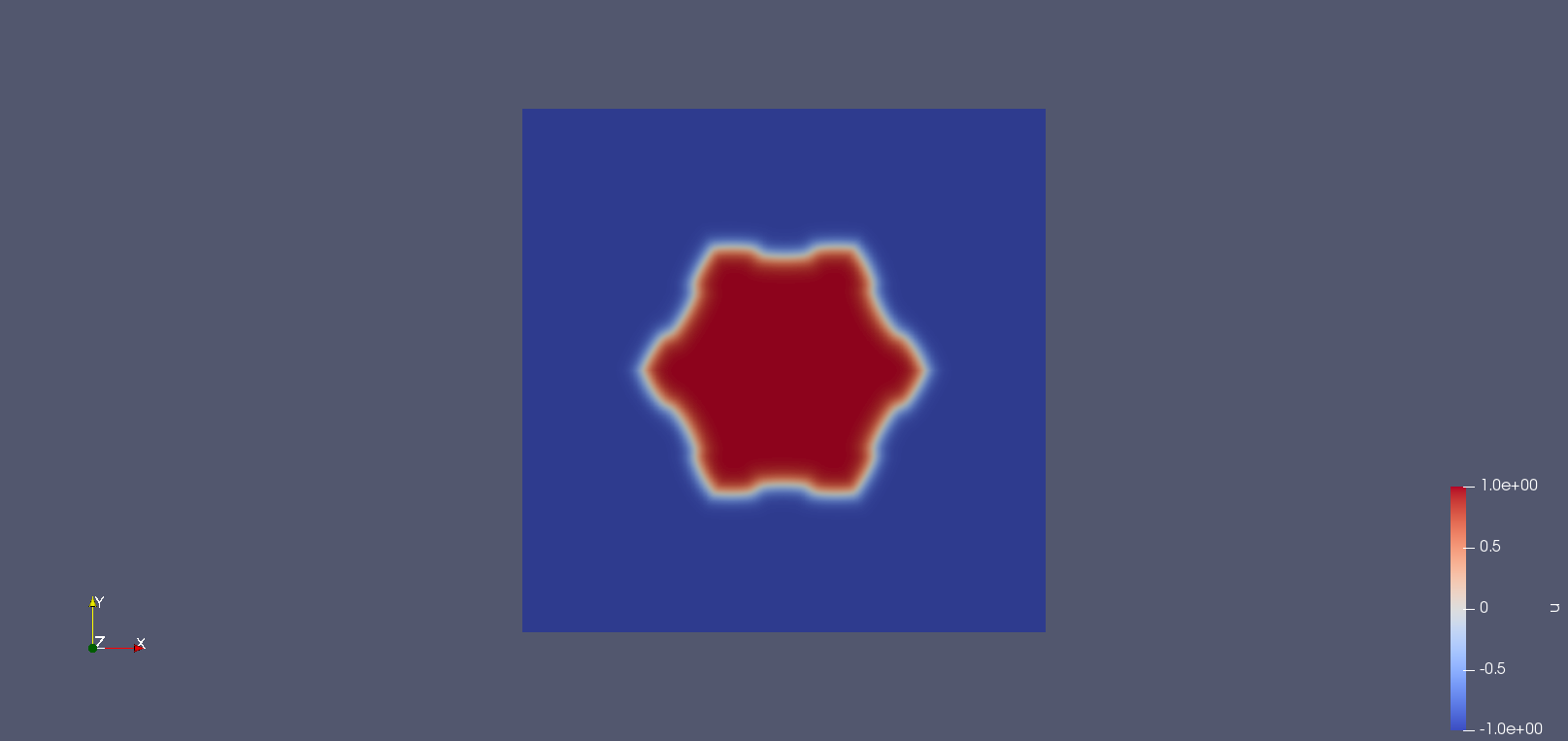}
	\includegraphics[trim={540px 112px 540px 112px}, clip, scale=.08]{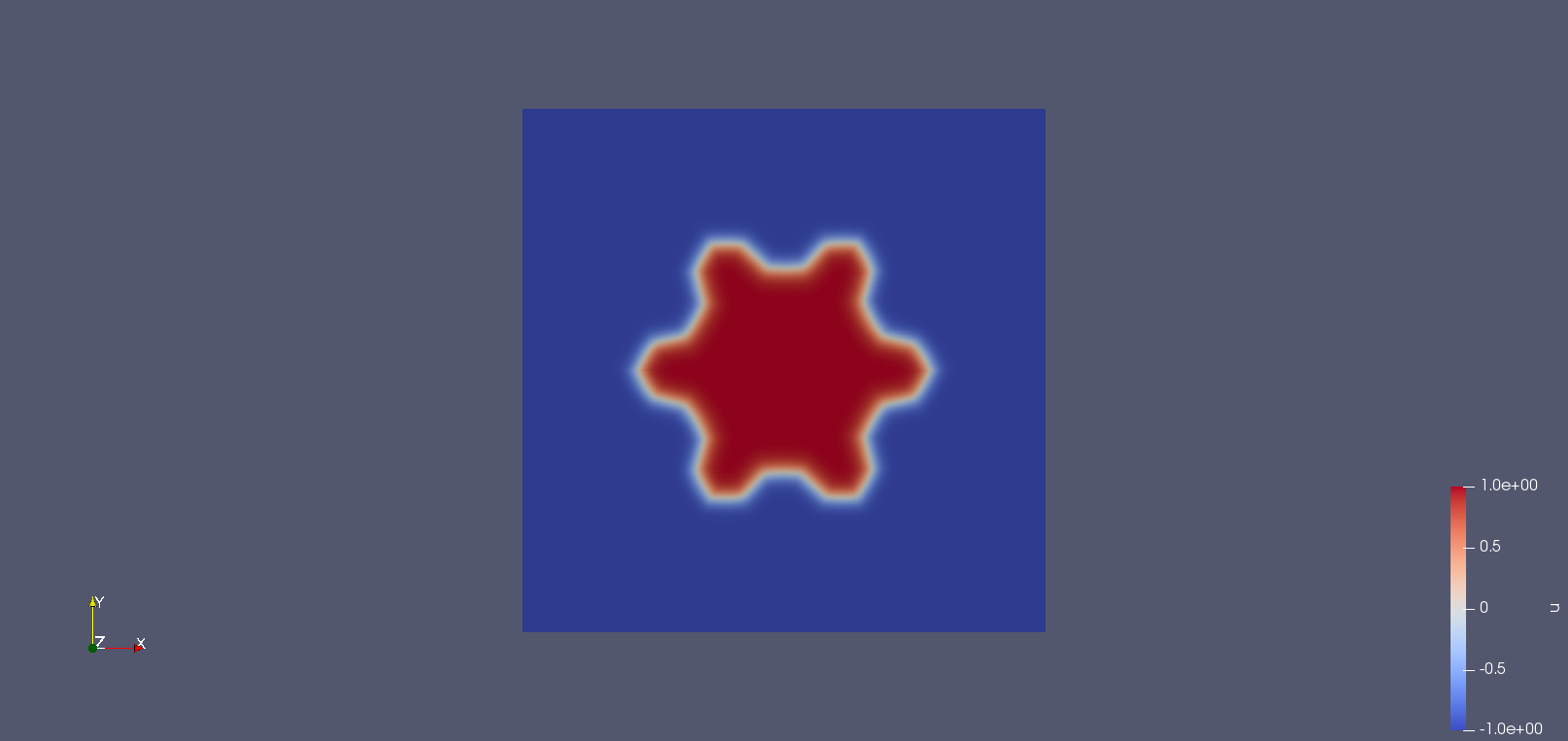}
	\includegraphics[trim={540px 112px 540px 112px}, clip, scale=.08]{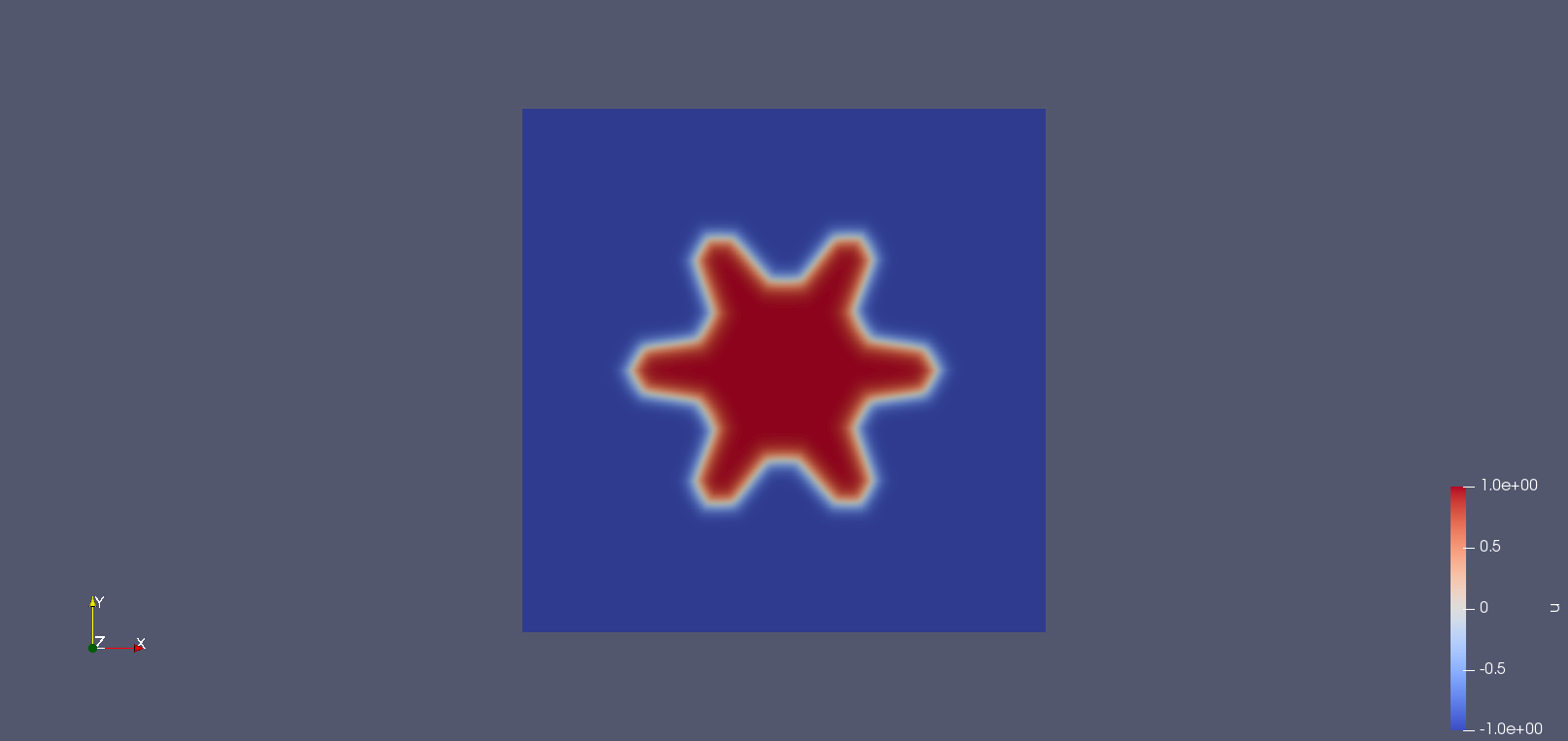}
	\includegraphics[trim={540px 112px 540px 112px}, clip, scale=.08]{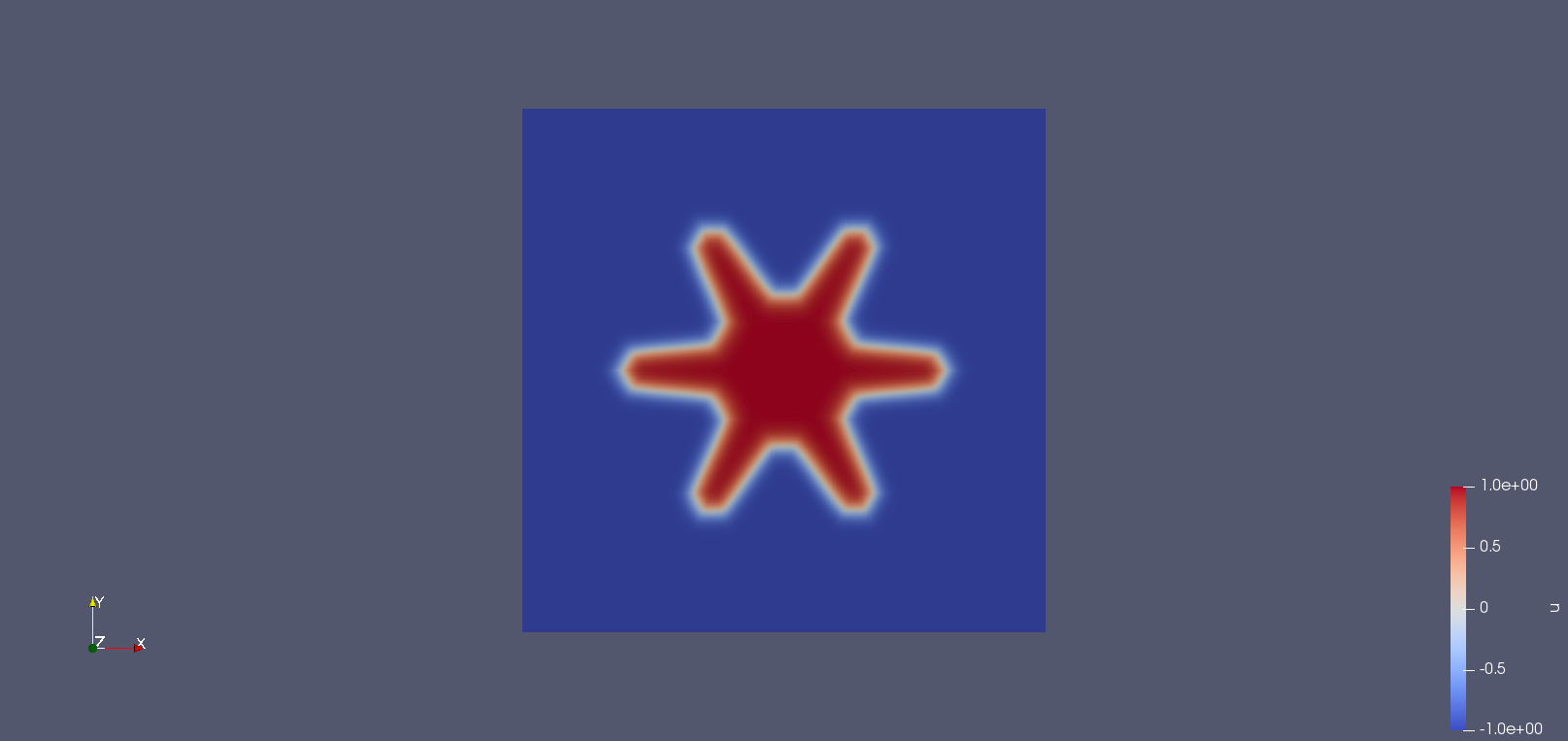}
	\includegraphics[trim={540px 112px 540px 112px}, clip, scale=.08]{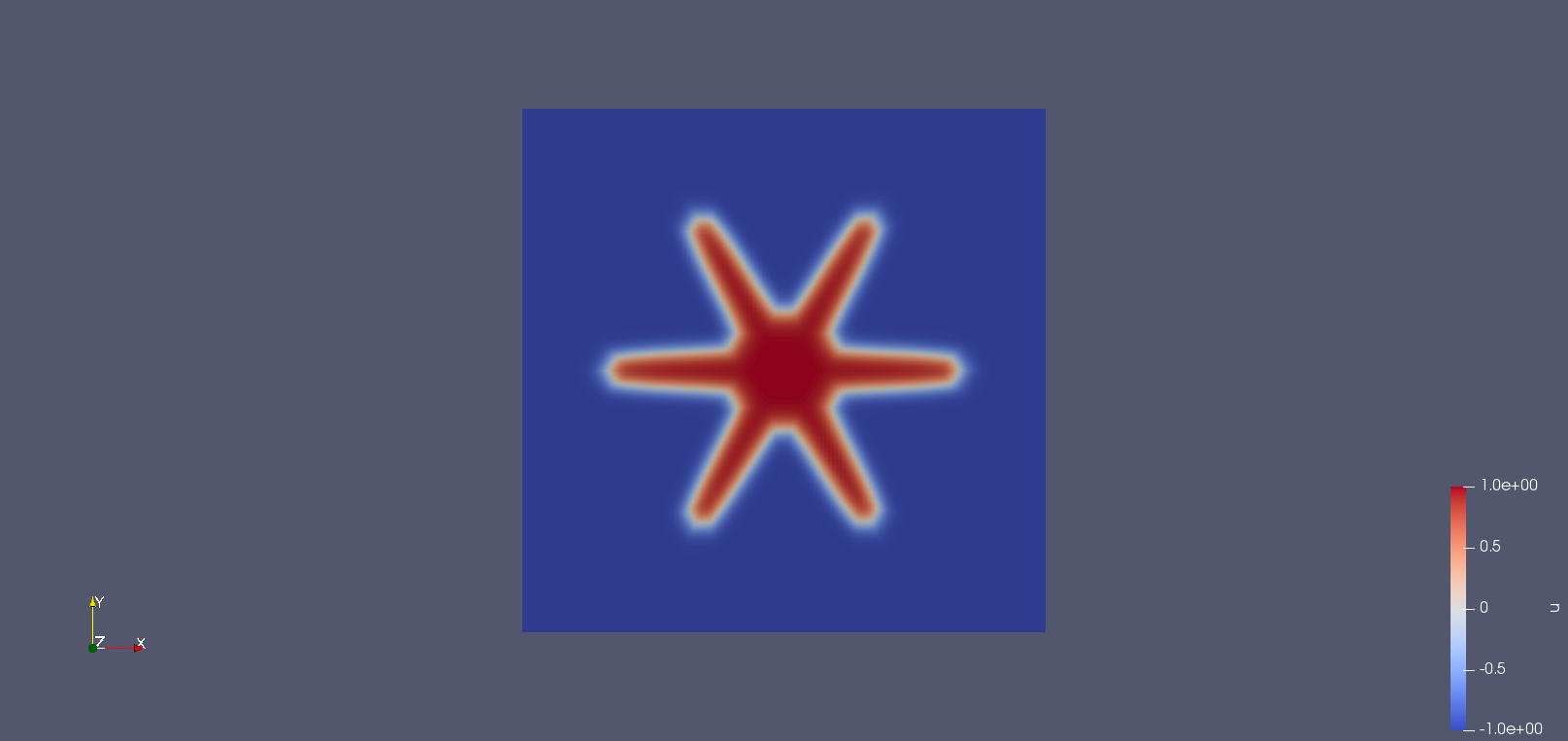}
	\includegraphics[trim={540px 112px 540px 112px}, clip, scale=.08]{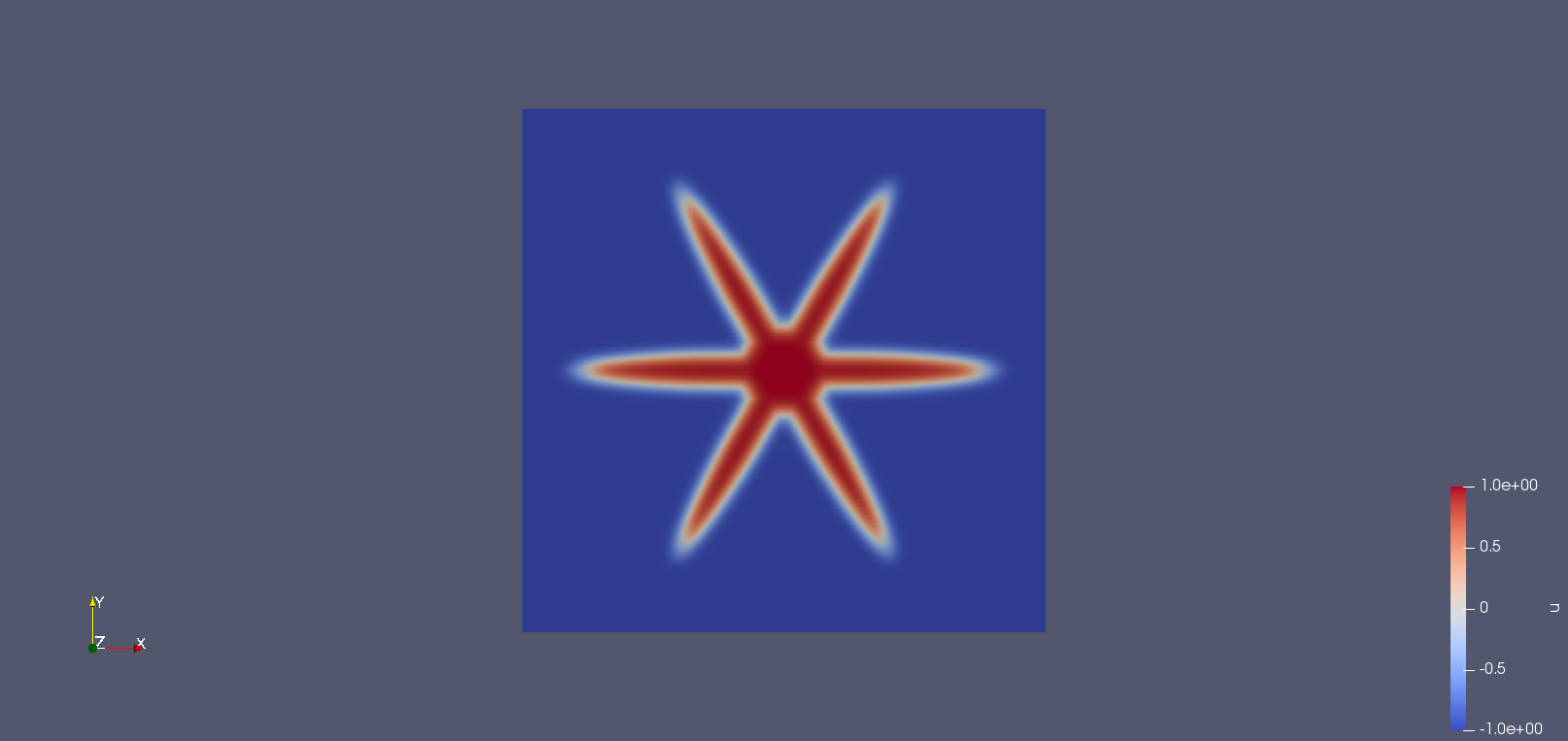}~~~
	\includegraphics[trim={1490px 0px 20px 490px}, clip, scale=.16]{{hexa_star_state.0000}.png}~\\~\\
	\includegraphics[trim={540px 112px 540px 112px}, clip, scale=.08]{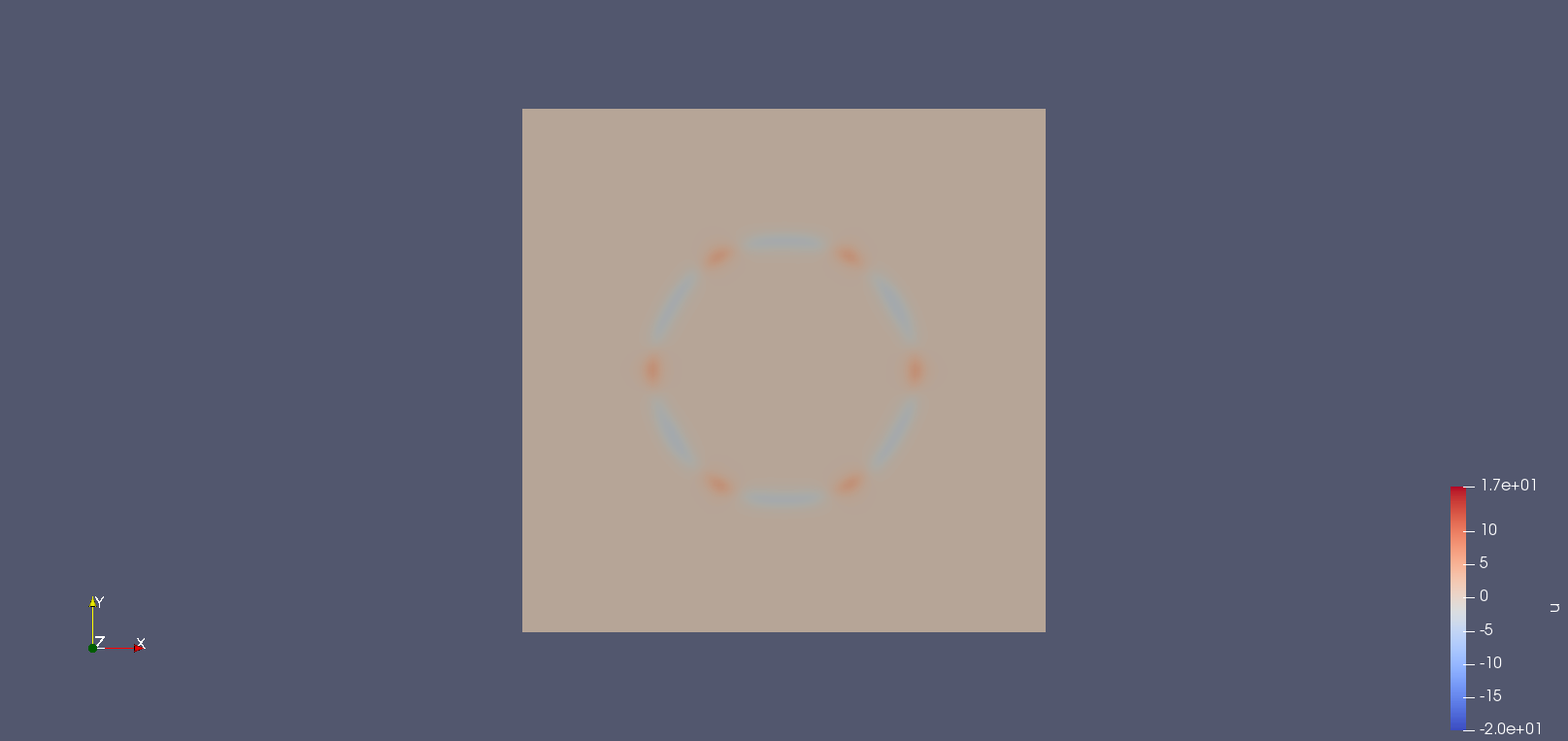}
	\includegraphics[trim={540px 112px 540px 112px}, clip, scale=.08]{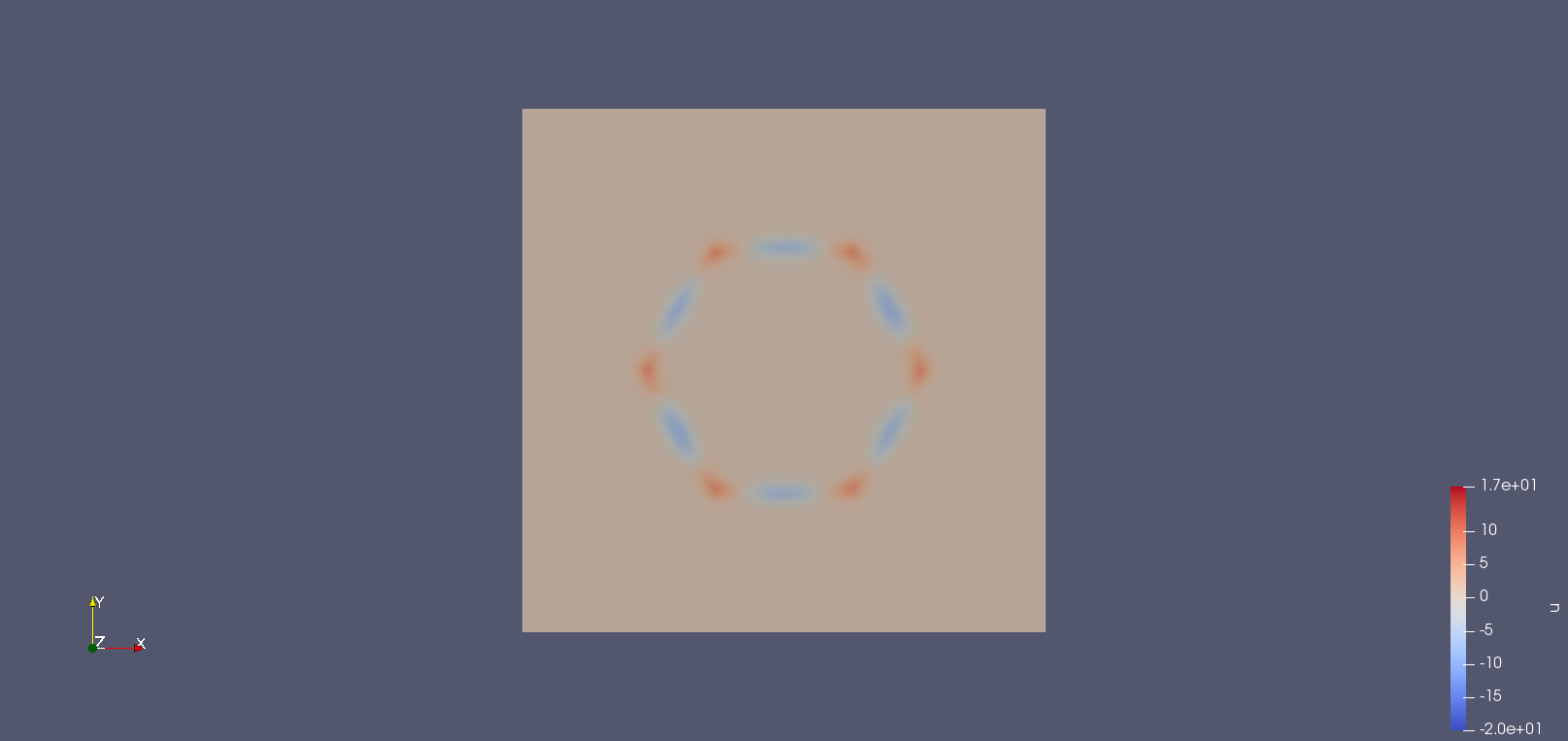}
	\includegraphics[trim={540px 112px 540px 112px}, clip, scale=.08]{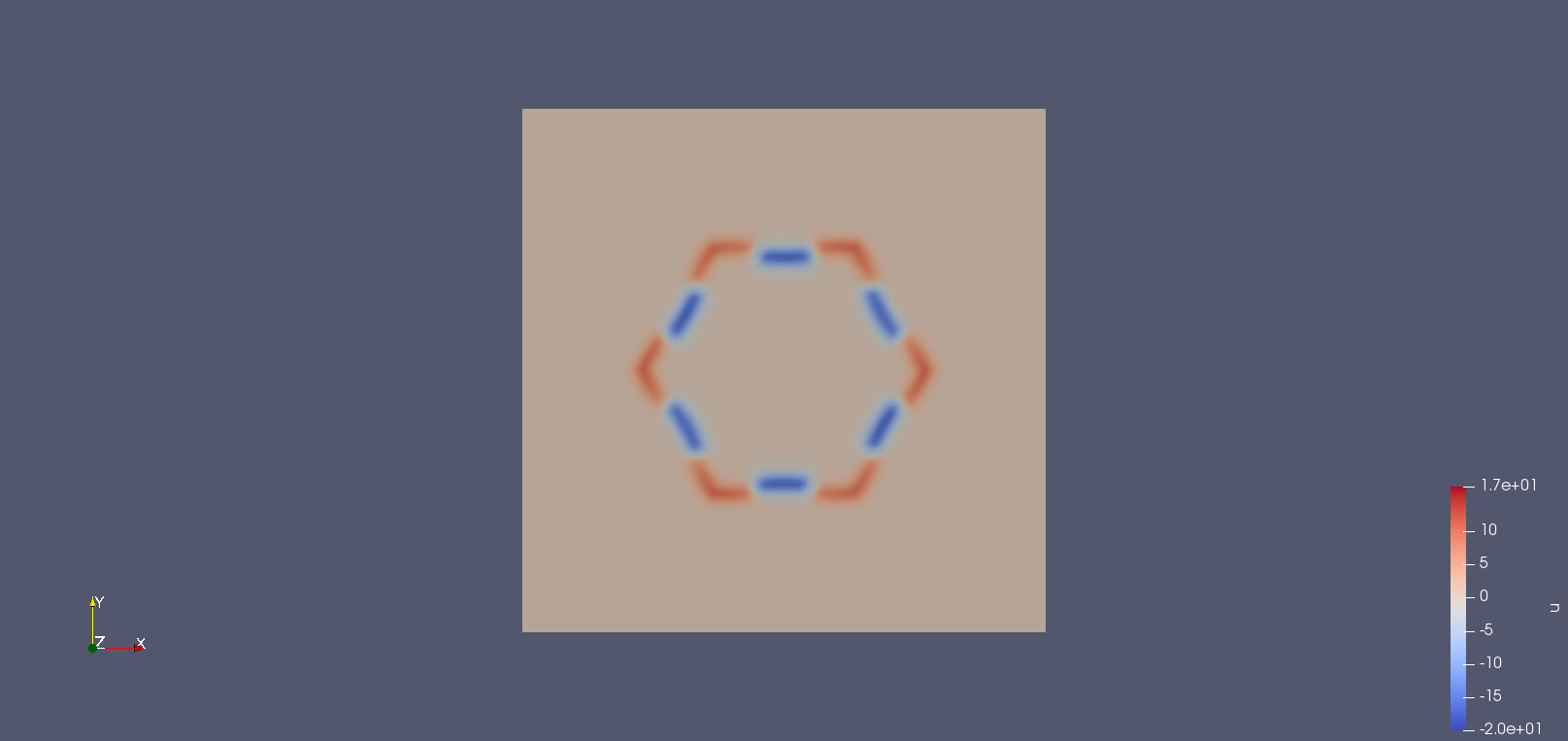}
	\includegraphics[trim={540px 112px 540px 112px}, clip, scale=.08]{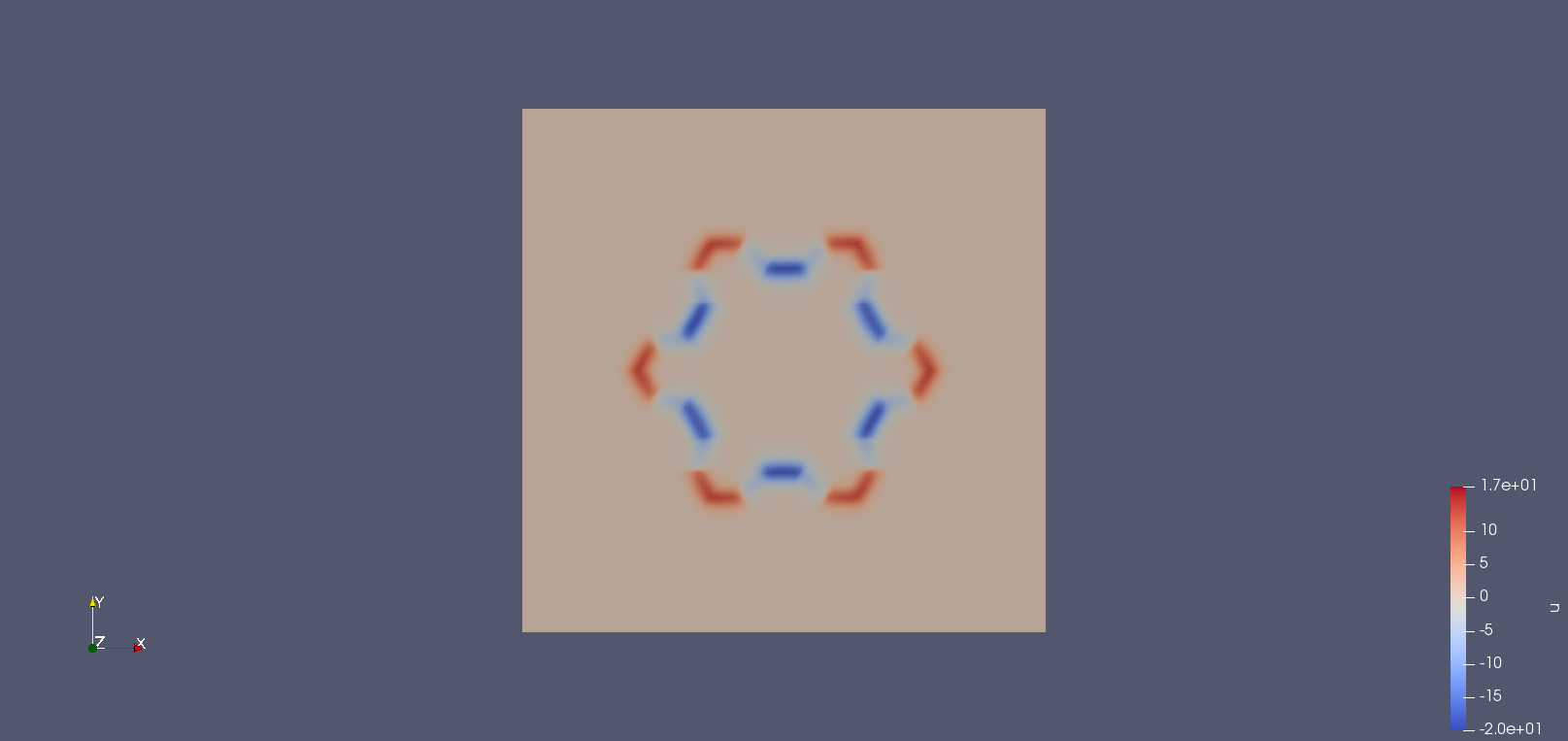}
	\includegraphics[trim={540px 112px 540px 112px}, clip, scale=.08]{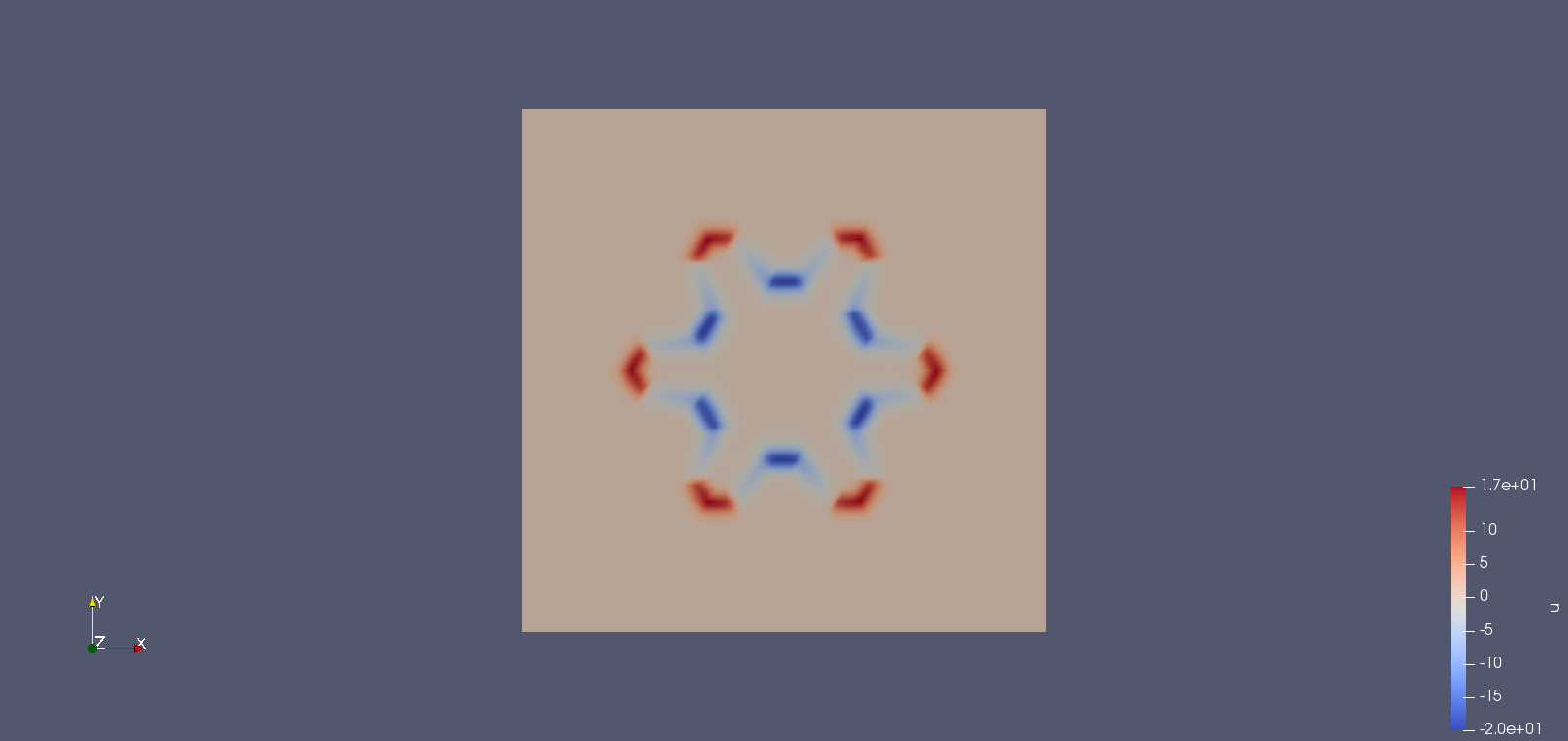}
	\includegraphics[trim={540px 112px 540px 112px}, clip, scale=.08]{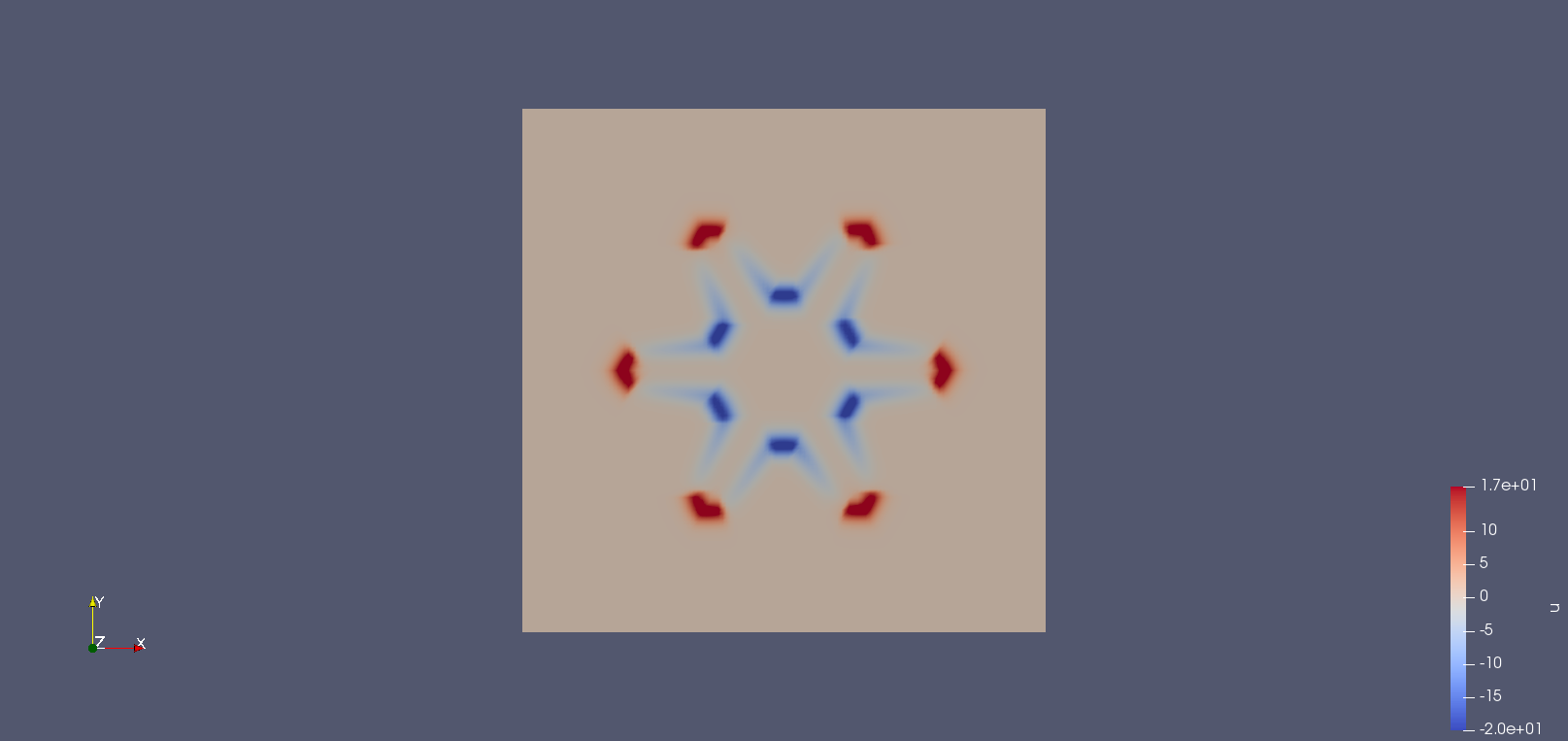}
	\includegraphics[trim={540px 112px 540px 112px}, clip, scale=.08]{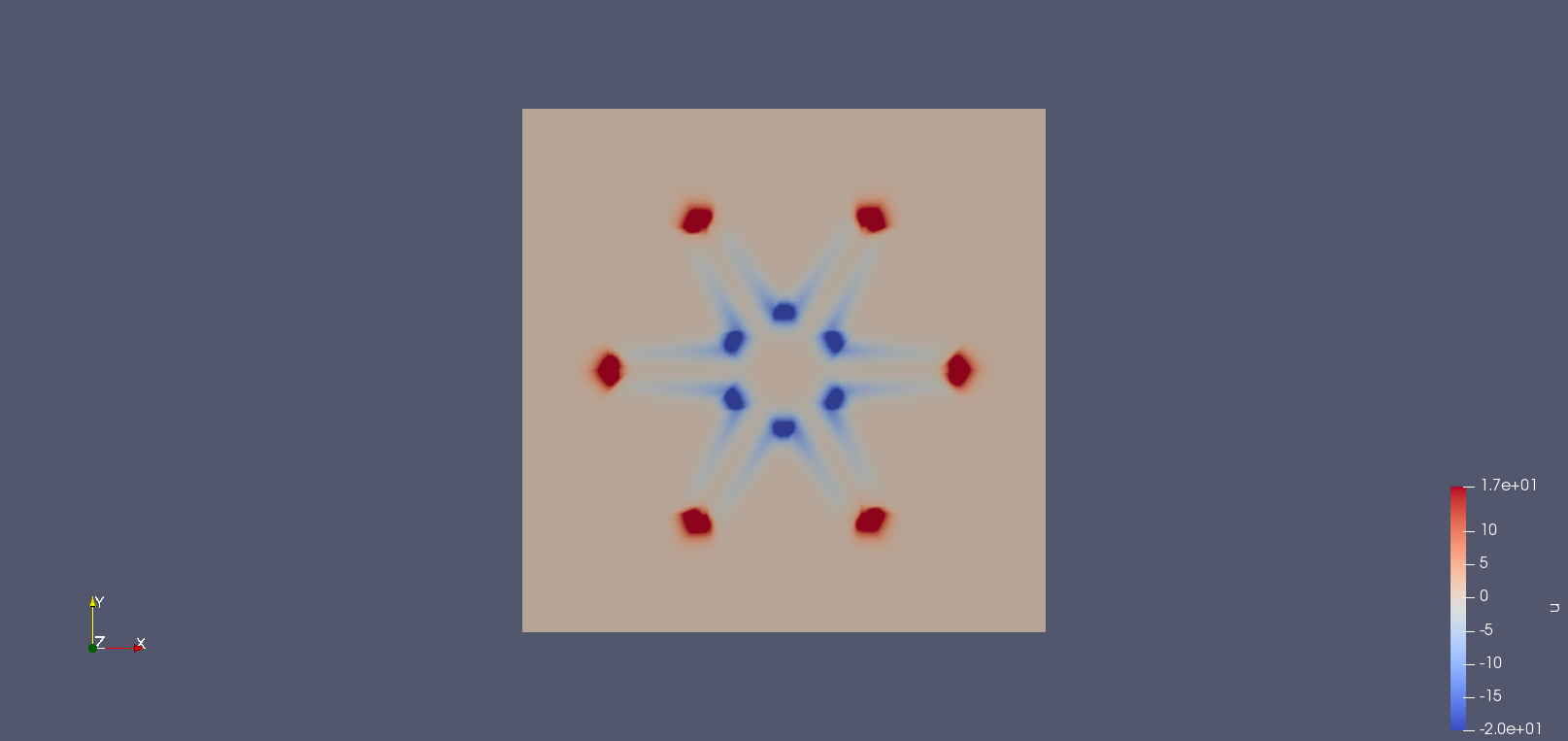}
	\includegraphics[trim={540px 112px 540px 112px}, clip, scale=.08]{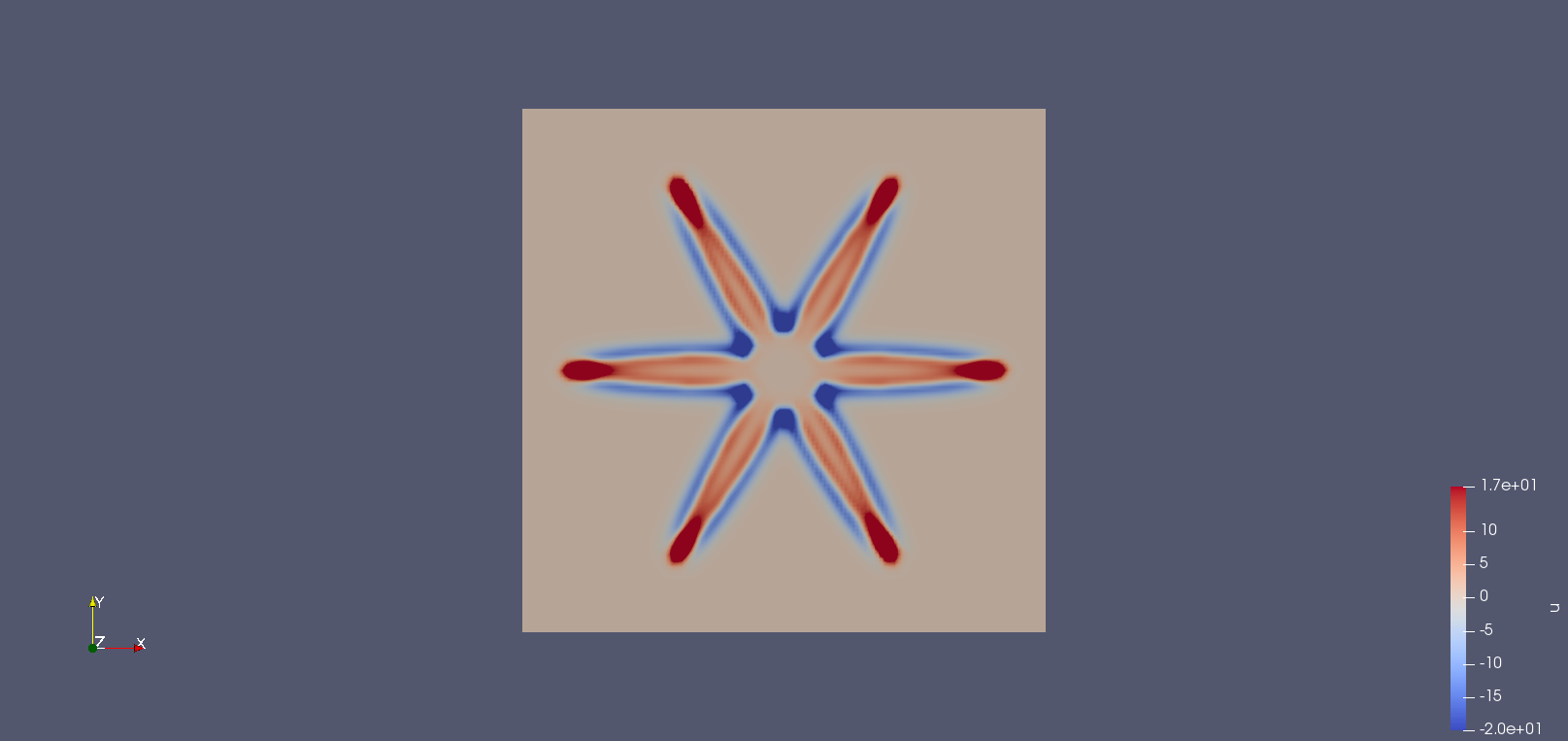}~~~
	\includegraphics[trim={1490px 0px 20px 490px}, clip, scale=.16]{{hexa_star_control.0000}.png}
\caption{Results for the hexagon case.} 
\label{fig:hexa_star6}
\end{center}
\end{subfigure}
\begin{subfigure}{1.\textwidth}
  \begin{center}
    \includegraphics[trim={540px 112px 540px 112px}, clip, scale=.08]{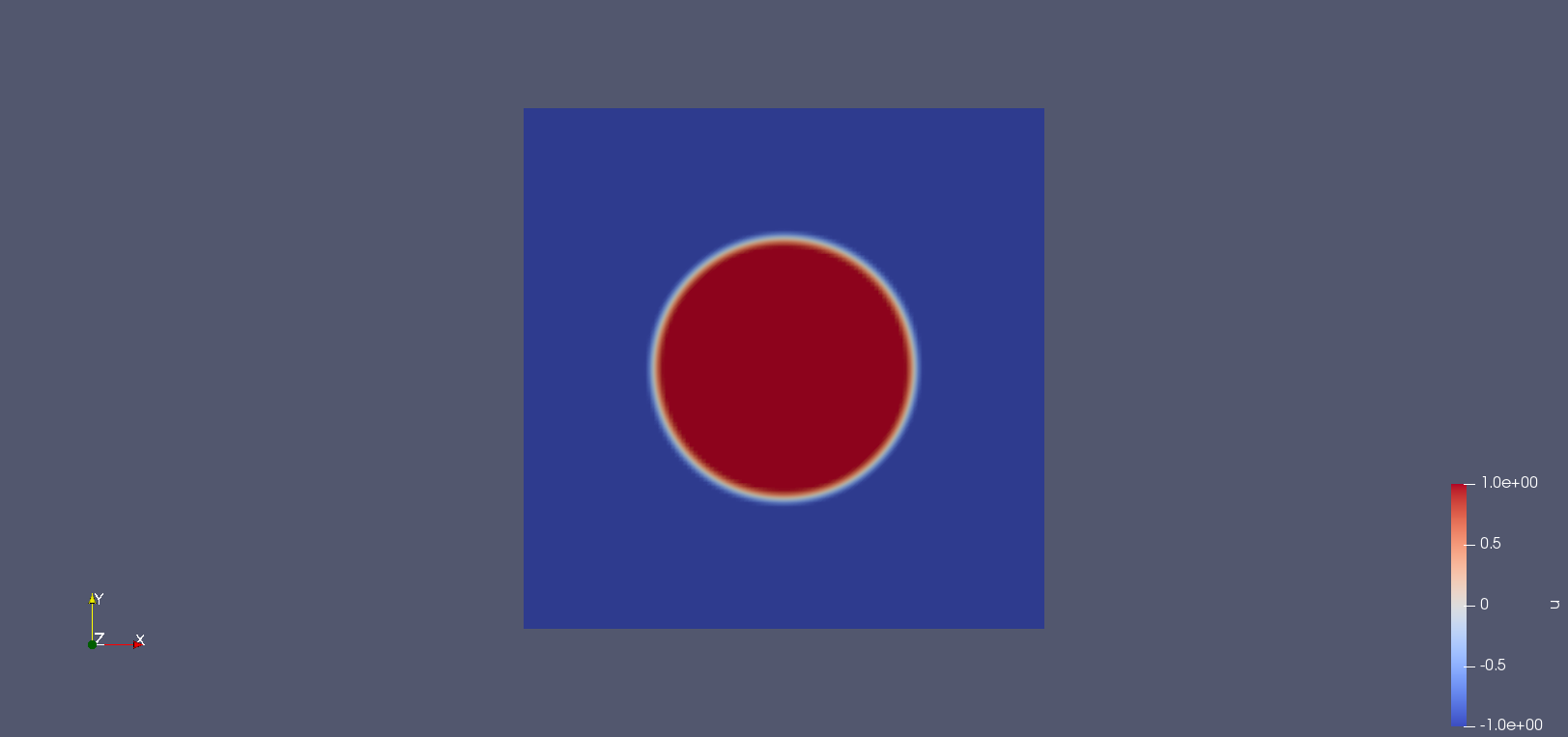}
	\includegraphics[trim={540px 112px 540px 112px}, clip, scale=.08]{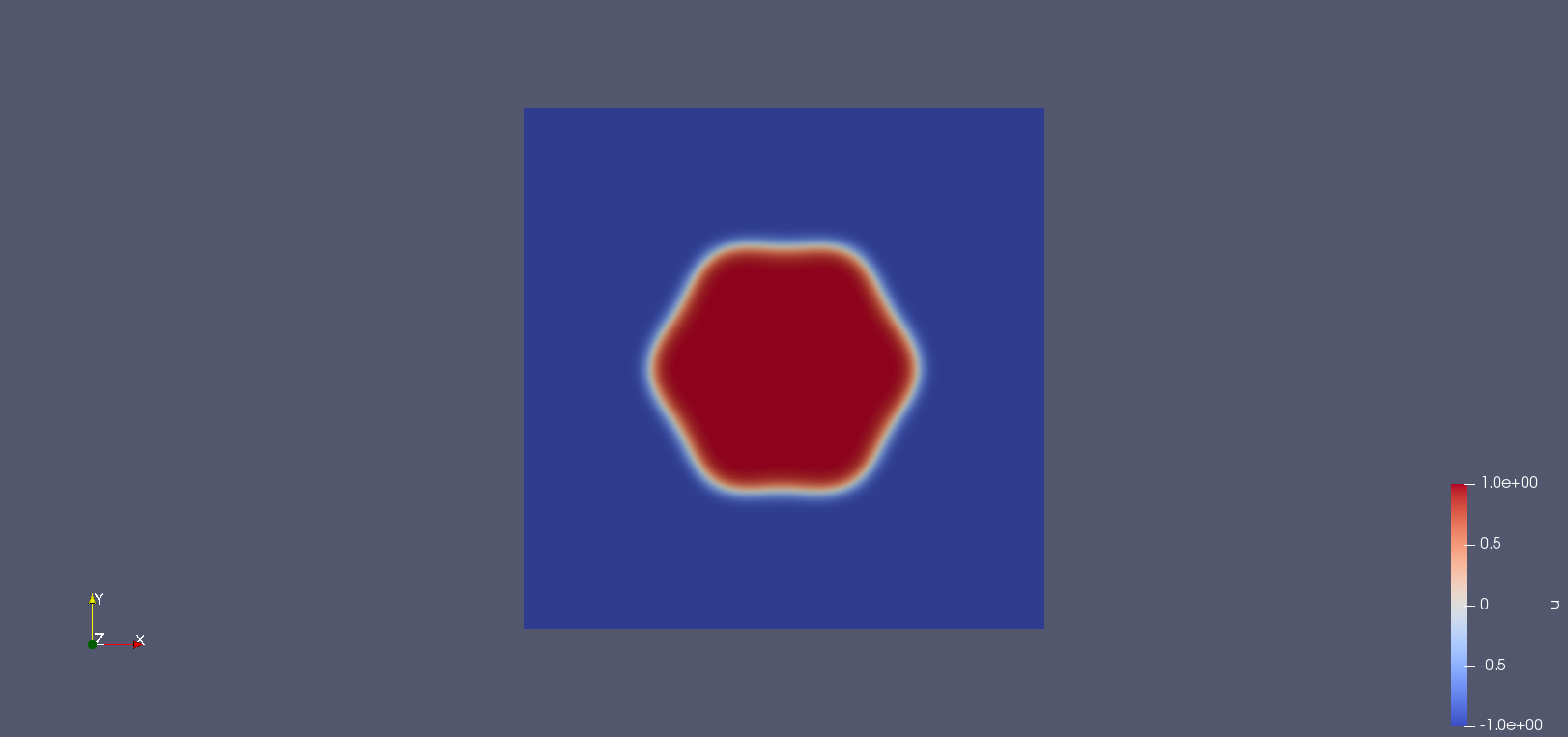}
	\includegraphics[trim={540px 112px 540px 112px}, clip, scale=.08]{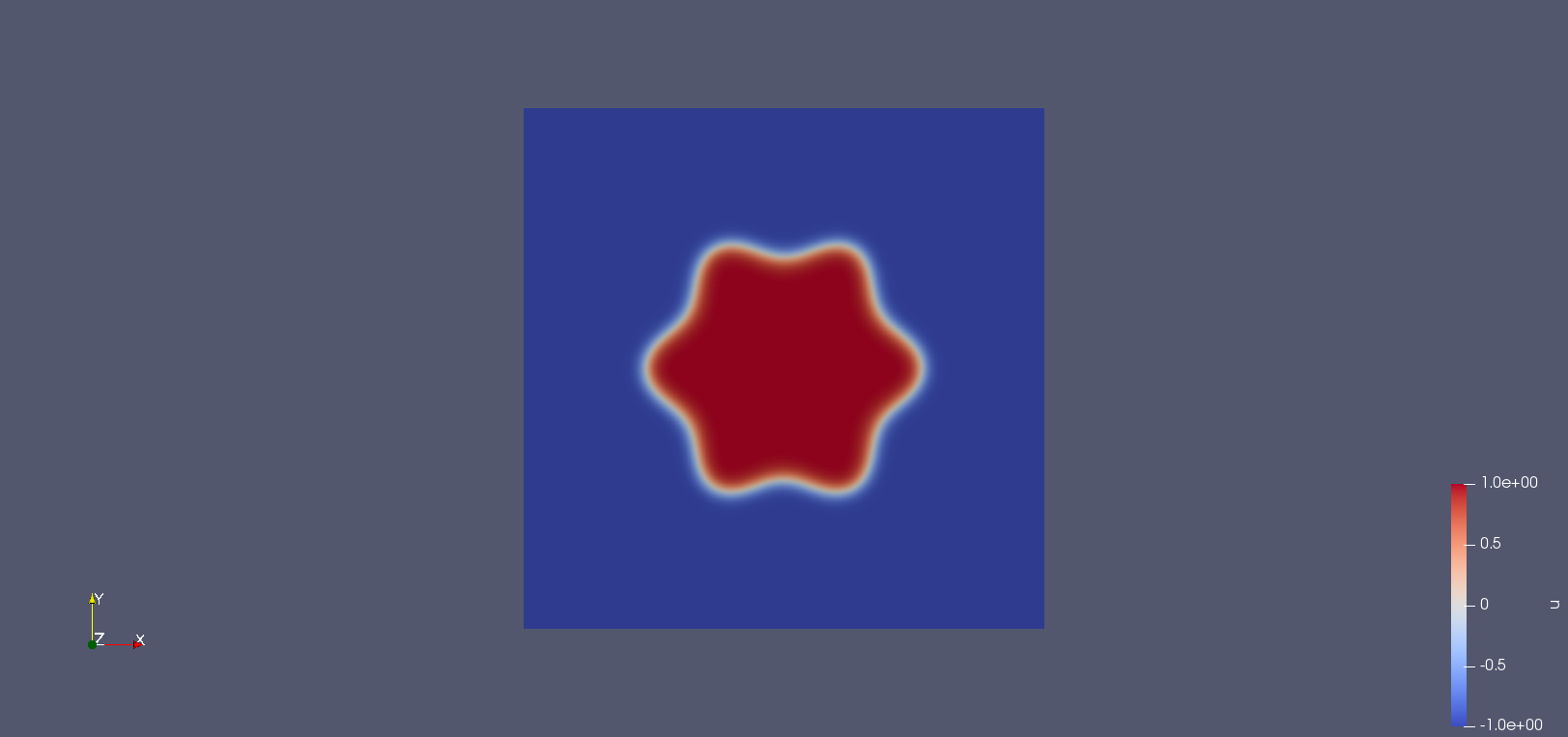}
	\includegraphics[trim={540px 112px 540px 112px}, clip, scale=.08]{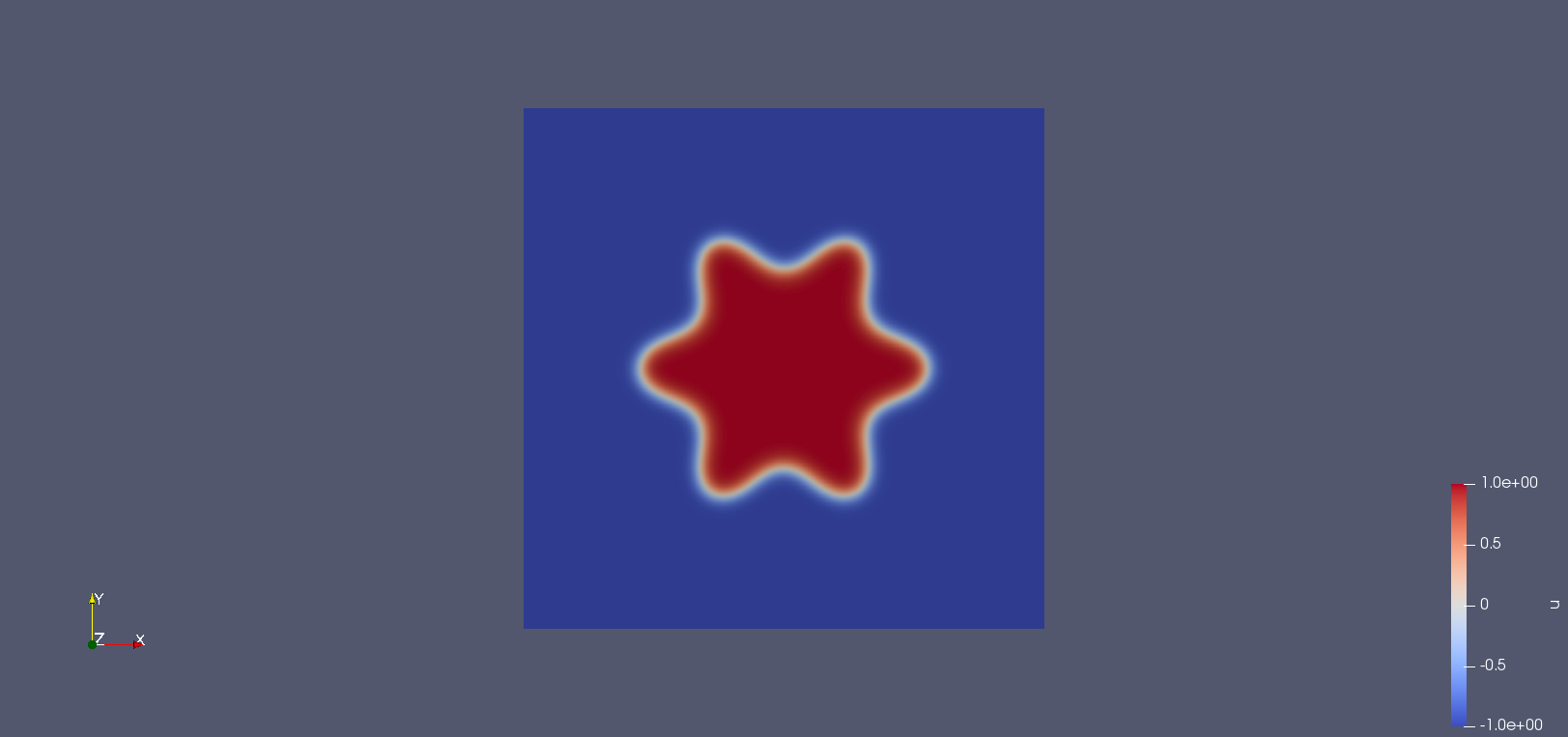}
	\includegraphics[trim={540px 112px 540px 112px}, clip, scale=.08]{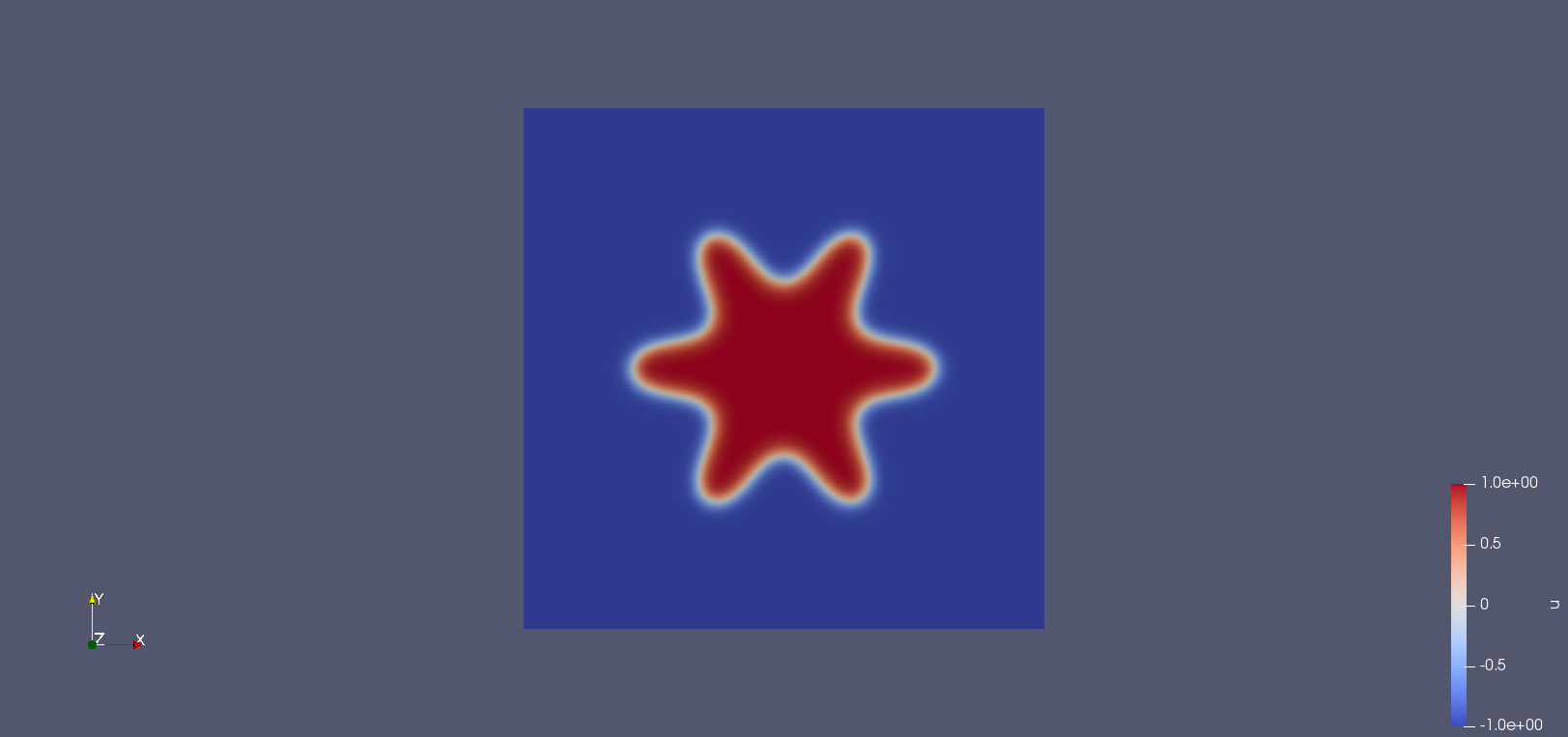}
	\includegraphics[trim={540px 112px 540px 112px}, clip, scale=.08]{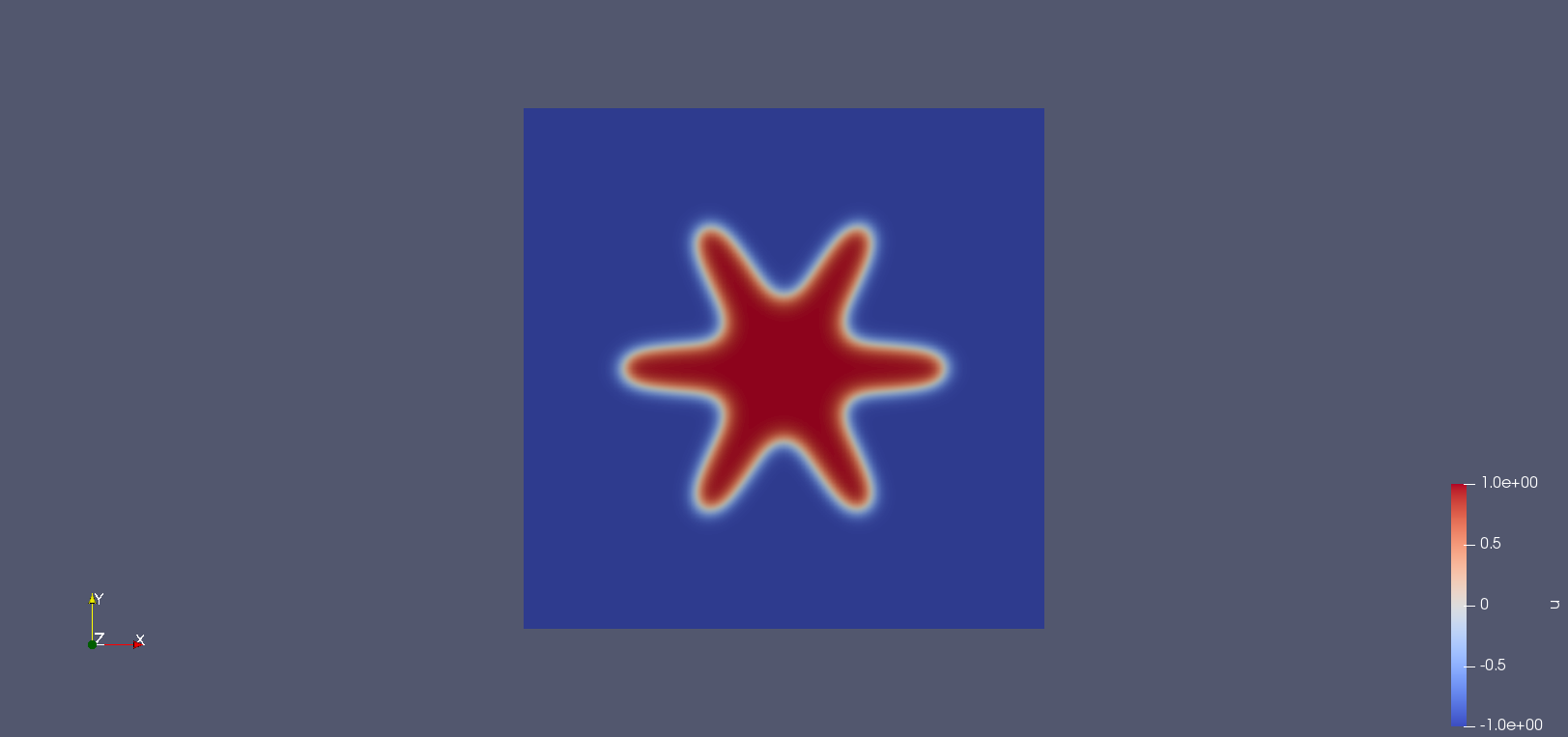}
	\includegraphics[trim={540px 112px 540px 112px}, clip, scale=.08]{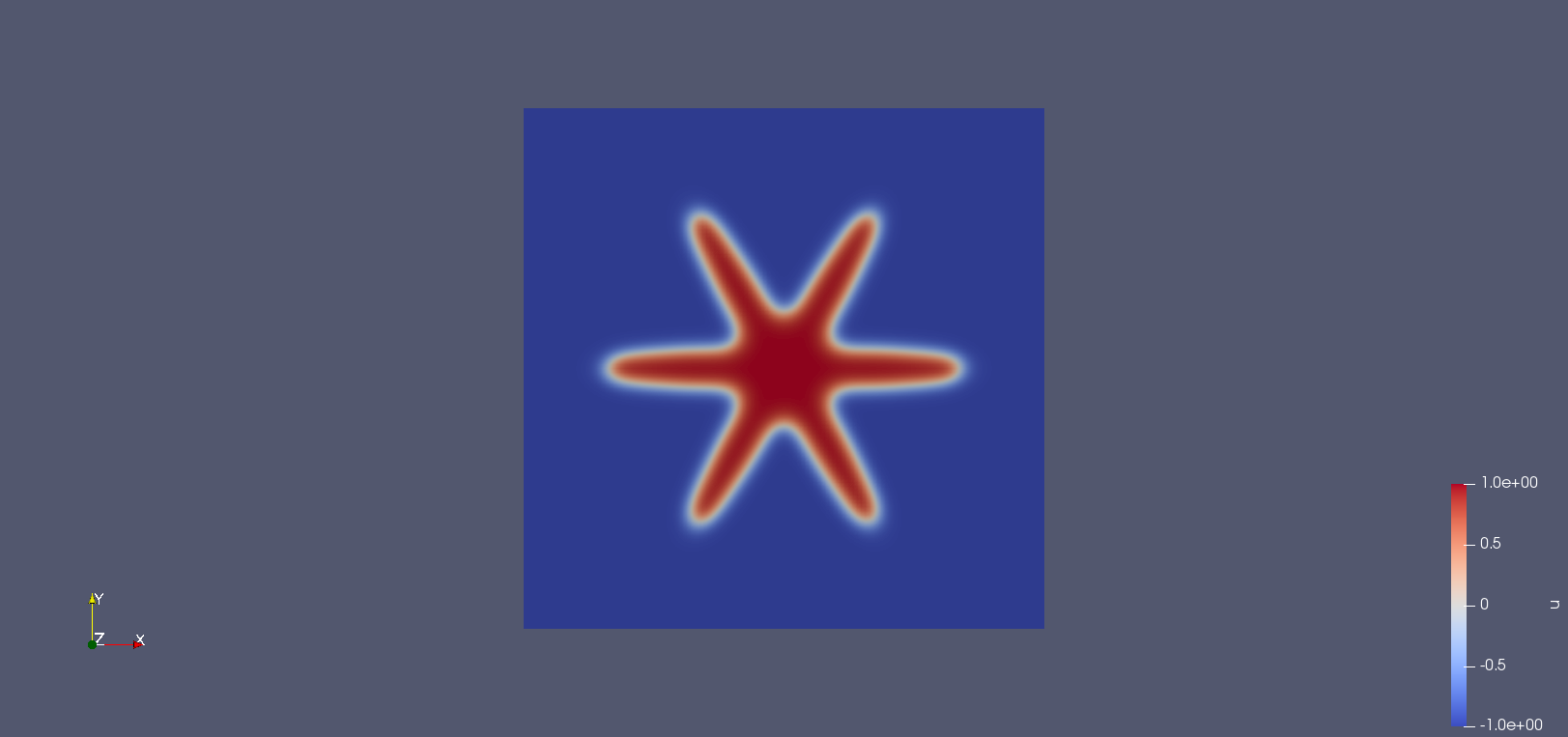}
	\includegraphics[trim={540px 112px 540px 112px}, clip, scale=.08]{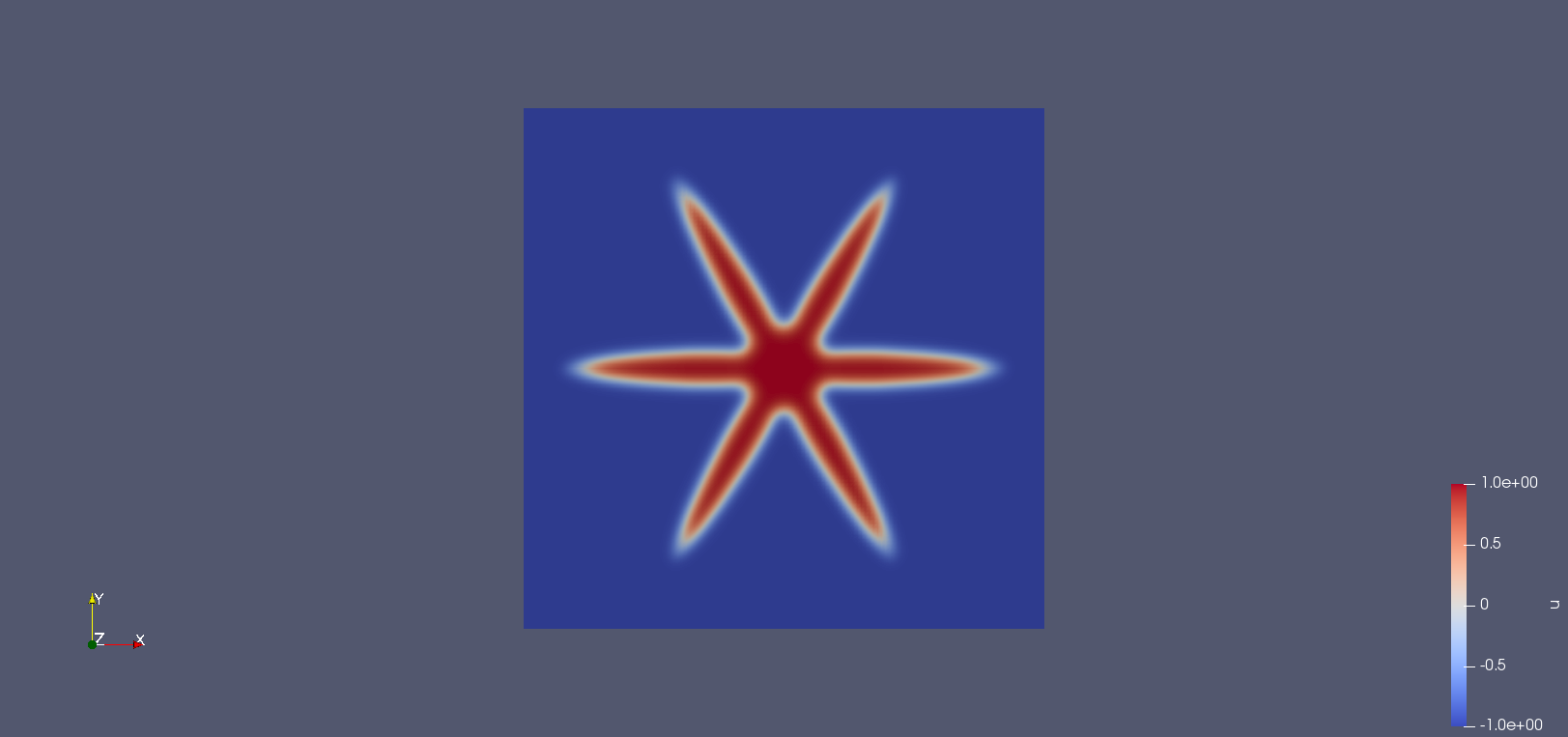}~~~
	\includegraphics[trim={1490px 0px 20px 490px}, clip, scale=.16]{{iso_starhexa_state.0000}.png}~\\~\\
	\includegraphics[trim={540px 112px 540px 112px}, clip, scale=.08]{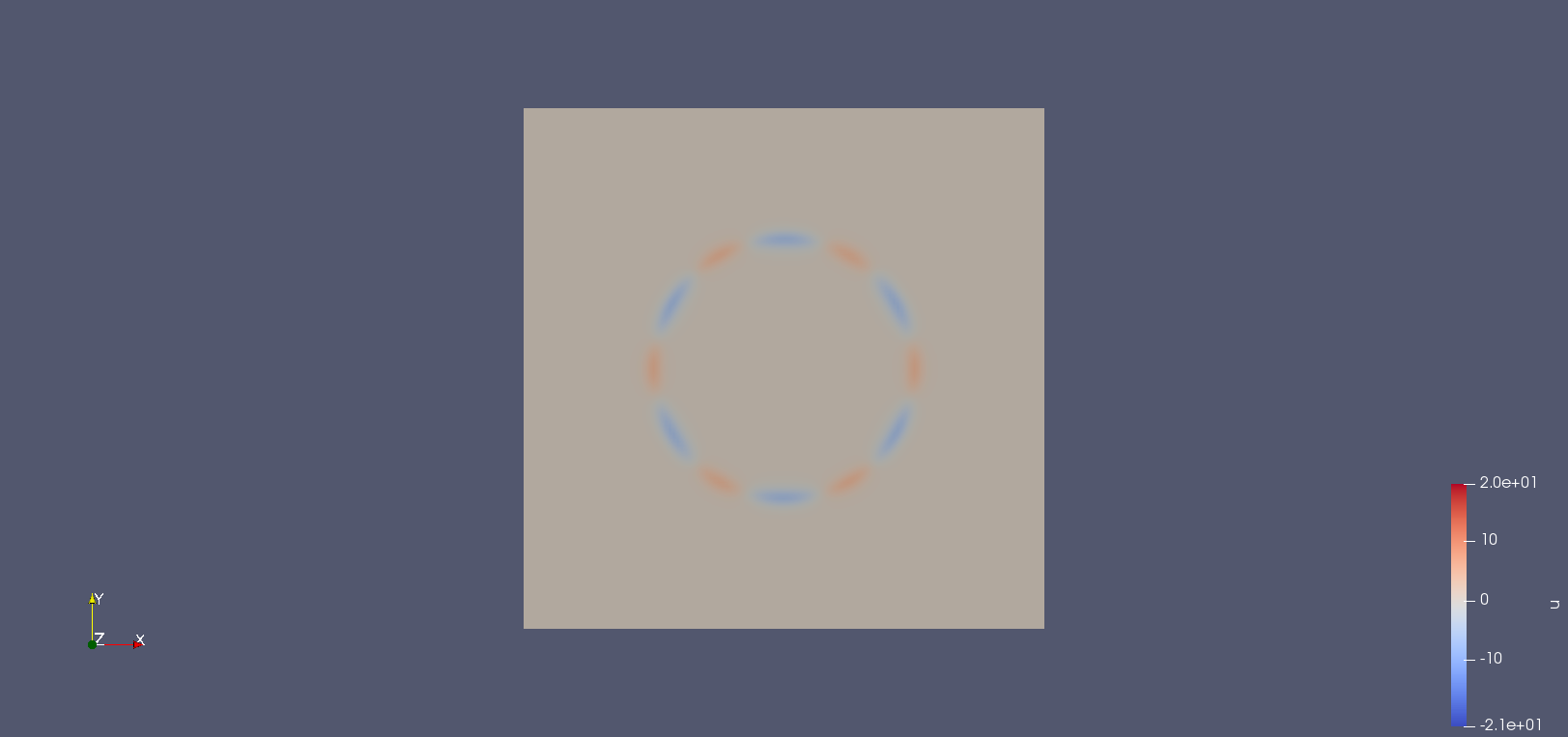}
	\includegraphics[trim={540px 112px 540px 112px}, clip, scale=.08]{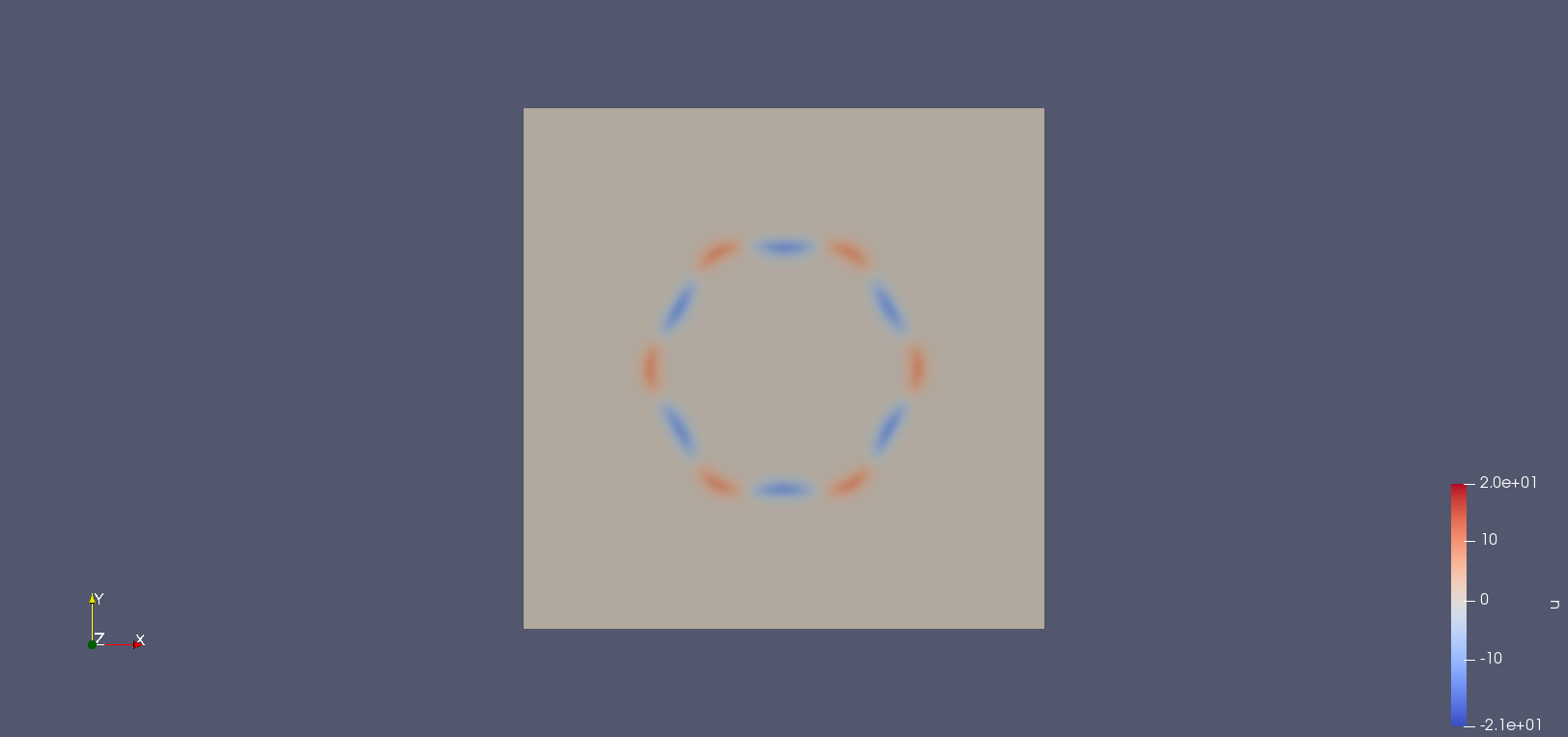}
	\includegraphics[trim={540px 112px 540px 112px}, clip, scale=.08]{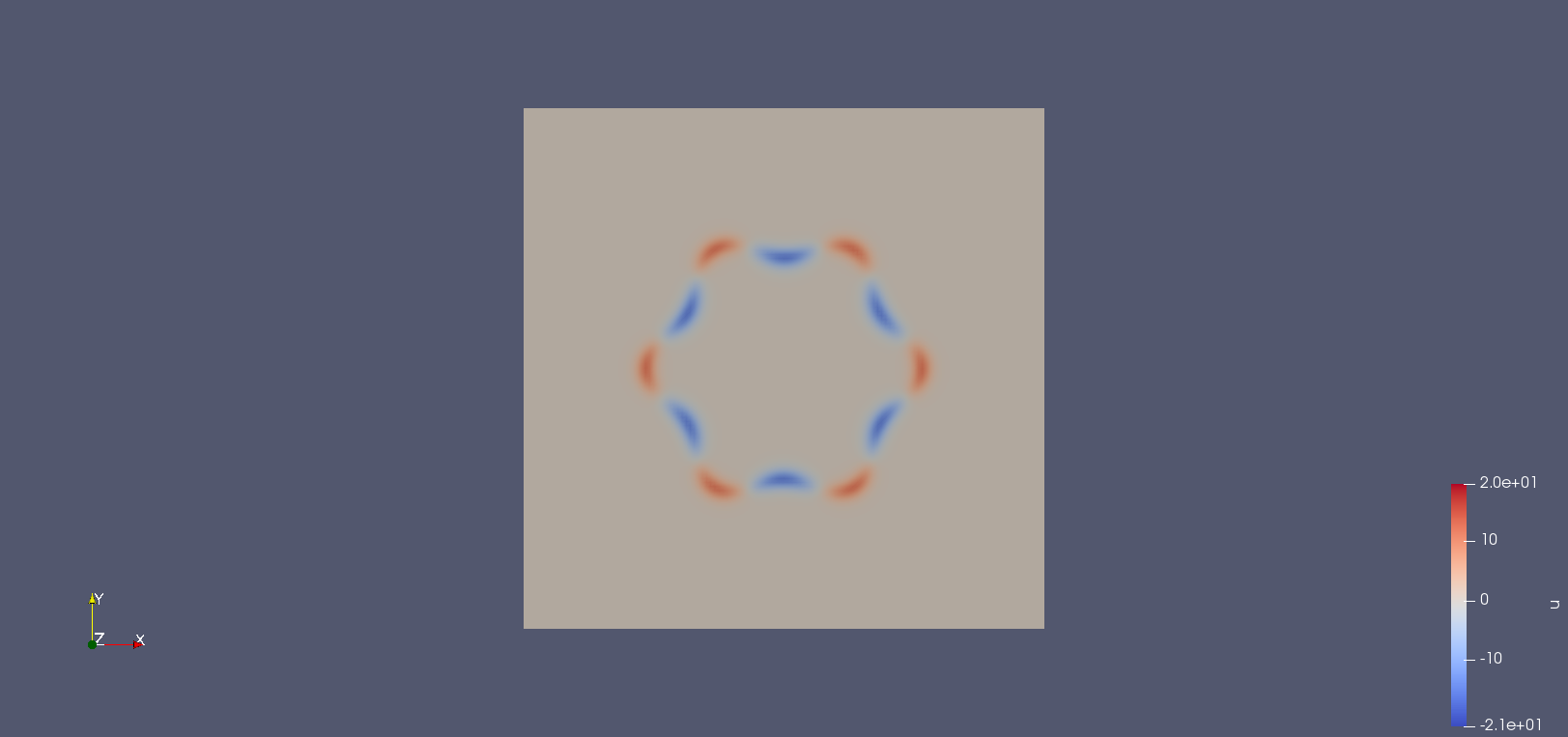}
	\includegraphics[trim={540px 112px 540px 112px}, clip, scale=.08]{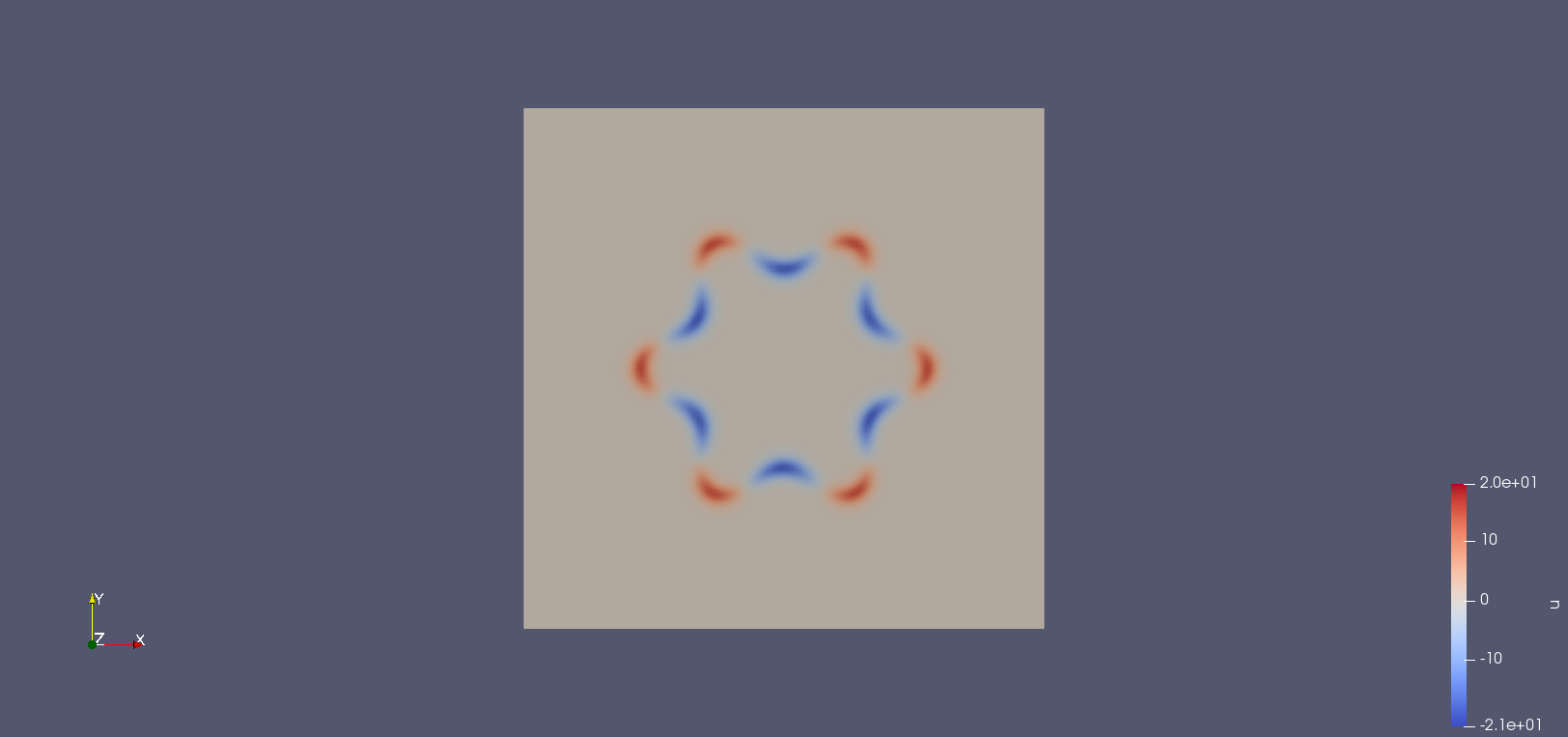}
	\includegraphics[trim={540px 112px 540px 112px}, clip, scale=.08]{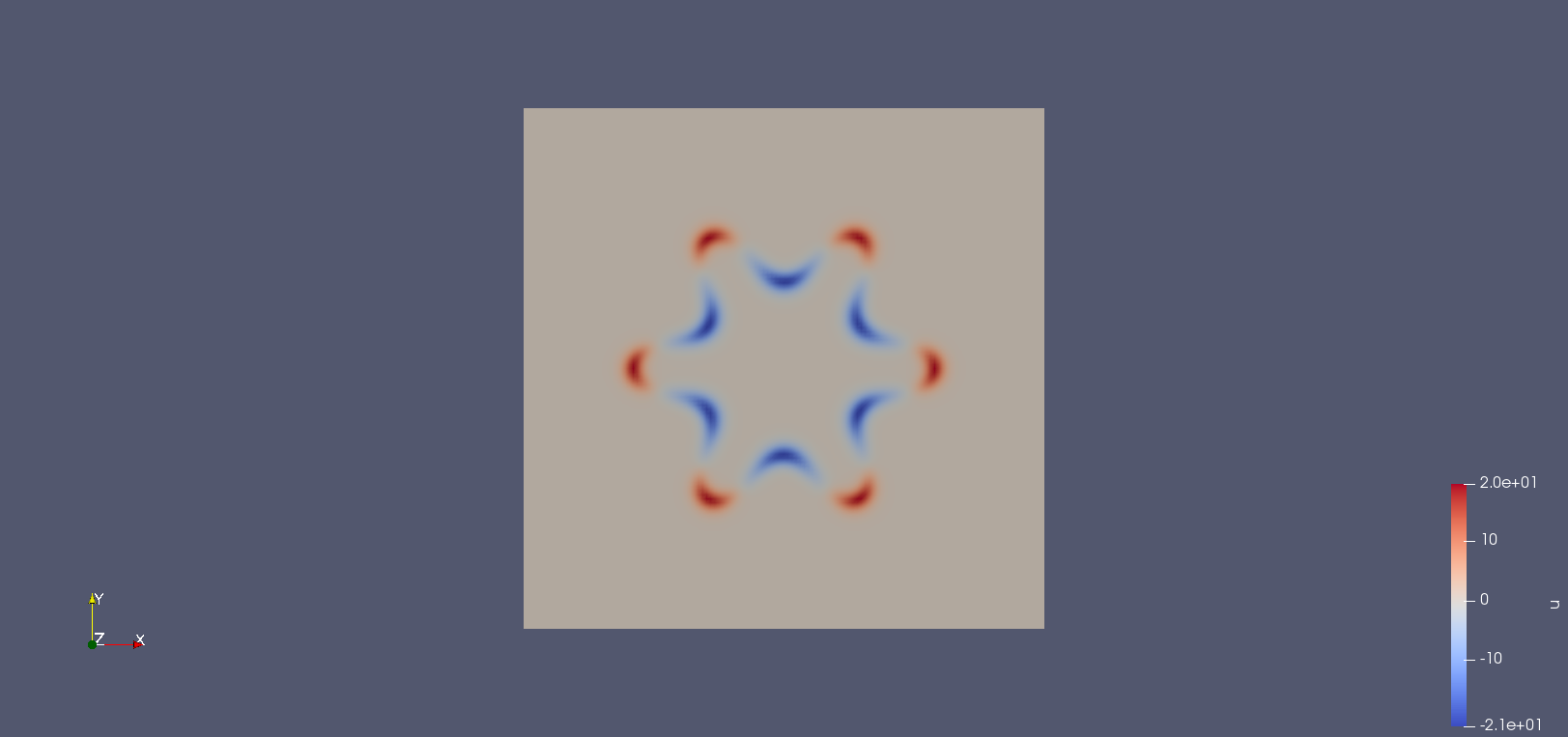}
	\includegraphics[trim={540px 112px 540px 112px}, clip, scale=.08]{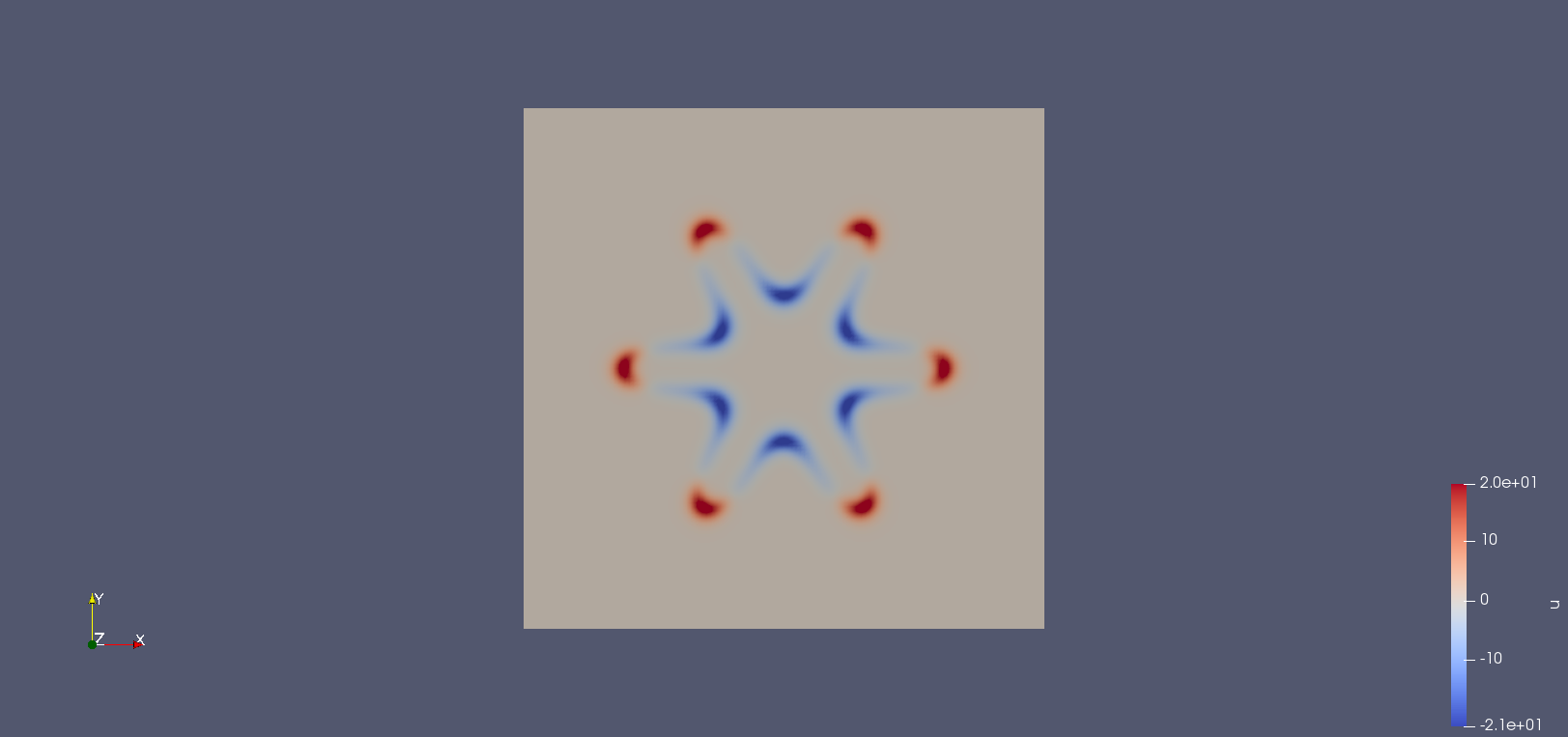}
	\includegraphics[trim={540px 112px 540px 112px}, clip, scale=.08]{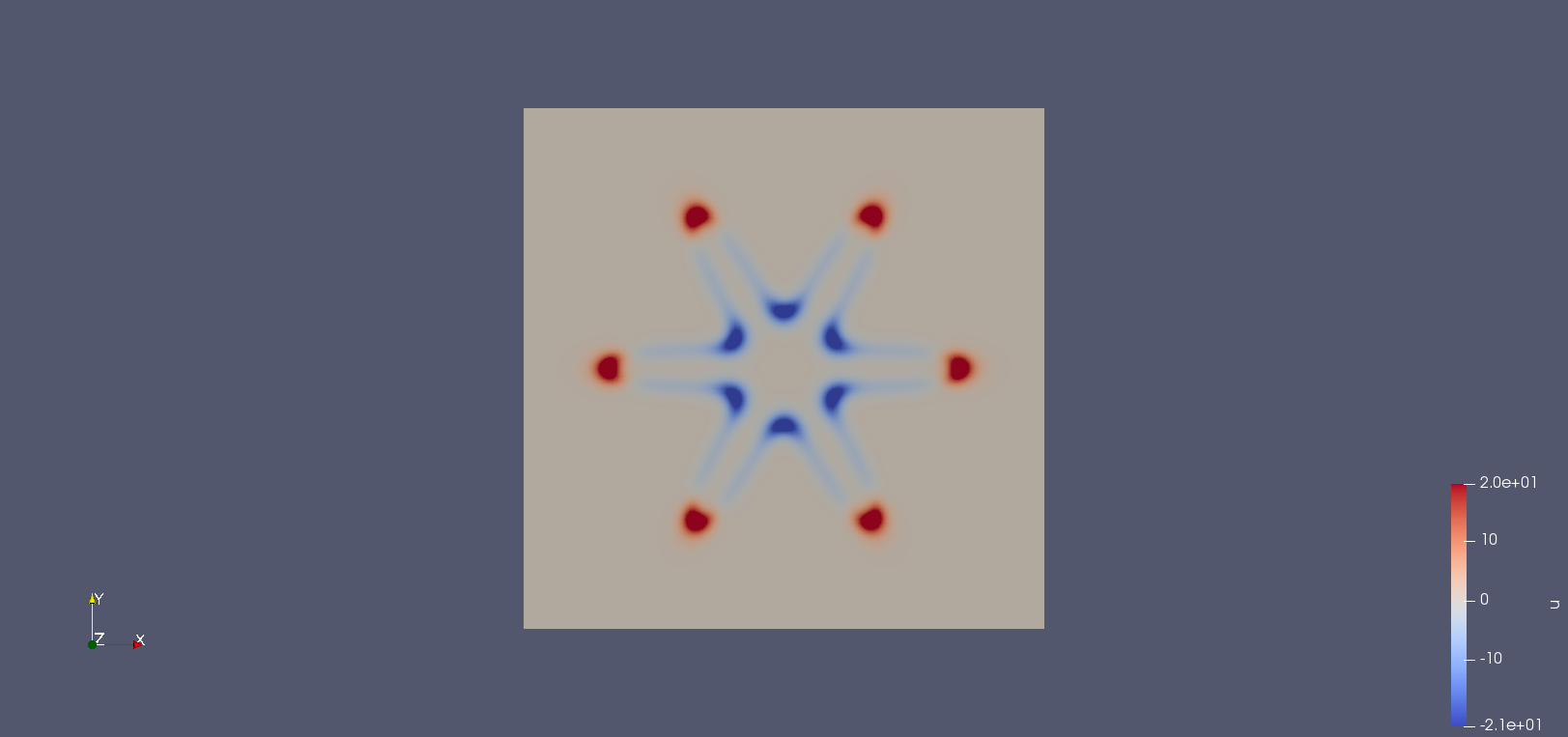}
	\includegraphics[trim={540px 112px 540px 112px}, clip, scale=.08]{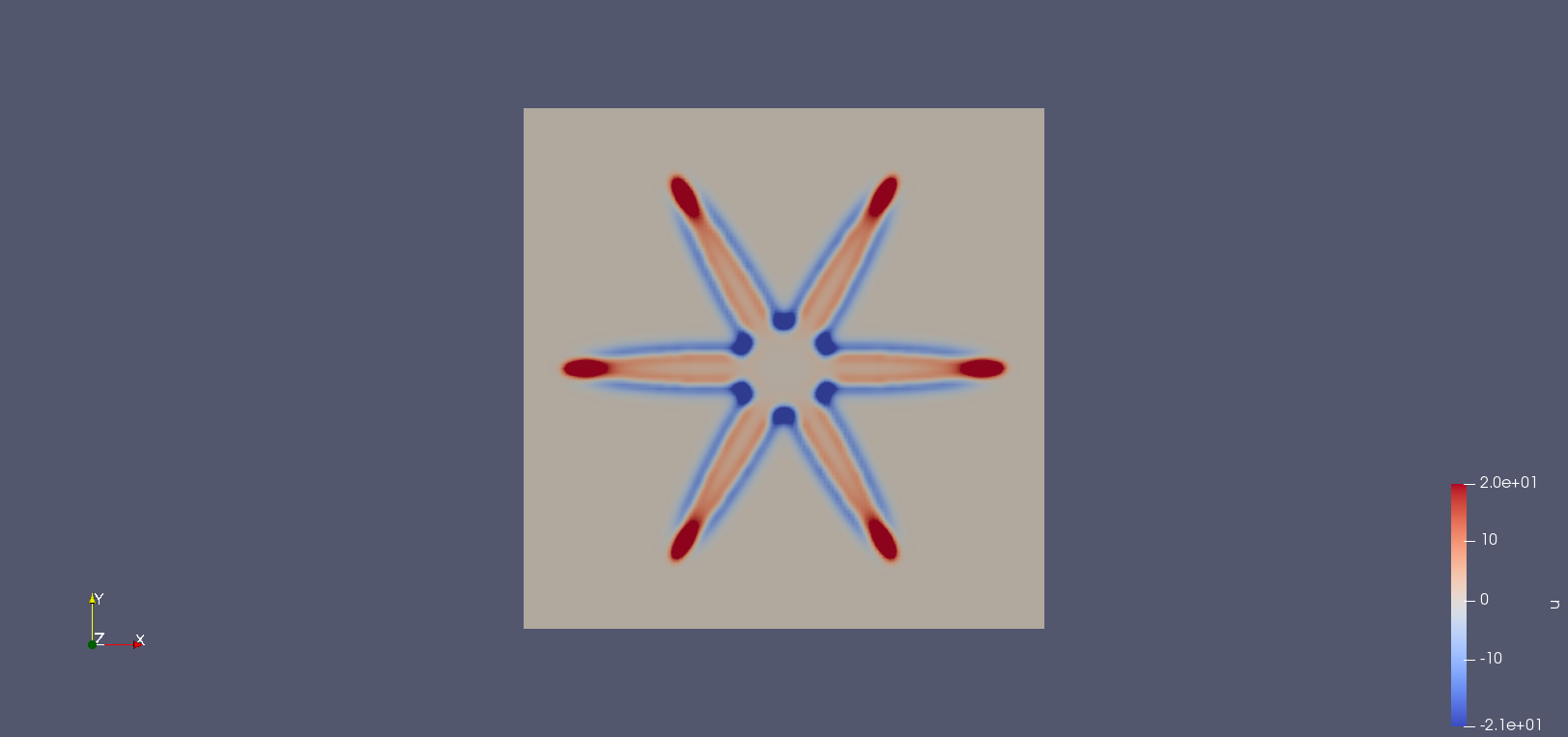}~~~
	\includegraphics[trim={1490px 0px 20px 490px}, clip, scale=.16]{{iso_starhexa_control.0000}.png}
	\caption{Results for the isotropic case.} 
	\label{fig:iso_star6}
\end{center}
\end{subfigure}
\caption{`circle to 6-star' solutions: states in the first and third row and controls in the second and fourth row.
\label{fig:isohex_star6}}
%  \johannes{    $t = timesnotlookedupsofar$.}
% Hinweis: die Zahlen entsprechen ungefähr diesen Zeitschritten, ich hab die Bilder von 0-7 darauf umgerechnet und auf "schöne" Zahlen gerundet
\end{center}
\end{figure}
\fi

%\ifgraphics
%\begin{figure}[htbp]
%\label{fig:uover time}
%\end{figure}
%\fi

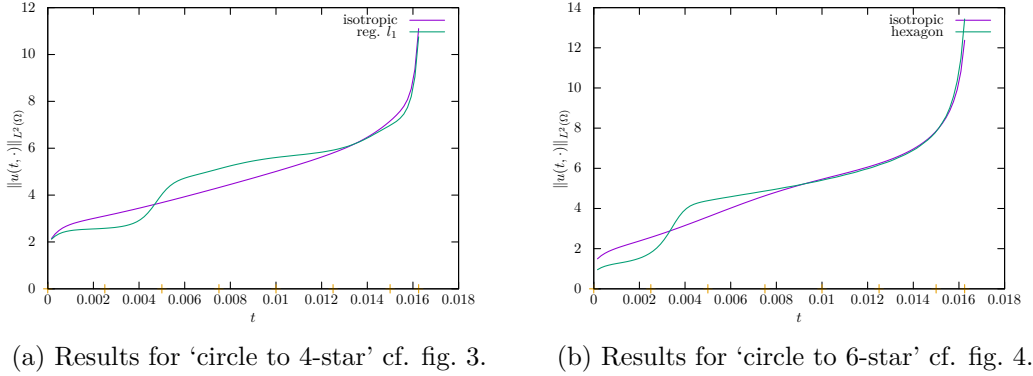
\begin{figure}
  \begin{center}
\begin{subfigure}{.48\textwidth}
  \begin{center}
    \begin{tikzpicture}[gnuplot, scale=0.5, every node/.style={scale=0.5}]
%% generated with GNUPLOT 5.2p2 (Lua 5.3; terminal rev. 99, script rev. 102)
%% Fr 07 Mai 2021 18:12:07 CEST
\path (0.000,0.000) rectangle (12.500,8.750);
\gpcolor{color=gp lt color border}
\gpsetlinetype{gp lt border}
\gpsetdashtype{gp dt solid}
\gpsetlinewidth{1.00}
\draw[gp path] (1.136,0.985)--(1.316,0.985);
\draw[gp path] (11.947,0.985)--(11.767,0.985);
\node[gp node right] at (0.952,0.985) {$0$};
\draw[gp path] (1.136,2.228)--(1.316,2.228);
\draw[gp path] (11.947,2.228)--(11.767,2.228);
\node[gp node right] at (0.952,2.228) {$2$};
\draw[gp path] (1.136,3.470)--(1.316,3.470);
\draw[gp path] (11.947,3.470)--(11.767,3.470);
\node[gp node right] at (0.952,3.470) {$4$};
\draw[gp path] (1.136,4.713)--(1.316,4.713);
\draw[gp path] (11.947,4.713)--(11.767,4.713);
\node[gp node right] at (0.952,4.713) {$6$};
\draw[gp path] (1.136,5.956)--(1.316,5.956);
\draw[gp path] (11.947,5.956)--(11.767,5.956);
\node[gp node right] at (0.952,5.956) {$8$};
\draw[gp path] (1.136,7.198)--(1.316,7.198);
\draw[gp path] (11.947,7.198)--(11.767,7.198);
\node[gp node right] at (0.952,7.198) {$10$};
\draw[gp path] (1.136,8.441)--(1.316,8.441);
\draw[gp path] (11.947,8.441)--(11.767,8.441);
\node[gp node right] at (0.952,8.441) {$12$};
\draw[gp path] (1.136,0.985)--(1.136,1.165);
\draw[gp path] (1.136,8.441)--(1.136,8.261);
\node[gp node center] at (1.136,0.677) {$0$};
\draw[gp path] (2.337,0.985)--(2.337,1.165);
\draw[gp path] (2.337,8.441)--(2.337,8.261);
\node[gp node center] at (2.337,0.677) {$0.002$};
\draw[gp path] (3.538,0.985)--(3.538,1.165);
\draw[gp path] (3.538,8.441)--(3.538,8.261);
\node[gp node center] at (3.538,0.677) {$0.004$};
\draw[gp path] (4.740,0.985)--(4.740,1.165);
\draw[gp path] (4.740,8.441)--(4.740,8.261);
\node[gp node center] at (4.740,0.677) {$0.006$};
\draw[gp path] (5.941,0.985)--(5.941,1.165);
\draw[gp path] (5.941,8.441)--(5.941,8.261);
\node[gp node center] at (5.941,0.677) {$0.008$};
\draw[gp path] (7.142,0.985)--(7.142,1.165);
\draw[gp path] (7.142,8.441)--(7.142,8.261);
\node[gp node center] at (7.142,0.677) {$0.01$};
\draw[gp path] (8.343,0.985)--(8.343,1.165);
\draw[gp path] (8.343,8.441)--(8.343,8.261);
\node[gp node center] at (8.343,0.677) {$0.012$};
\draw[gp path] (9.545,0.985)--(9.545,1.165);
\draw[gp path] (9.545,8.441)--(9.545,8.261);
\node[gp node center] at (9.545,0.677) {$0.014$};
\draw[gp path] (10.746,0.985)--(10.746,1.165);
\draw[gp path] (10.746,8.441)--(10.746,8.261);
\node[gp node center] at (10.746,0.677) {$0.016$};
\draw[gp path] (11.947,0.985)--(11.947,1.165);
\draw[gp path] (11.947,8.441)--(11.947,8.261);
\node[gp node center] at (11.947,0.677) {$0.018$};
\draw[gp path] (1.136,8.441)--(1.136,0.985)--(11.947,0.985)--(11.947,8.441)--cycle;
\node[gp node center,rotate=-270] at (0.276,4.713) {$\|u(t,\cdot)\|_{L^2(\Omega)}$};
\node[gp node center] at (6.541,0.215) {$t$};
\node[gp node right] at (10.479,8.107) {isotropic};
\gpcolor{rgb color={0.580,0.000,0.827}}
\draw[gp path] (10.663,8.107)--(11.579,8.107);
\draw[gp path] (1.234,2.318)--(1.331,2.441)--(1.429,2.530)--(1.526,2.597)--(1.624,2.647)%
  --(1.722,2.688)--(1.819,2.722)--(1.917,2.750)--(2.014,2.776)--(2.112,2.800)--(2.210,2.823)%
  --(2.307,2.845)--(2.405,2.866)--(2.502,2.888)--(2.600,2.909)--(2.698,2.930)--(2.795,2.952)%
  --(2.893,2.973)--(2.990,2.995)--(3.088,3.018)--(3.186,3.040)--(3.283,3.063)--(3.381,3.085)%
  --(3.478,3.109)--(3.576,3.132)--(3.674,3.156)--(3.771,3.180)--(3.869,3.204)--(3.966,3.228)%
  --(4.064,3.252)--(4.162,3.277)--(4.259,3.302)--(4.357,3.327)--(4.454,3.352)--(4.552,3.378)%
  --(4.650,3.403)--(4.747,3.429)--(4.845,3.455)--(4.942,3.481)--(5.040,3.507)--(5.138,3.533)%
  --(5.235,3.560)--(5.333,3.586)--(5.430,3.613)--(5.528,3.639)--(5.626,3.666)--(5.723,3.693)%
  --(5.821,3.720)--(5.918,3.747)--(6.016,3.774)--(6.114,3.802)--(6.211,3.829)--(6.309,3.857)%
  --(6.406,3.884)--(6.504,3.912)--(6.602,3.940)--(6.699,3.968)--(6.797,3.997)--(6.894,4.025)%
  --(6.992,4.054)--(7.090,4.082)--(7.187,4.111)--(7.285,4.141)--(7.382,4.170)--(7.480,4.200)%
  --(7.578,4.230)--(7.675,4.260)--(7.773,4.291)--(7.870,4.322)--(7.968,4.353)--(8.066,4.385)%
  --(8.163,4.418)--(8.261,4.451)--(8.358,4.485)--(8.456,4.520)--(8.554,4.555)--(8.651,4.592)%
  --(8.749,4.629)--(8.846,4.668)--(8.944,4.708)--(9.042,4.750)--(9.139,4.794)--(9.237,4.841)%
  --(9.334,4.889)--(9.432,4.941)--(9.530,4.996)--(9.627,5.054)--(9.725,5.116)--(9.822,5.183)%
  --(9.920,5.253)--(10.018,5.329)--(10.115,5.409)--(10.213,5.495)--(10.310,5.589)--(10.408,5.695)%
  --(10.506,5.825)--(10.603,6.004)--(10.701,6.288)--(10.798,6.807)--(10.896,7.885);
\gpcolor{color=gp lt color border}
\node[gp node right] at (10.479,7.799) {reg. $l_1$};
\gpcolor{rgb color={0.000,0.620,0.451}}
\draw[gp path] (10.663,7.799)--(11.579,7.799);
\draw[gp path] (1.234,2.300)--(1.331,2.383)--(1.429,2.441)--(1.526,2.480)--(1.624,2.507)%
  --(1.722,2.526)--(1.819,2.540)--(1.917,2.549)--(2.014,2.557)--(2.112,2.562)--(2.210,2.567)%
  --(2.307,2.572)--(2.405,2.576)--(2.502,2.581)--(2.600,2.587)--(2.698,2.594)--(2.795,2.602)%
  --(2.893,2.611)--(2.990,2.624)--(3.088,2.639)--(3.186,2.659)--(3.283,2.685)--(3.381,2.718)%
  --(3.478,2.762)--(3.576,2.821)--(3.674,2.898)--(3.771,2.997)--(3.869,3.116)--(3.966,3.250)%
  --(4.064,3.388)--(4.162,3.517)--(4.259,3.628)--(4.357,3.718)--(4.454,3.789)--(4.552,3.844)%
  --(4.650,3.888)--(4.747,3.923)--(4.845,3.954)--(4.942,3.981)--(5.040,4.007)--(5.138,4.033)%
  --(5.235,4.058)--(5.333,4.084)--(5.430,4.110)--(5.528,4.136)--(5.626,4.163)--(5.723,4.189)%
  --(5.821,4.215)--(5.918,4.240)--(6.016,4.264)--(6.114,4.287)--(6.211,4.308)--(6.309,4.330)%
  --(6.406,4.350)--(6.504,4.369)--(6.602,4.387)--(6.699,4.403)--(6.797,4.419)--(6.894,4.434)%
  --(6.992,4.448)--(7.090,4.461)--(7.187,4.473)--(7.285,4.485)--(7.382,4.495)--(7.480,4.506)%
  --(7.578,4.517)--(7.675,4.528)--(7.773,4.538)--(7.870,4.548)--(7.968,4.560)--(8.066,4.572)%
  --(8.163,4.585)--(8.261,4.597)--(8.358,4.611)--(8.456,4.628)--(8.554,4.646)--(8.651,4.664)%
  --(8.749,4.686)--(8.846,4.711)--(8.944,4.738)--(9.042,4.767)--(9.139,4.802)--(9.237,4.840)%
  --(9.334,4.880)--(9.432,4.926)--(9.530,4.976)--(9.627,5.026)--(9.725,5.083)--(9.822,5.136)%
  --(9.920,5.192)--(10.018,5.248)--(10.115,5.302)--(10.213,5.362)--(10.310,5.424)--(10.408,5.507)%
  --(10.506,5.619)--(10.603,5.794)--(10.701,6.083)--(10.798,6.604)--(10.896,7.667);
\gpcolor{rgb color={0.902,0.624,0.000}}
\gpsetpointsize{4.00}
\gppoint{gp mark 1}{(1.136,0.985)}
\gppoint{gp mark 1}{(2.638,0.985)}
\gppoint{gp mark 1}{(4.139,0.985)}
\gppoint{gp mark 1}{(5.641,0.985)}
\gppoint{gp mark 1}{(7.142,0.985)}
\gppoint{gp mark 1}{(8.644,0.985)}
\gppoint{gp mark 1}{(10.145,0.985)}
\gppoint{gp mark 1}{(10.896,0.985)}
\gpcolor{color=gp lt color border}
\draw[gp path] (1.136,8.441)--(1.136,0.985)--(11.947,0.985)--(11.947,8.441)--cycle;
%% coordinates of the plot area
\gpdefrectangularnode{gp plot 1}{\pgfpoint{1.136cm}{0.985cm}}{\pgfpoint{11.947cm}{8.441cm}}
\end{tikzpicture}
%% gnuplot variables
    \label{fig:4star_compare}\caption{Results for `circle to 4-star' cf. \cref{fig:l1iso_star4}. }
  \end{center}
\end{subfigure}
%\hspace{.05\textwidth}
\begin{subfigure}{.48\textwidth}
  \begin{center}
    \begin{tikzpicture}[gnuplot, scale=0.5, every node/.style={scale=0.5}]
%% generated with GNUPLOT 5.2p2 (Lua 5.3; terminal rev. 99, script rev. 102)
%% Fr 07 Mai 2021 18:12:07 CEST
\path (0.000,0.000) rectangle (12.500,8.750);
\gpcolor{color=gp lt color border}
\gpsetlinetype{gp lt border}
\gpsetdashtype{gp dt solid}
\gpsetlinewidth{1.00}
\draw[gp path] (1.136,0.985)--(1.316,0.985);
\draw[gp path] (11.947,0.985)--(11.767,0.985);
\node[gp node right] at (0.952,0.985) {$0$};
\draw[gp path] (1.136,2.050)--(1.316,2.050);
\draw[gp path] (11.947,2.050)--(11.767,2.050);
\node[gp node right] at (0.952,2.050) {$2$};
\draw[gp path] (1.136,3.115)--(1.316,3.115);
\draw[gp path] (11.947,3.115)--(11.767,3.115);
\node[gp node right] at (0.952,3.115) {$4$};
\draw[gp path] (1.136,4.180)--(1.316,4.180);
\draw[gp path] (11.947,4.180)--(11.767,4.180);
\node[gp node right] at (0.952,4.180) {$6$};
\draw[gp path] (1.136,5.246)--(1.316,5.246);
\draw[gp path] (11.947,5.246)--(11.767,5.246);
\node[gp node right] at (0.952,5.246) {$8$};
\draw[gp path] (1.136,6.311)--(1.316,6.311);
\draw[gp path] (11.947,6.311)--(11.767,6.311);
\node[gp node right] at (0.952,6.311) {$10$};
\draw[gp path] (1.136,7.376)--(1.316,7.376);
\draw[gp path] (11.947,7.376)--(11.767,7.376);
\node[gp node right] at (0.952,7.376) {$12$};
\draw[gp path] (1.136,8.441)--(1.316,8.441);
\draw[gp path] (11.947,8.441)--(11.767,8.441);
\node[gp node right] at (0.952,8.441) {$14$};
\draw[gp path] (1.136,0.985)--(1.136,1.165);
\draw[gp path] (1.136,8.441)--(1.136,8.261);
\node[gp node center] at (1.136,0.677) {$0$};
\draw[gp path] (2.337,0.985)--(2.337,1.165);
\draw[gp path] (2.337,8.441)--(2.337,8.261);
\node[gp node center] at (2.337,0.677) {$0.002$};
\draw[gp path] (3.538,0.985)--(3.538,1.165);
\draw[gp path] (3.538,8.441)--(3.538,8.261);
\node[gp node center] at (3.538,0.677) {$0.004$};
\draw[gp path] (4.740,0.985)--(4.740,1.165);
\draw[gp path] (4.740,8.441)--(4.740,8.261);
\node[gp node center] at (4.740,0.677) {$0.006$};
\draw[gp path] (5.941,0.985)--(5.941,1.165);
\draw[gp path] (5.941,8.441)--(5.941,8.261);
\node[gp node center] at (5.941,0.677) {$0.008$};
\draw[gp path] (7.142,0.985)--(7.142,1.165);
\draw[gp path] (7.142,8.441)--(7.142,8.261);
\node[gp node center] at (7.142,0.677) {$0.01$};
\draw[gp path] (8.343,0.985)--(8.343,1.165);
\draw[gp path] (8.343,8.441)--(8.343,8.261);
\node[gp node center] at (8.343,0.677) {$0.012$};
\draw[gp path] (9.545,0.985)--(9.545,1.165);
\draw[gp path] (9.545,8.441)--(9.545,8.261);
\node[gp node center] at (9.545,0.677) {$0.014$};
\draw[gp path] (10.746,0.985)--(10.746,1.165);
\draw[gp path] (10.746,8.441)--(10.746,8.261);
\node[gp node center] at (10.746,0.677) {$0.016$};
\draw[gp path] (11.947,0.985)--(11.947,1.165);
\draw[gp path] (11.947,8.441)--(11.947,8.261);
\node[gp node center] at (11.947,0.677) {$0.018$};
\draw[gp path] (1.136,8.441)--(1.136,0.985)--(11.947,0.985)--(11.947,8.441)--cycle;
\node[gp node center,rotate=-270] at (0.276,4.713) {$\|u(t,\cdot)\|_{L^2(\Omega)}$};
\node[gp node center] at (6.541,0.215) {$t$};
\node[gp node right] at (10.479,8.107) {isotropic};
\gpcolor{rgb color={0.580,0.000,0.827}}
\draw[gp path] (10.663,8.107)--(11.579,8.107);
\draw[gp path] (1.234,1.781)--(1.331,1.863)--(1.429,1.926)--(1.526,1.977)--(1.624,2.019)%
  --(1.722,2.057)--(1.819,2.091)--(1.917,2.123)--(2.014,2.154)--(2.112,2.184)--(2.210,2.214)%
  --(2.307,2.244)--(2.405,2.275)--(2.502,2.306)--(2.600,2.337)--(2.698,2.369)--(2.795,2.401)%
  --(2.893,2.434)--(2.990,2.467)--(3.088,2.501)--(3.186,2.536)--(3.283,2.571)--(3.381,2.606)%
  --(3.478,2.642)--(3.576,2.679)--(3.674,2.715)--(3.771,2.752)--(3.869,2.789)--(3.966,2.827)%
  --(4.064,2.865)--(4.162,2.902)--(4.259,2.940)--(4.357,2.978)--(4.454,3.016)--(4.552,3.054)%
  --(4.650,3.091)--(4.747,3.128)--(4.845,3.165)--(4.942,3.202)--(5.040,3.239)--(5.138,3.275)%
  --(5.235,3.310)--(5.333,3.345)--(5.430,3.380)--(5.528,3.414)--(5.626,3.447)--(5.723,3.480)%
  --(5.821,3.512)--(5.918,3.543)--(6.016,3.574)--(6.114,3.605)--(6.211,3.634)--(6.309,3.663)%
  --(6.406,3.692)--(6.504,3.720)--(6.602,3.747)--(6.699,3.774)--(6.797,3.801)--(6.894,3.827)%
  --(6.992,3.853)--(7.090,3.878)--(7.187,3.903)--(7.285,3.929)--(7.382,3.954)--(7.480,3.979)%
  --(7.578,4.004)--(7.675,4.029)--(7.773,4.055)--(7.870,4.080)--(7.968,4.107)--(8.066,4.133)%
  --(8.163,4.161)--(8.261,4.189)--(8.358,4.217)--(8.456,4.247)--(8.554,4.278)--(8.651,4.310)%
  --(8.749,4.343)--(8.846,4.377)--(8.944,4.414)--(9.042,4.452)--(9.139,4.493)--(9.237,4.536)%
  --(9.334,4.582)--(9.432,4.631)--(9.530,4.684)--(9.627,4.742)--(9.725,4.805)--(9.822,4.874)%
  --(9.920,4.950)--(10.018,5.036)--(10.115,5.133)--(10.213,5.244)--(10.310,5.372)--(10.408,5.525)%
  --(10.506,5.712)--(10.603,5.950)--(10.701,6.271)--(10.798,6.750)--(10.896,7.576);
\gpcolor{color=gp lt color border}
\node[gp node right] at (10.479,7.799) {hexagon};
\gpcolor{rgb color={0.000,0.620,0.451}}
\draw[gp path] (10.663,7.799)--(11.579,7.799);
\draw[gp path] (1.234,1.492)--(1.331,1.546)--(1.429,1.584)--(1.526,1.611)--(1.624,1.634)%
  --(1.722,1.653)--(1.819,1.670)--(1.917,1.688)--(2.014,1.707)--(2.112,1.729)--(2.210,1.755)%
  --(2.307,1.785)--(2.405,1.821)--(2.502,1.866)--(2.600,1.920)--(2.698,1.987)--(2.795,2.069)%
  --(2.893,2.170)--(2.990,2.292)--(3.088,2.435)--(3.186,2.593)--(3.283,2.755)--(3.381,2.902)%
  --(3.478,3.023)--(3.576,3.115)--(3.674,3.182)--(3.771,3.227)--(3.869,3.260)--(3.966,3.286)%
  --(4.064,3.308)--(4.162,3.327)--(4.259,3.345)--(4.357,3.363)--(4.454,3.380)--(4.552,3.397)%
  --(4.650,3.413)--(4.747,3.429)--(4.845,3.446)--(4.942,3.462)--(5.040,3.479)--(5.138,3.495)%
  --(5.235,3.510)--(5.333,3.526)--(5.430,3.543)--(5.528,3.560)--(5.626,3.576)--(5.723,3.592)%
  --(5.821,3.608)--(5.918,3.625)--(6.016,3.643)--(6.114,3.662)--(6.211,3.680)--(6.309,3.696)%
  --(6.406,3.713)--(6.504,3.732)--(6.602,3.753)--(6.699,3.774)--(6.797,3.792)--(6.894,3.811)%
  --(6.992,3.830)--(7.090,3.852)--(7.187,3.875)--(7.285,3.898)--(7.382,3.919)--(7.480,3.939)%
  --(7.578,3.963)--(7.675,3.989)--(7.773,4.015)--(7.870,4.038)--(7.968,4.062)--(8.066,4.089)%
  --(8.163,4.120)--(8.261,4.149)--(8.358,4.174)--(8.456,4.203)--(8.554,4.237)--(8.651,4.274)%
  --(8.749,4.306)--(8.846,4.338)--(8.944,4.373)--(9.042,4.415)--(9.139,4.460)--(9.237,4.502)%
  --(9.334,4.546)--(9.432,4.598)--(9.530,4.652)--(9.627,4.708)--(9.725,4.778)--(9.822,4.853)%
  --(9.920,4.928)--(10.018,5.016)--(10.115,5.117)--(10.213,5.244)--(10.310,5.394)--(10.408,5.560)%
  --(10.506,5.787)--(10.603,6.088)--(10.701,6.497)--(10.798,7.094)--(10.896,8.140);
\gpcolor{rgb color={0.902,0.624,0.000}}
\gpsetpointsize{4.00}
\gppoint{gp mark 1}{(1.136,0.985)}
\gppoint{gp mark 1}{(2.638,0.985)}
\gppoint{gp mark 1}{(4.139,0.985)}
\gppoint{gp mark 1}{(5.641,0.985)}
\gppoint{gp mark 1}{(7.142,0.985)}
\gppoint{gp mark 1}{(8.644,0.985)}
\gppoint{gp mark 1}{(10.145,0.985)}
\gppoint{gp mark 1}{(10.896,0.985)}
\gpcolor{color=gp lt color border}
\draw[gp path] (1.136,8.441)--(1.136,0.985)--(11.947,0.985)--(11.947,8.441)--cycle;
%% coordinates of the plot area
\gpdefrectangularnode{gp plot 1}{\pgfpoint{1.136cm}{0.985cm}}{\pgfpoint{11.947cm}{8.441cm}}
\end{tikzpicture}
%% gnuplot variables
    \label{fig:6star_compare}\caption{Results for `circle to 6-star' cf. \cref{fig:isohex_star6}. }
  \end{center}
  \end{subfigure}
\caption{Time evolution $\|u(t,\cdot)\|_{L^2(\Omega)}$.
	\label{fig:star_iso_compare}}
\end{center}
\end{figure}

In \cref{tab:comp_star} the computed (local) minima are listed together with their single constituents---the difference of the optimal state to the desired state $j_1=\|y-y_\Omega\|^2_{L^2(\Omega)}$ as well as the contribution of the control
$j_2= \tfrac\lambda{2\varepsilon} \|u\|^2_{L^2(Q)}$.
We note that the (local) optima for the isotropic case are slightly below their anisotropic counterparts.
%https://www.tablesgenerator.com/
\begin{table}[h]
	\scriptsize
	\begin{tabular}{|l|l|l|l|l|}
		\hline
		& \multicolumn{2}{c|}{4 star} & \multicolumn{2}{c|}{6 star} \\ \cline{2-5} 
		& iso          & $l1$         & iso          & hexa         \\ \hline
		$j(u)$      & 0.102184            & 0.107378            & 0.115034            & 0.12248            \\ \hline
		$j_1 + j_2$ & 0.0115987 + 0.0905854            & 0.0108916 + 0.0964865            & 0.0122366 + 0.102798            & 0.0156595 + 0.106821            \\ \hline
	\end{tabular}
\caption{Values of the cost functional.}
\label{tab:comp_star}
\end{table}

\kommentar{
Also we compare the numerical effort by comparing the needed trust region Newton steps.
\luise{Fehlt}
\johannes{
  l1: 
hexa:
\\
eingehen auf TR Schritte und Kontrollcost im Vergleich zu iso
}
}

\subsubsection{Splitting and merging geometries}
Finally we consider examples where topology changes are necessary to aim at the target.
In \cref{fig:split} we present the results for splitting a circle, a square and a hexagon into two of such respectively. The underlying model equation uses the corresponding (an-)isotropy.
In \cref{fig:merge} the solutions of merging two of these objects into one are given.
Here the target objects are the initial states of the splitting examples and vice versa.
The norms of the corresponding controls over the time can be seen in \cref{fig:u_norm_others}.
To avoid potential confusion% for comparison with other data
, we point out that the scales of the ordinates are adapted to better fit the plots.

While the hexagon is splitted by squeezing it together vertically, the circle is controlled to develop first a hole in the middle and then to increase the hole until the split is present. The square is divided at the whole middle line simultaneously---as far as we could see visually. The controls are largest at times where they force topology changes as can be observed in \cref{fig:unorm_split}.
%{To avoid potential confusion% for comparison with other data
%	, we point out that the scale of the ordinate is adapted to better fit the plot.}

Considering the examples for `merging' (see \cref{fig:merge}) we observe a very similar behavior of the states as for the `splitting' solutions, but backwards in time.
There is less control necessary which is indicated by the values for $j_2$ in \cref{tab:comp_splitting,tab:comp_merge}
%where $j_1$ and $j_2$ are given as for the star like examples in Table \ref{tab:comp_star}.
{where the cost functionals are given as for the star like examples before.}
For the isotropic and hexagon example the splitting cost is higher by a factor of approximately $1.5$.
%Considering the Figures \ref{fig:unorm_split} and
%\ref{fig:unorm_merge} one has to keep in mind the employed limits of the axis for $\|u\|_{L^2(\Omega)}$.
That the difference is not bigger is probably due to the short time interval that forces the evolution of the gradient flow to be accelerated to obtain the target in time---a phenomenon present for splitting as well as for merging, with comparable impact. This also leads to the nearly constant time behavior for a long period { as can be seen in \cref{fig:unorm_merge}}.
% Finally this might explain the small difference for the $l_1$ example between \cref{tab:comp_splitting,tab:comp_merge} as the initial and final squares differ in size.
\kommentar{
The reduction at the end {of \cref{fig:unorm_merge}}, that is most notable for the square example, is due to the agreement of the merging with the uncontrolled Allen-Cahn evolution.
%Nevertheless the merging of the square at the end is mainly due to the uncontrolled Allen-Cahn evolution.
The steep peak at the end is explained by the adjustment of small details as is also the case in the earlier examples.
Also notable is the almost linear increase for the isotropic example in the beginning. Comparing with \cref{fig:merge} lets us conclude that this is due to the deformation of the circle to build two prongs at the top and bottom, which strongly opposes the natural Allen-Cahn flow. As soon as the two circles touch at these two points (see the fourth snapshot) the control cost decreases again. }

Altogether these examples have demonstrated that it is possible to steer to a variety of shapes even if they have a different topology as the initial state. Also targeting strong crystal-like structures is possible which might find use in material science or chemical applications.

\ifgraphics
\begin{figure}[htbp]
  \begin{center}
\begin{subfigure}{1.0\textwidth}
  \begin{center}
	\includegraphics[trim={540px 112px 540px 112px}, clip, scale=.08]{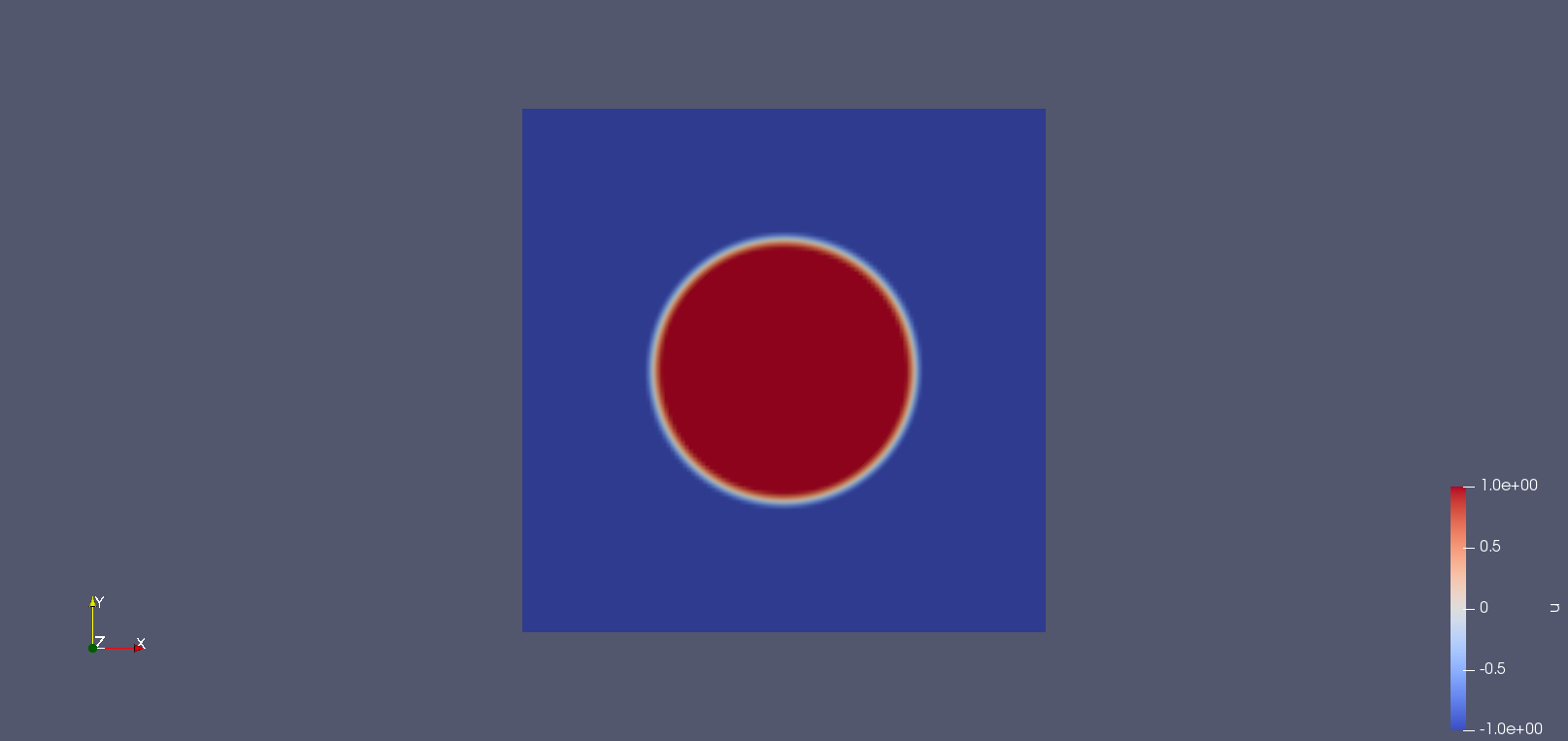}
	\includegraphics[trim={540px 112px 540px 112px}, clip, scale=.08]{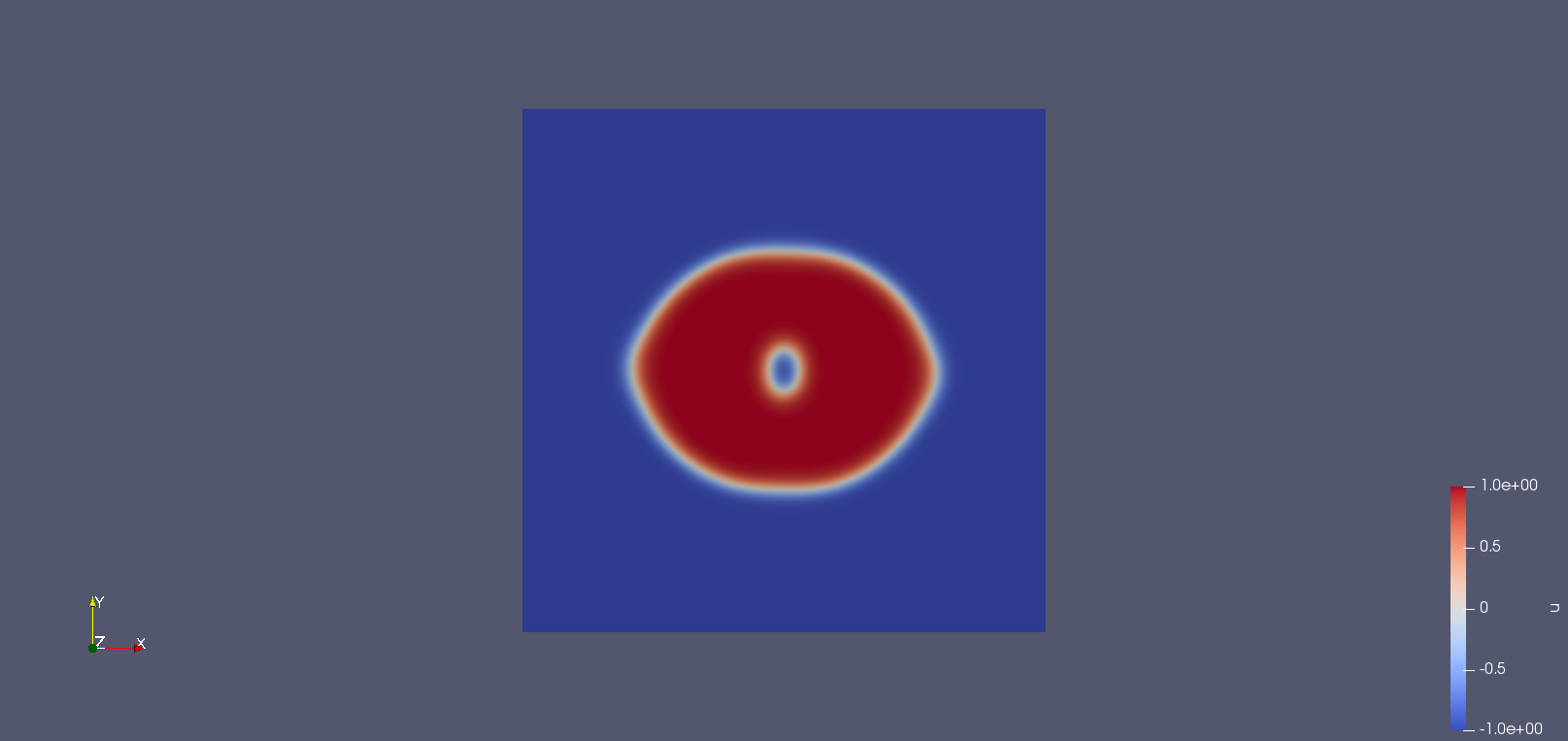}
	\includegraphics[trim={540px 112px 540px 112px}, clip, scale=.08]{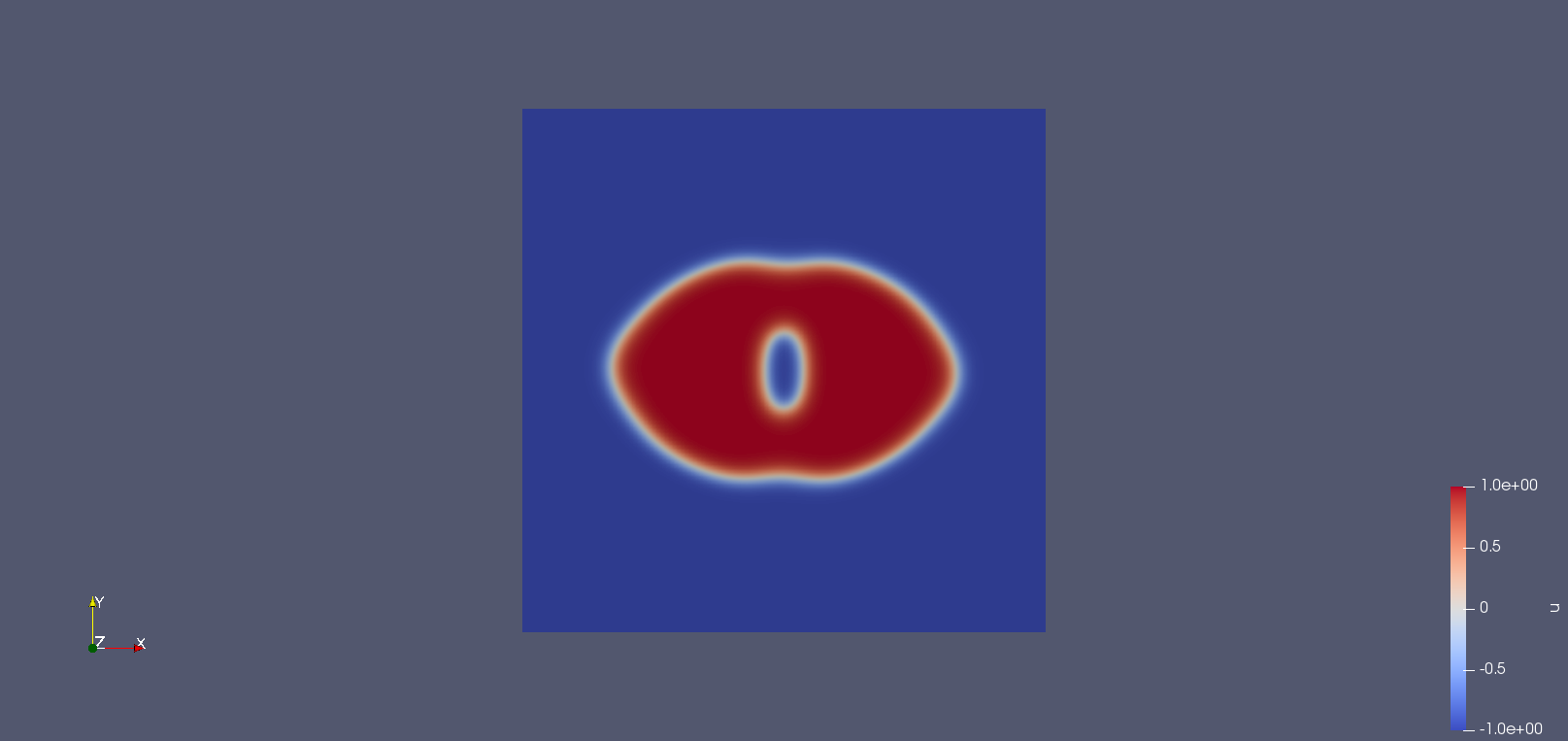}
	\includegraphics[trim={540px 112px 540px 112px}, clip, scale=.08]{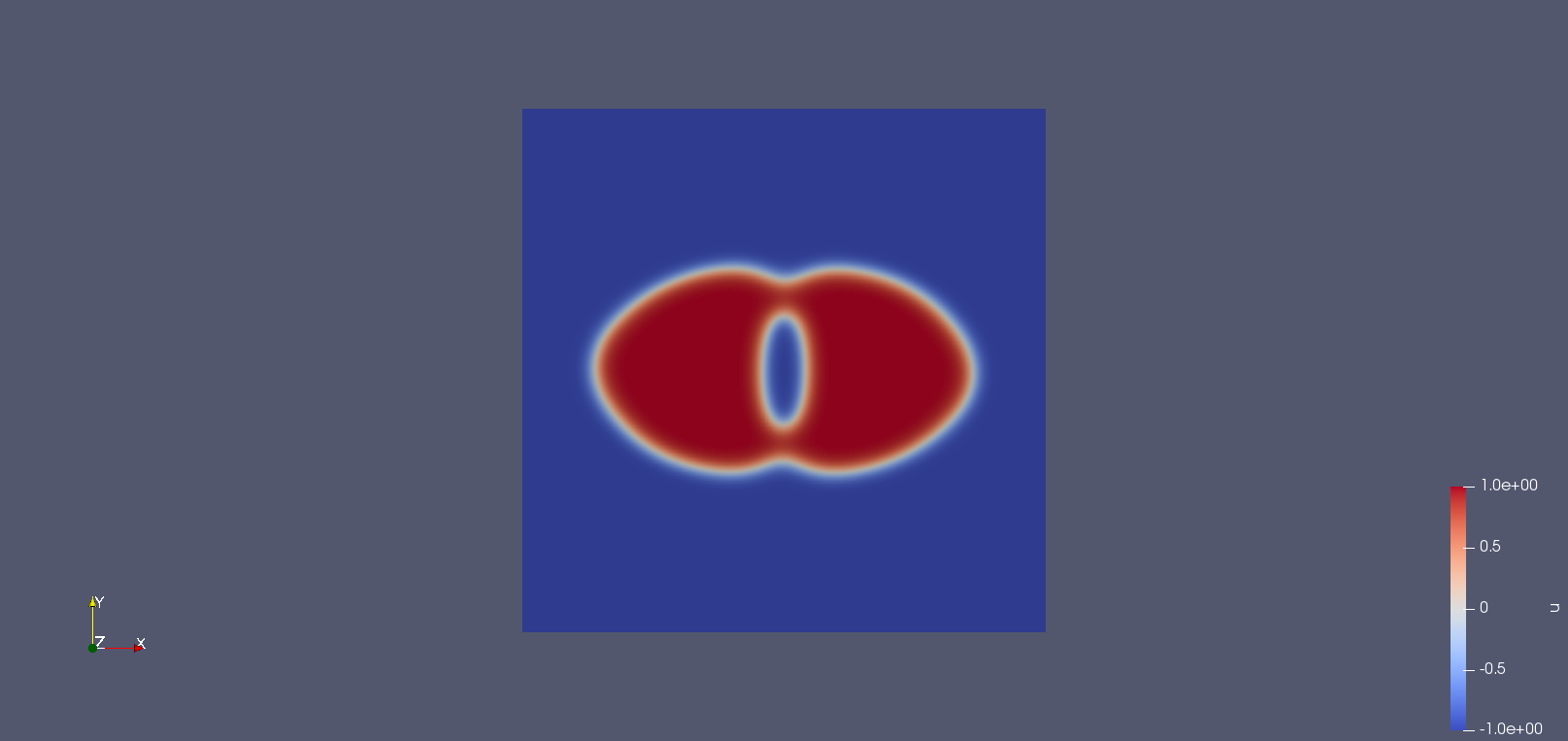}
	\includegraphics[trim={540px 112px 540px 112px}, clip, scale=.08]{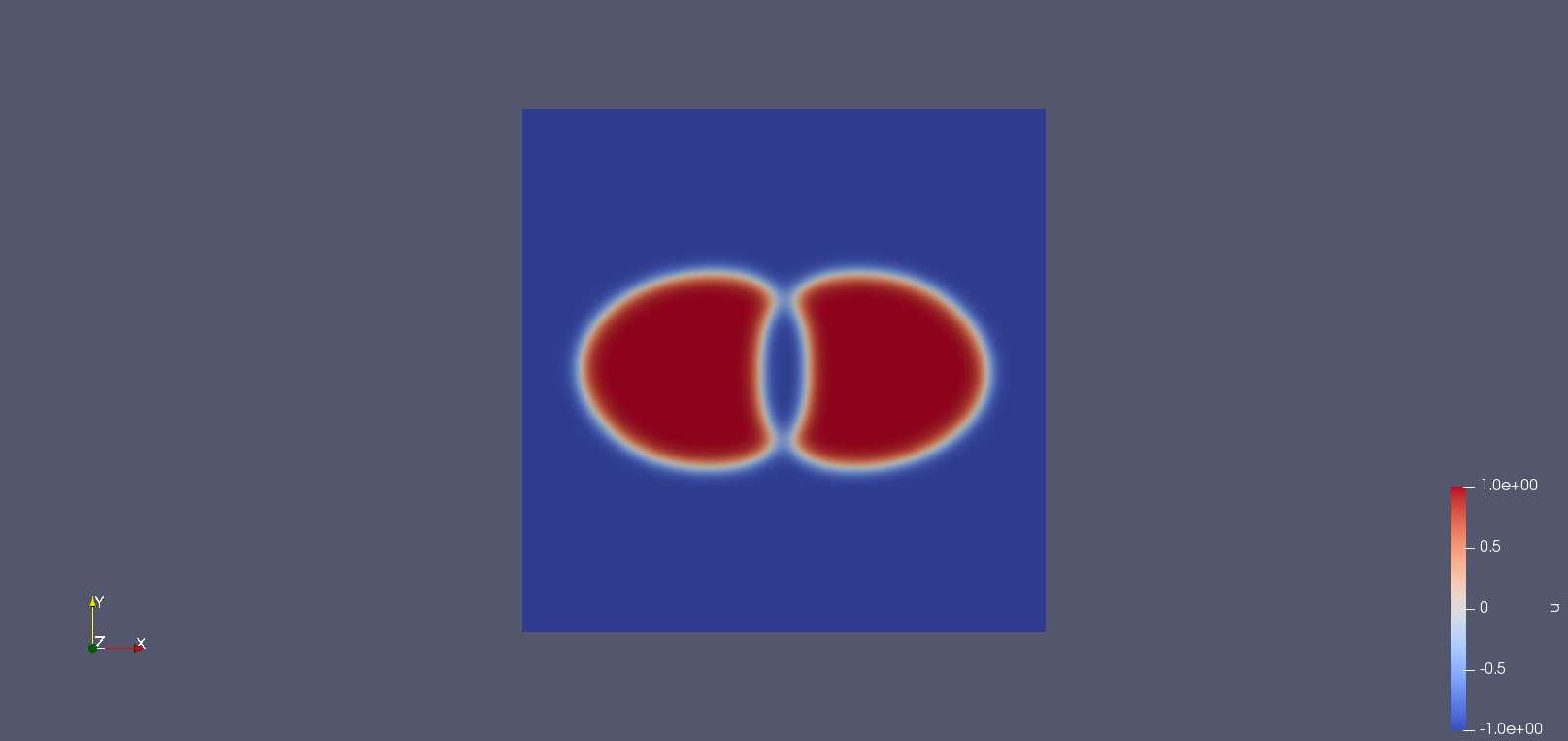}
	\includegraphics[trim={540px 112px 540px 112px}, clip, scale=.08]{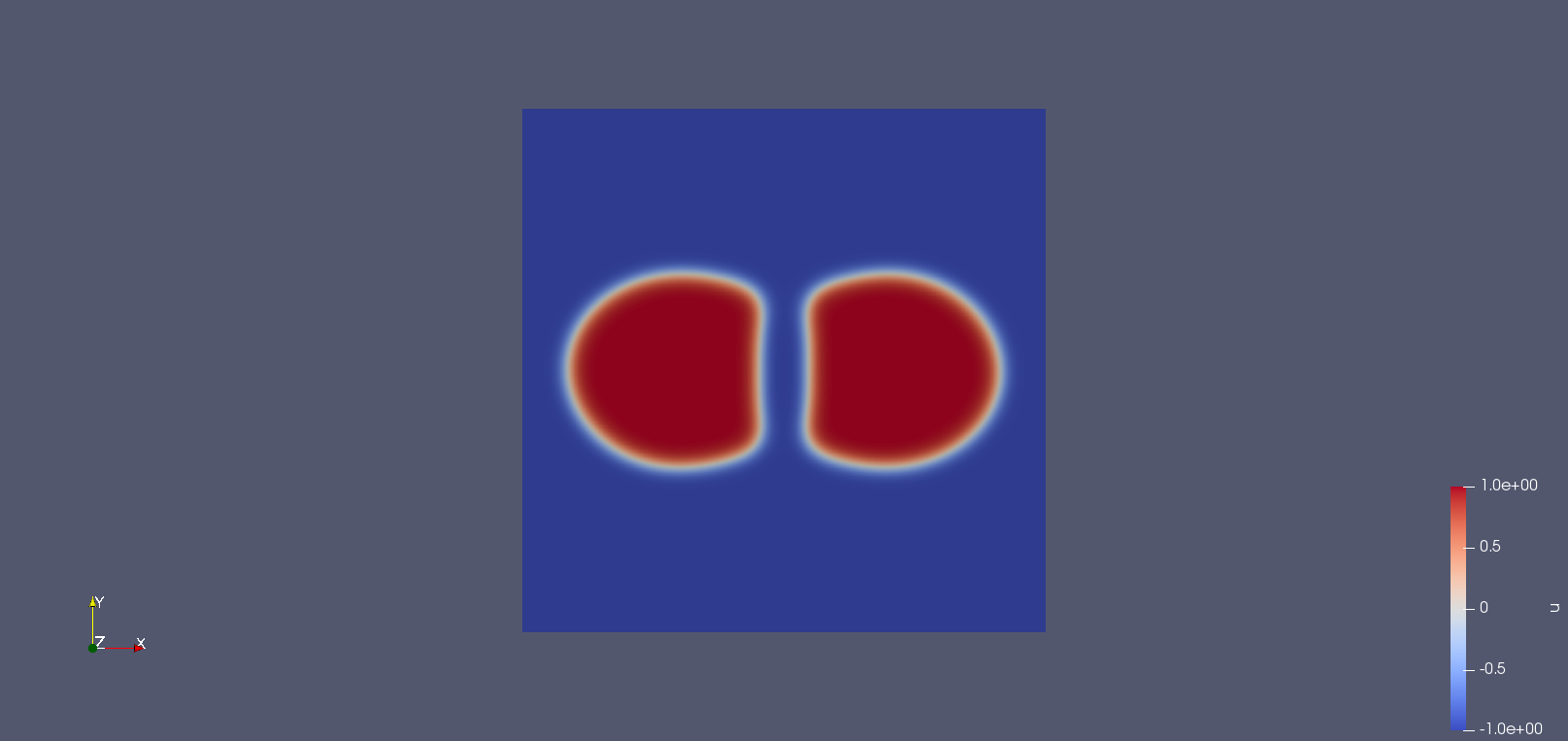}
	\includegraphics[trim={540px 112px 540px 112px}, clip, scale=.08]{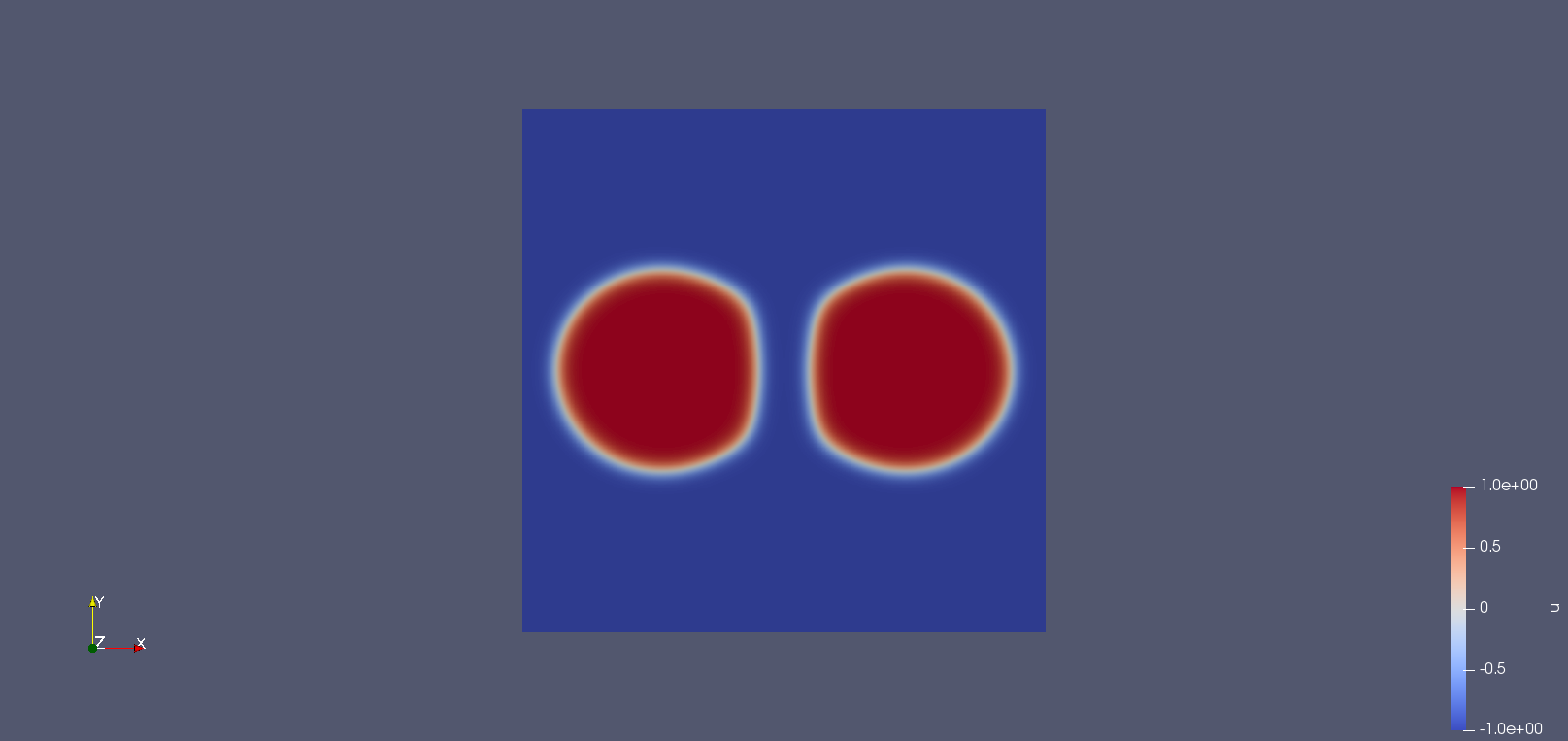}
	\includegraphics[trim={540px 112px 540px 112px}, clip, scale=.08]{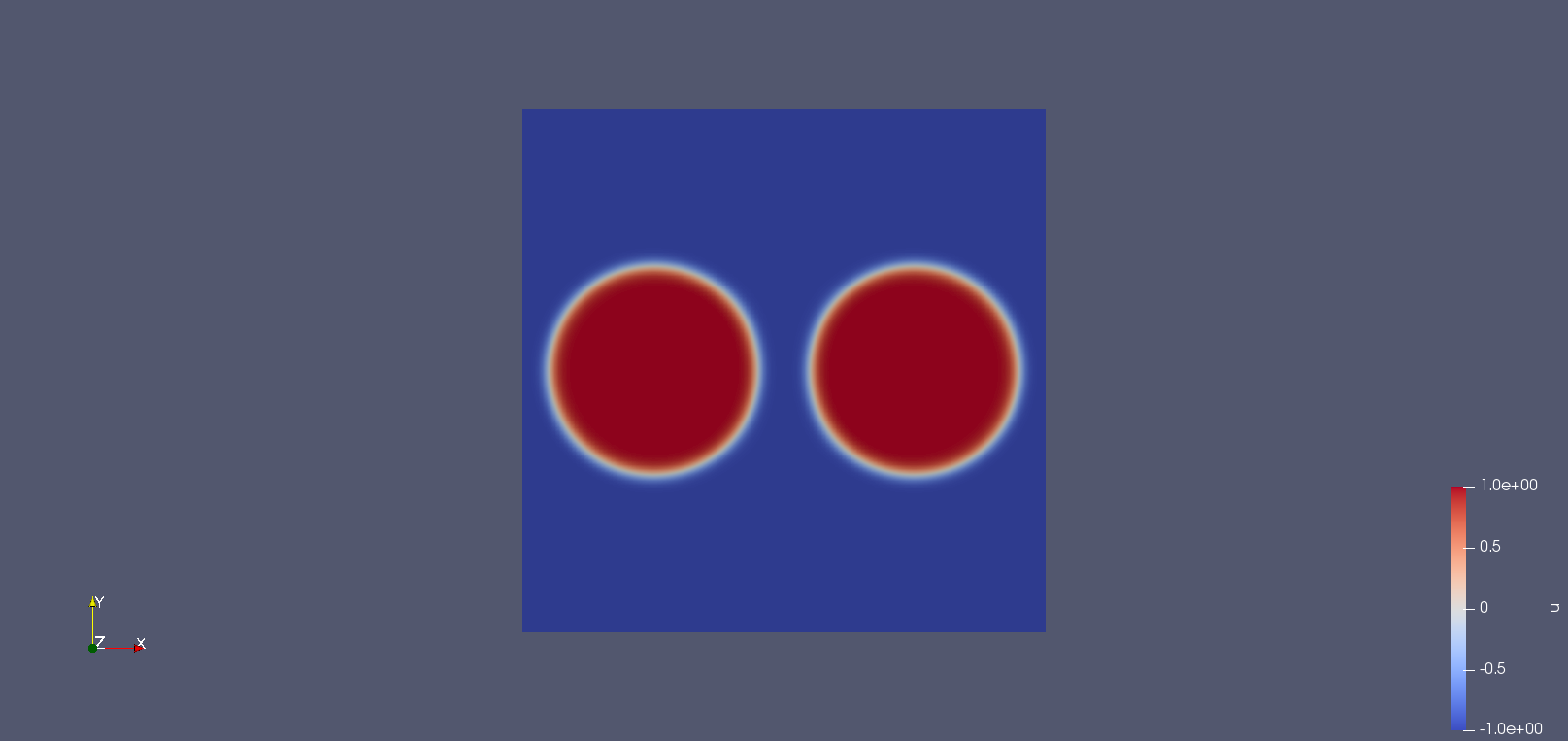}~~~
	\includegraphics[trim={1490px 0px 20px 490px}, clip, scale=.16]{{iso_splitting_state.0000}.png}~\\~\\
	\includegraphics[trim={540px 112px 540px 112px}, clip, scale=.08]{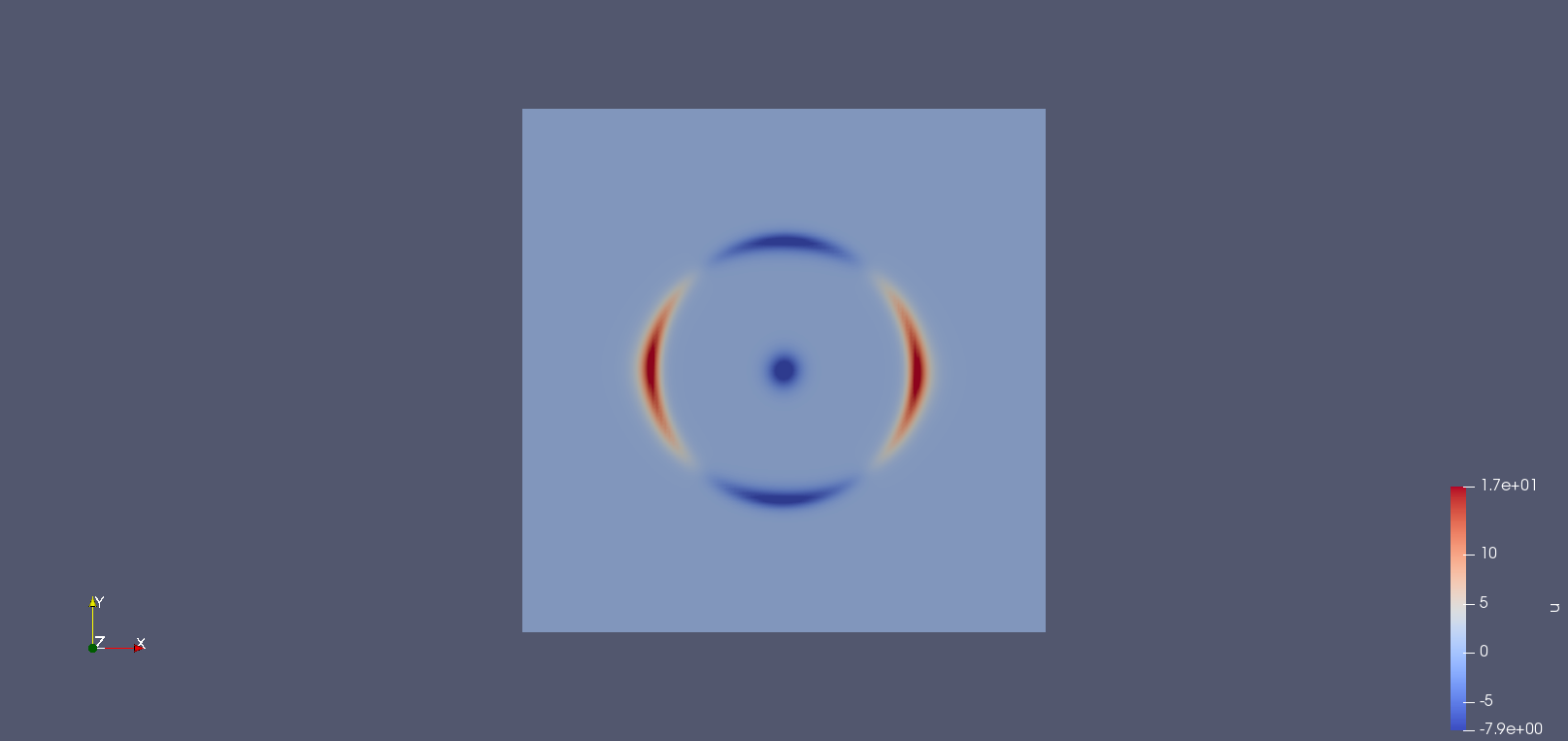}
	\includegraphics[trim={540px 112px 540px 112px}, clip, scale=.08]{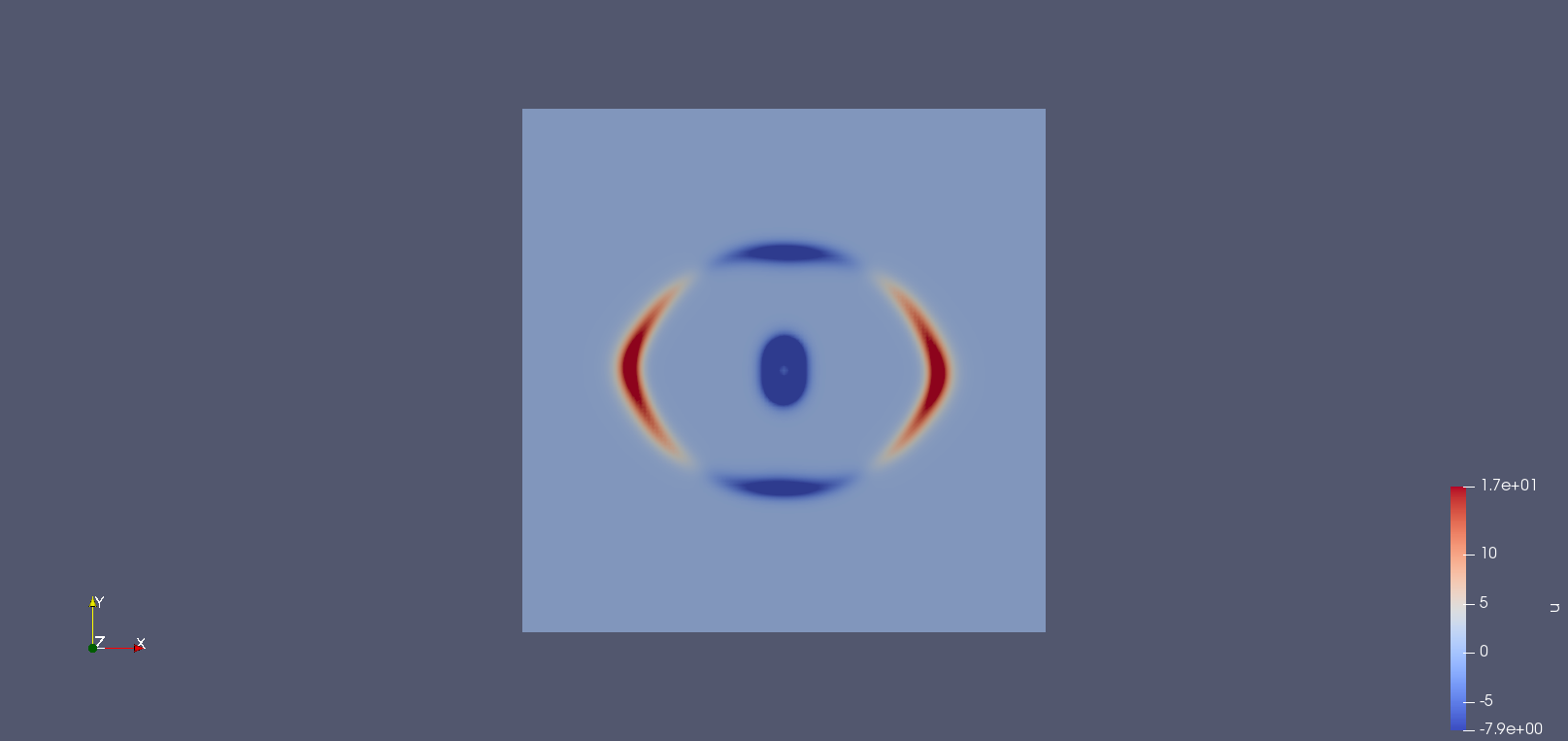}
	\includegraphics[trim={540px 112px 540px 112px}, clip, scale=.08]{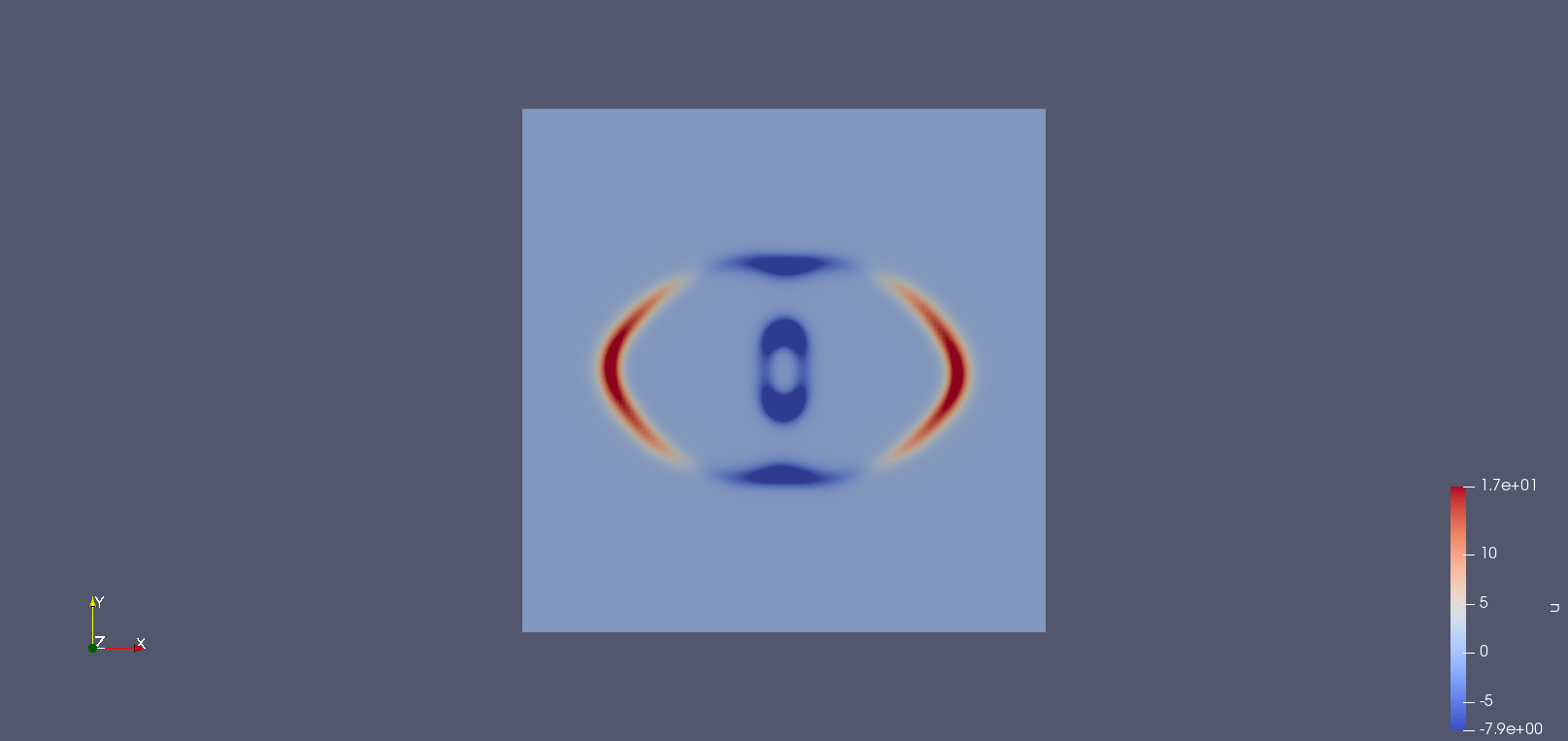}
	\includegraphics[trim={540px 112px 540px 112px}, clip, scale=.08]{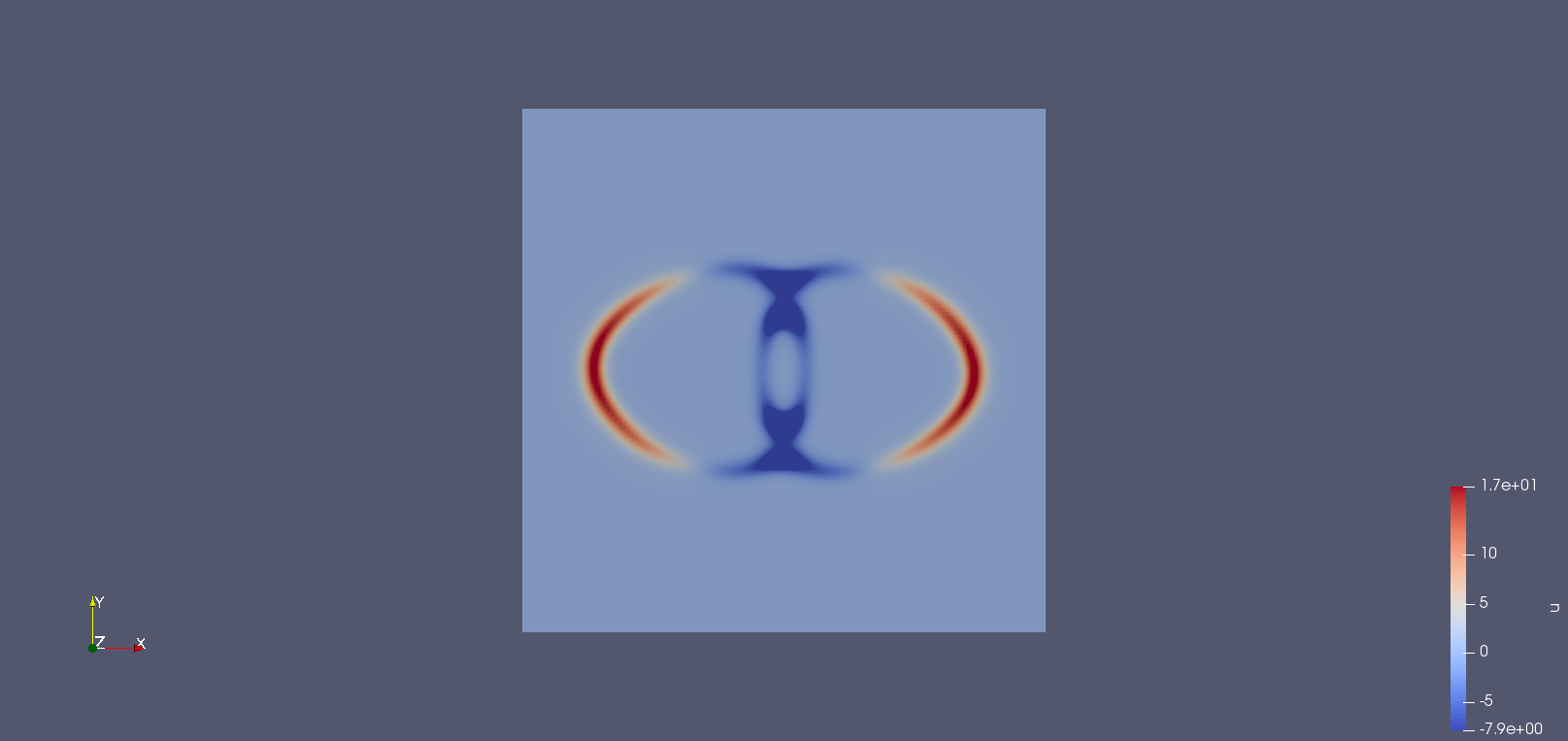}
	\includegraphics[trim={540px 112px 540px 112px}, clip, scale=.08]{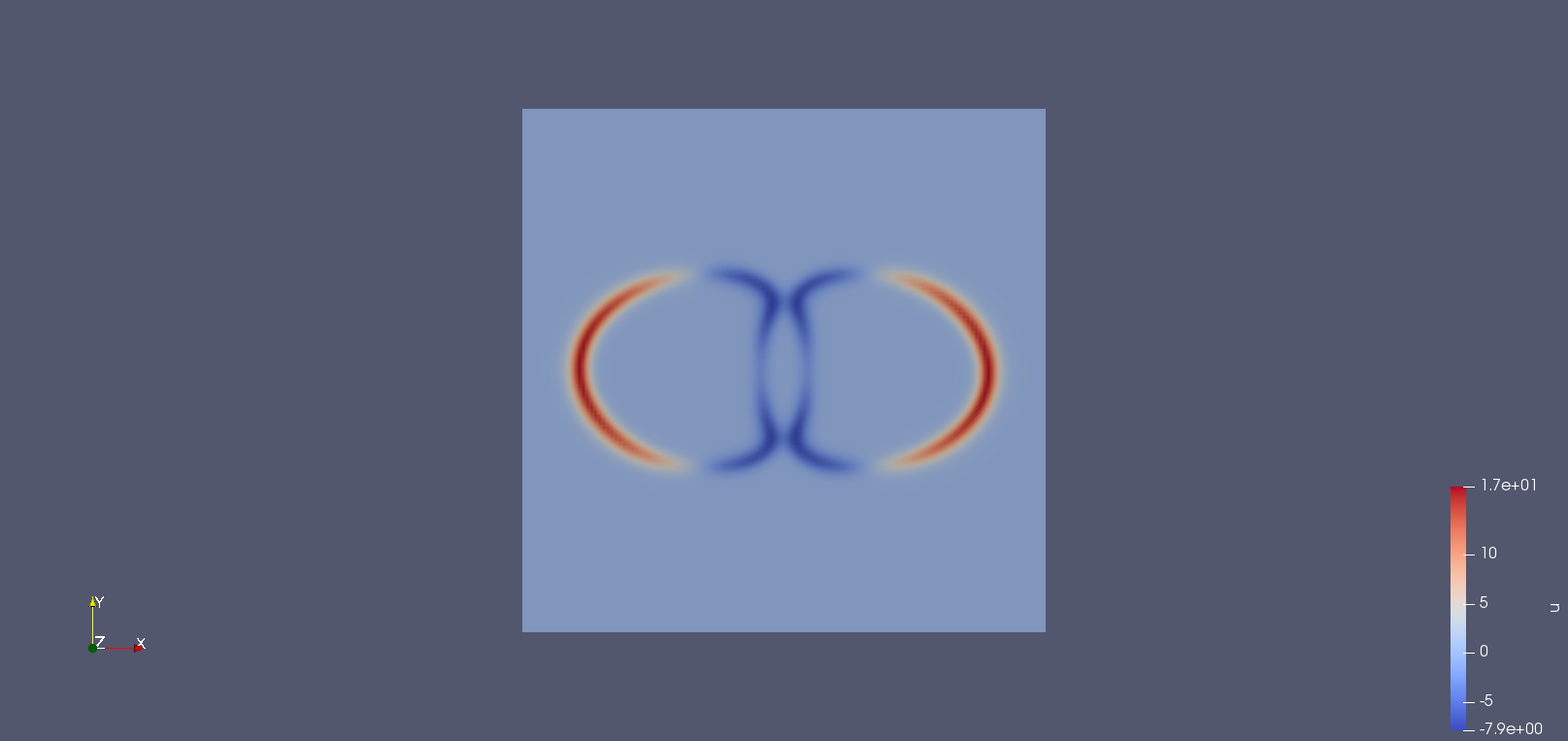}
	\includegraphics[trim={540px 112px 540px 112px}, clip, scale=.08]{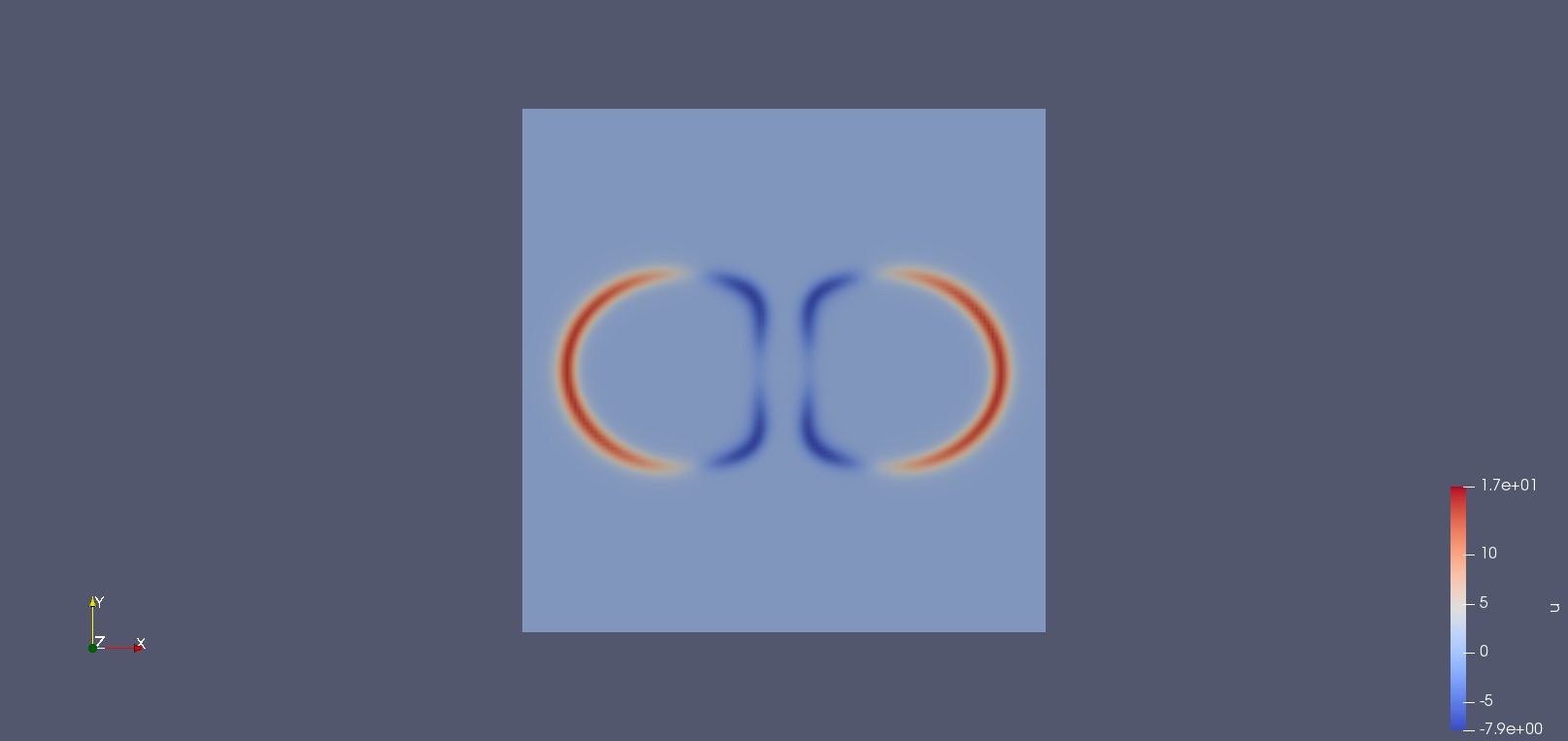}
	\includegraphics[trim={540px 112px 540px 112px}, clip, scale=.08]{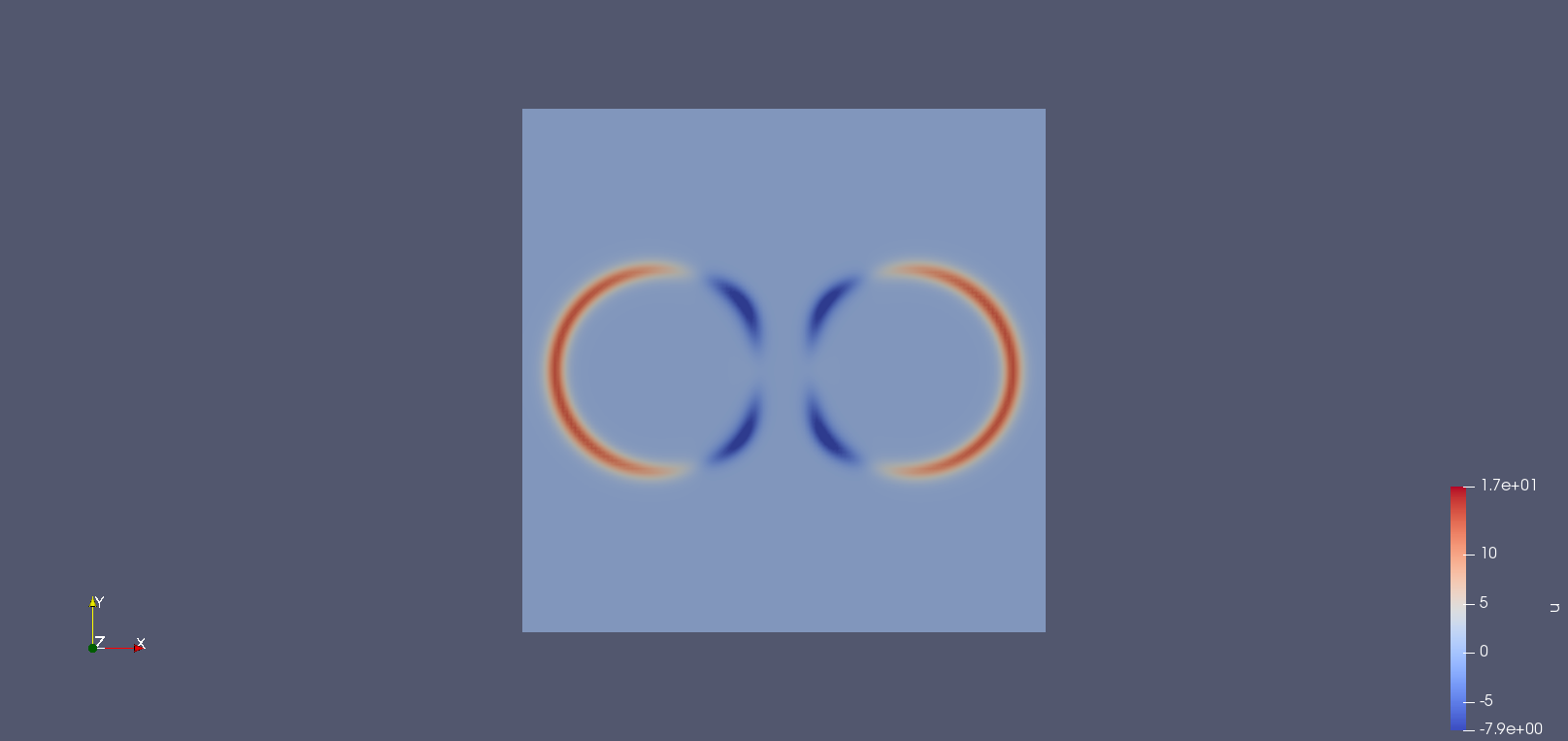}
	\includegraphics[trim={540px 112px 540px 112px}, clip, scale=.08]{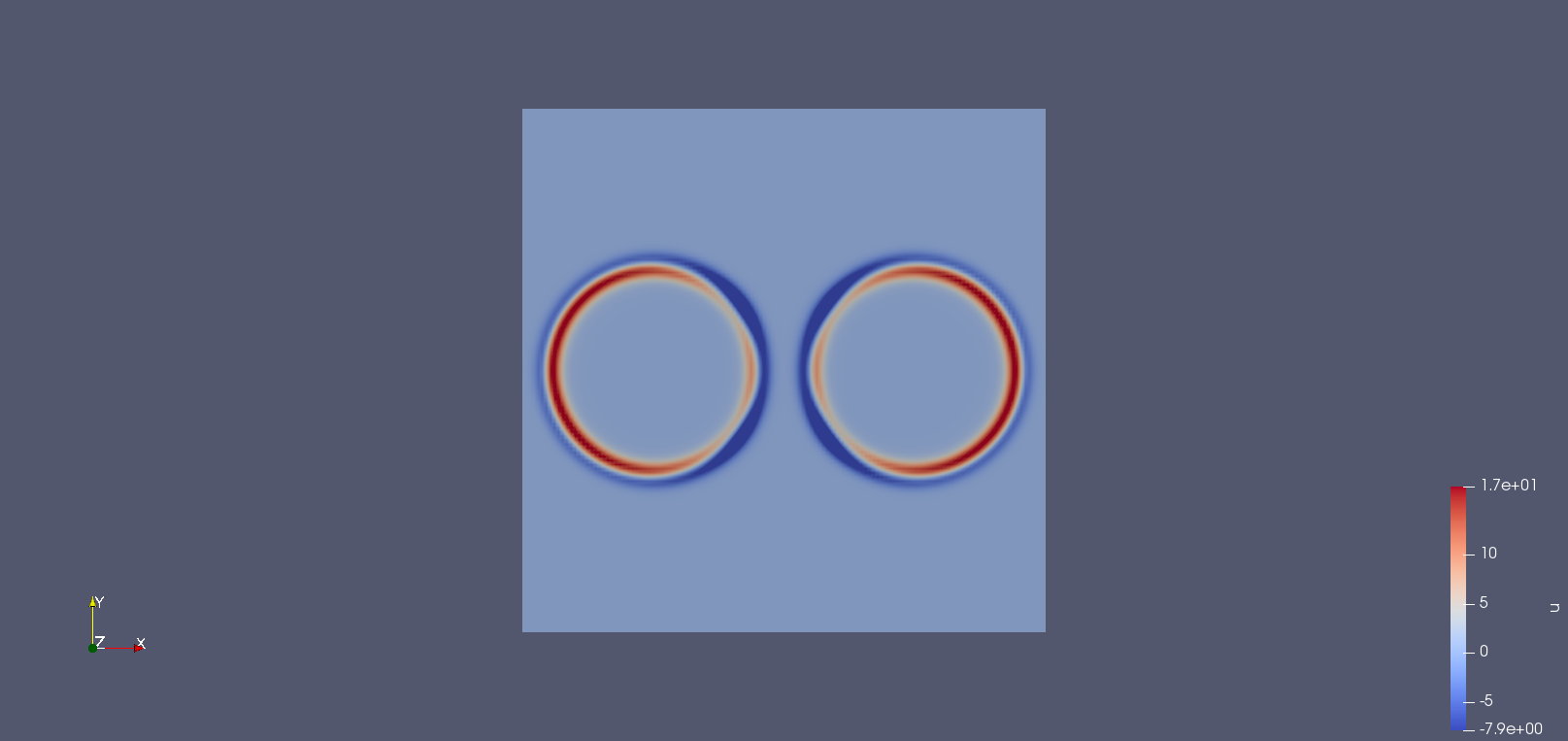}~~~
	\includegraphics[trim={1490px 0px 20px 490px}, clip, scale=.16]{{iso_splitting_control.0000}.png}
\caption{Result for the `splitting circle' in the isotropic case.	\label{fig:iso_splitting}} 
\end{center}    
\end{subfigure}
\medskip

\begin{subfigure}{1.0\textwidth}
\begin{center}   
	\includegraphics[trim={540px 112px 540px 112px}, clip, scale=.08]{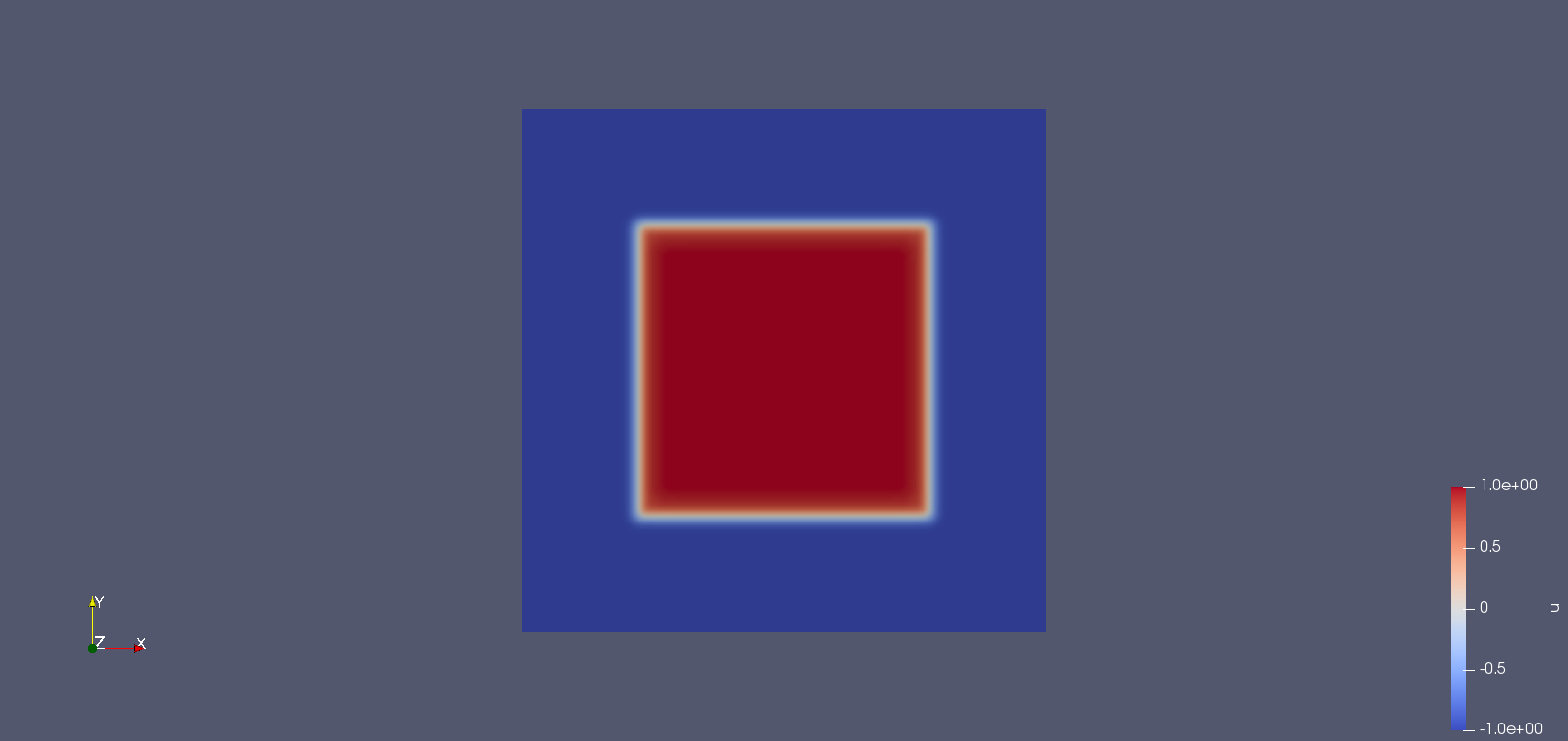}
	\includegraphics[trim={540px 112px 540px 112px}, clip, scale=.08]{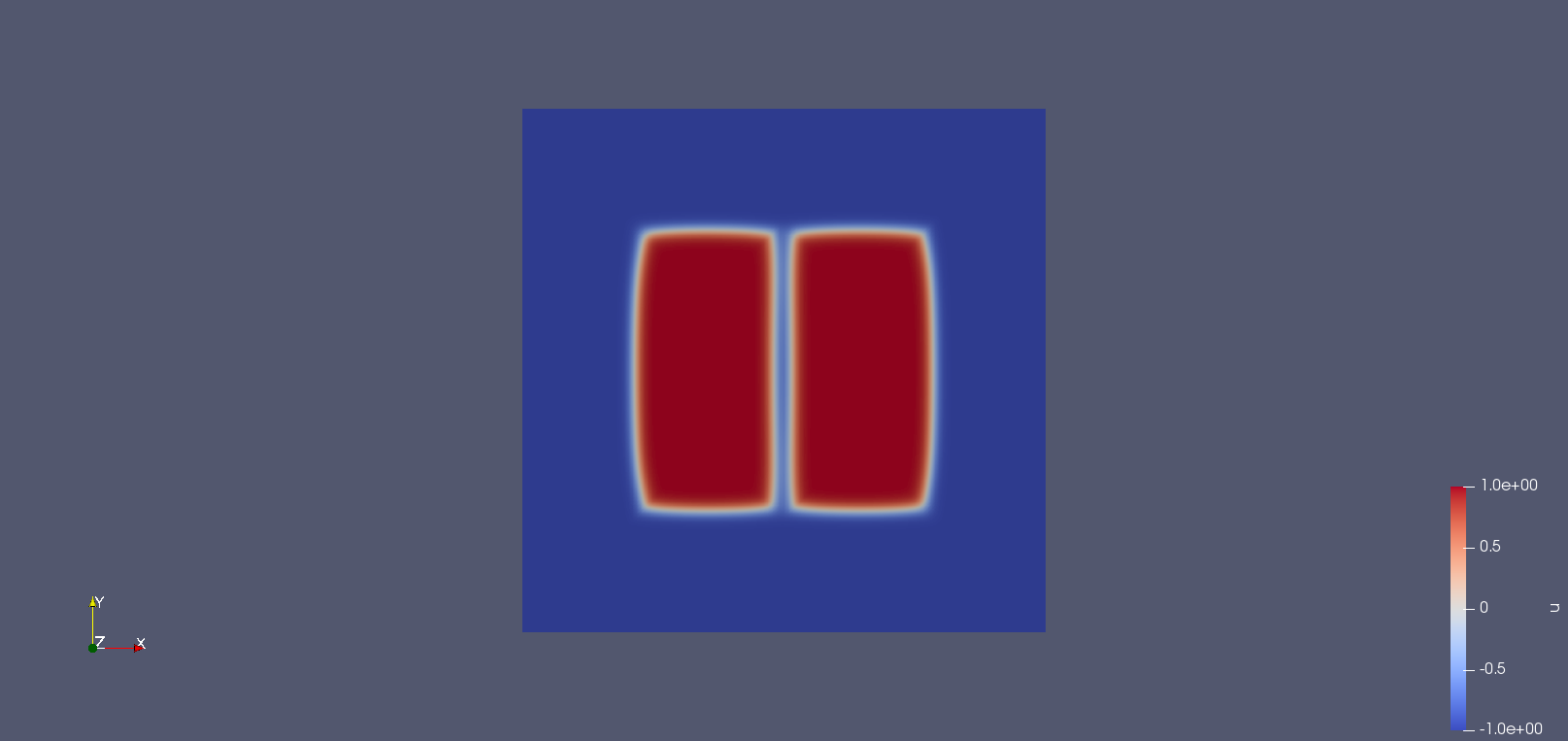}
	\includegraphics[trim={540px 112px 540px 112px}, clip, scale=.08]{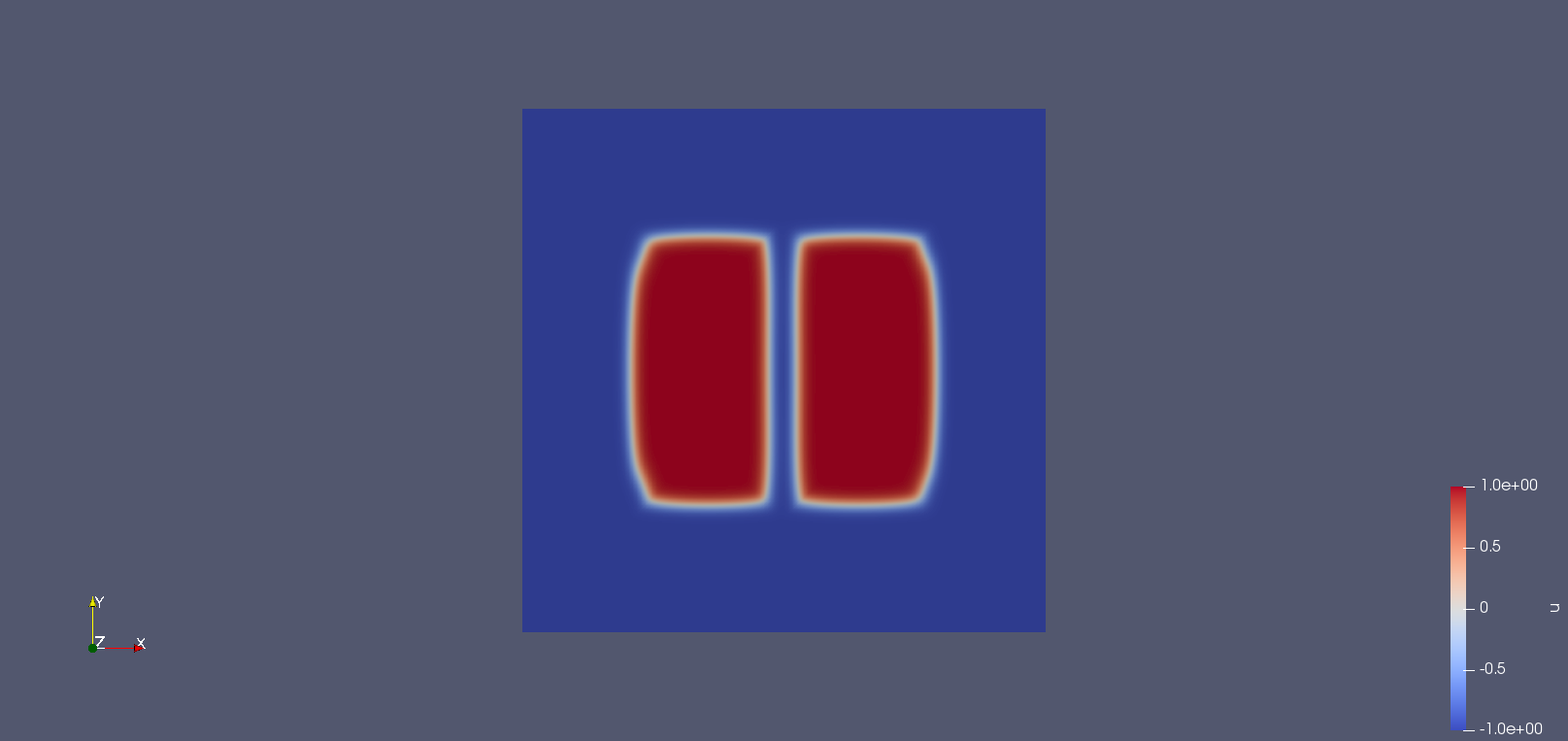}
	\includegraphics[trim={540px 112px 540px 112px}, clip, scale=.08]{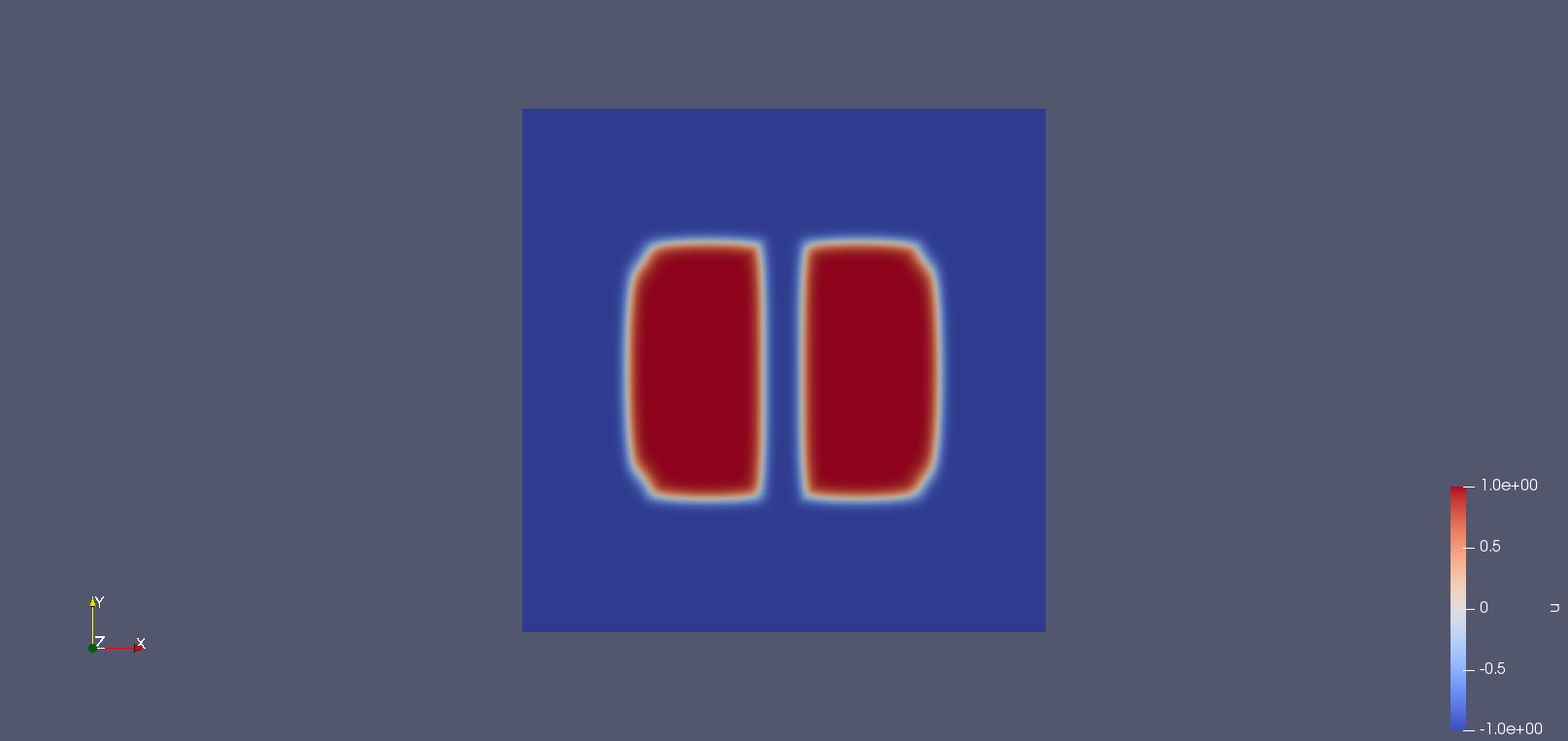}
	\includegraphics[trim={540px 112px 540px 112px}, clip, scale=.08]{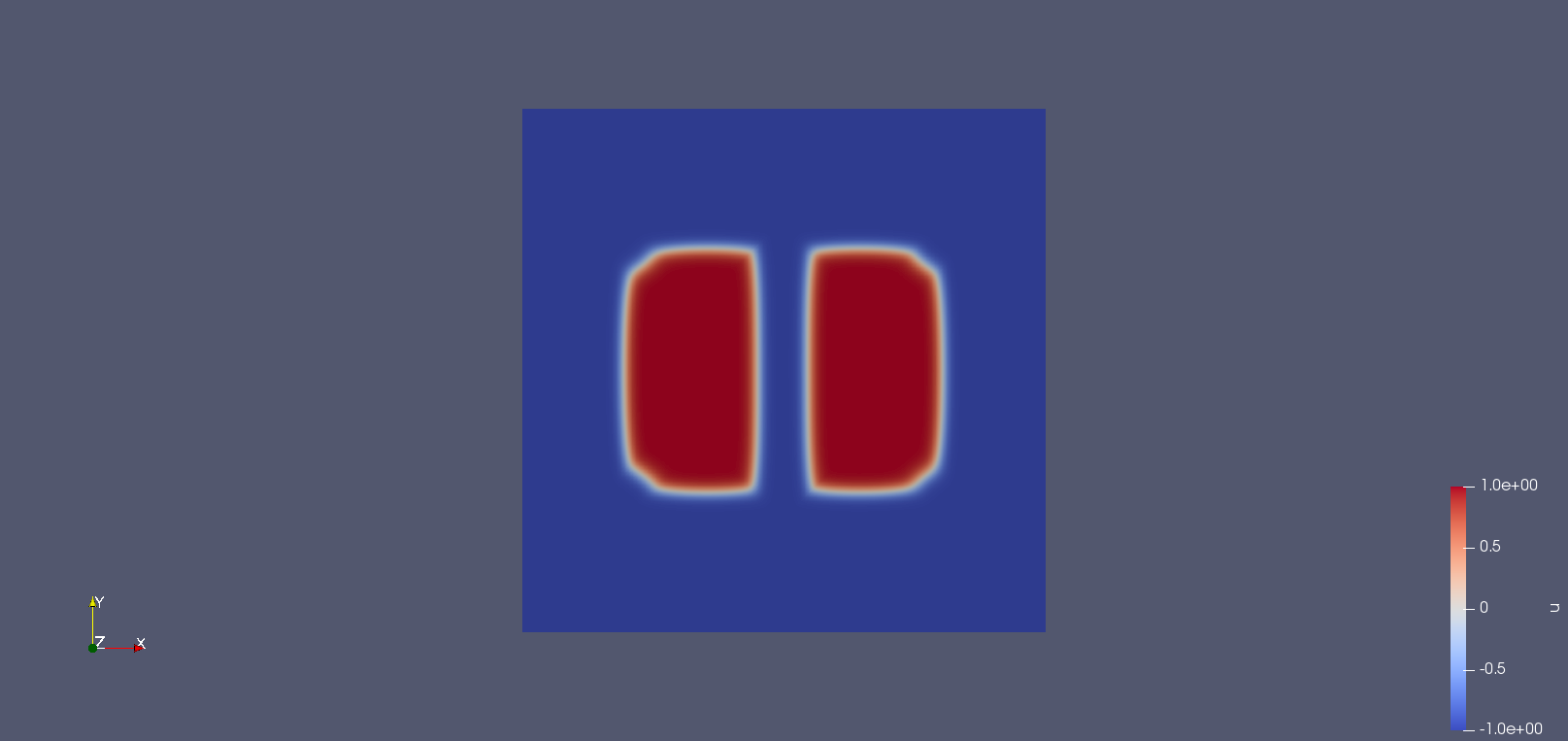}
	\includegraphics[trim={540px 112px 540px 112px}, clip, scale=.08]{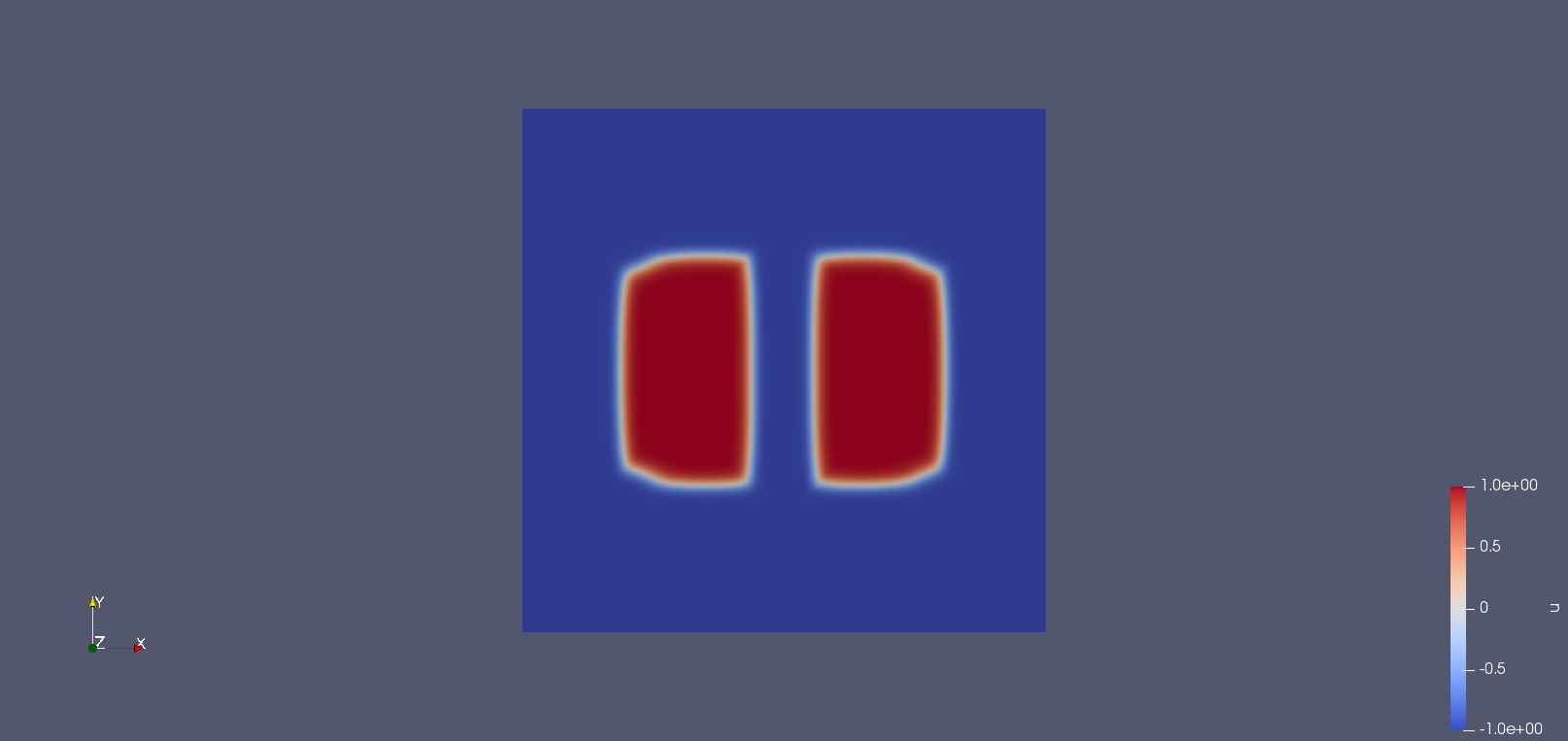}
	\includegraphics[trim={540px 112px 540px 112px}, clip, scale=.08]{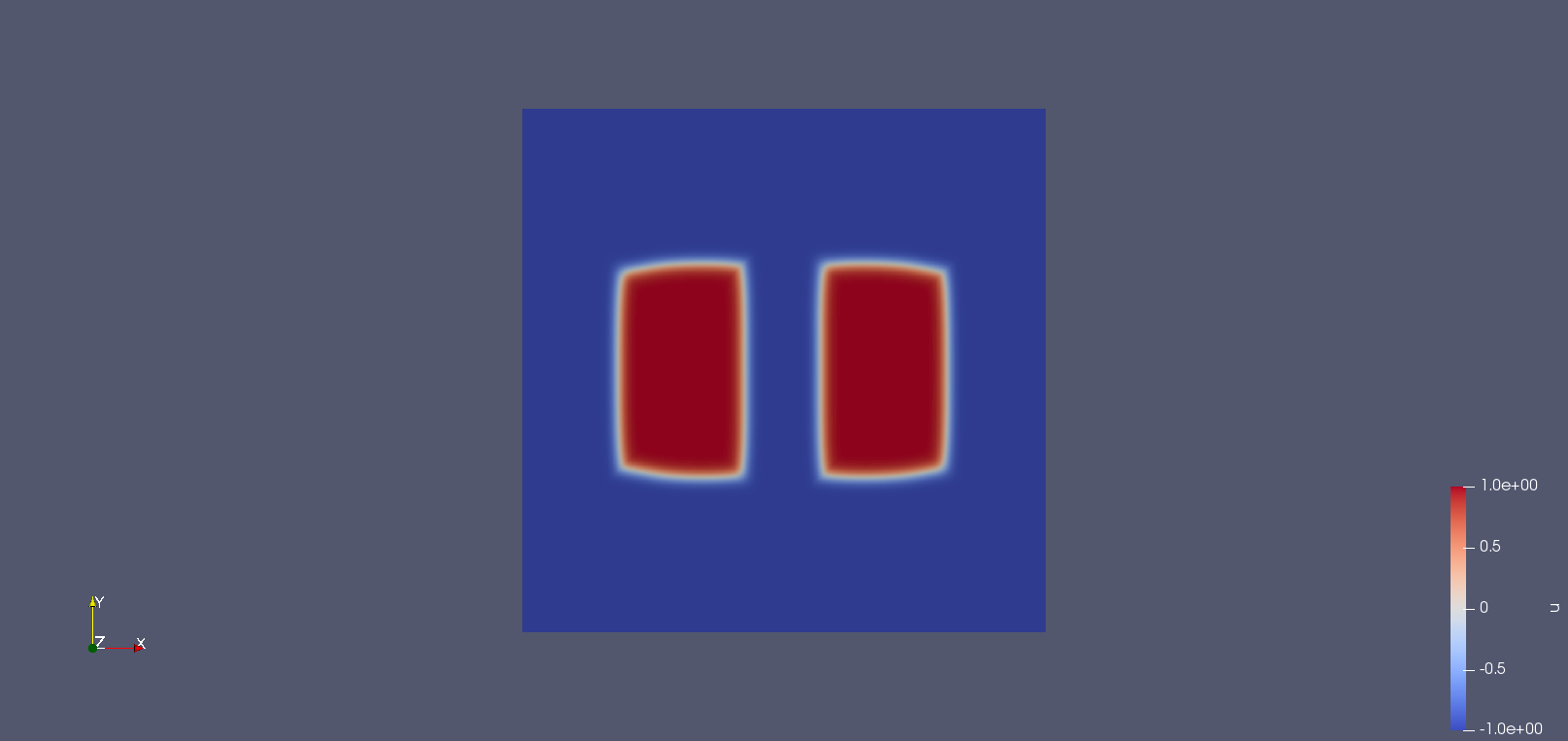}
	\includegraphics[trim={540px 112px 540px 112px}, clip, scale=.08]{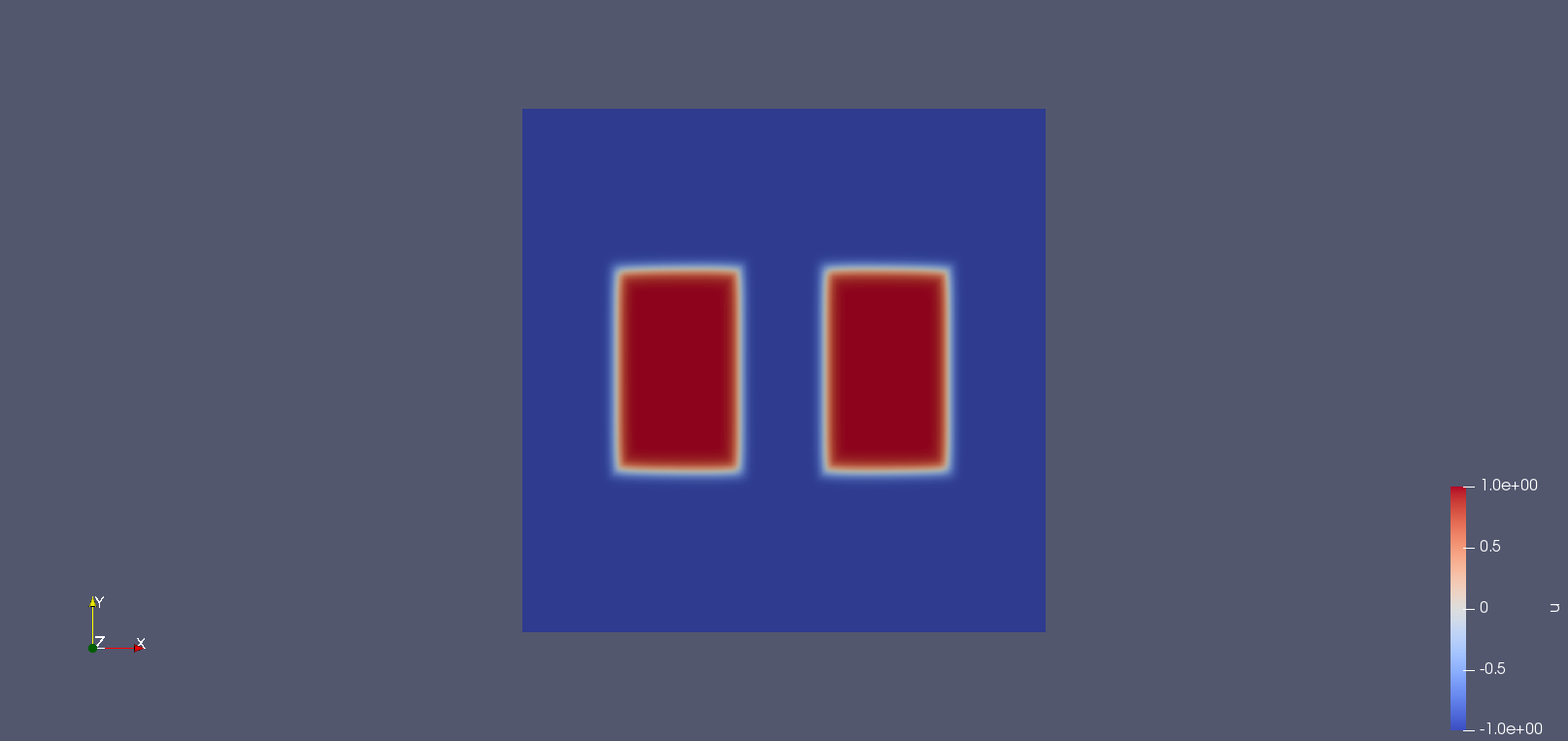}~~~
	\includegraphics[trim={1490px 0px 20px 490px}, clip, scale=.16]{{l1_splitting_state.0000}.png}~\\~\\
	\includegraphics[trim={540px 112px 540px 112px}, clip, scale=.08]{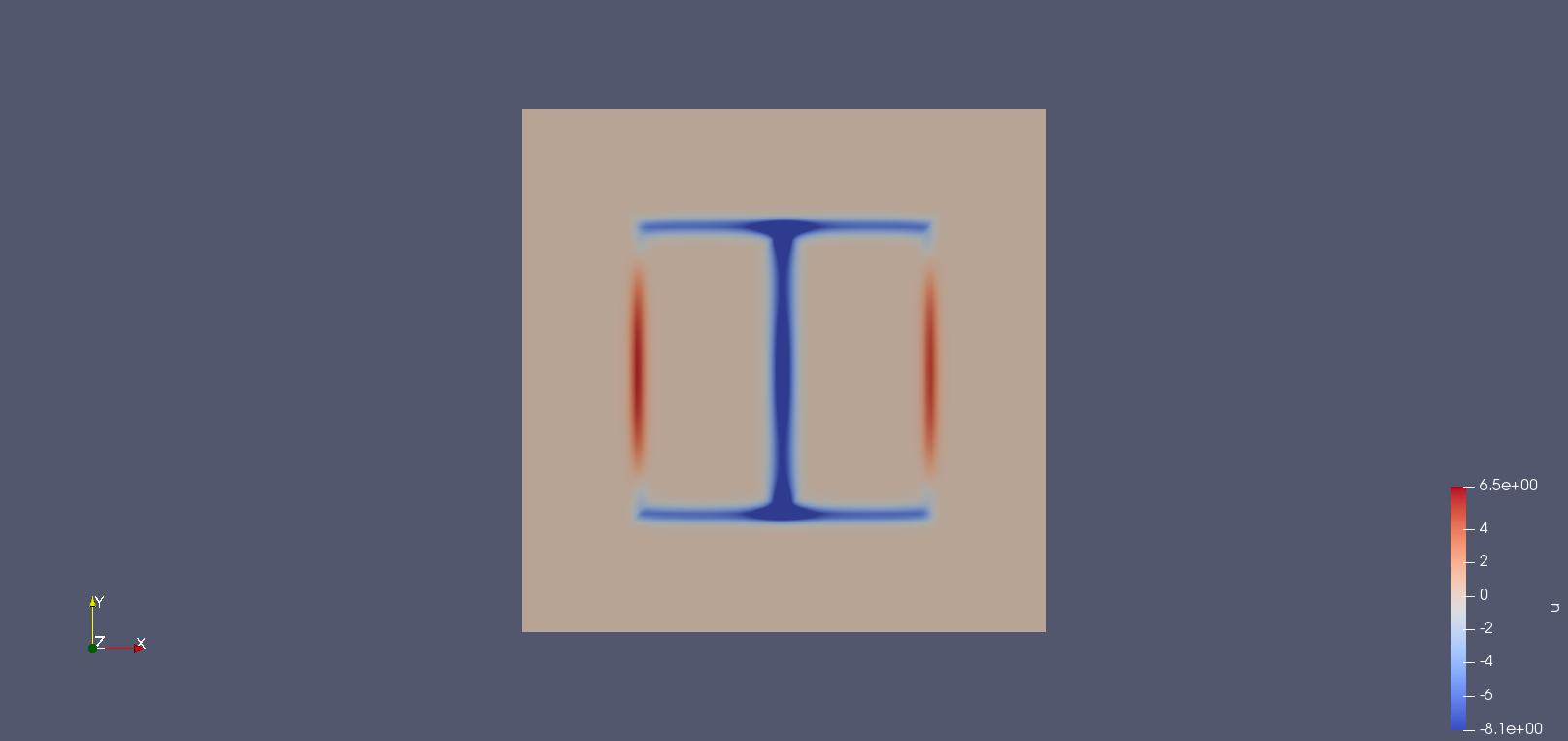}
	\includegraphics[trim={540px 112px 540px 112px}, clip, scale=.08]{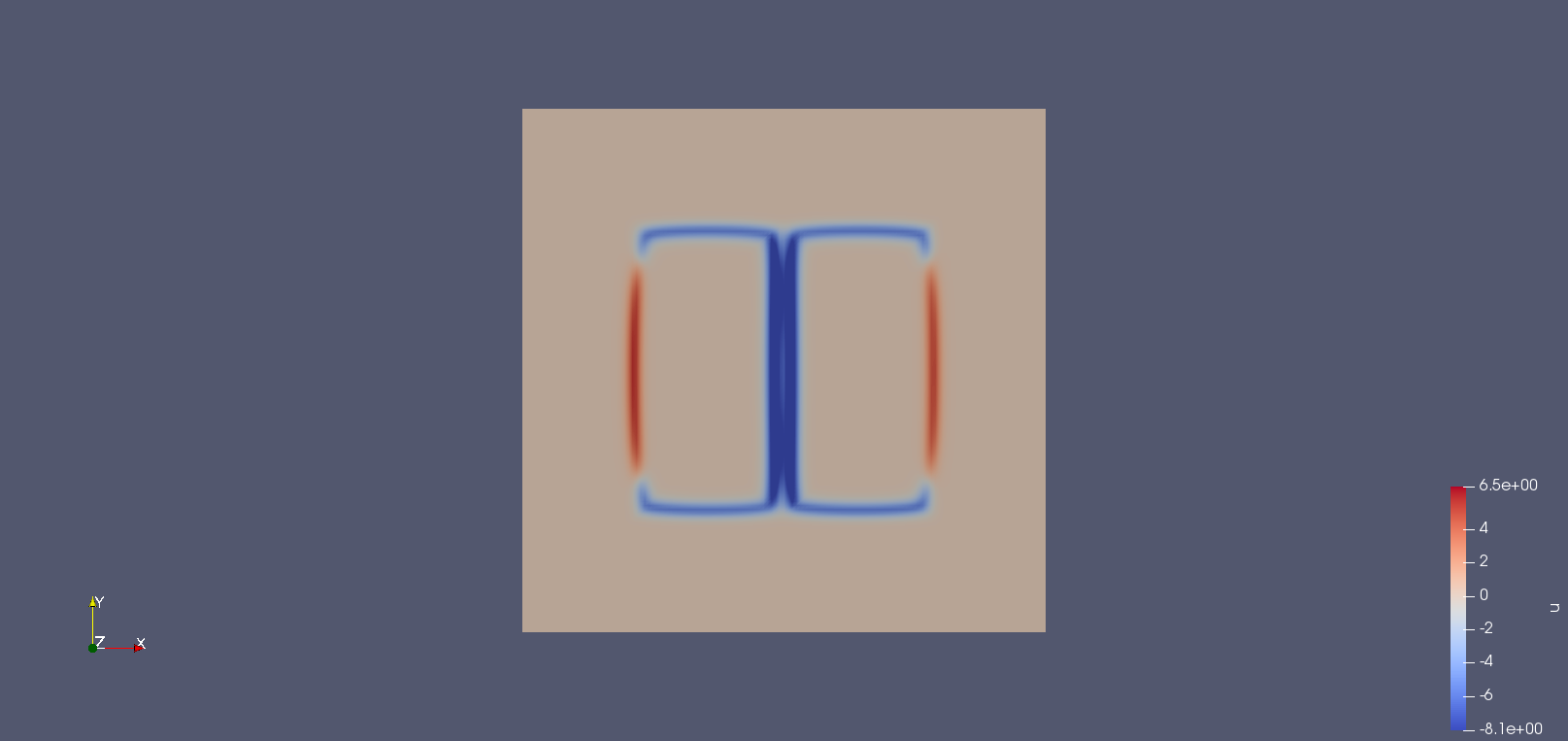}
	\includegraphics[trim={540px 112px 540px 112px}, clip, scale=.08]{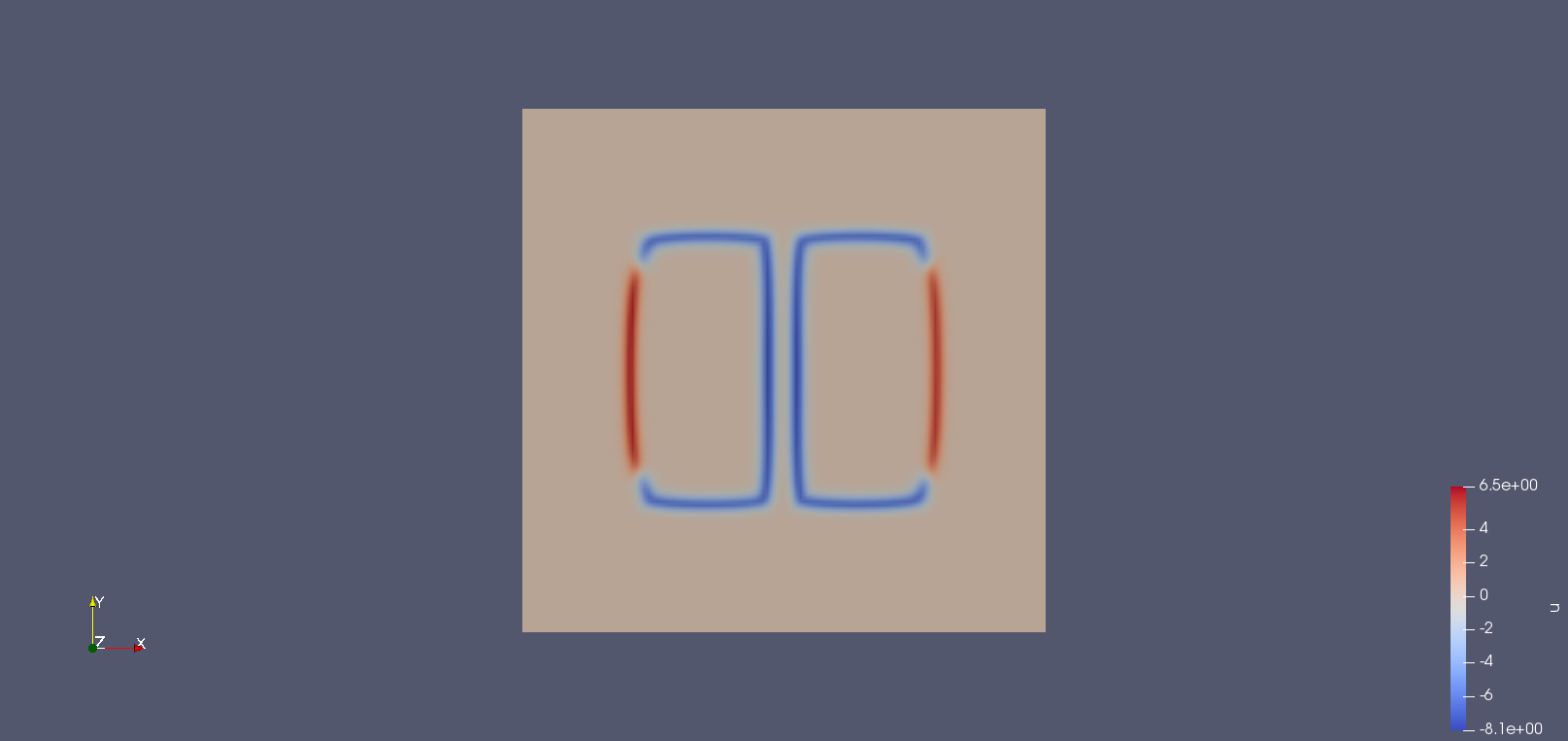}
	\includegraphics[trim={540px 112px 540px 112px}, clip, scale=.08]{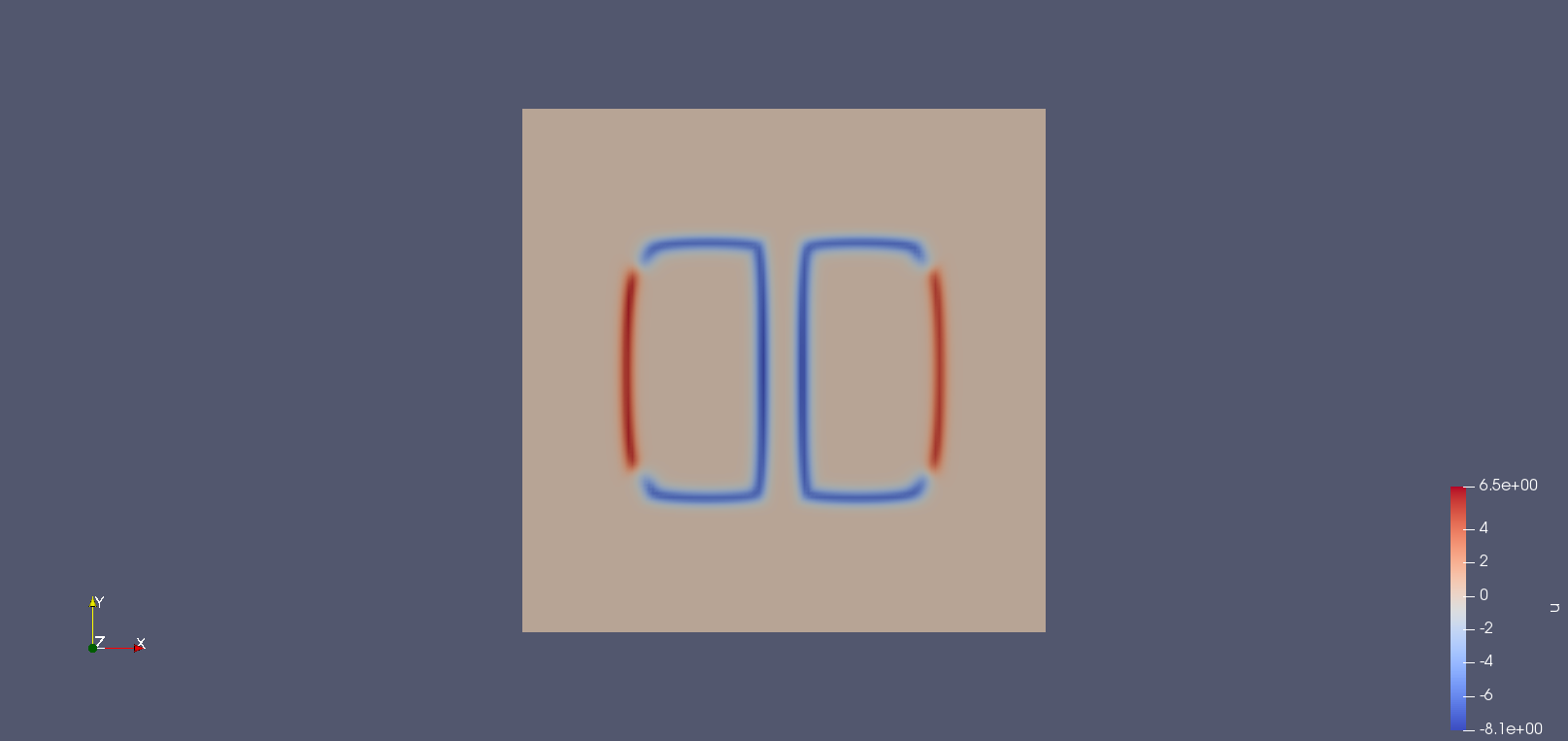}
	\includegraphics[trim={540px 112px 540px 112px}, clip, scale=.08]{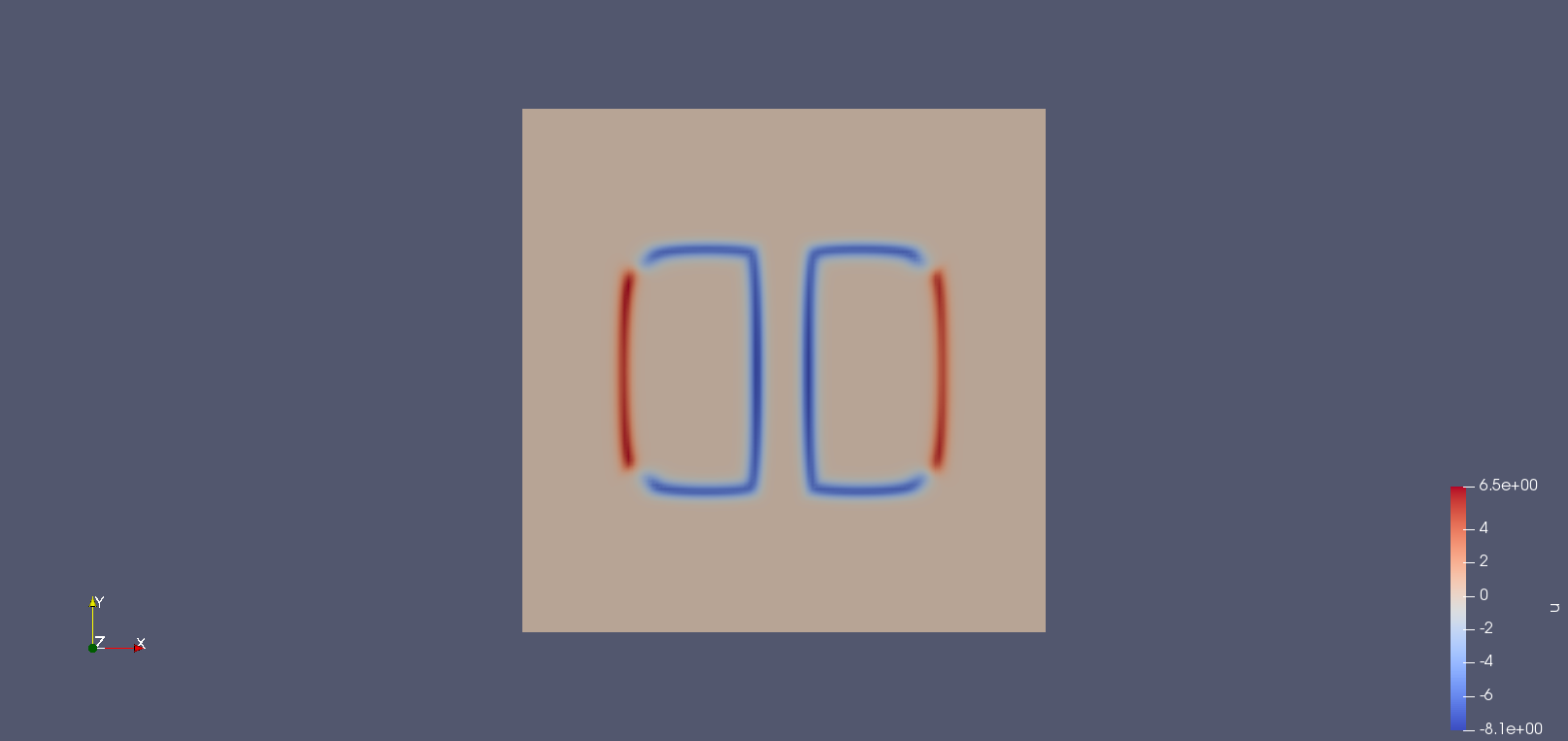}
	\includegraphics[trim={540px 112px 540px 112px}, clip, scale=.08]{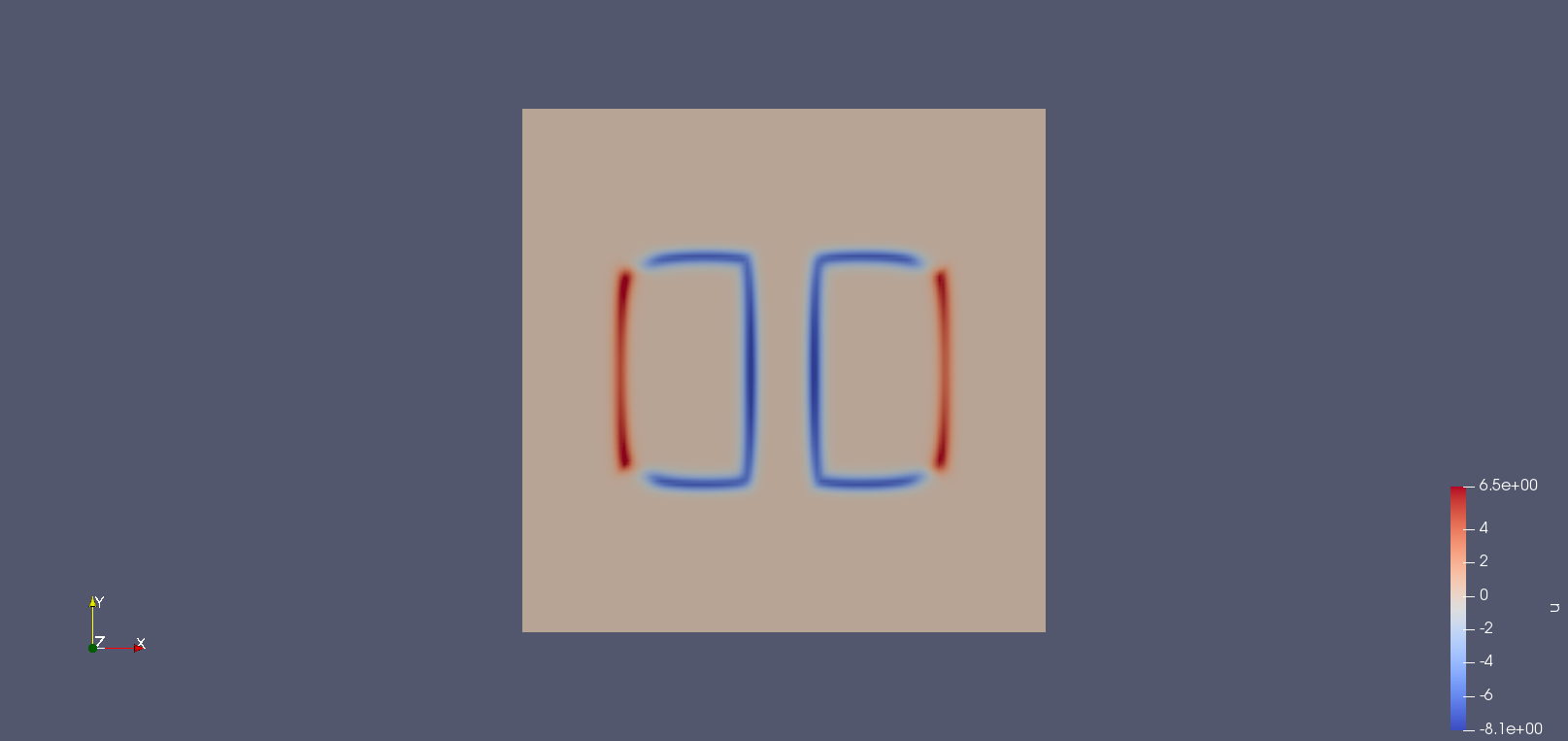}
	\includegraphics[trim={540px 112px 540px 112px}, clip, scale=.08]{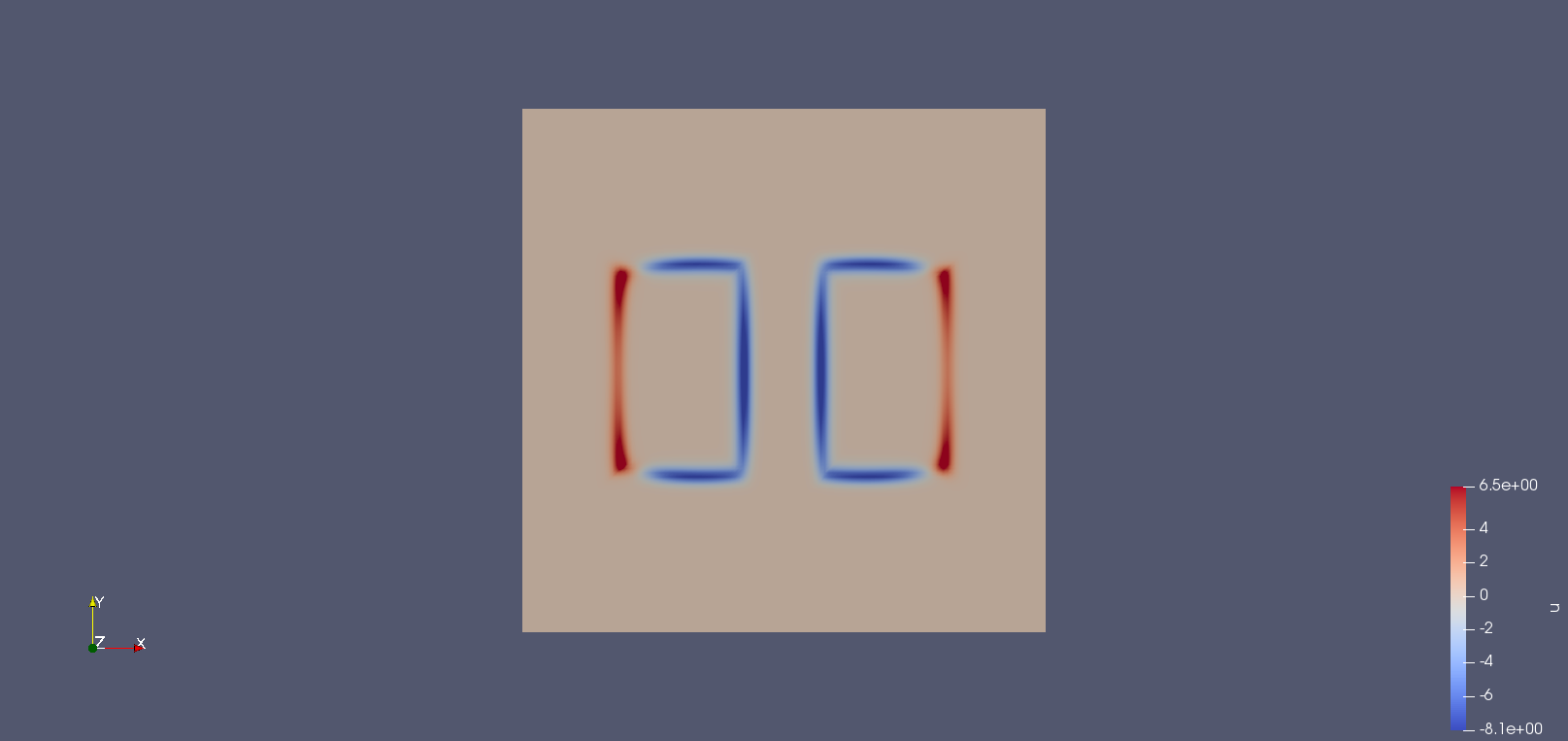}
	\includegraphics[trim={540px 112px 540px 112px}, clip, scale=.08]{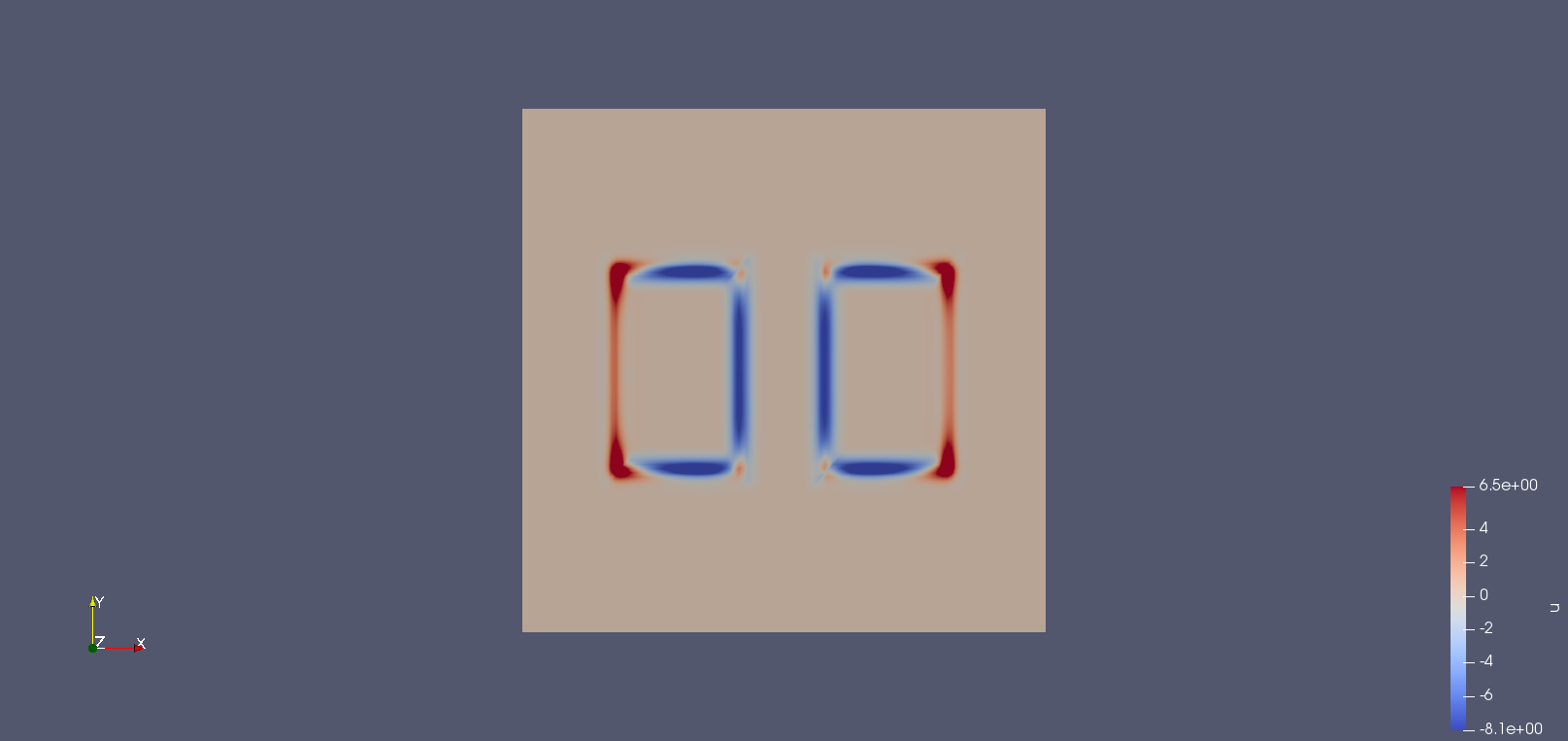}~~~
	\includegraphics[trim={1490px 0px 20px 490px}, clip, scale=.16]{{l1_splitting_control.0000}.png}
	\caption{Result for the `splitting square' with the regularized $l_1$-norm.
	\label{fig:l1_splitting}}
\end{center}    
\end{subfigure}
\medskip

\begin{subfigure}{1.0\textwidth}
\begin{center}  
	\includegraphics[trim={540px 112px 540px 112px}, clip, scale=.08]{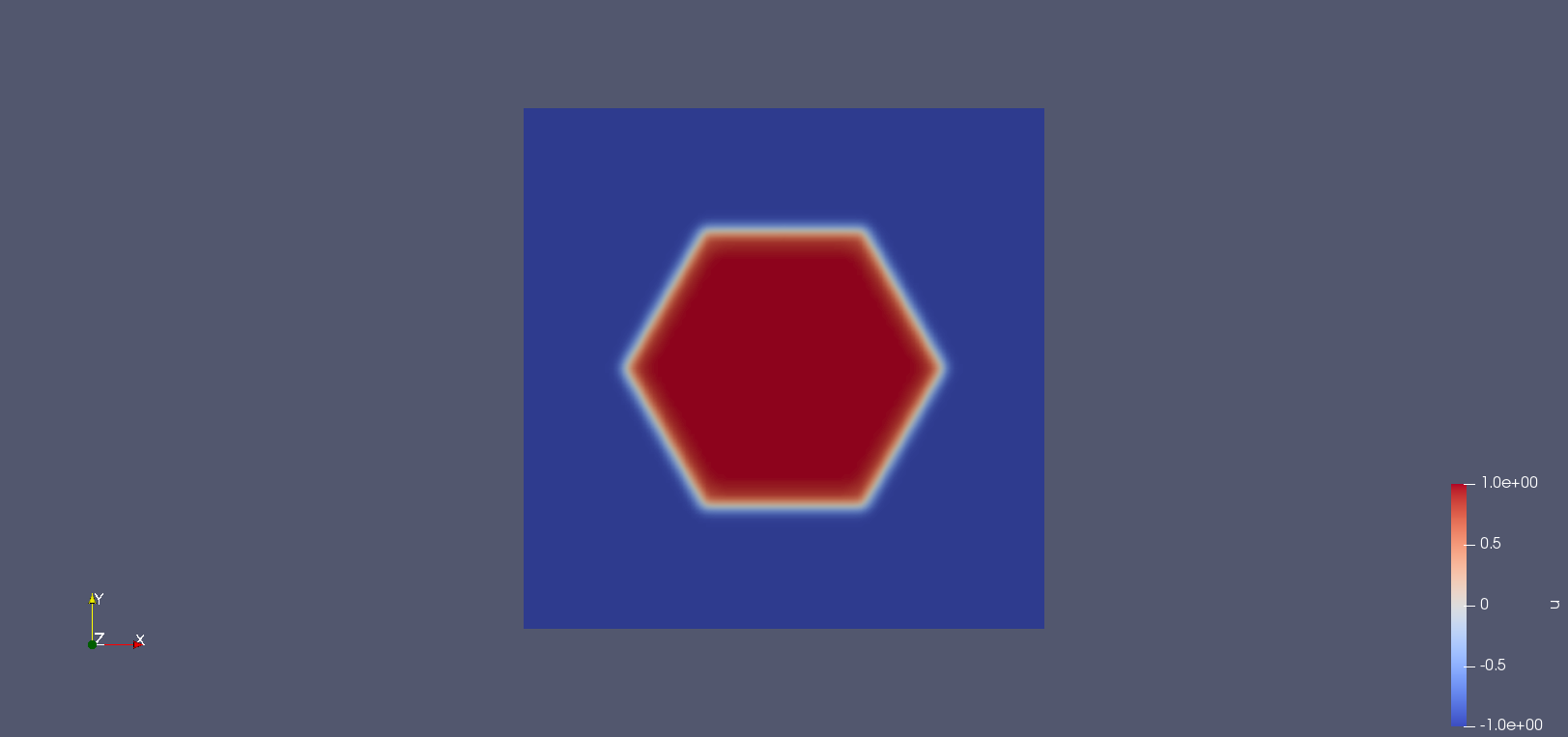}
	\includegraphics[trim={540px 112px 540px 112px}, clip, scale=.08]{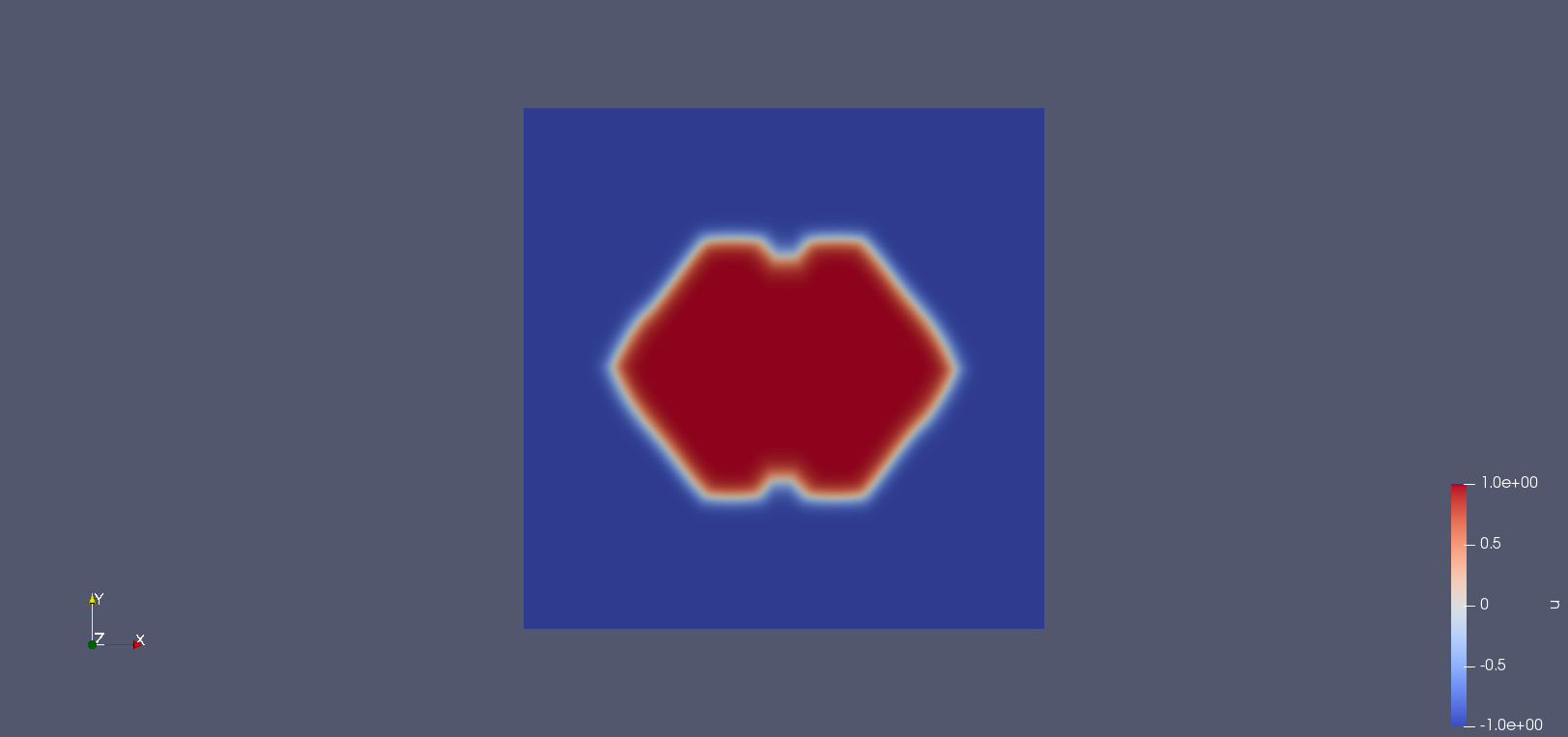}
	\includegraphics[trim={540px 112px 540px 112px}, clip, scale=.08]{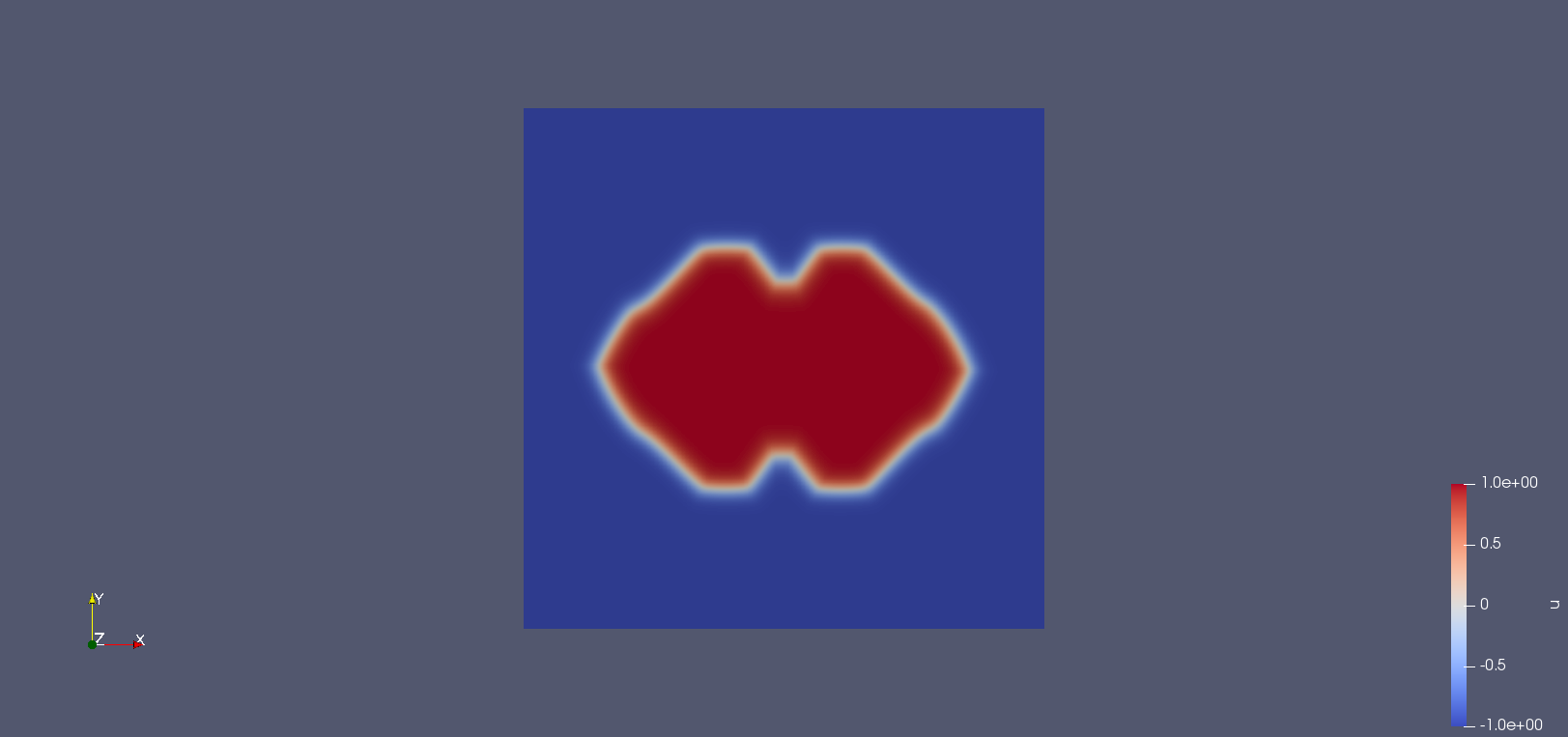}
	\includegraphics[trim={540px 112px 540px 112px}, clip, scale=.08]{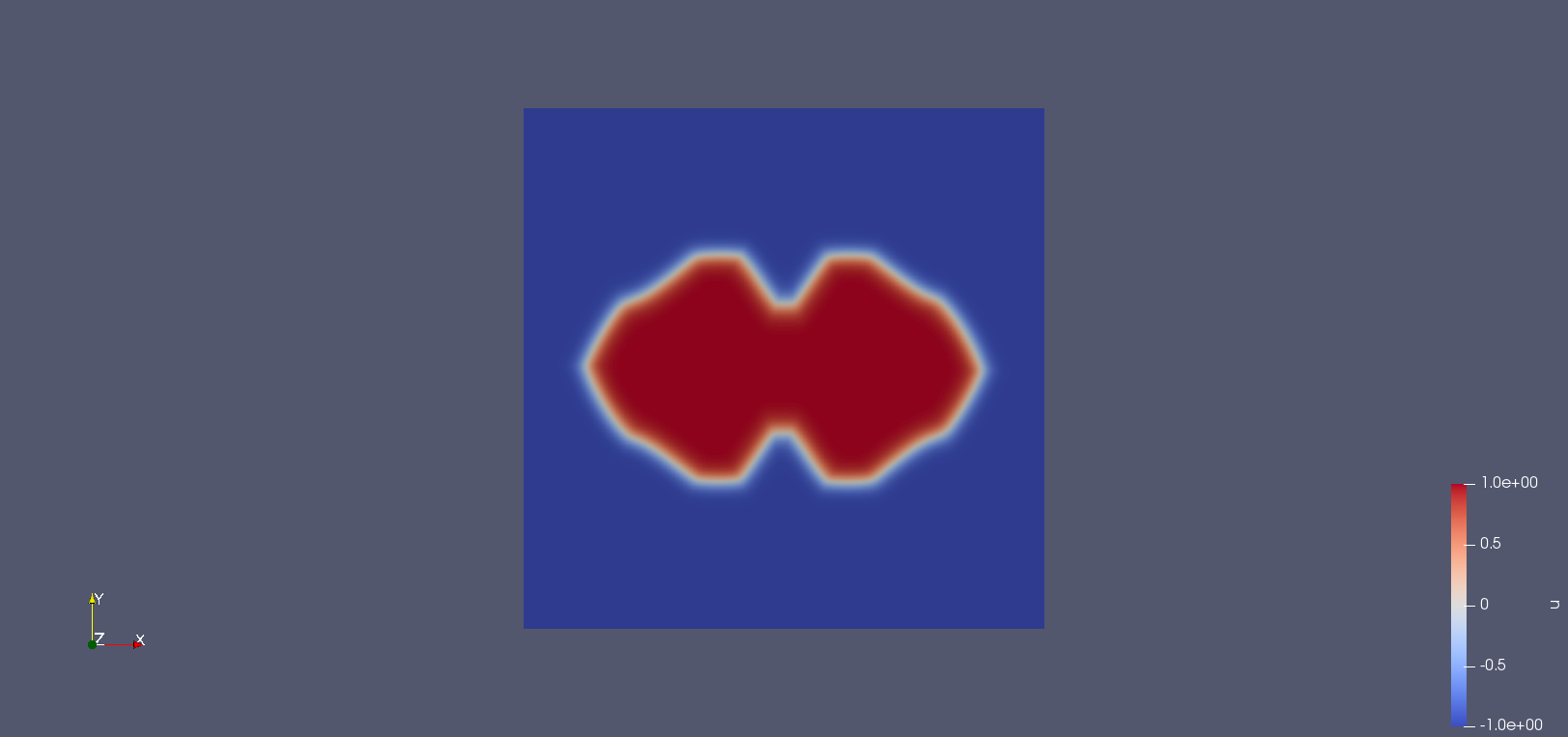}
	\includegraphics[trim={540px 112px 540px 112px}, clip, scale=.08]{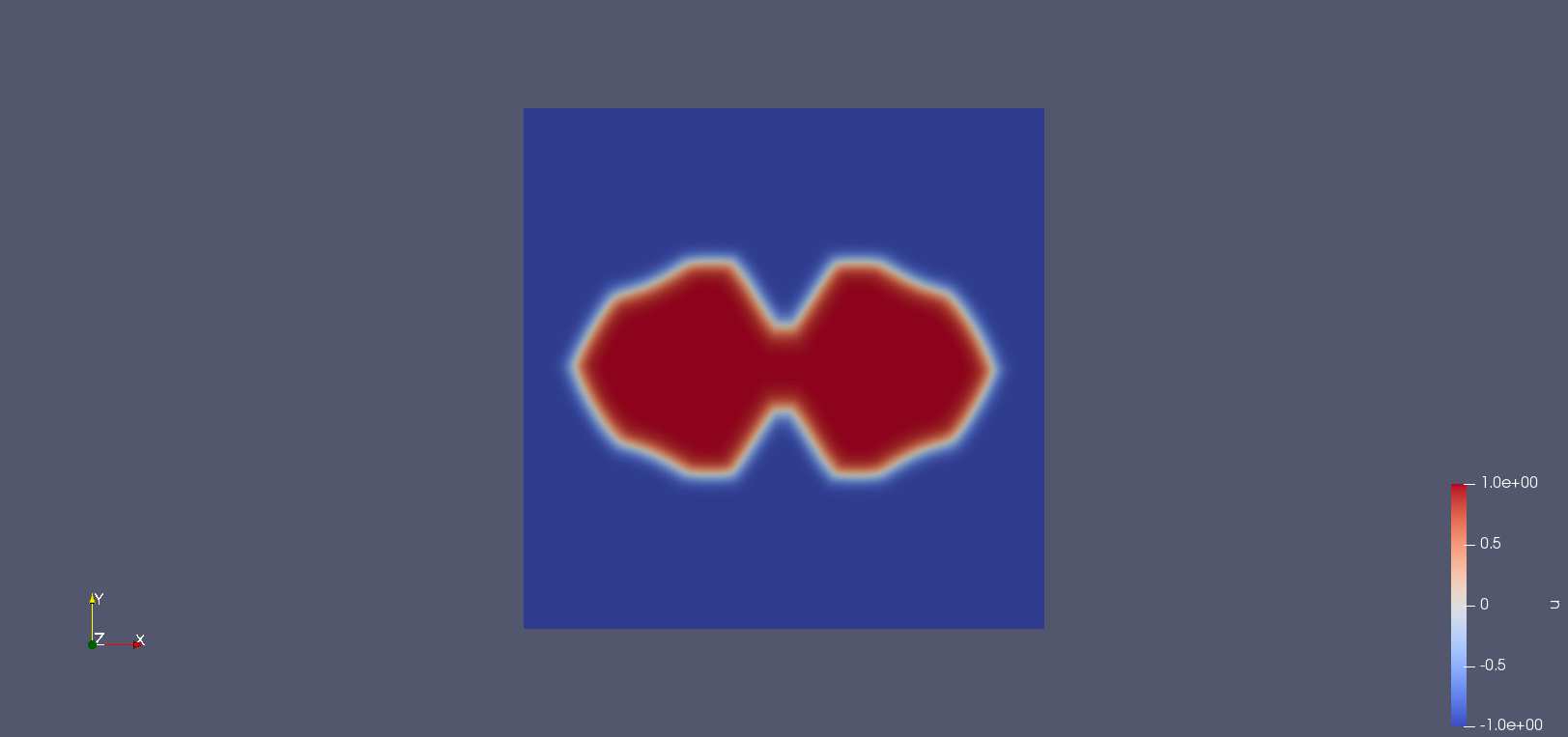}
	\includegraphics[trim={540px 112px 540px 112px}, clip, scale=.08]{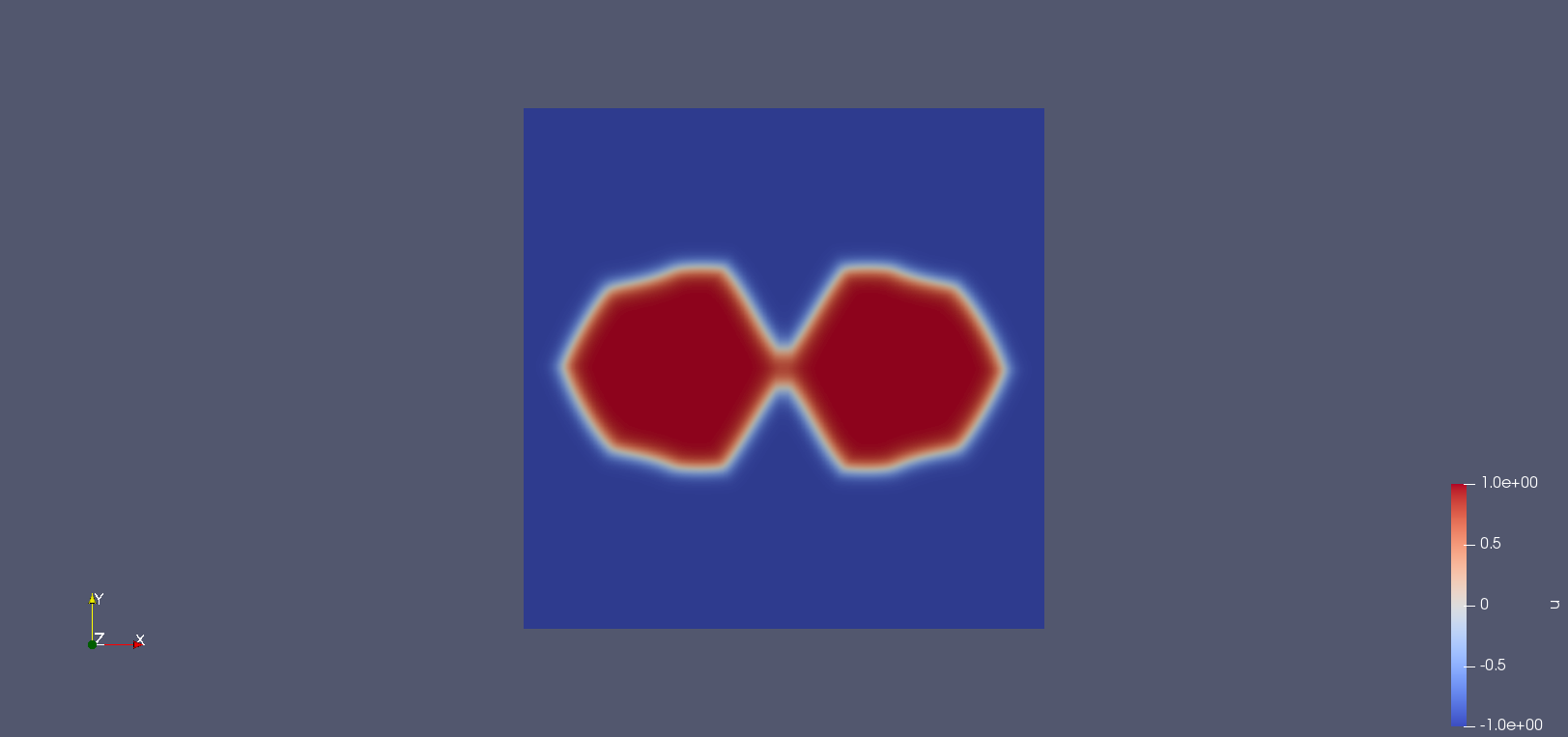}
	\includegraphics[trim={540px 112px 540px 112px}, clip, scale=.08]{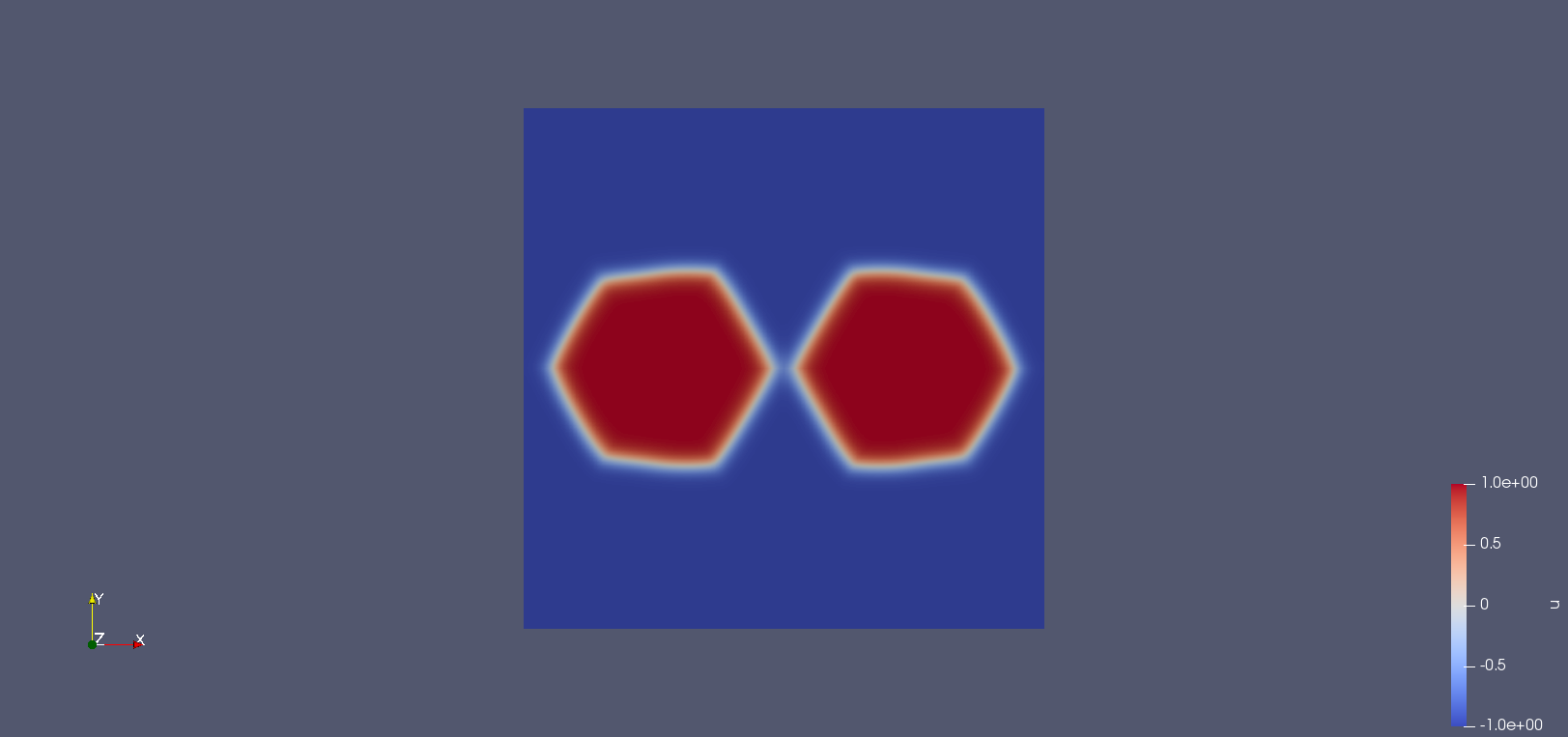}
	\includegraphics[trim={540px 112px 540px 112px}, clip, scale=.08]{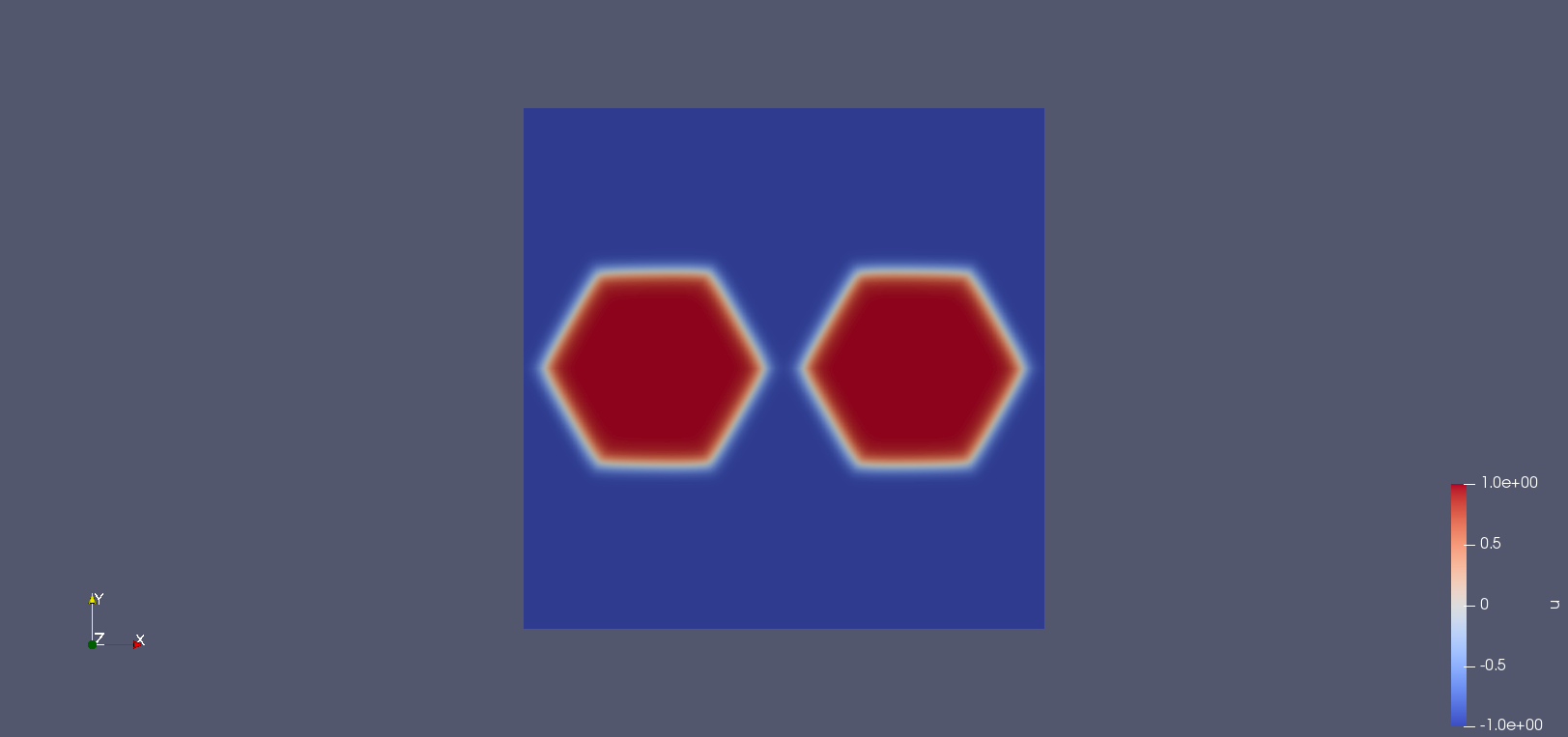}~~~
	\includegraphics[trim={1490px 0px 20px 490px}, clip, scale=.16]{{hexa_splitting_state.0000}.png}~\\~\\
	\includegraphics[trim={540px 112px 540px 112px}, clip, scale=.08]{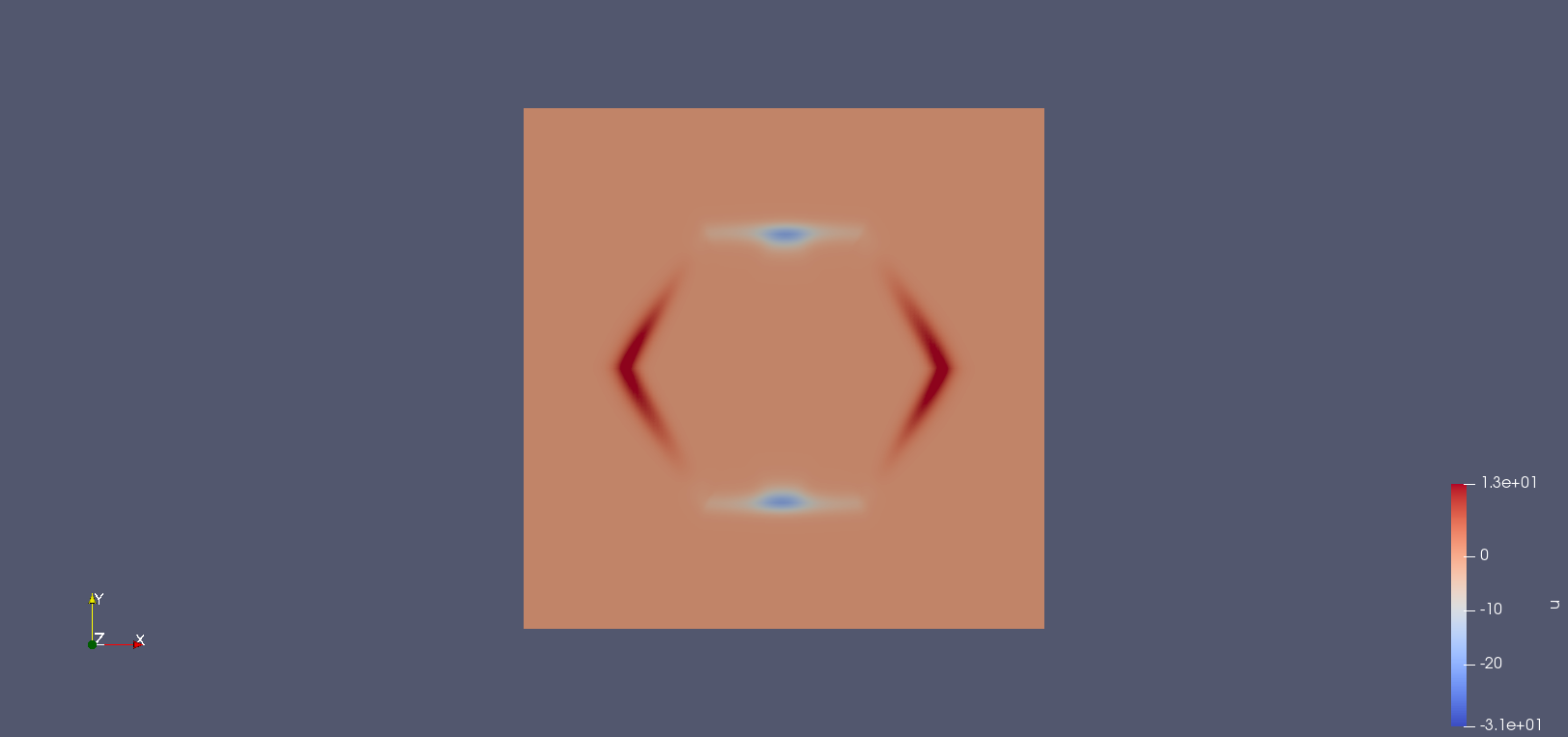}
	\includegraphics[trim={540px 112px 540px 112px}, clip, scale=.08]{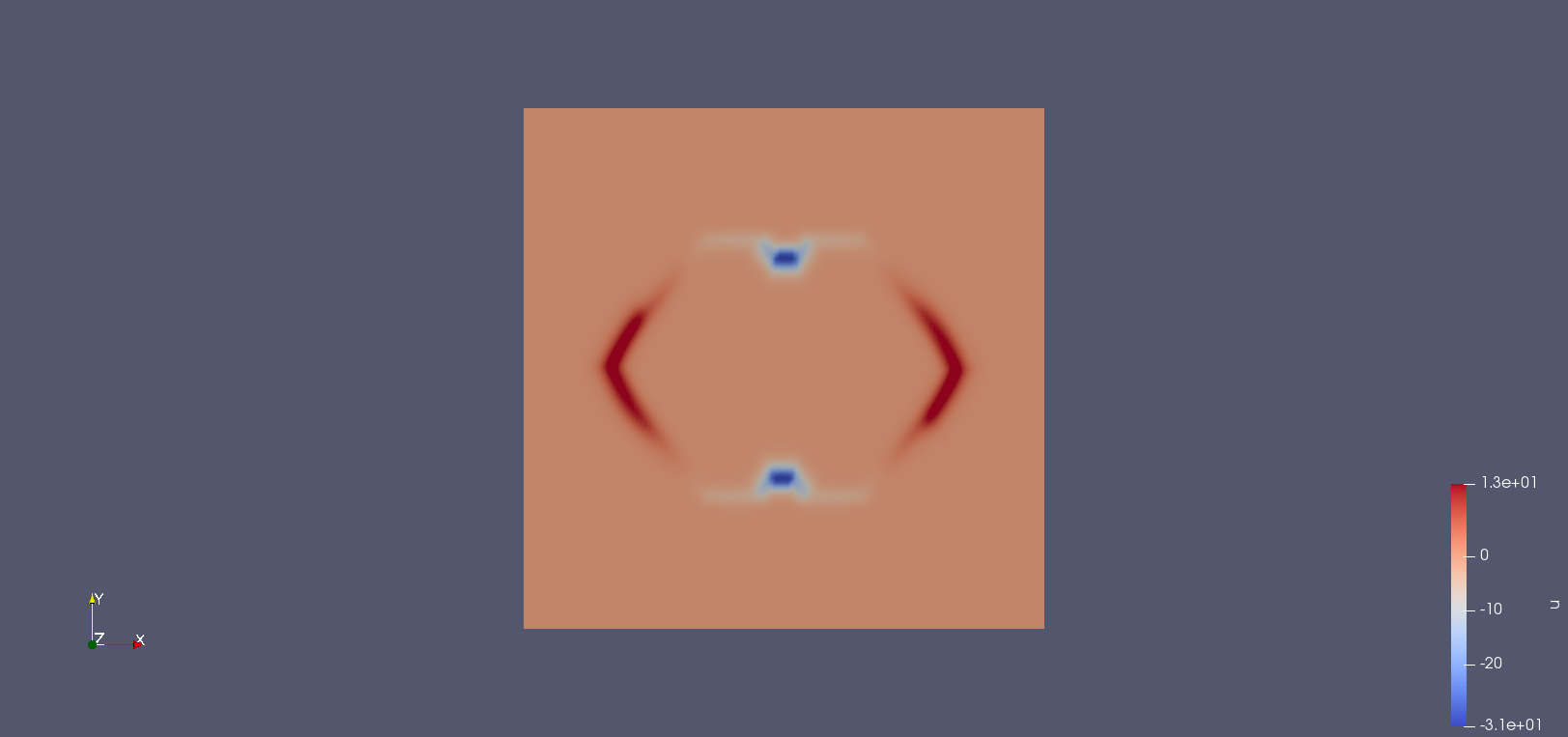}
	\includegraphics[trim={540px 112px 540px 112px}, clip, scale=.08]{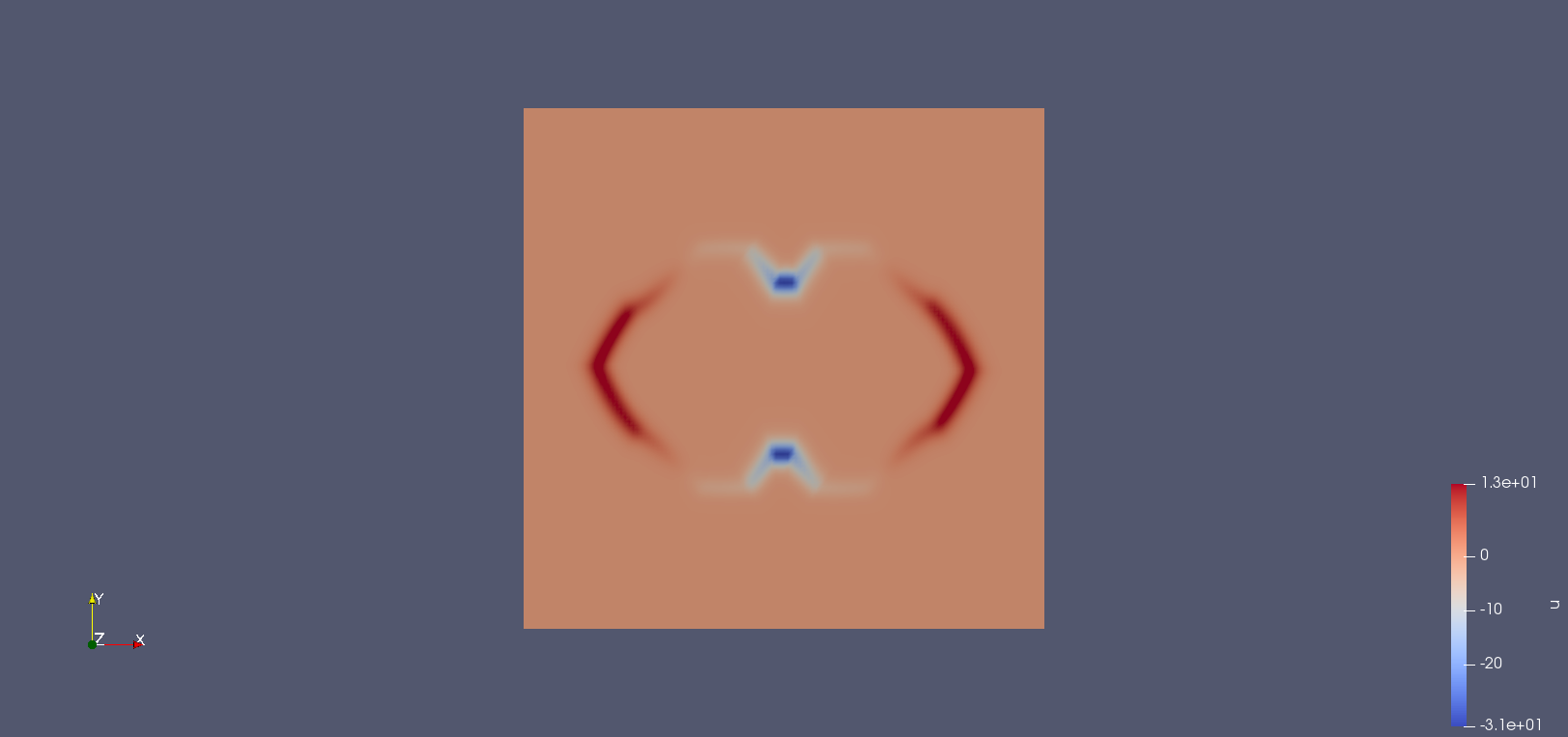}
	\includegraphics[trim={540px 112px 540px 112px}, clip, scale=.08]{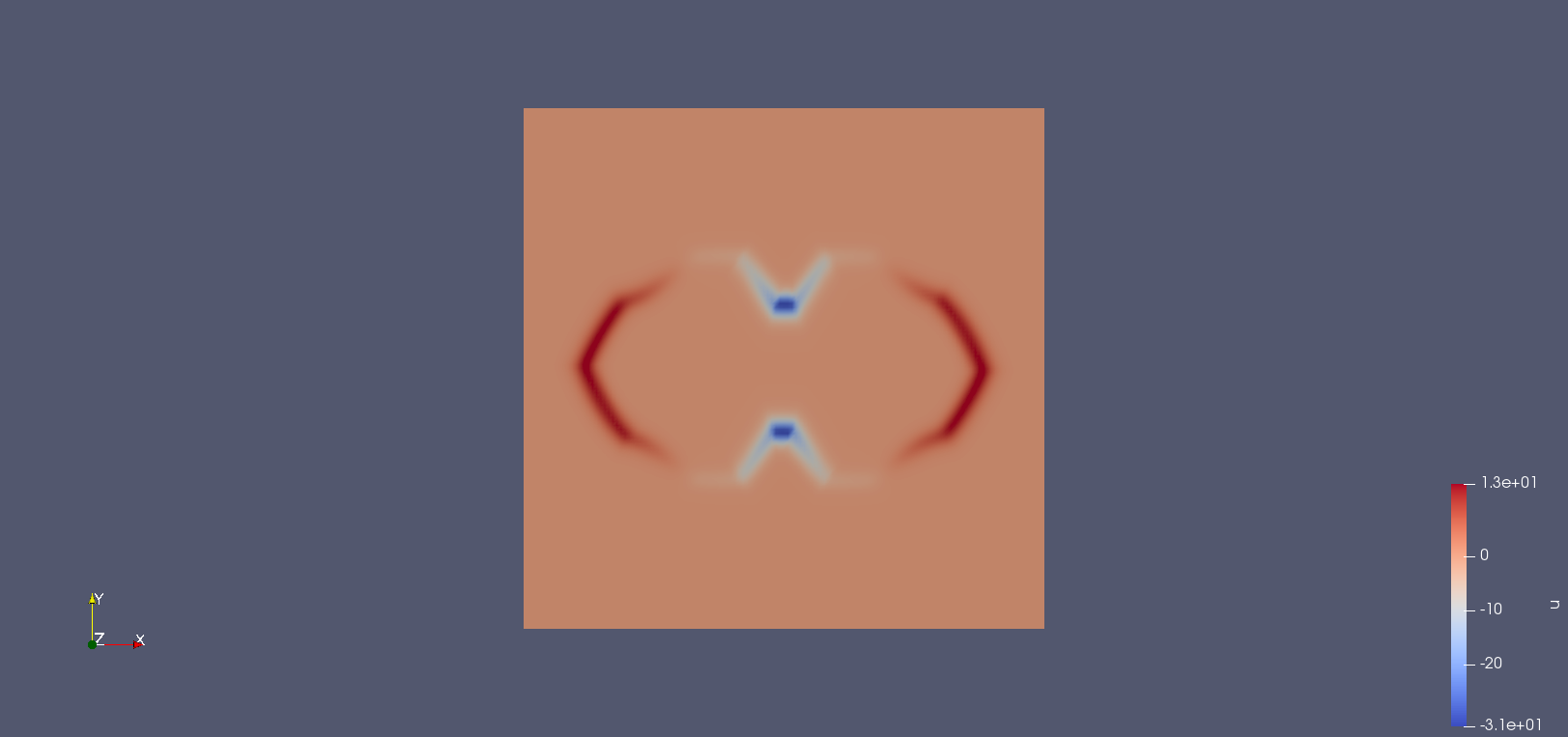}
	\includegraphics[trim={540px 112px 540px 112px}, clip, scale=.08]{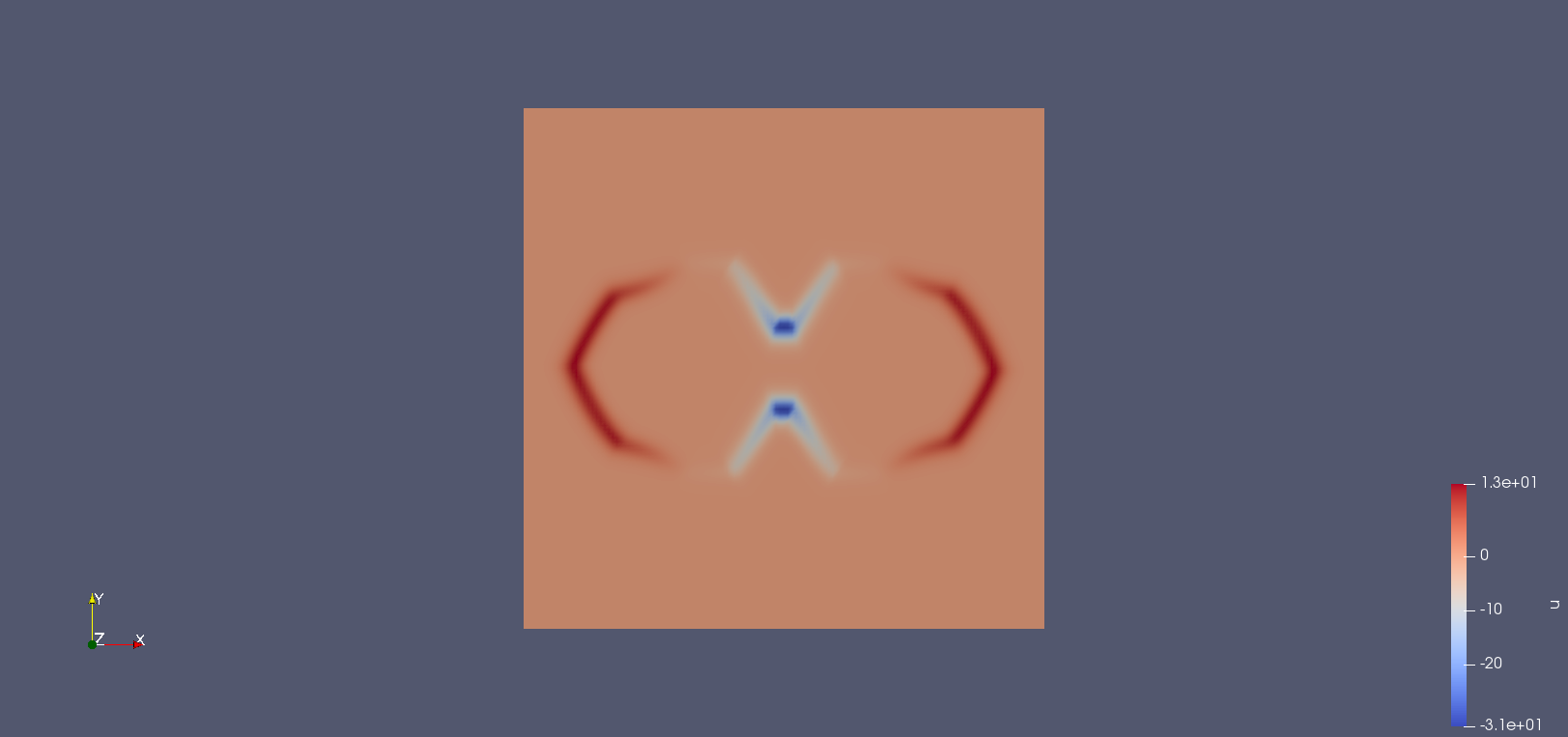}
	\includegraphics[trim={540px 112px 540px 112px}, clip, scale=.08]{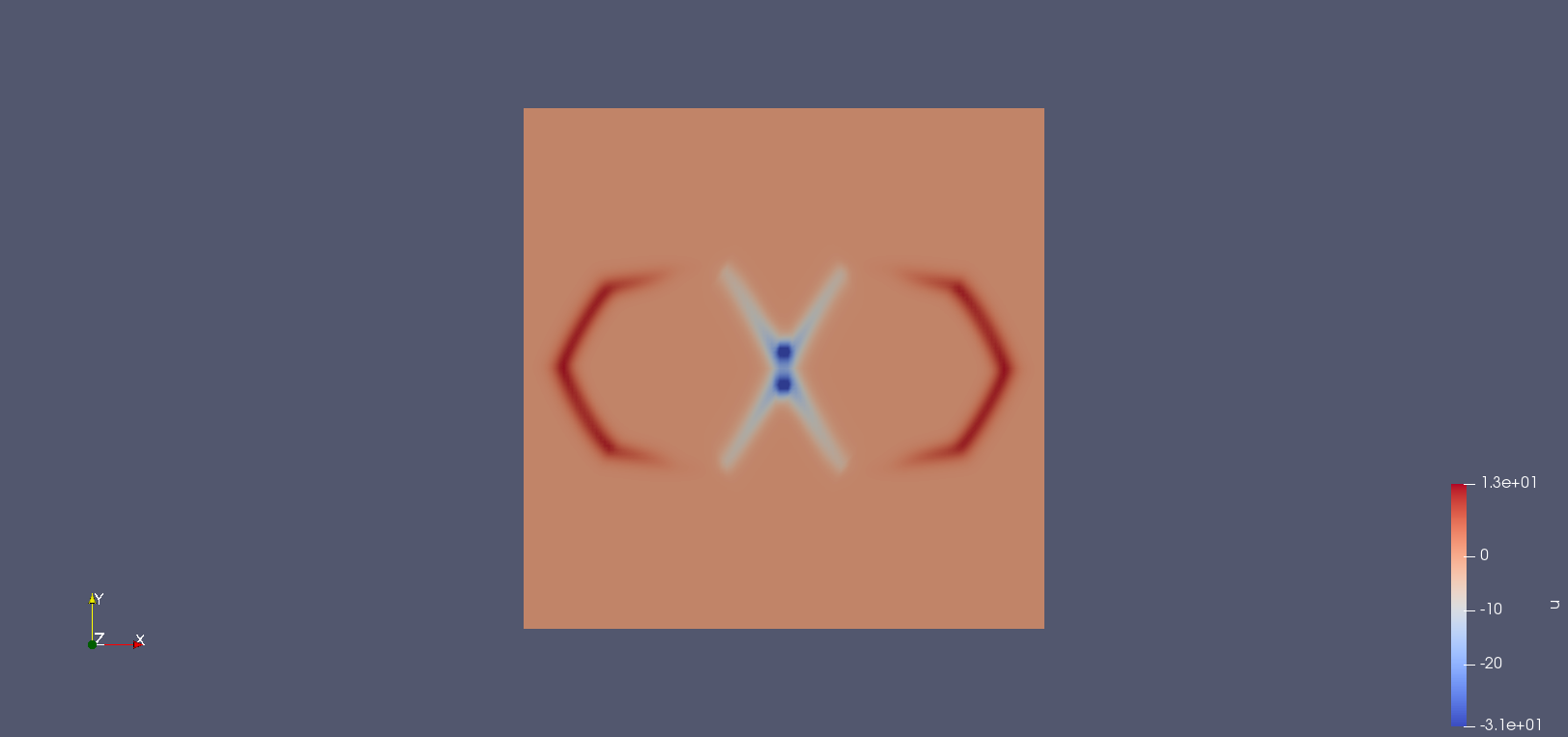}
	\includegraphics[trim={540px 112px 540px 112px}, clip, scale=.08]{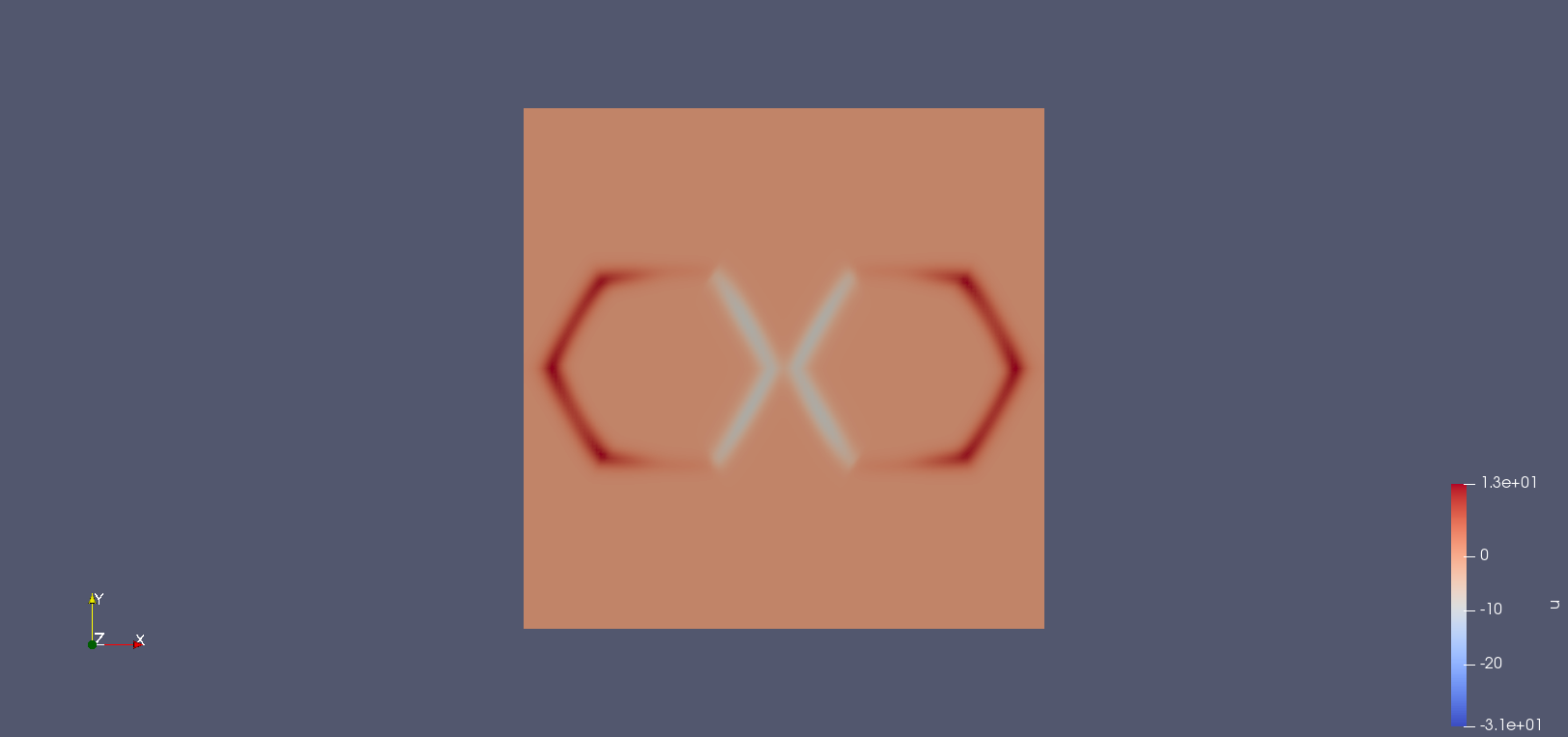}
	\includegraphics[trim={540px 112px 540px 112px}, clip, scale=.08]{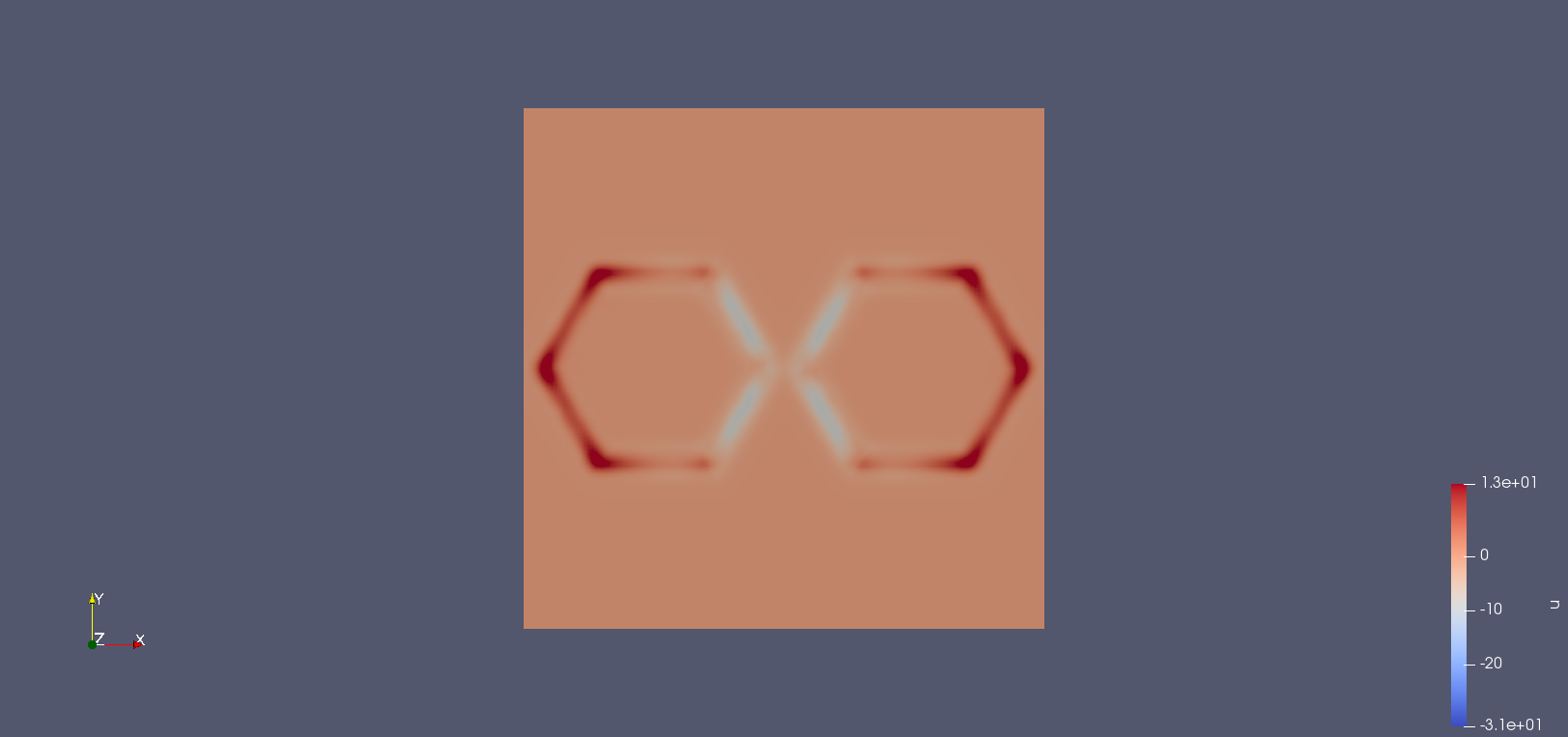}~~~
	\includegraphics[trim={1490px 0px 20px 490px}, clip, scale=.16]{{hexa_splitting_control.0000}.png}
\caption{Result for the `splitting hexagon' using the hexagon anisotropy.
	\label{fig:hexa_splitting}}
\end{center}   
\end{subfigure}
\caption{States (above) and corresponding controls (below) for the solution of splitting geometries.
  %at $t = timesnotlookedupsofar$.} % Hinweis: die Zahlen entsprechen ungefähr diesen Zeitschritten, ich hab die Bilder von 0-7 darauf umgerechnet und auf "schöne" Zahlen gerundet
  \label{fig:split}}
\end{center}  
\end{figure}
\fi

%\subsubsection{Merging geometries}	
\ifgraphics
\begin{figure}[htbp]
  \begin{center}
\begin{subfigure}{1.0\textwidth}
\begin{center}  
	\includegraphics[trim={540px 112px 540px 112px}, clip, scale=.08]{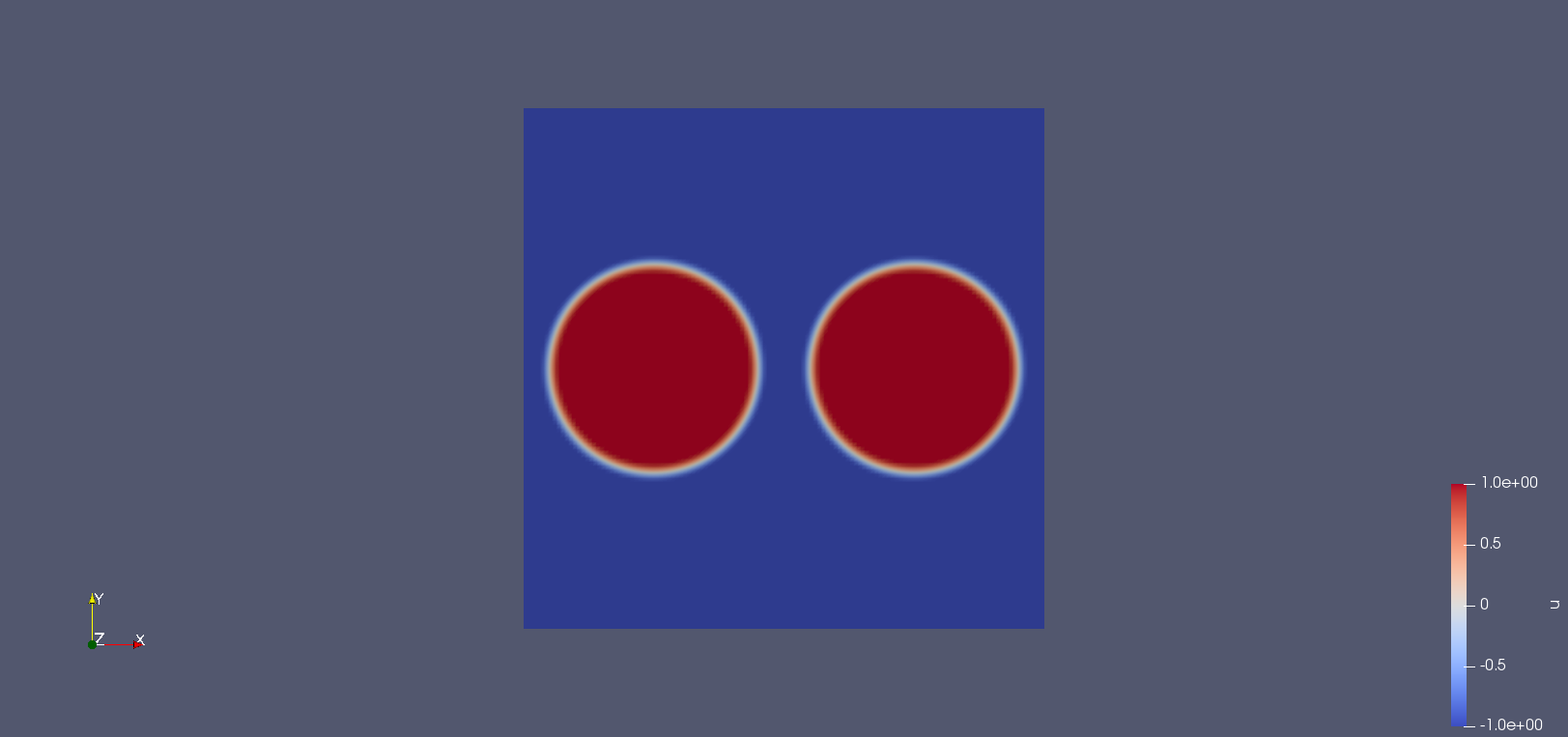}
	\includegraphics[trim={540px 112px 540px 112px}, clip, scale=.08]{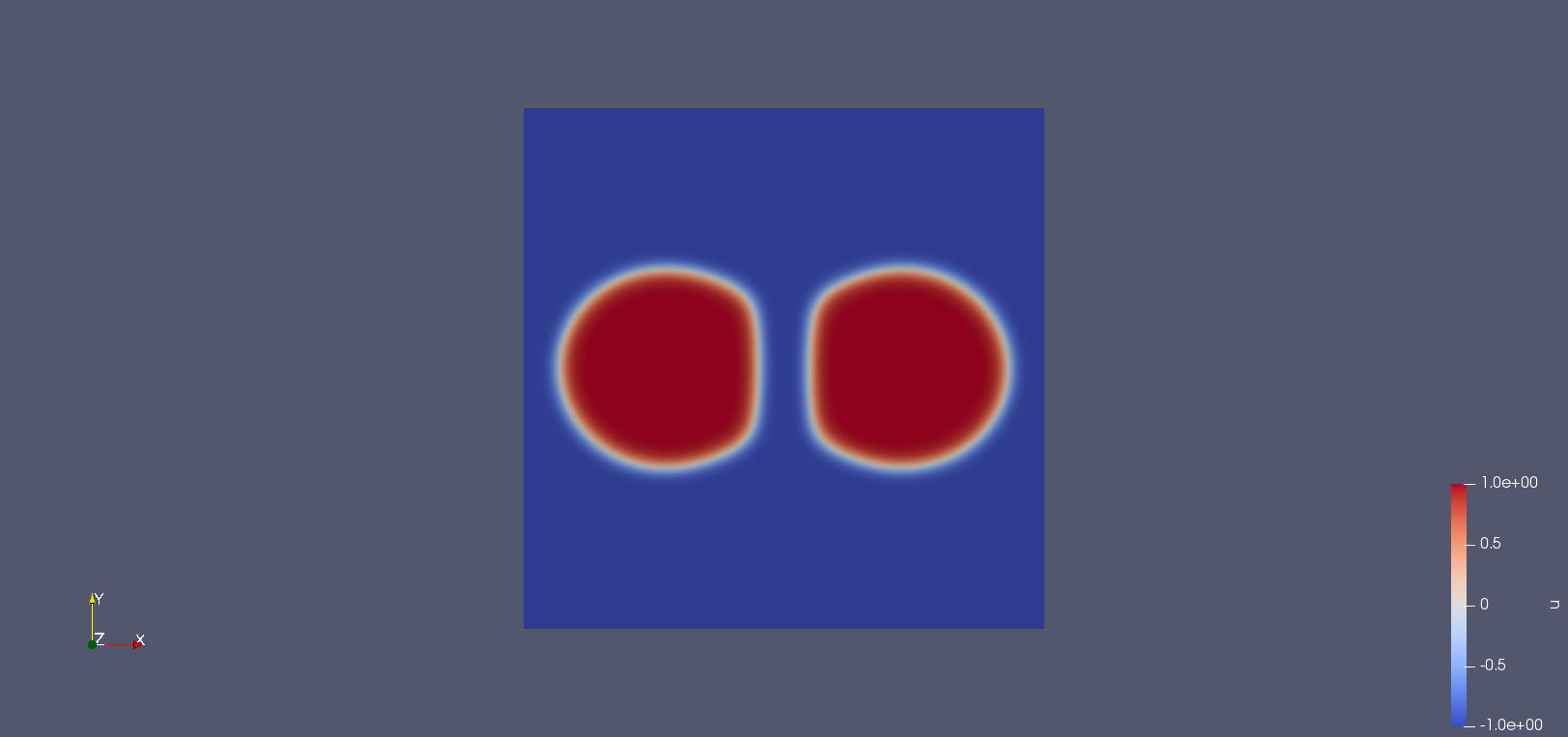}
	\includegraphics[trim={540px 112px 540px 112px}, clip, scale=.08]{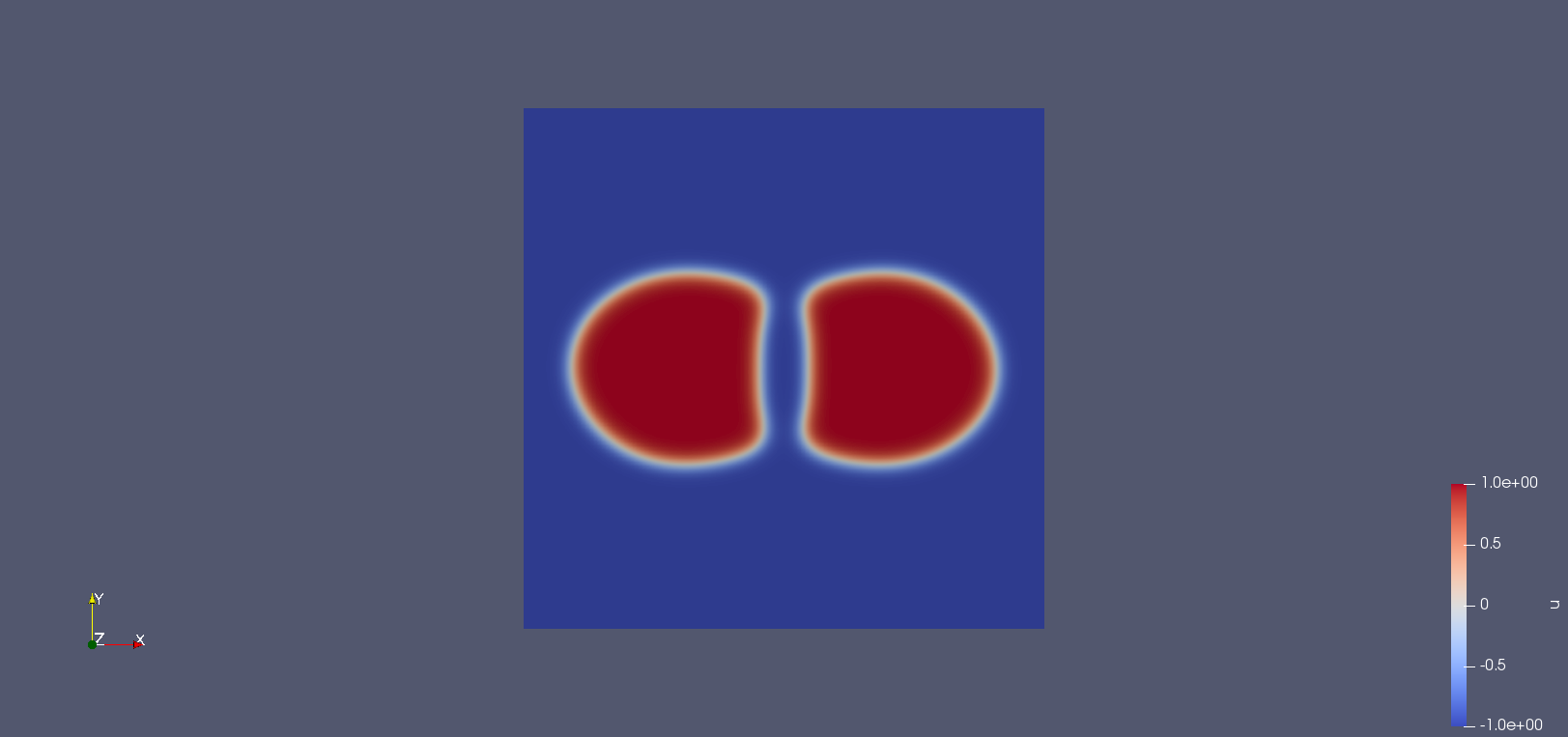}
	\includegraphics[trim={540px 112px 540px 112px}, clip, scale=.08]{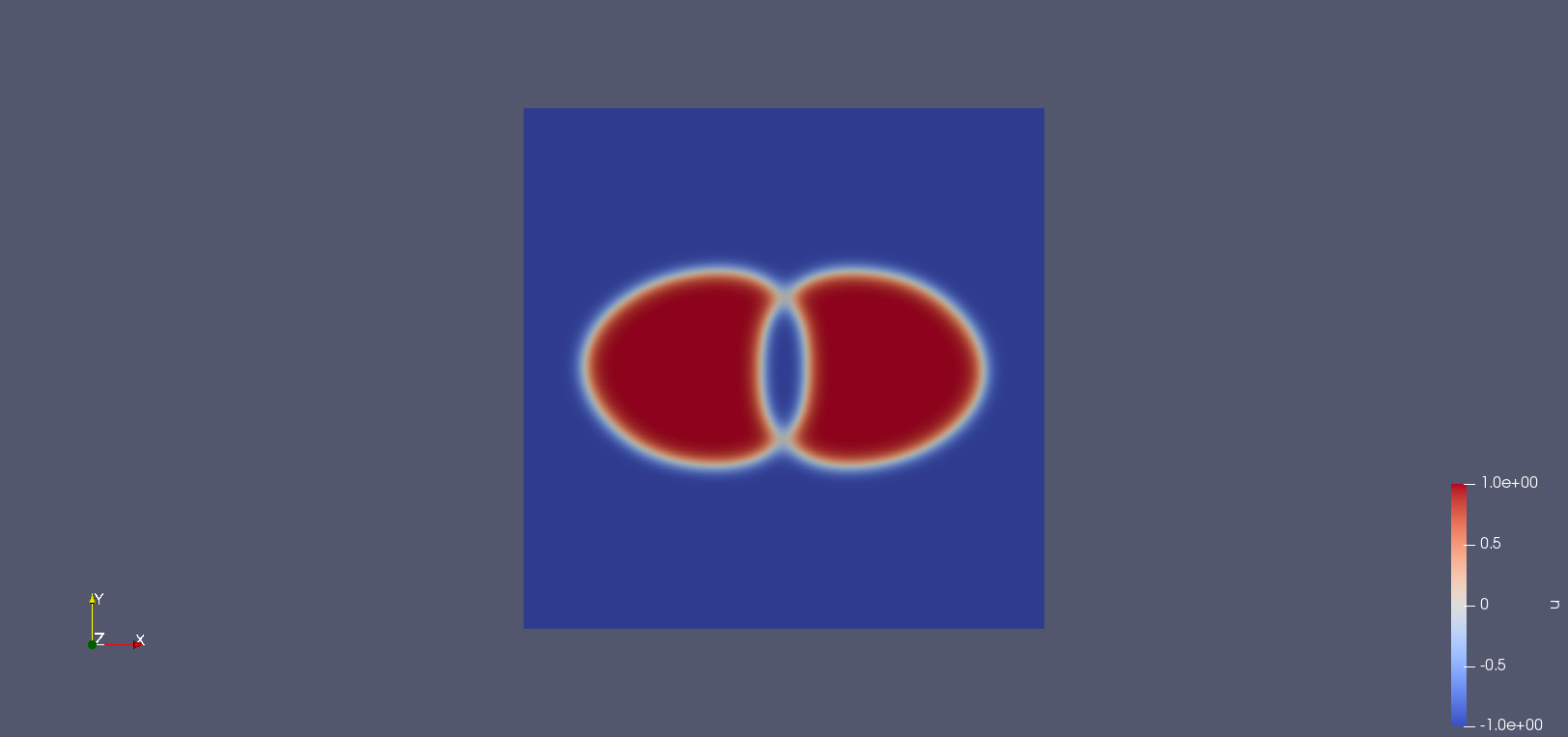}
	\includegraphics[trim={540px 112px 540px 112px}, clip, scale=.08]{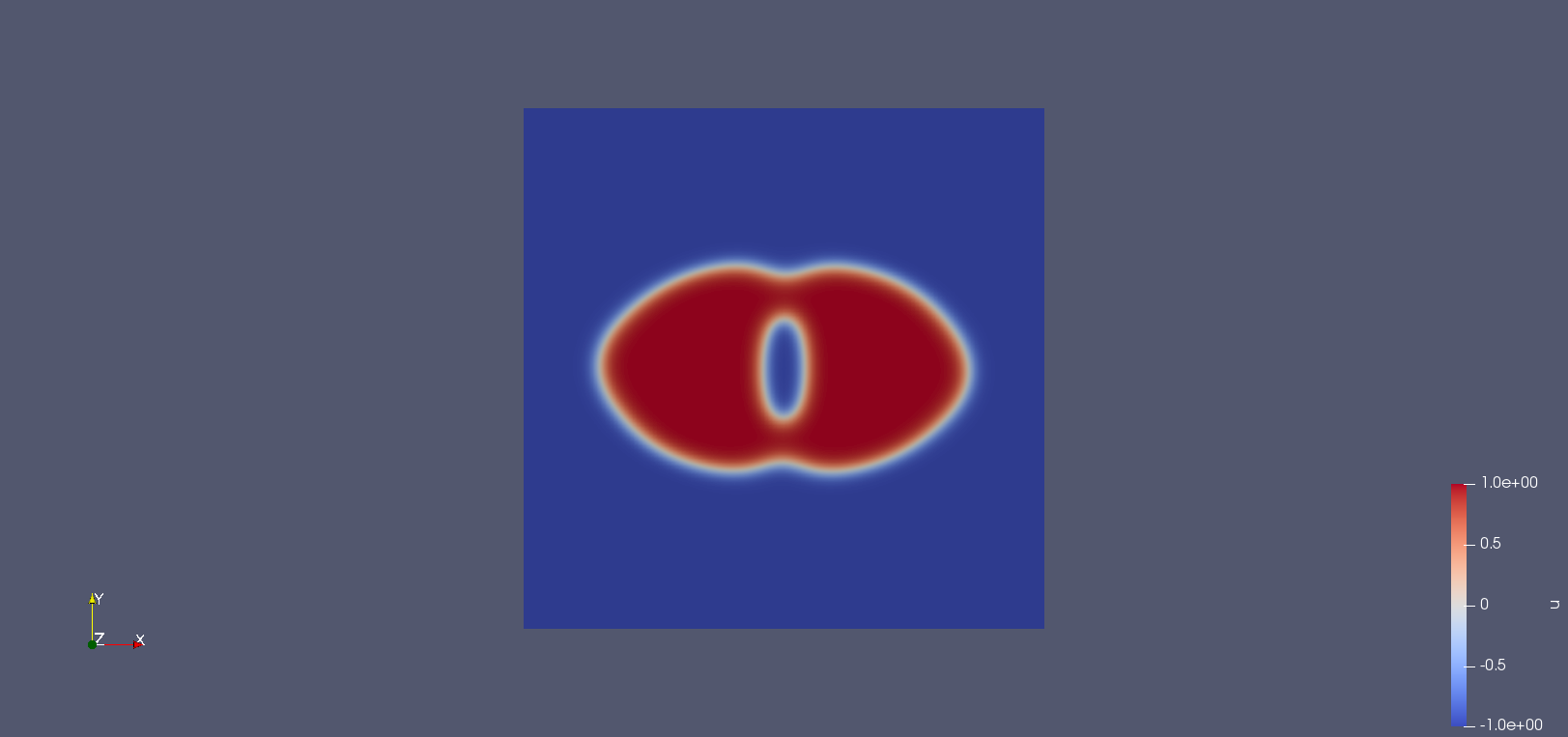}
	\includegraphics[trim={540px 112px 540px 112px}, clip, scale=.08]{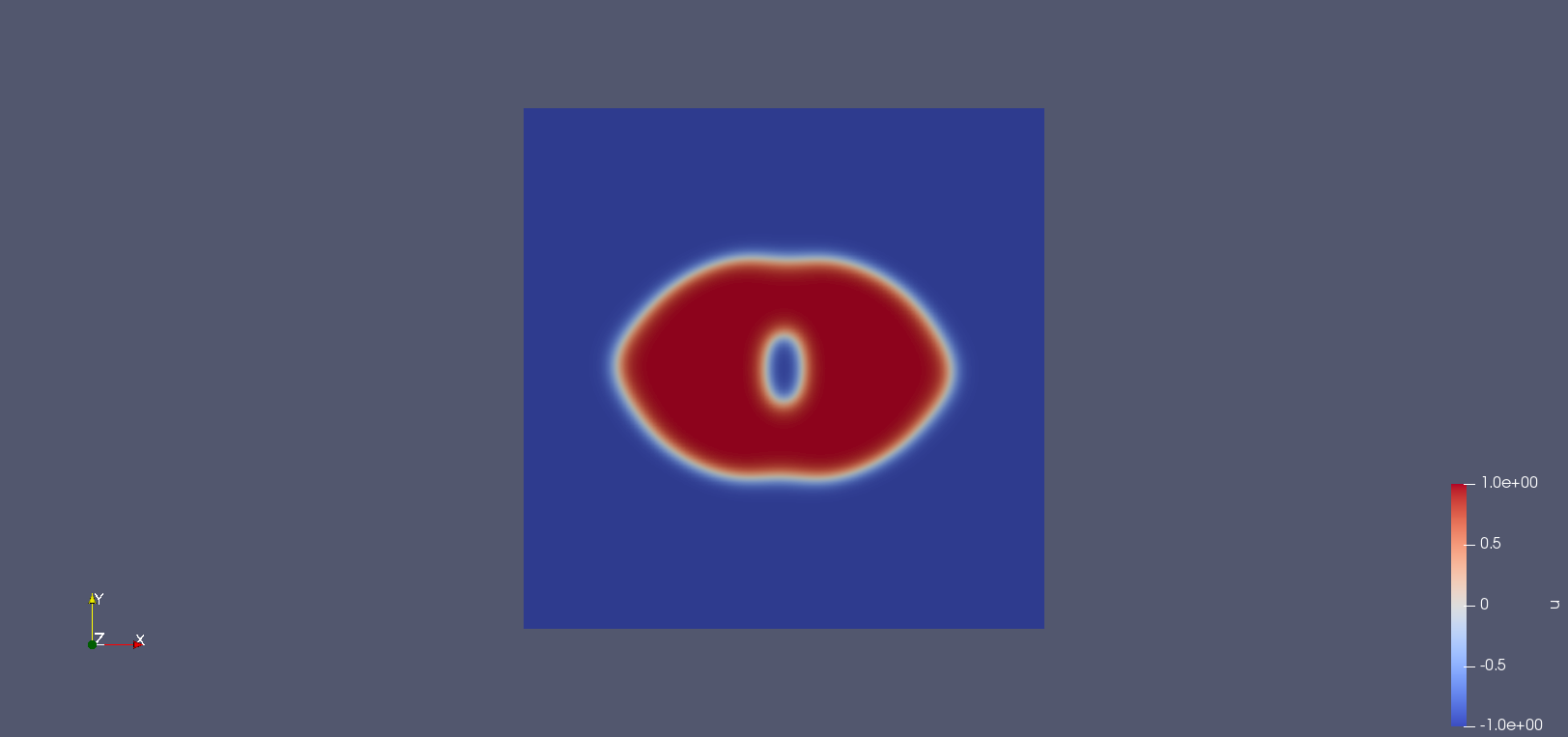}
	\includegraphics[trim={540px 112px 540px 112px}, clip, scale=.08]{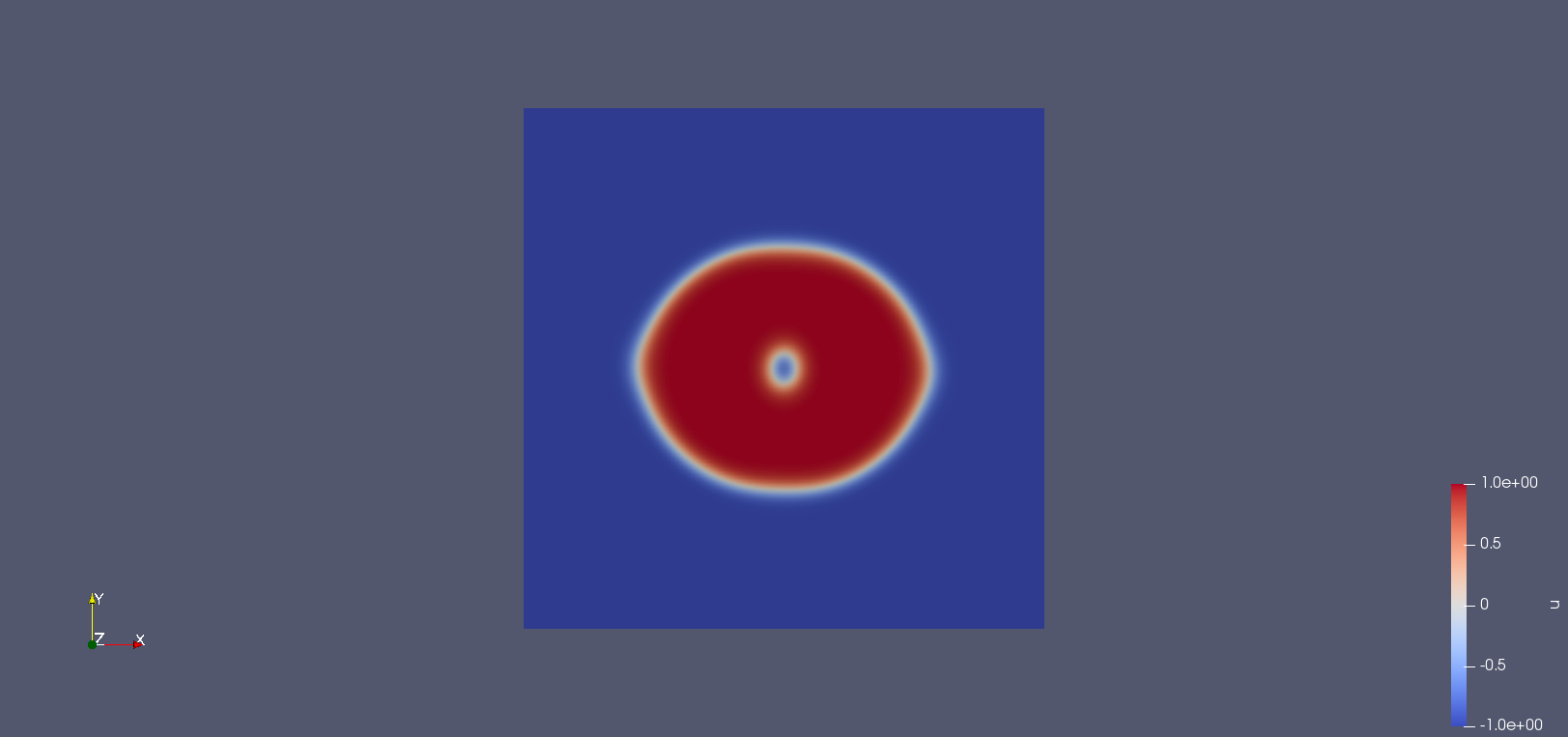}
	\includegraphics[trim={540px 112px 540px 112px}, clip, scale=.08]{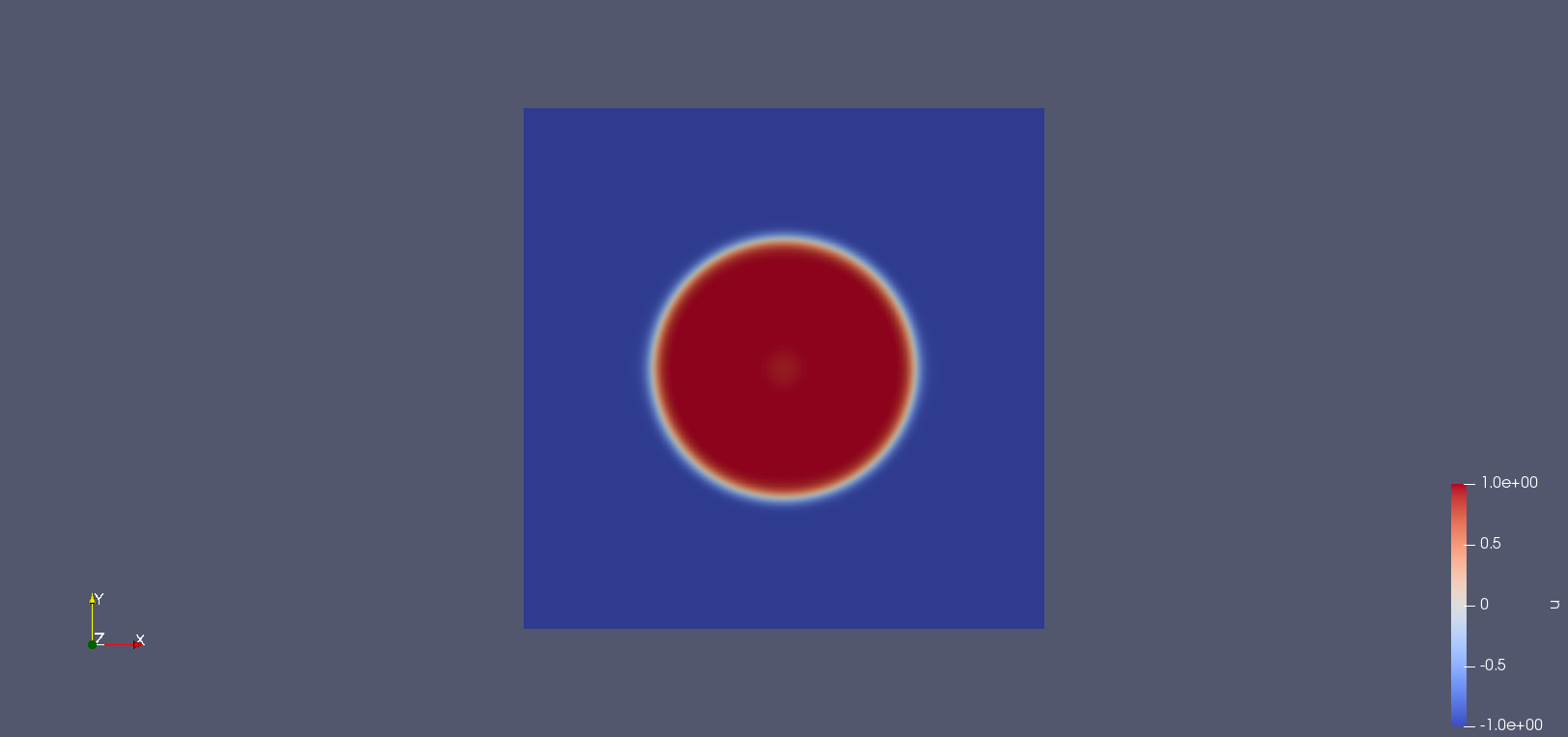}~~~
	\includegraphics[trim={1490px 0px 20px 490px}, clip, scale=.16]{{iso_merge_state.0000}.png}~\\~\\
	\includegraphics[trim={540px 112px 540px 112px}, clip, scale=.08]{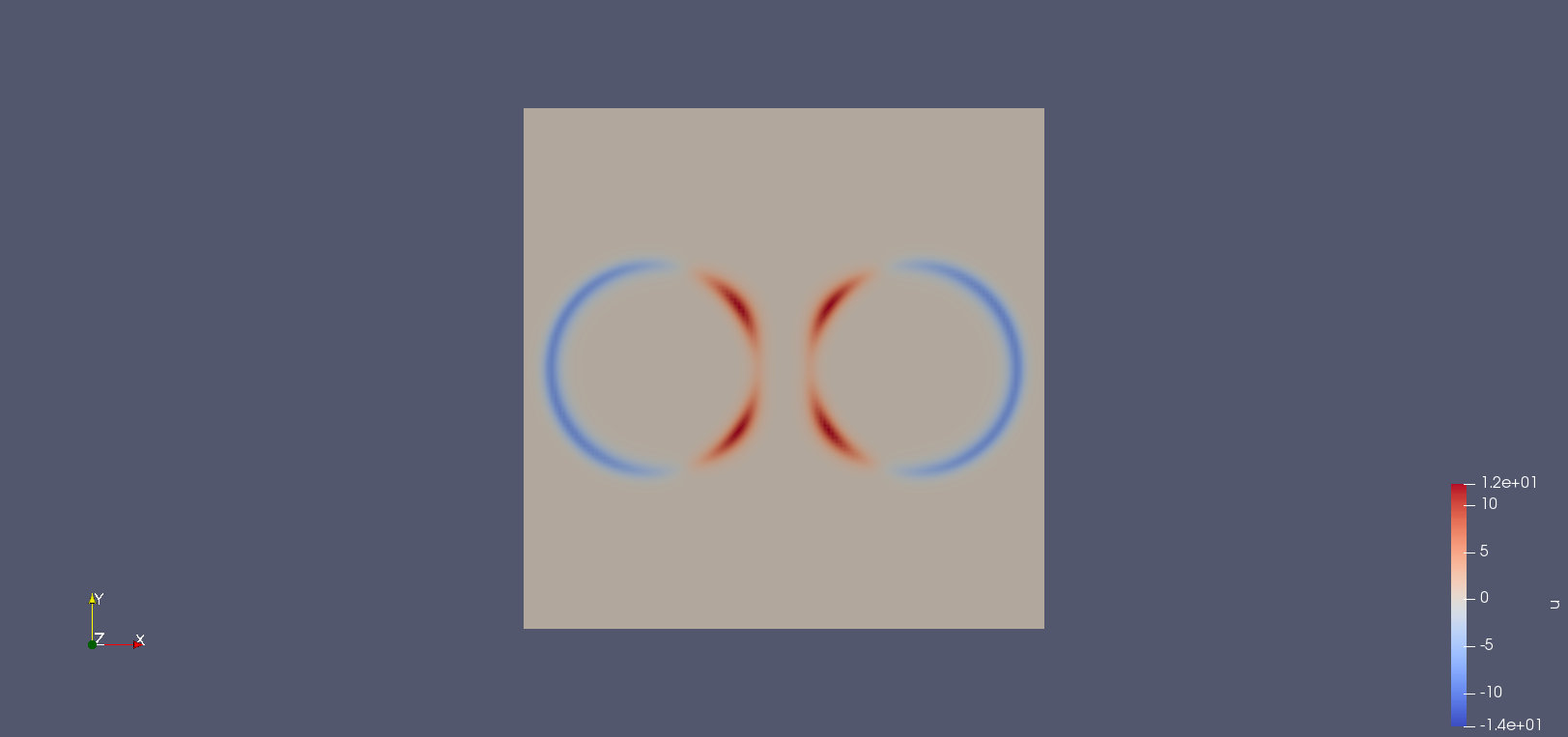}
	\includegraphics[trim={540px 112px 540px 112px}, clip, scale=.08]{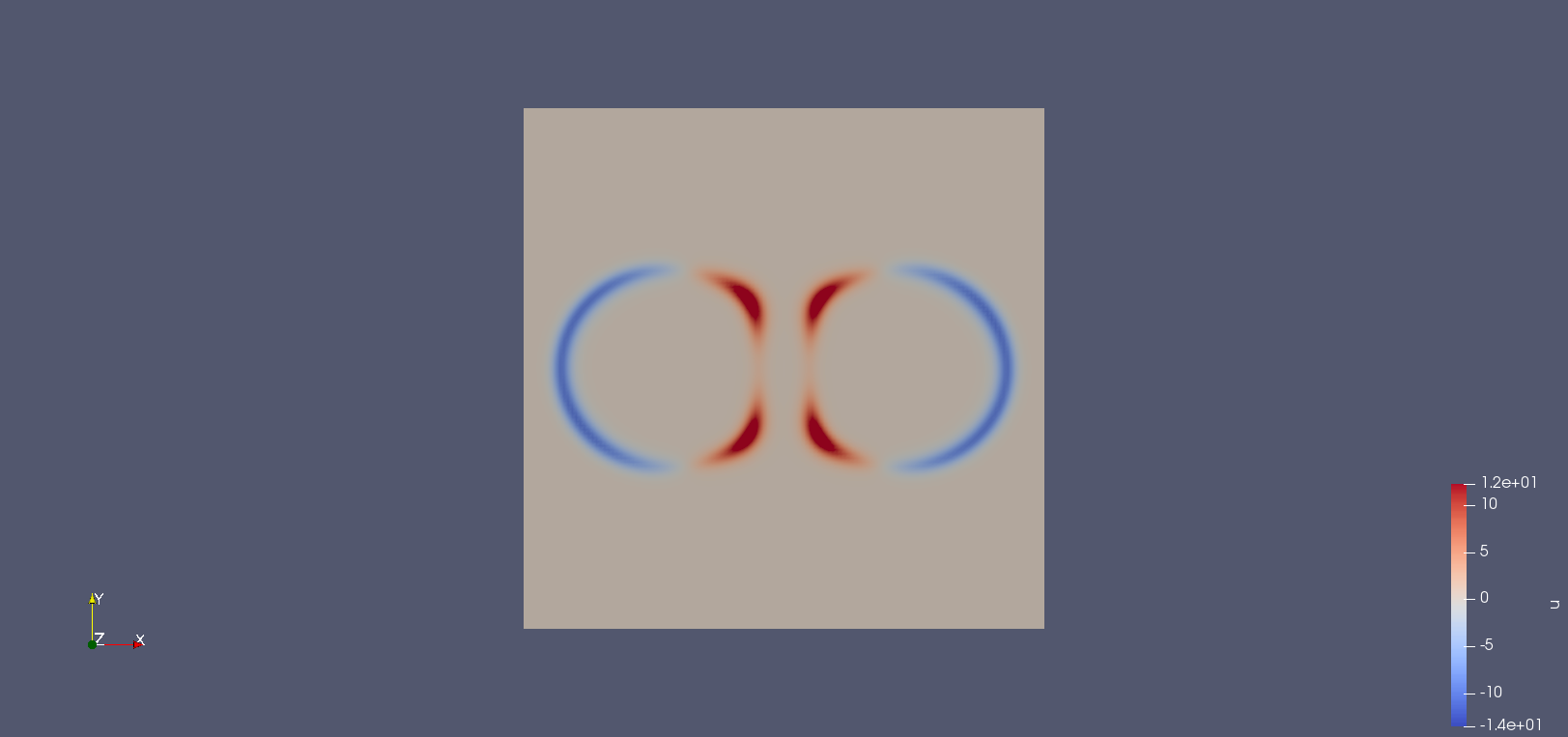}
	\includegraphics[trim={540px 112px 540px 112px}, clip, scale=.08]{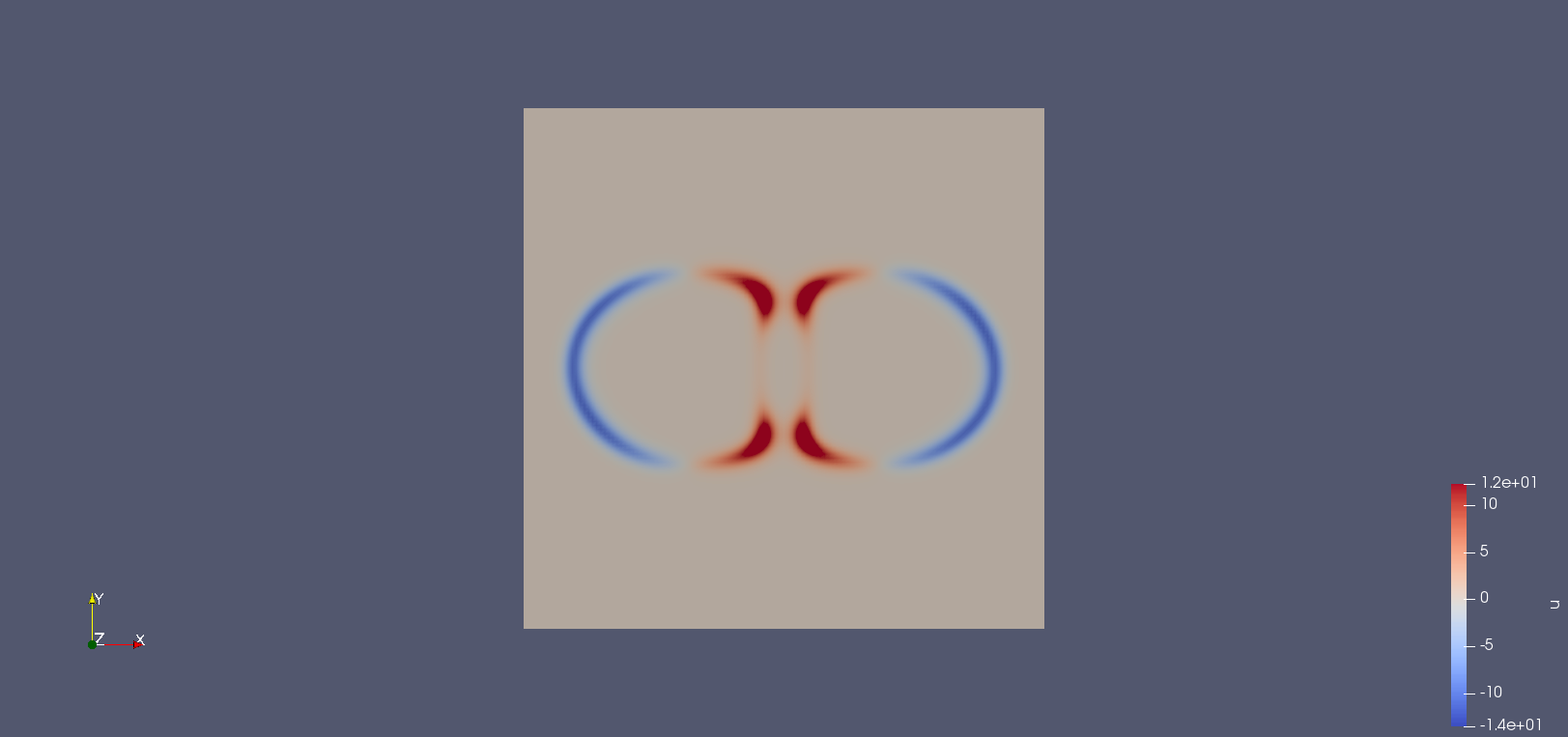}
	\includegraphics[trim={540px 112px 540px 112px}, clip, scale=.08]{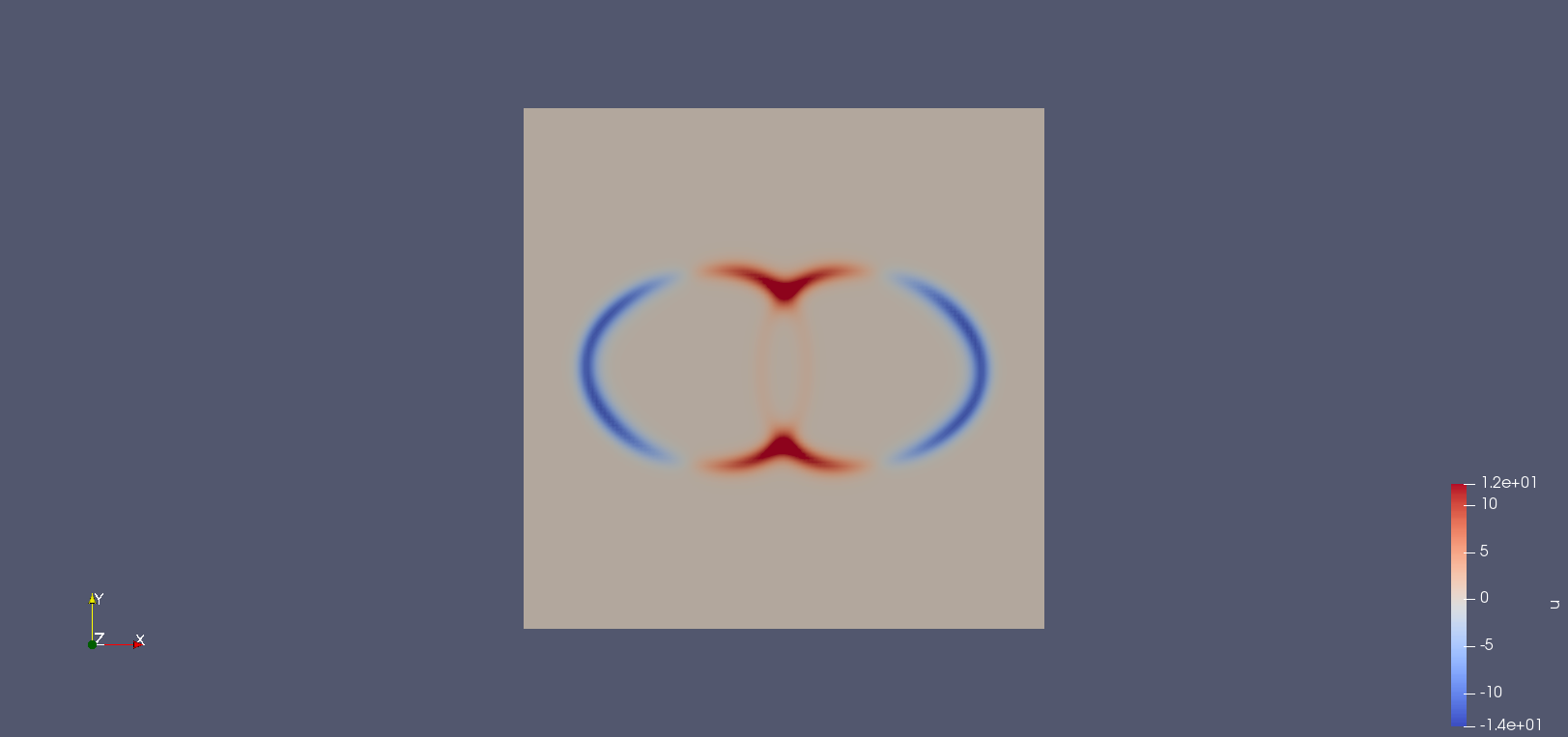}
	\includegraphics[trim={540px 112px 540px 112px}, clip, scale=.08]{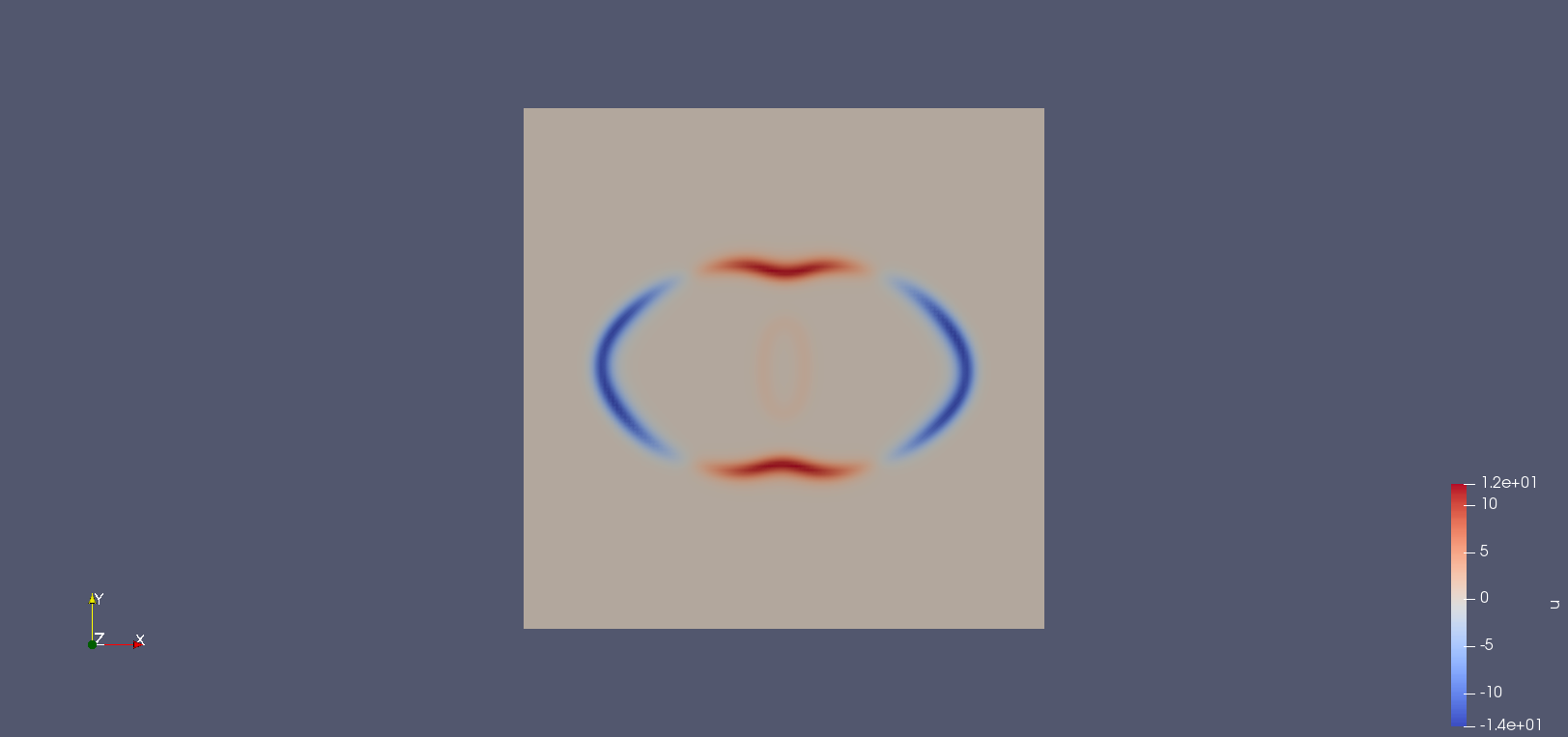}
	\includegraphics[trim={540px 112px 540px 112px}, clip, scale=.08]{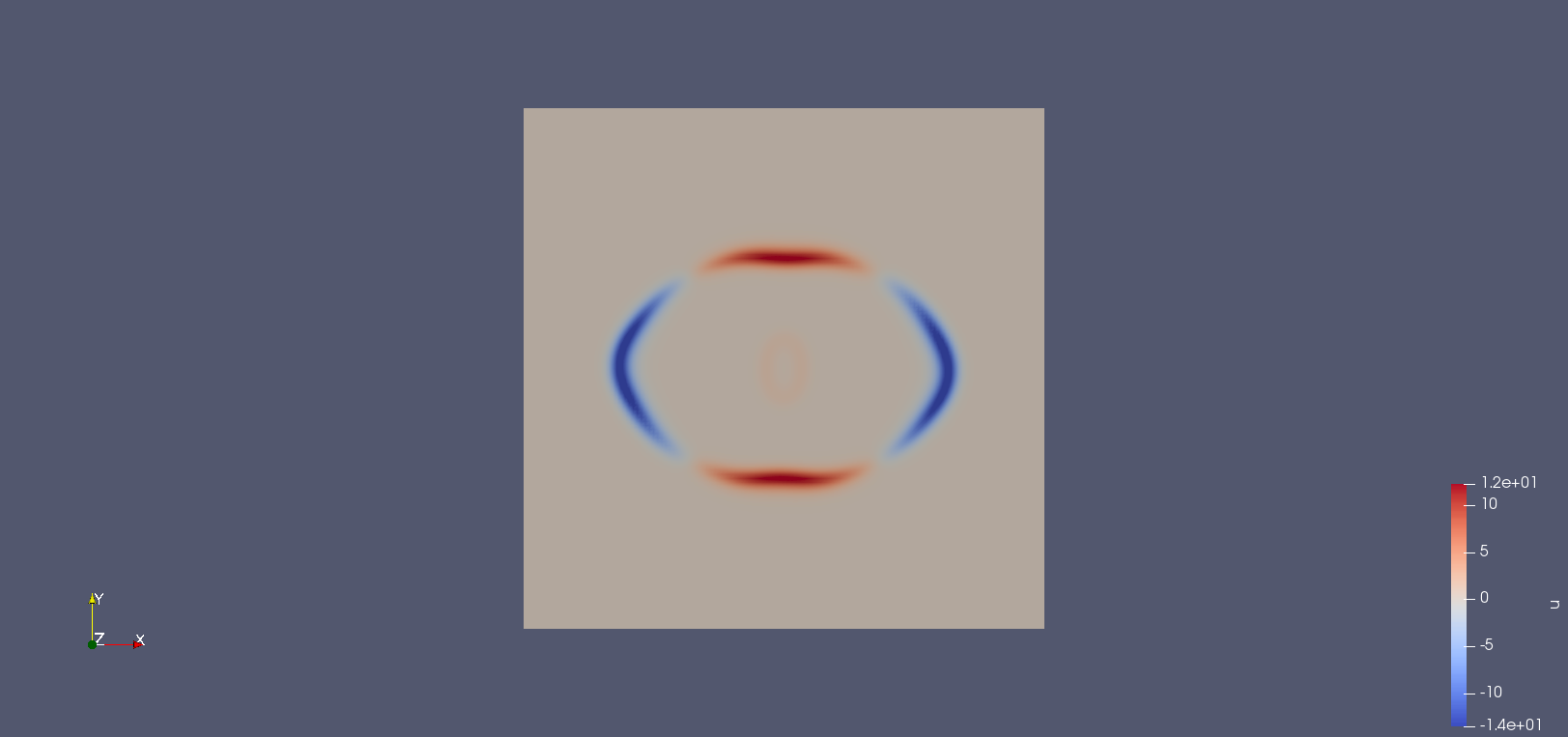}
	\includegraphics[trim={540px 112px 540px 112px}, clip, scale=.08]{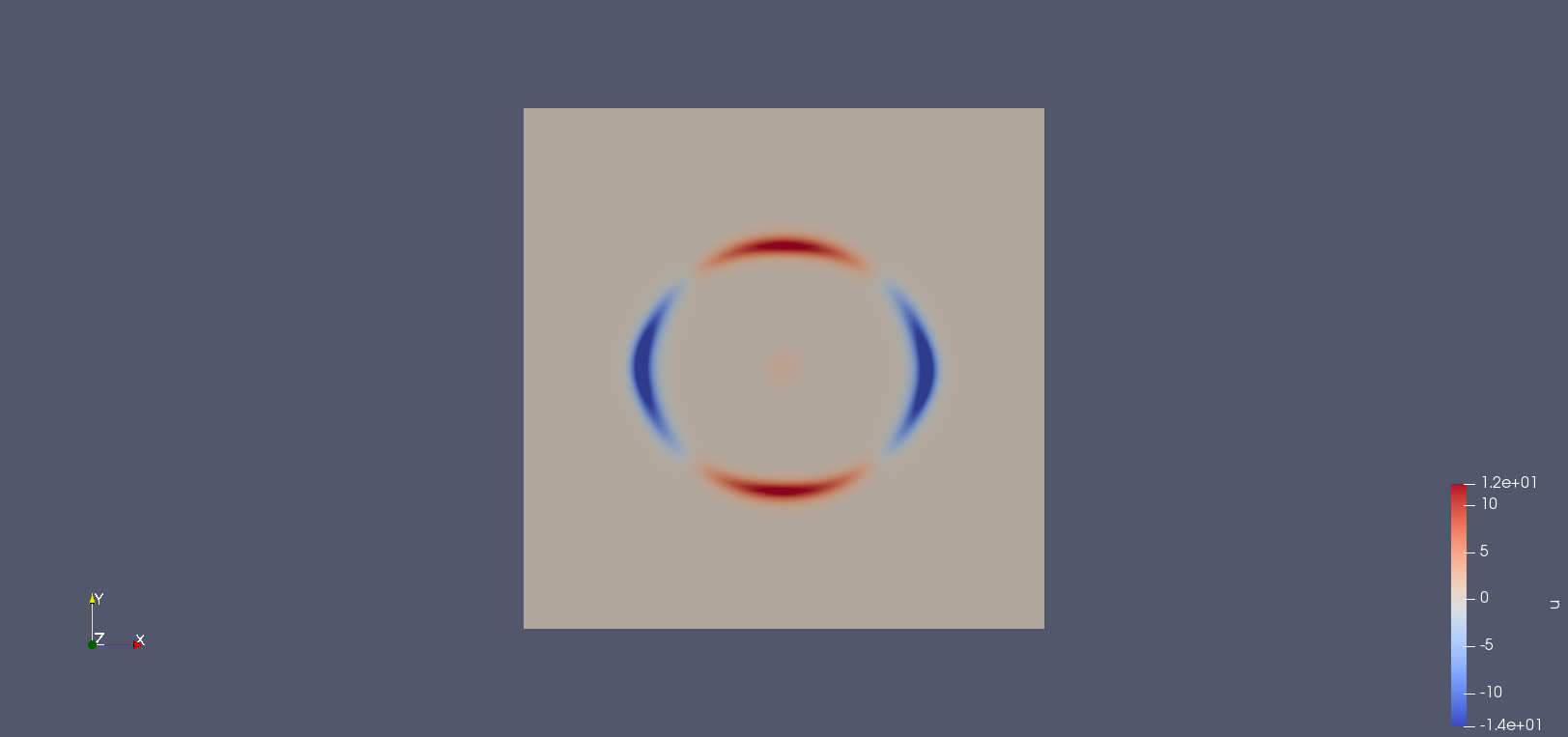}
	\includegraphics[trim={540px 112px 540px 112px}, clip, scale=.08]{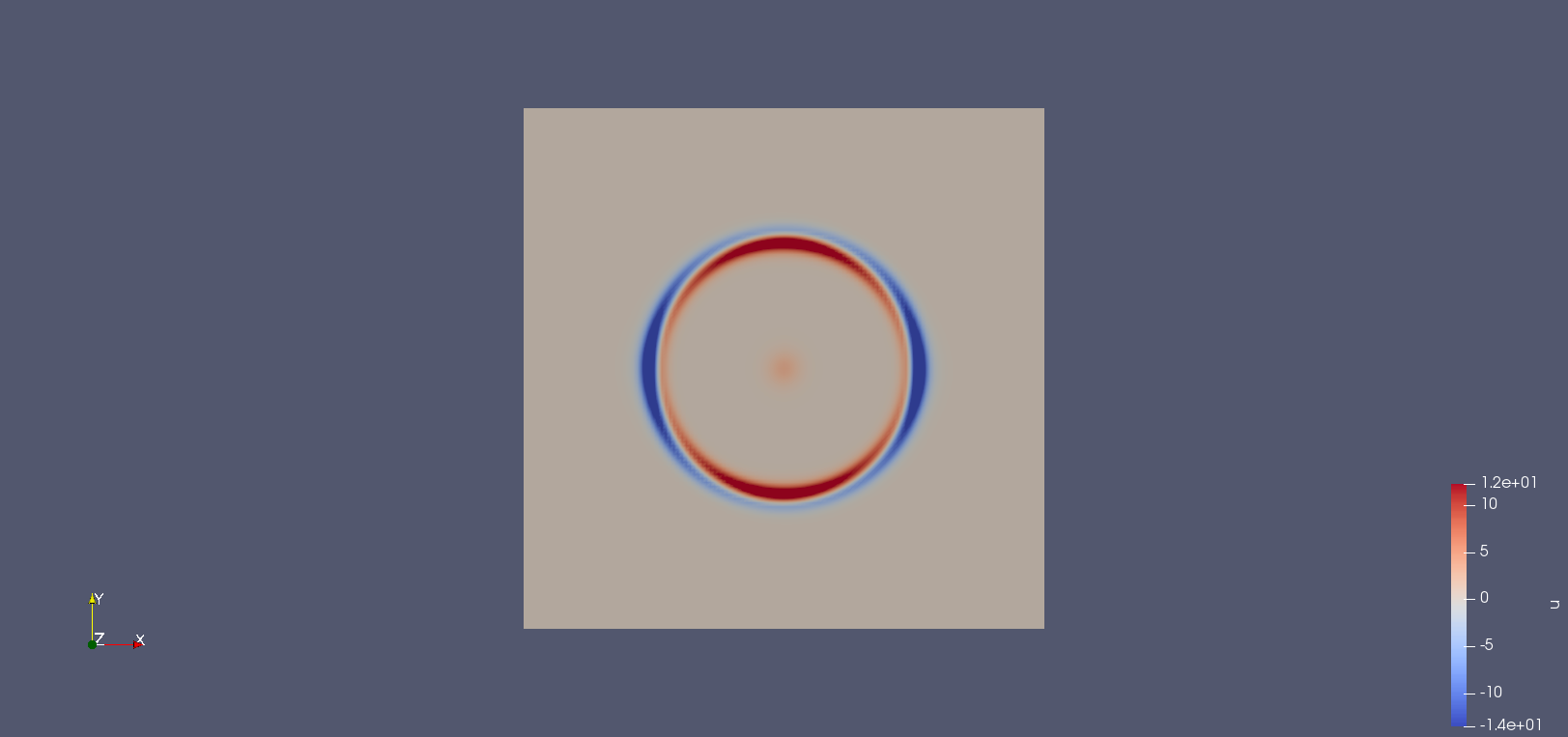}~~~
	\includegraphics[trim={1490px 0px 20px 490px}, clip, scale=.16]{{iso_merge_control.0000}.png}
	\caption{Result for `merge circle' in the isotropic case. 
	\label{fig:iso_merge}}
\end{center}    
\end{subfigure}
\medskip

\begin{subfigure}{1.0\textwidth}
\begin{center}  
		\includegraphics[trim={540px 112px 540px 112px}, clip, scale=.08]{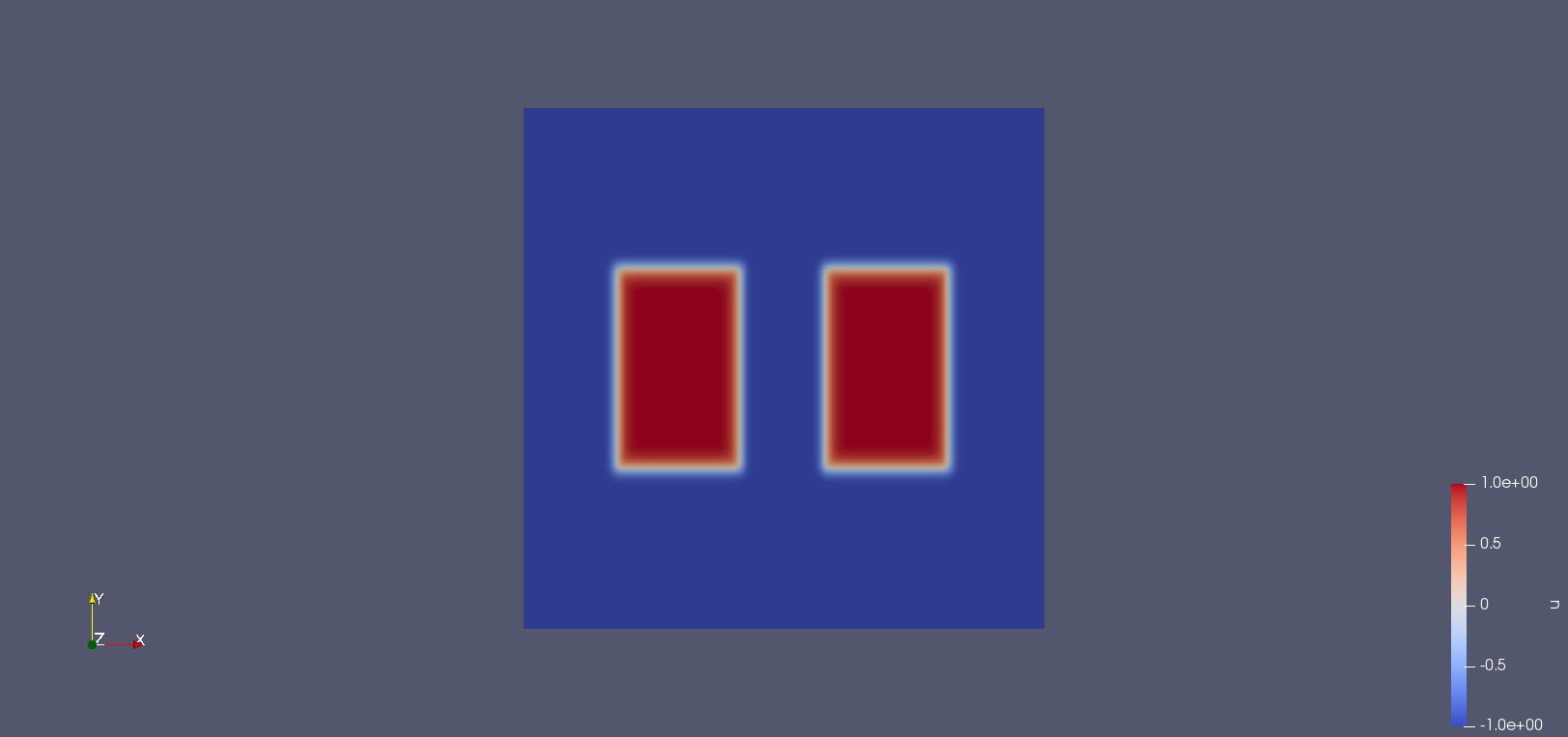}
		\includegraphics[trim={540px 112px 540px 112px}, clip, scale=.08]{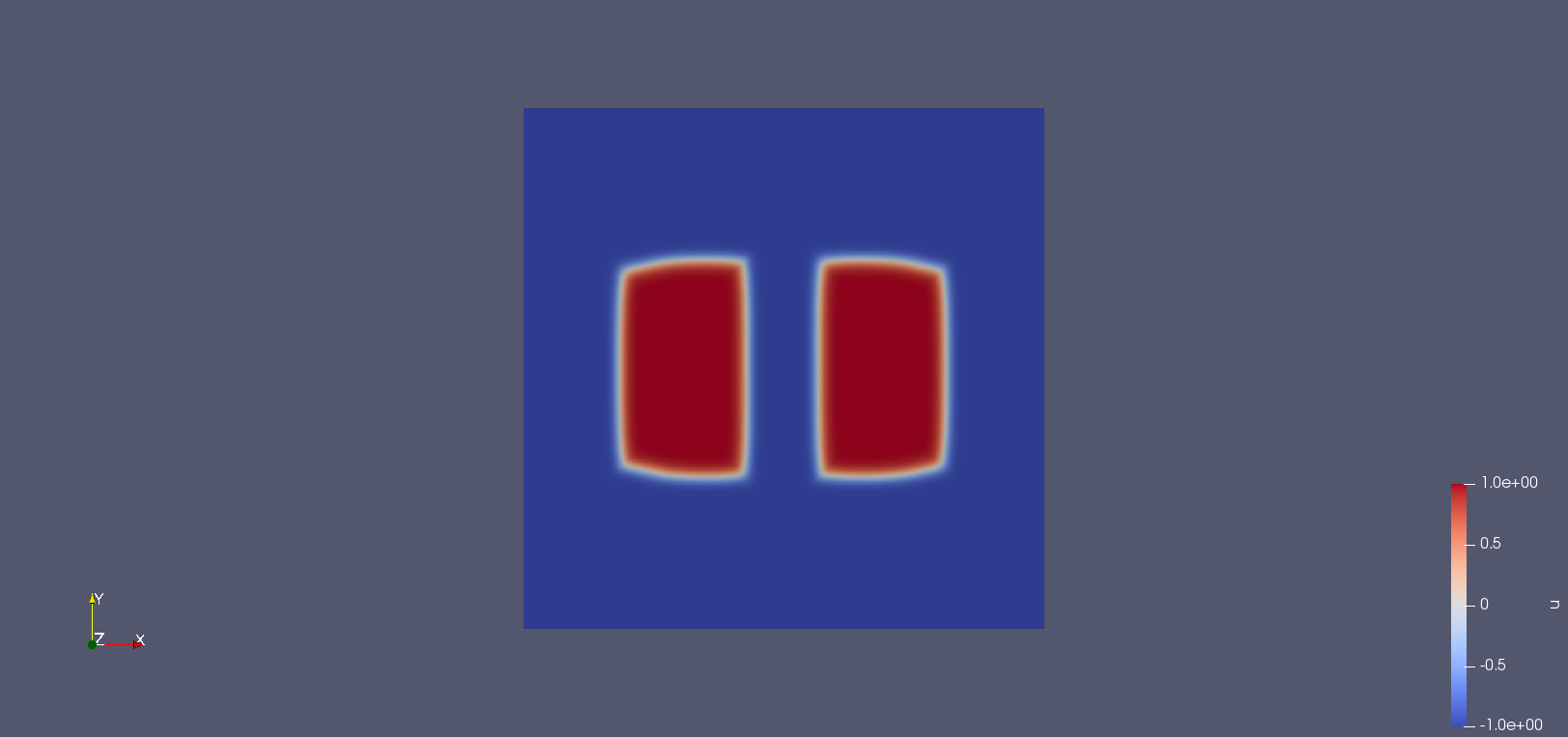}
		\includegraphics[trim={540px 112px 540px 112px}, clip, scale=.08]{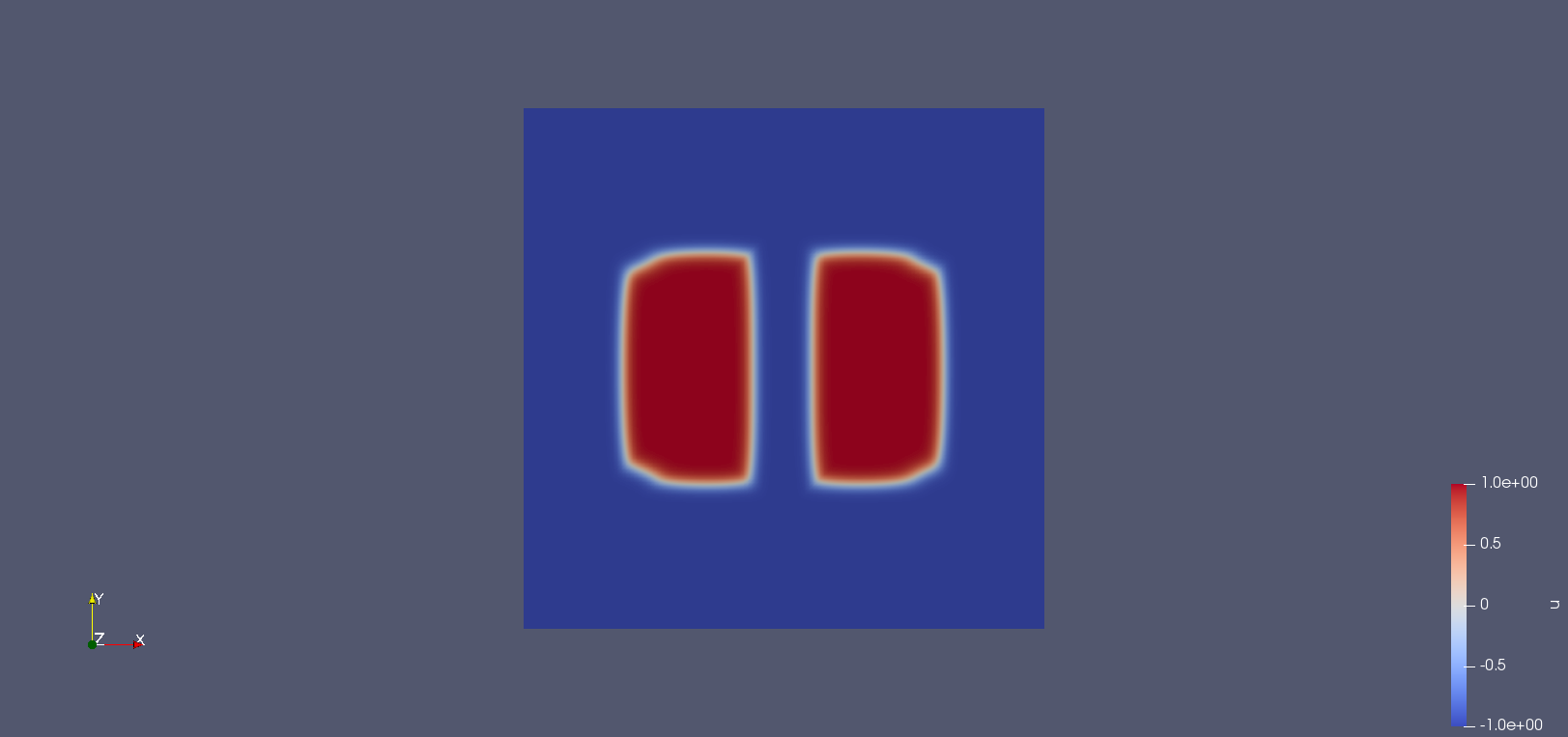}
		\includegraphics[trim={540px 112px 540px 112px}, clip, scale=.08]{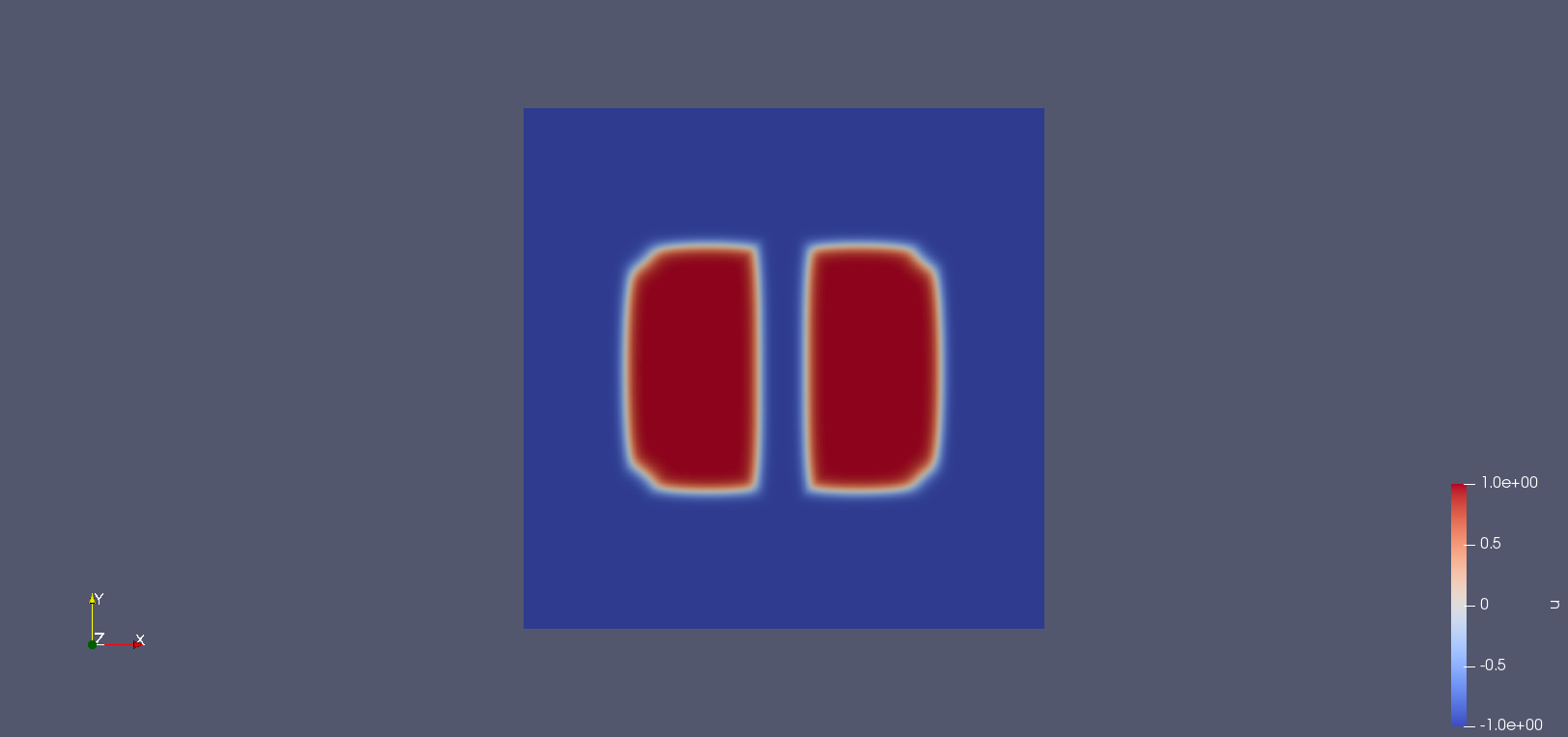}
		\includegraphics[trim={540px 112px 540px 112px}, clip, scale=.08]{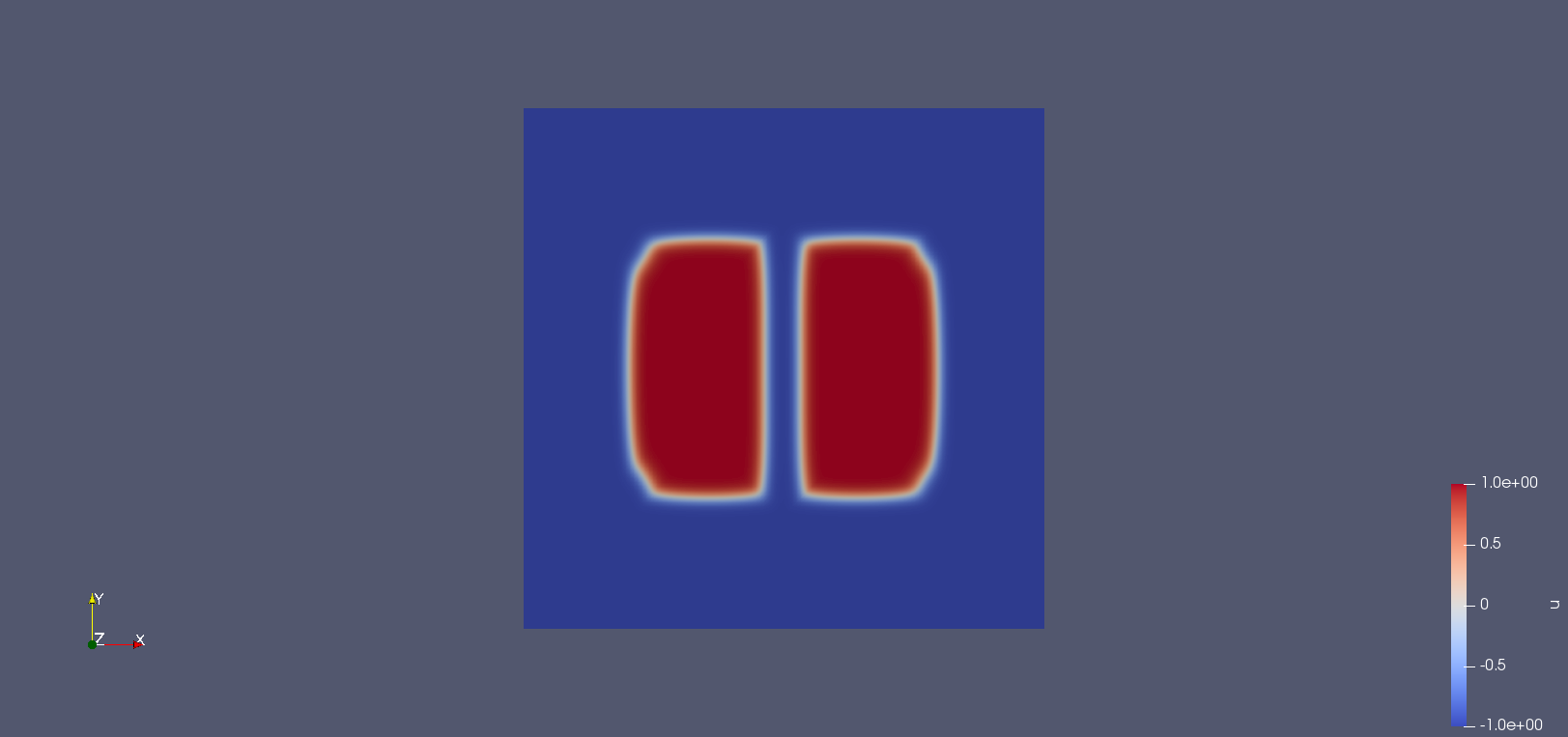}
		\includegraphics[trim={540px 112px 540px 112px}, clip, scale=.08]{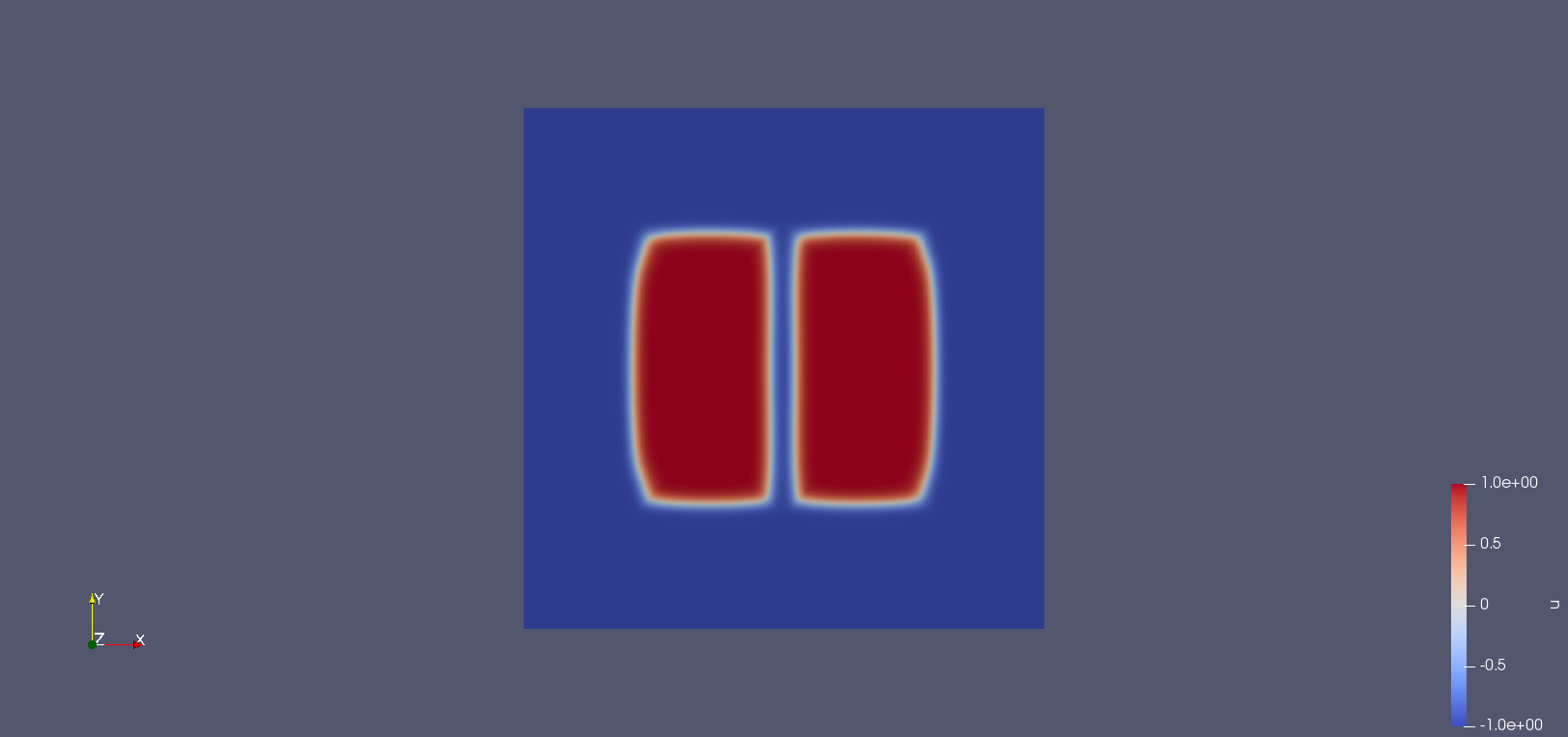}
		\includegraphics[trim={540px 112px 540px 112px}, clip, scale=.08]{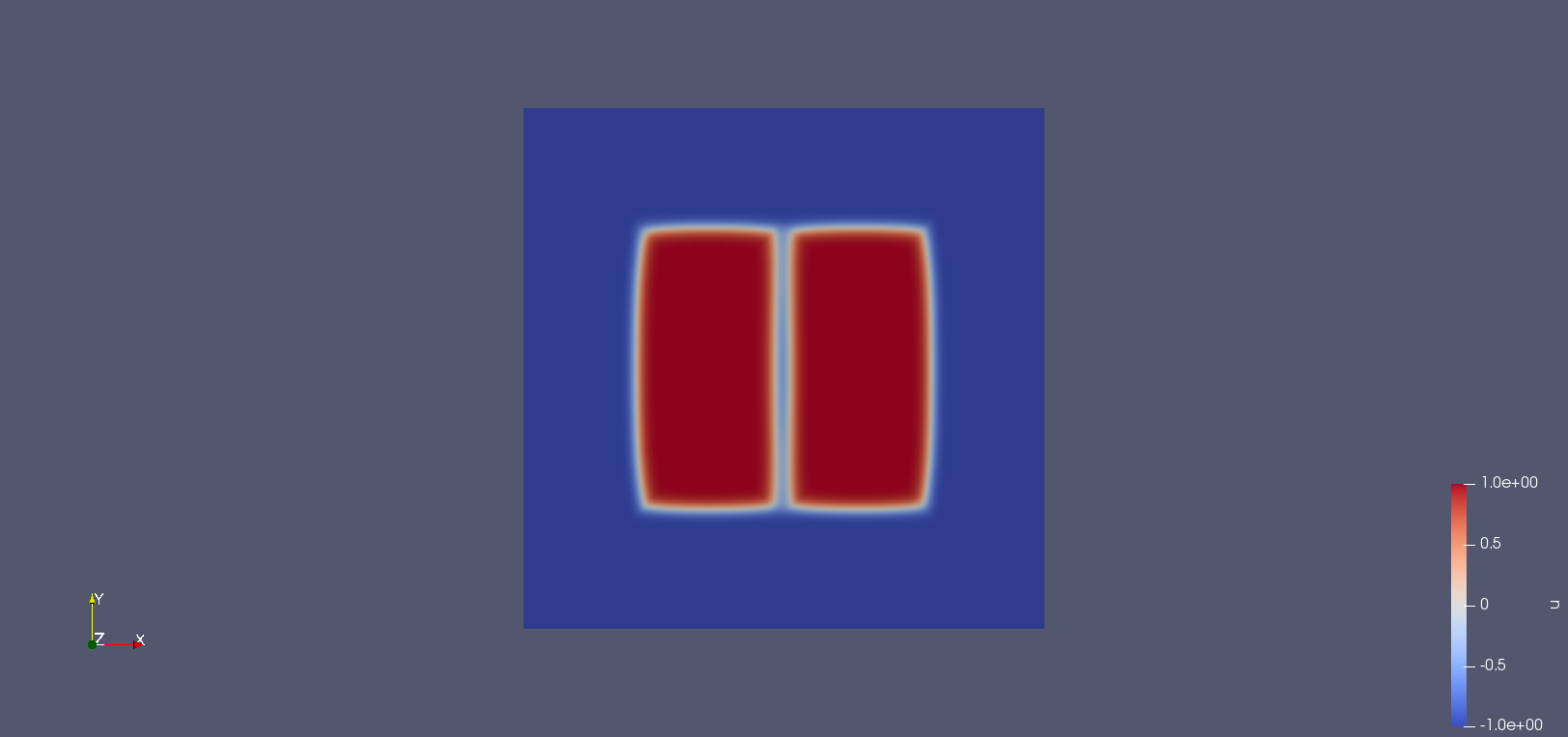}
		\includegraphics[trim={540px 112px 540px 112px}, clip, scale=.08]{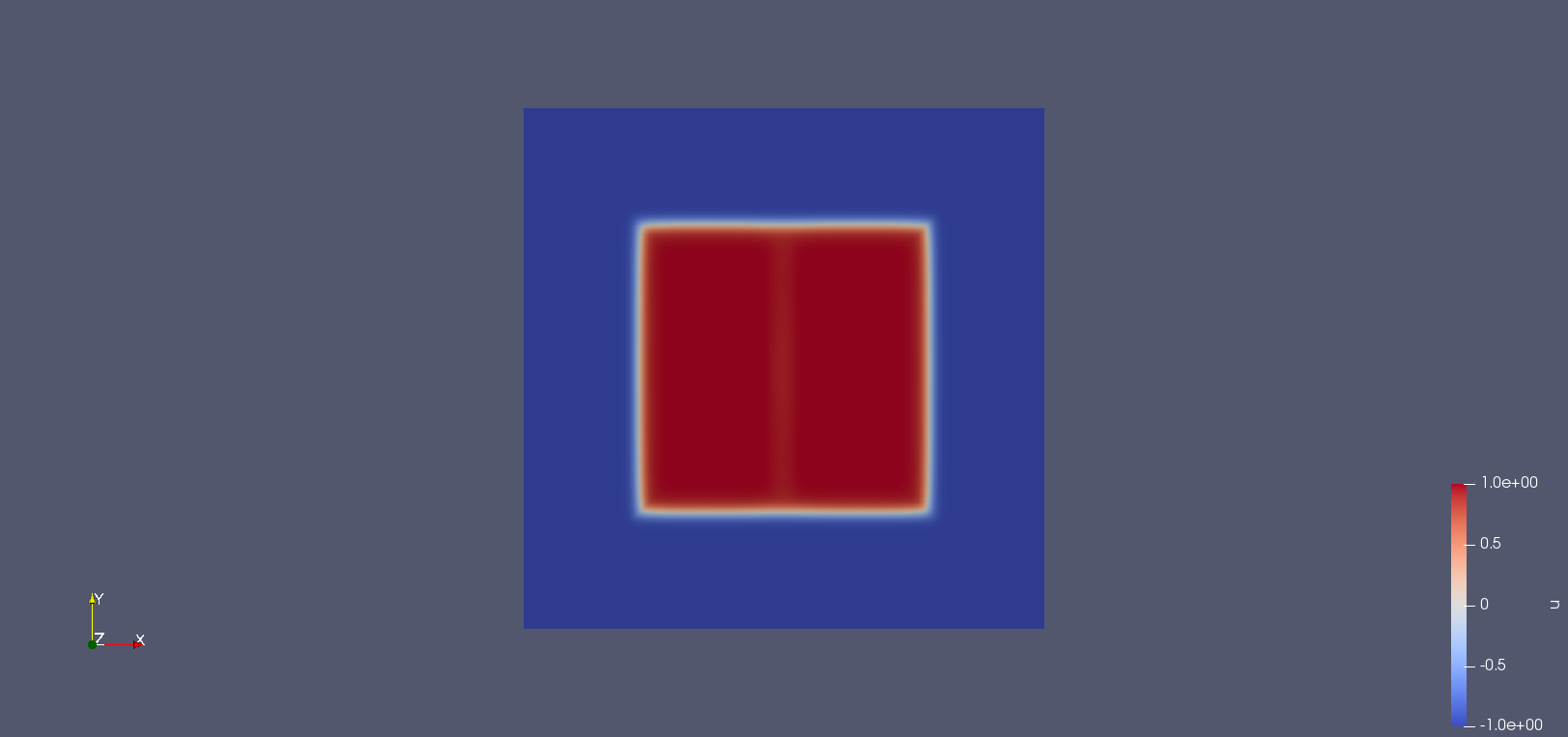}~~~
		\includegraphics[trim={1490px 0px 20px 490px}, clip, scale=.16]{{l1_merge_state.0000}.png}~\\~\\
		\includegraphics[trim={540px 112px 540px 112px}, clip, scale=.08]{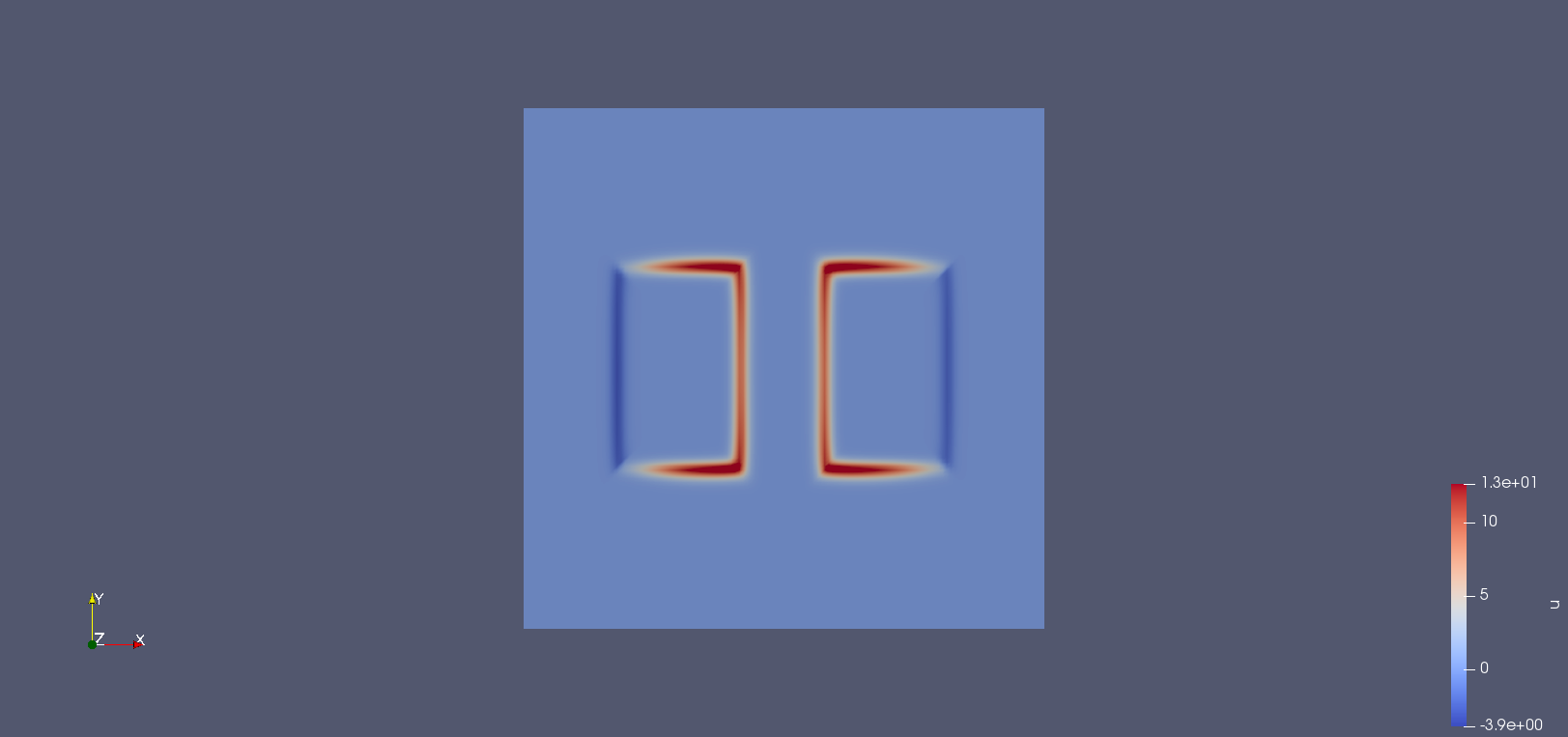}
		\includegraphics[trim={540px 112px 540px 112px}, clip, scale=.08]{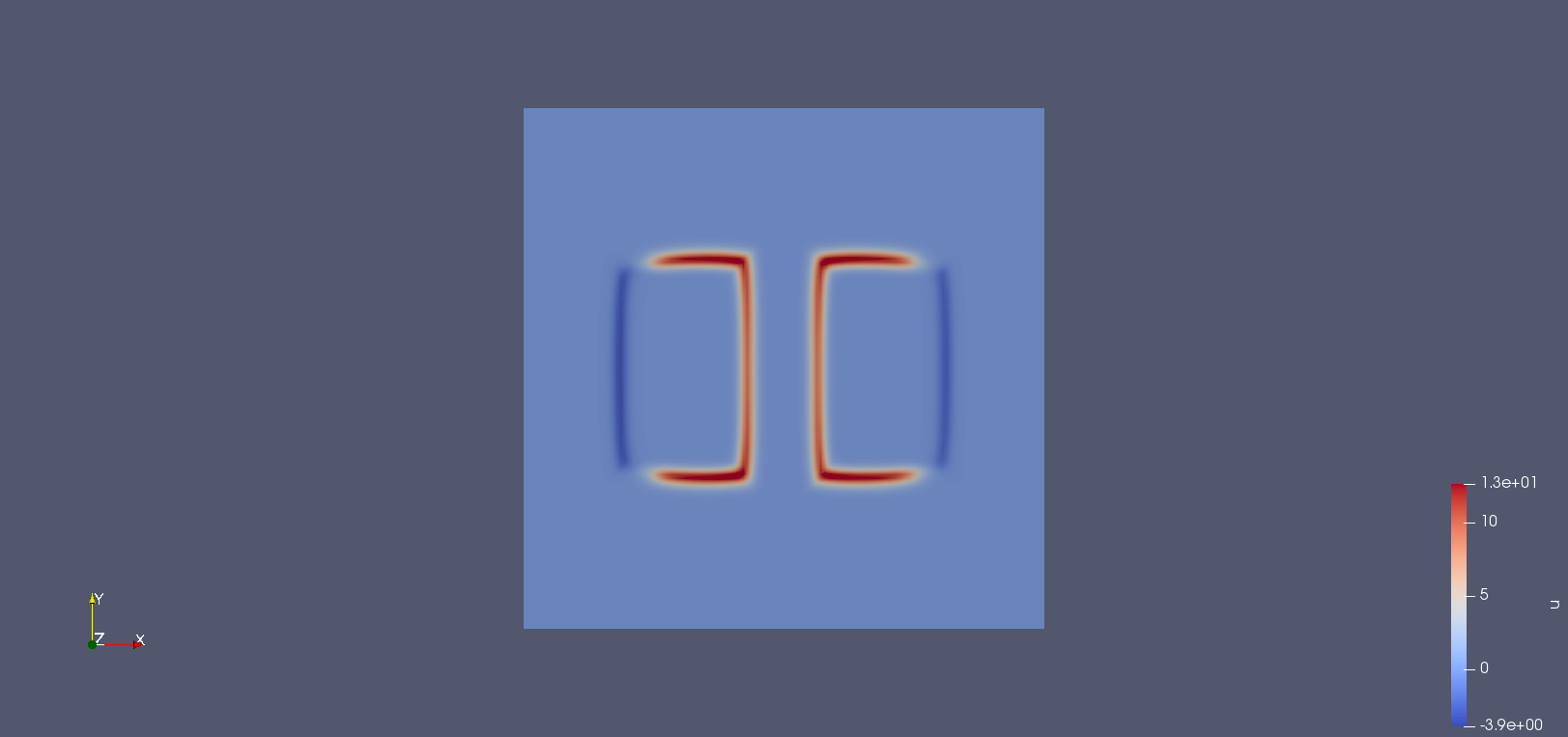}
		\includegraphics[trim={540px 112px 540px 112px}, clip, scale=.08]{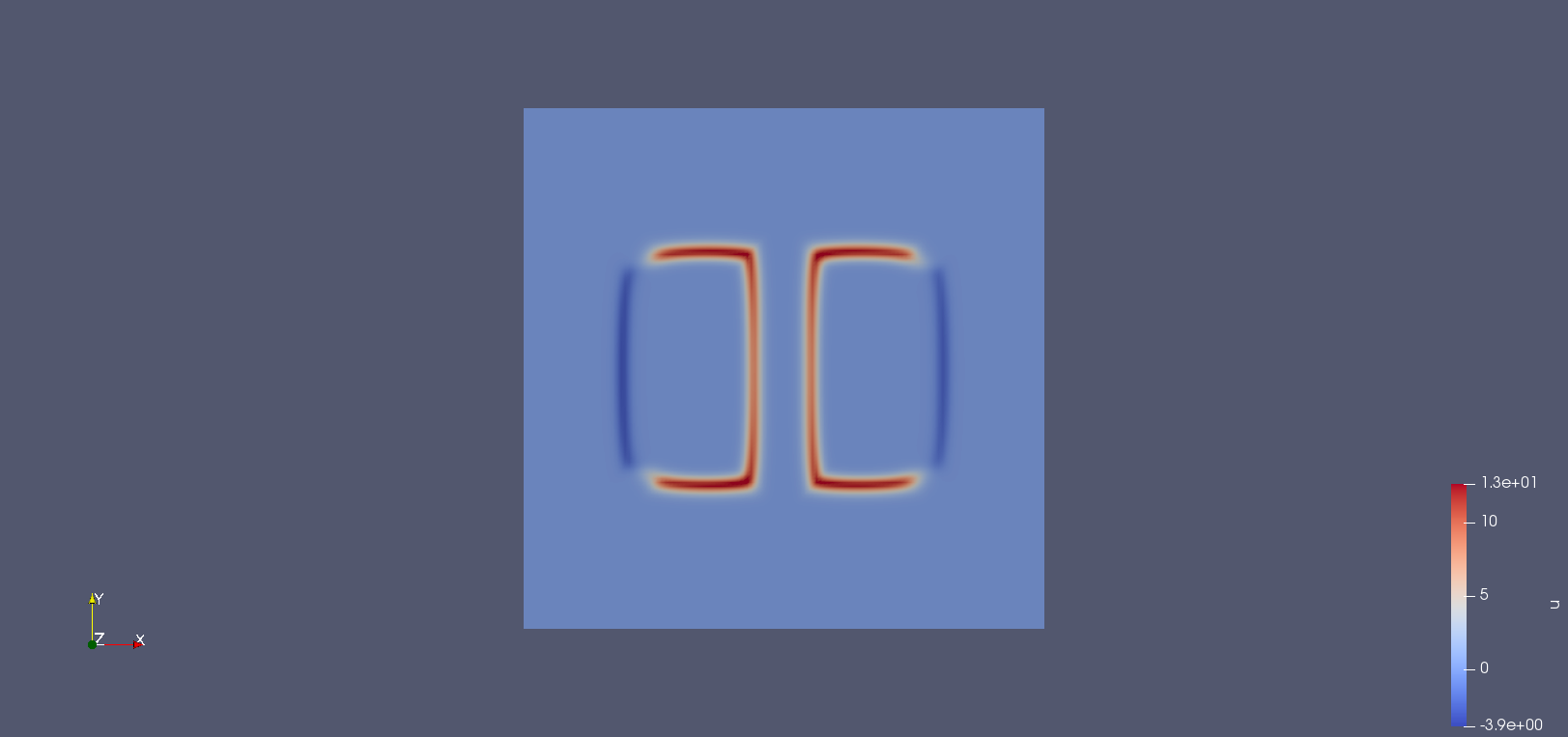}
		\includegraphics[trim={540px 112px 540px 112px}, clip, scale=.08]{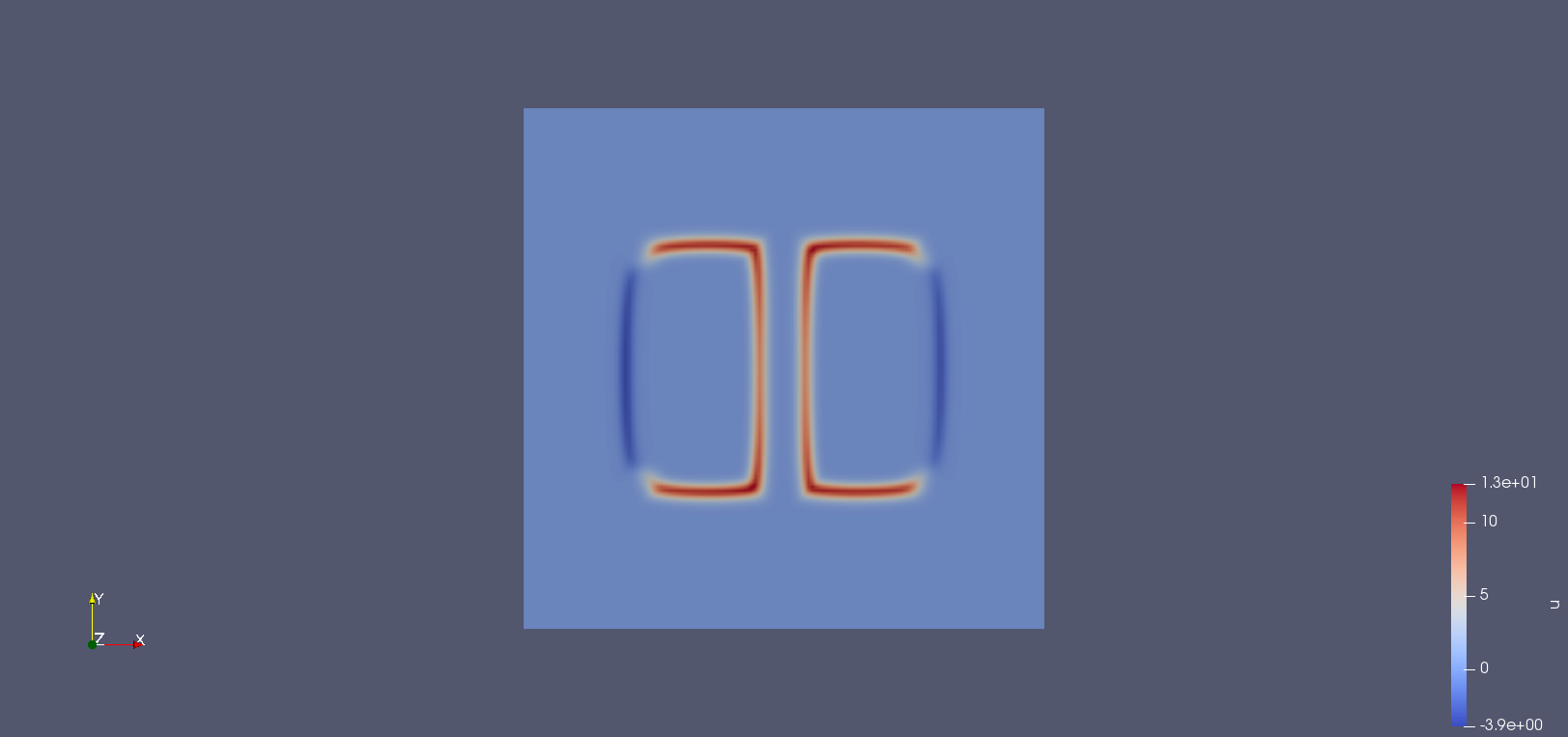}
		\includegraphics[trim={540px 112px 540px 112px}, clip, scale=.08]{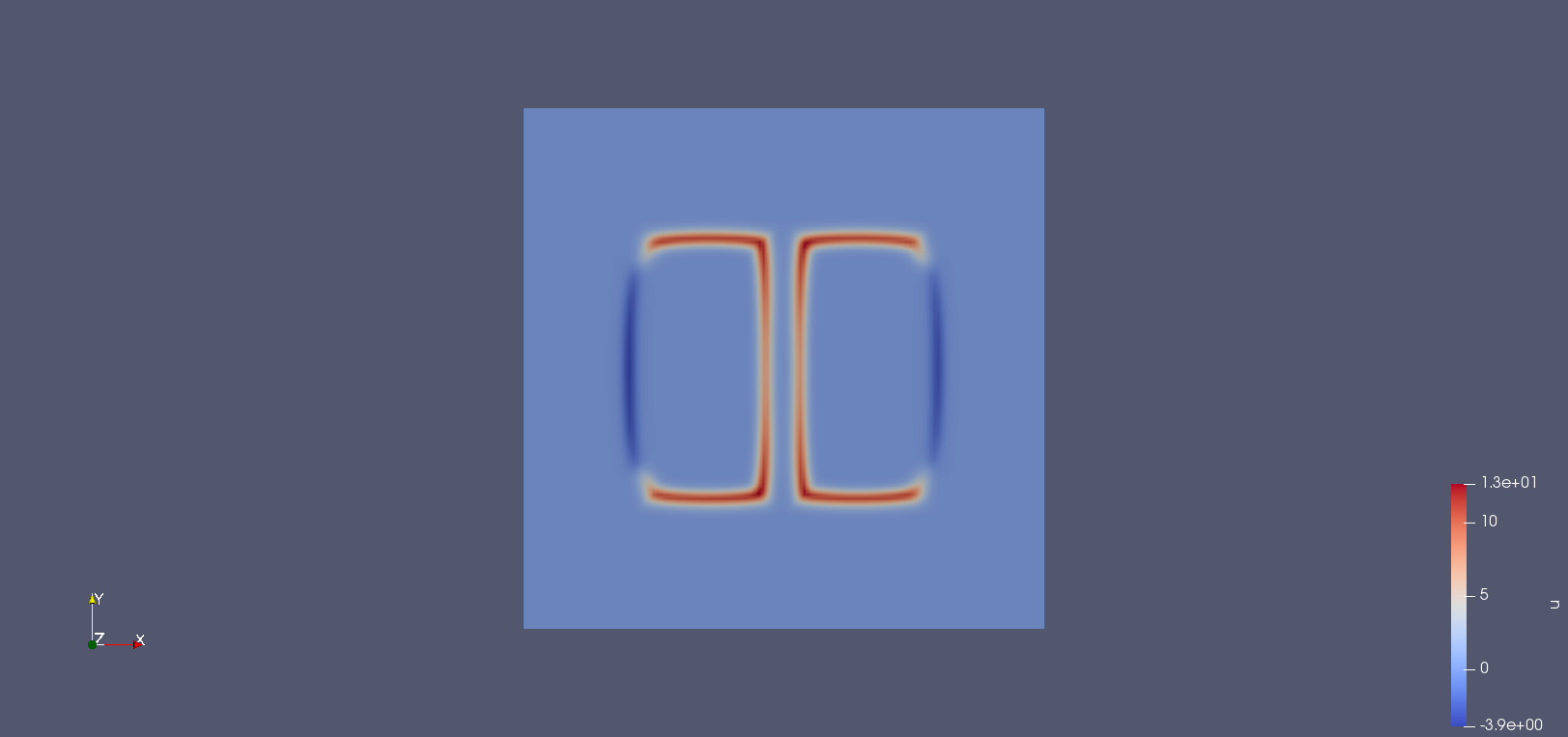}
		\includegraphics[trim={540px 112px 540px 112px}, clip, scale=.08]{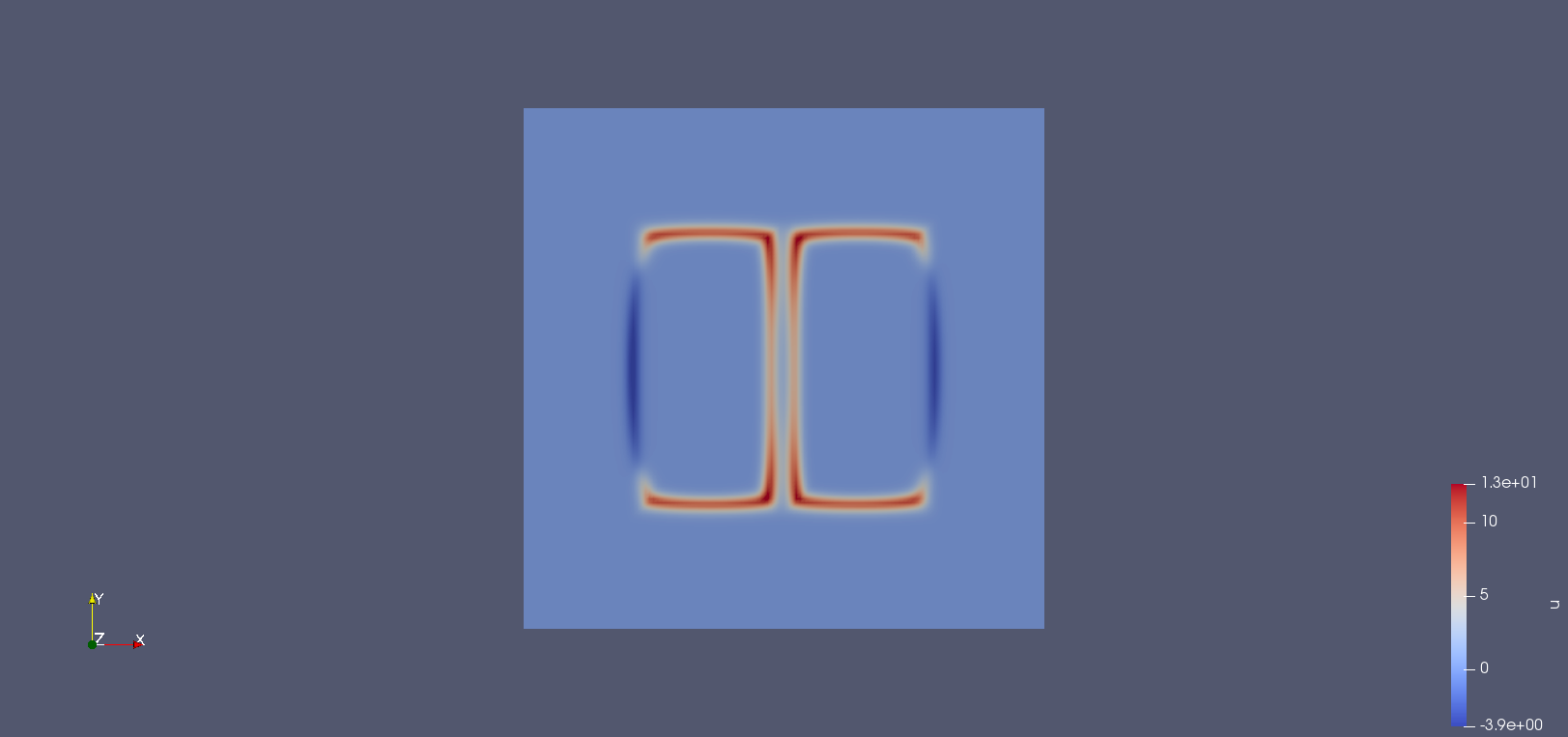}
		\includegraphics[trim={540px 112px 540px 112px}, clip, scale=.08]{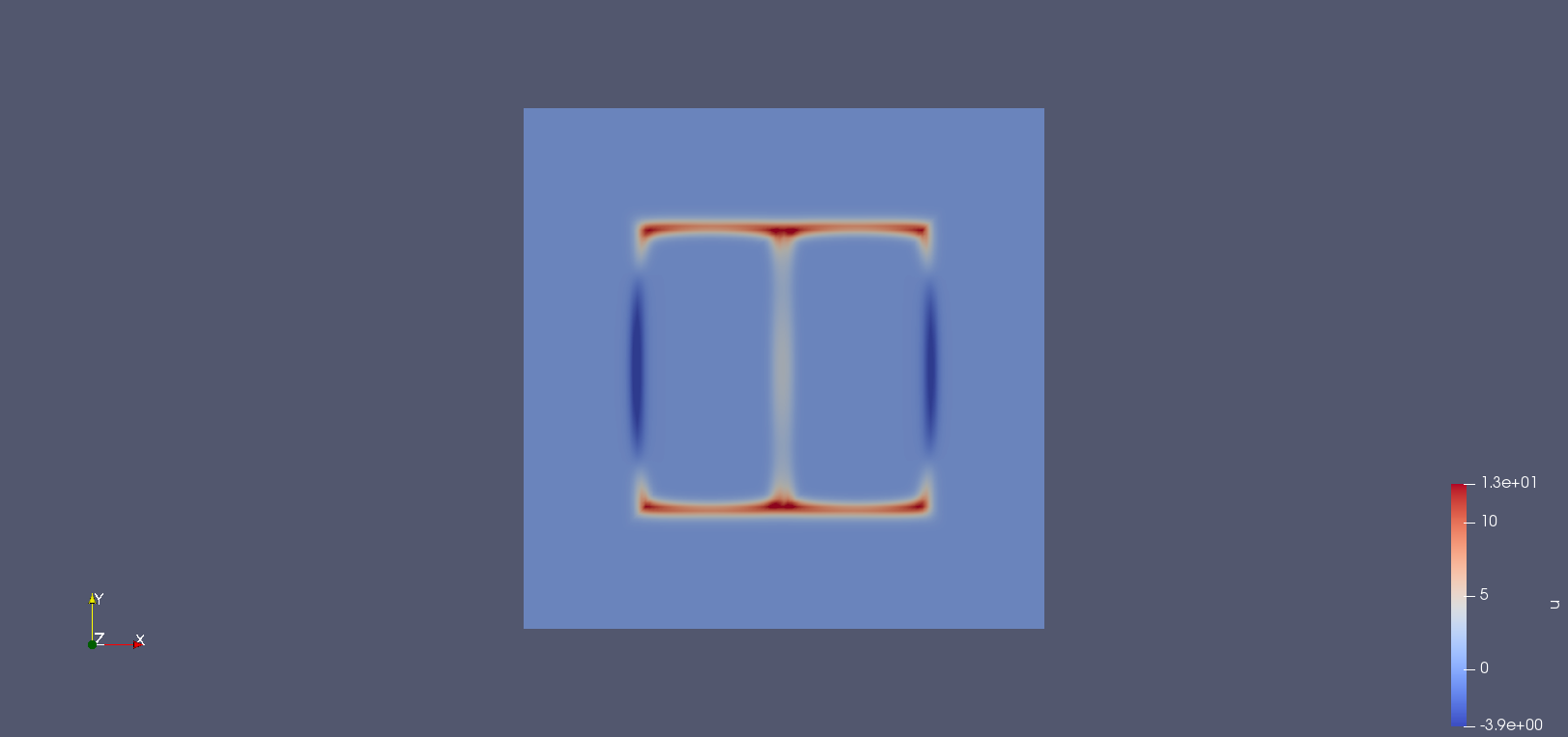}
		\includegraphics[trim={540px 112px 540px 112px}, clip, scale=.08]{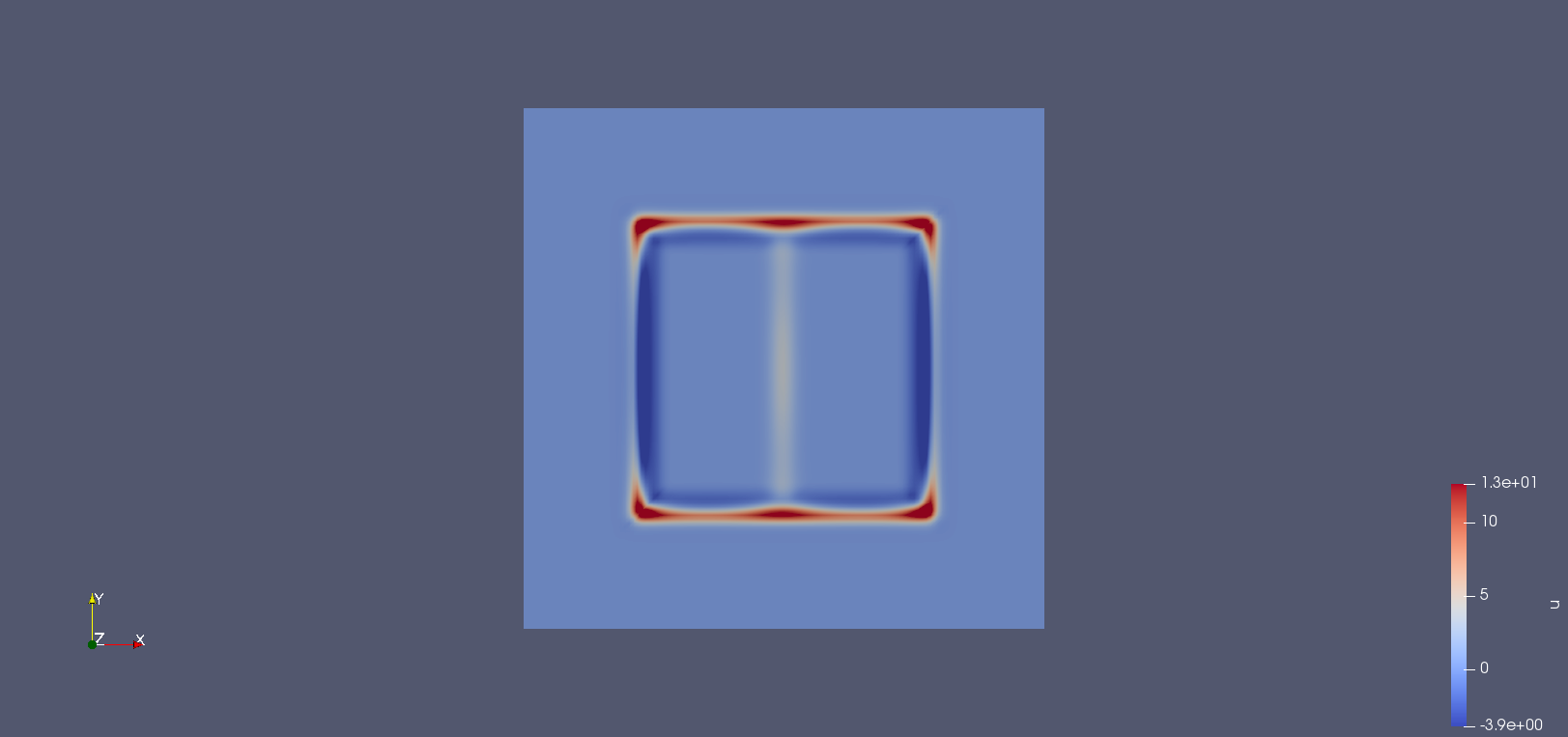}~~~
		\includegraphics[trim={1490px 0px 20px 490px}, clip, scale=.16]{{l1_merge_control.0000}.png}
	\caption{Result for `merge square' with the regularized $l_1$-norm.
	\label{fig:l1_merge}}
\end{center}    
\end{subfigure}
\medskip

\begin{subfigure}{1.0\textwidth}
\begin{center}  
	\includegraphics[trim={540px 112px 540px 112px}, clip, scale=.08]{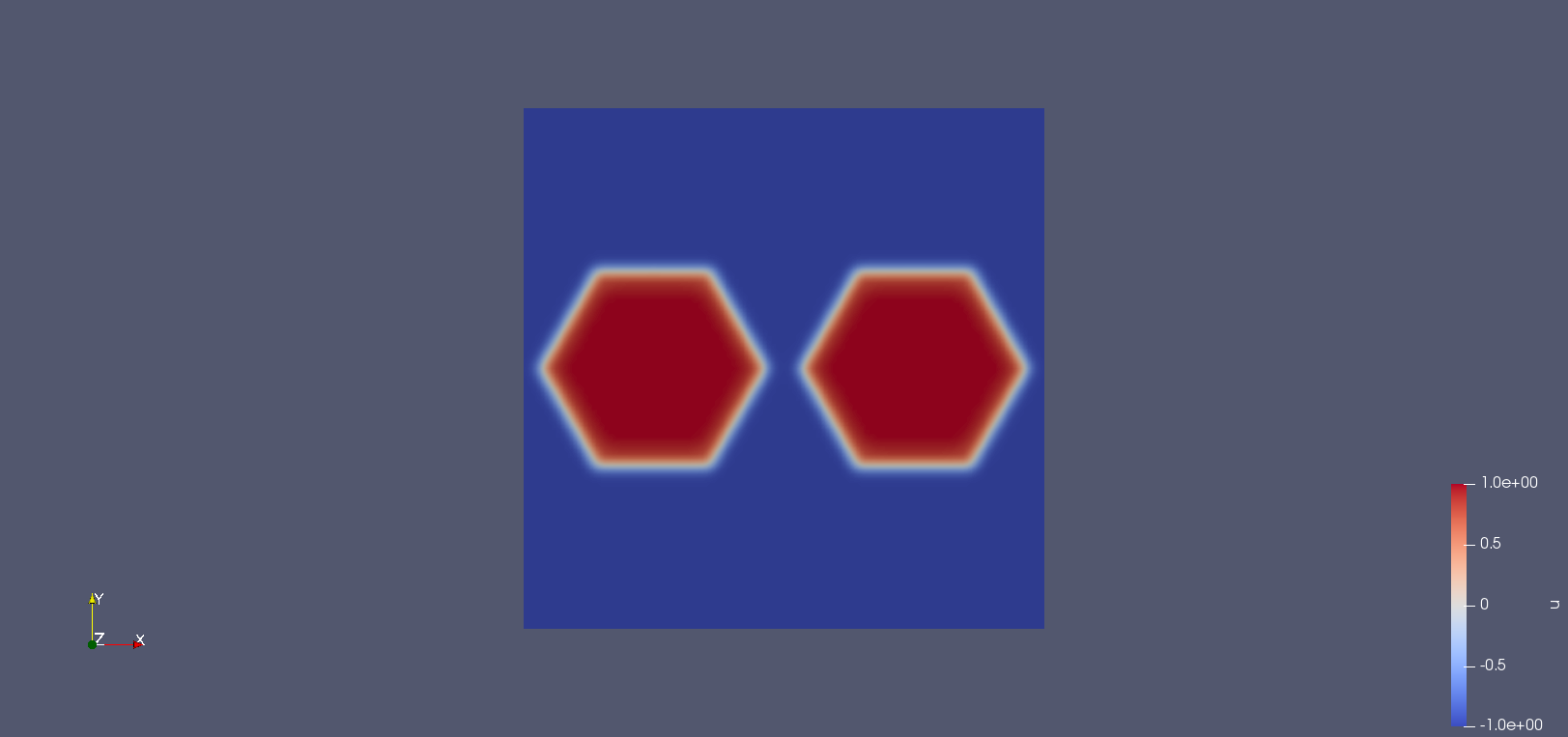}
	\includegraphics[trim={540px 112px 540px 112px}, clip, scale=.08]{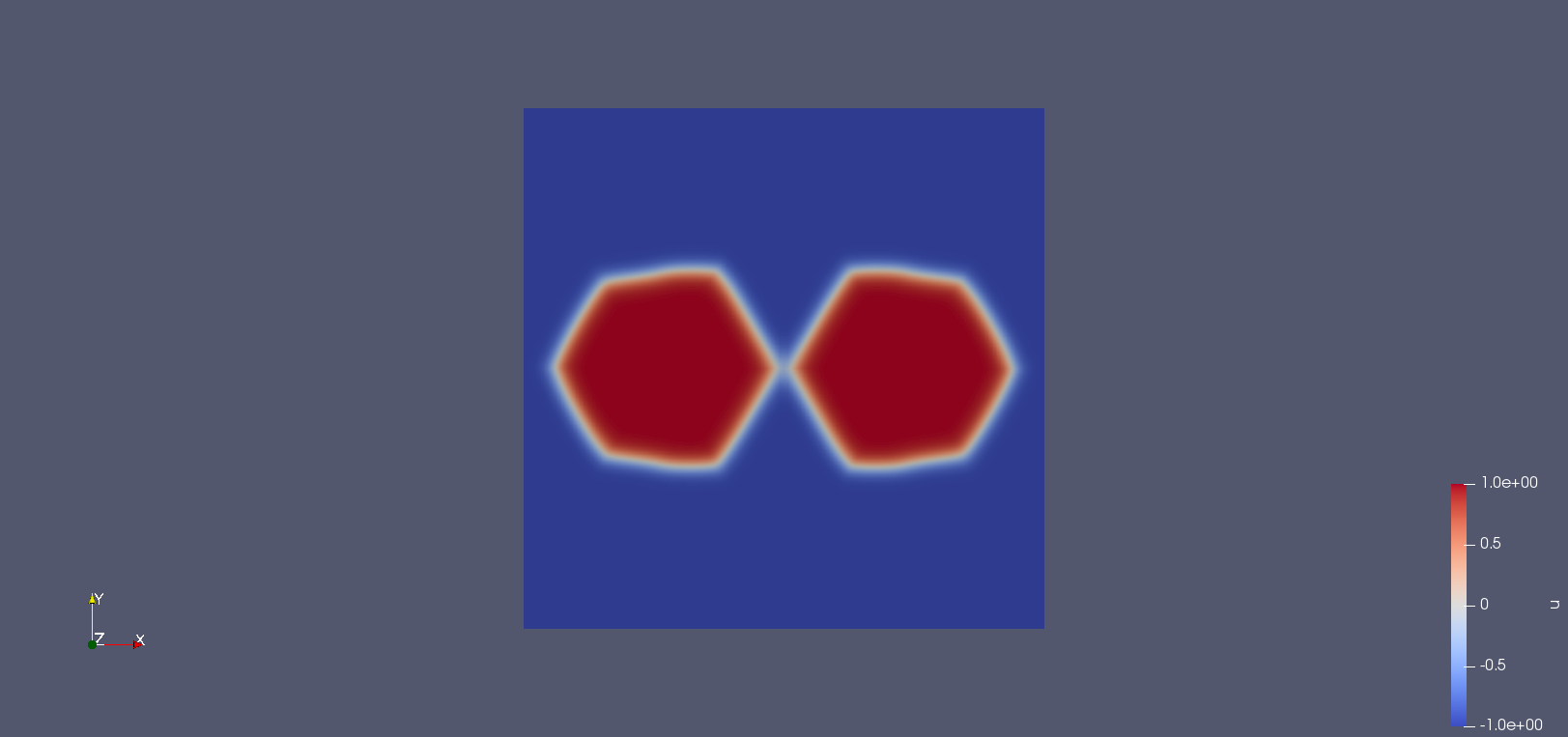}
	\includegraphics[trim={540px 112px 540px 112px}, clip, scale=.08]{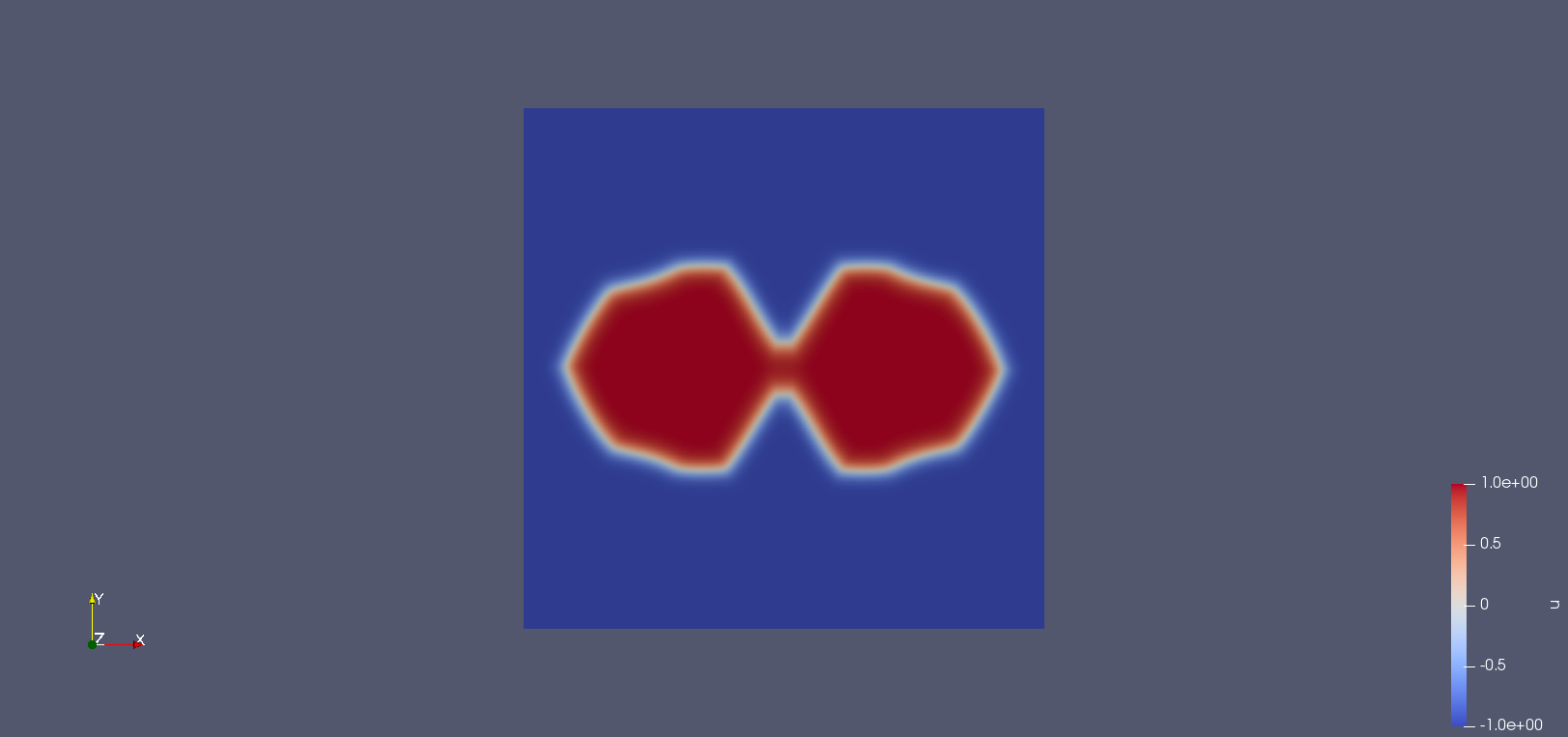}
	\includegraphics[trim={540px 112px 540px 112px}, clip, scale=.08]{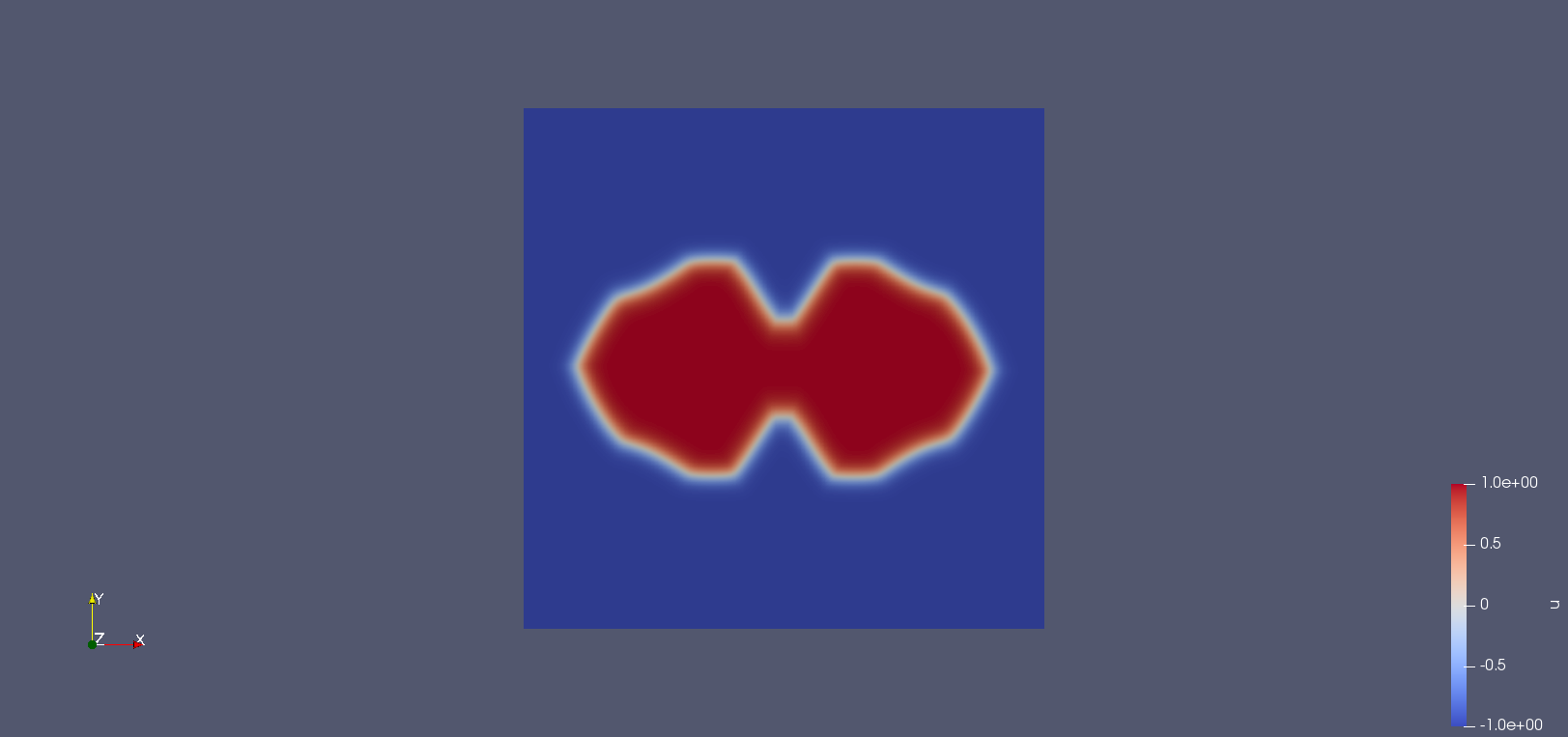}
	\includegraphics[trim={540px 112px 540px 112px}, clip, scale=.08]{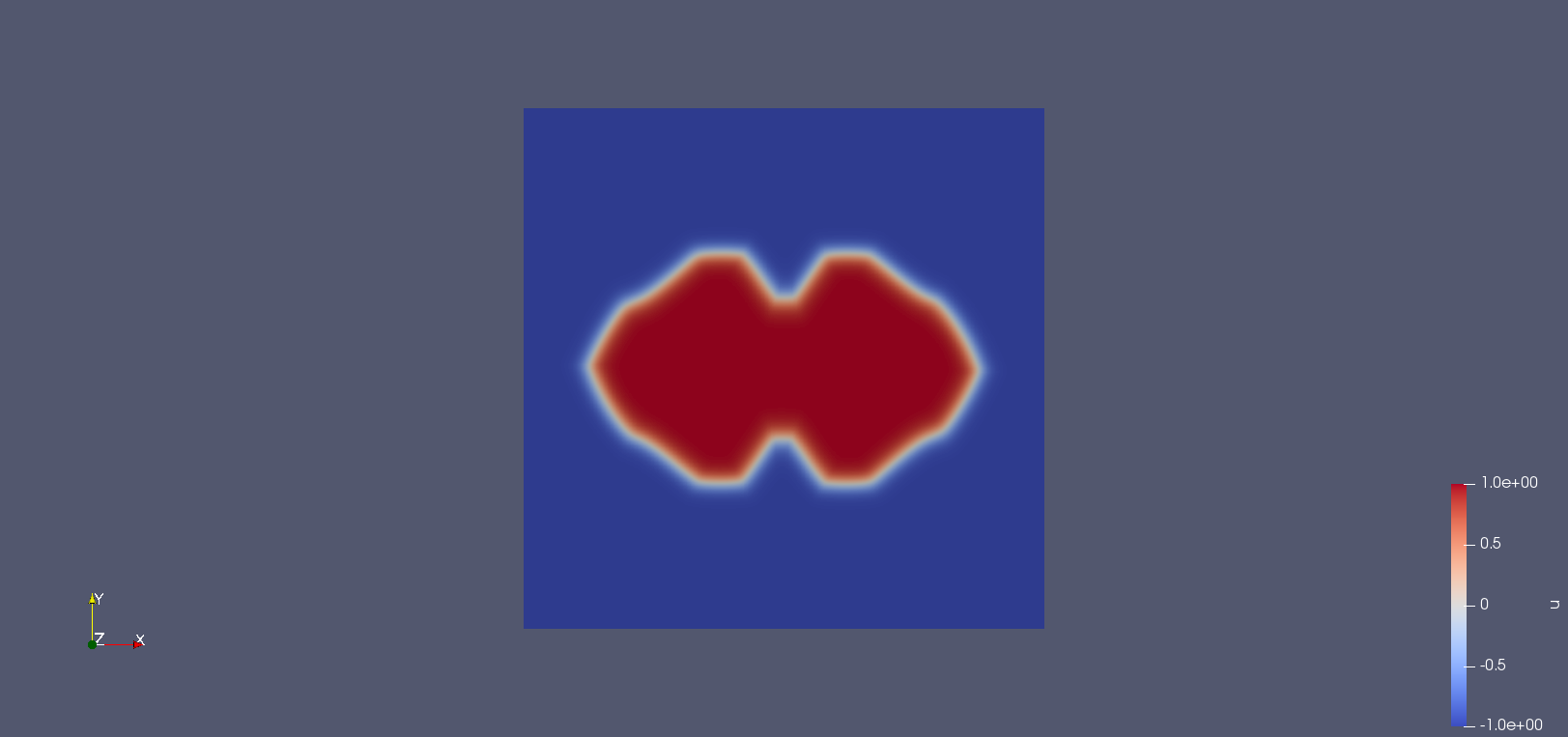}
	\includegraphics[trim={540px 112px 540px 112px}, clip, scale=.08]{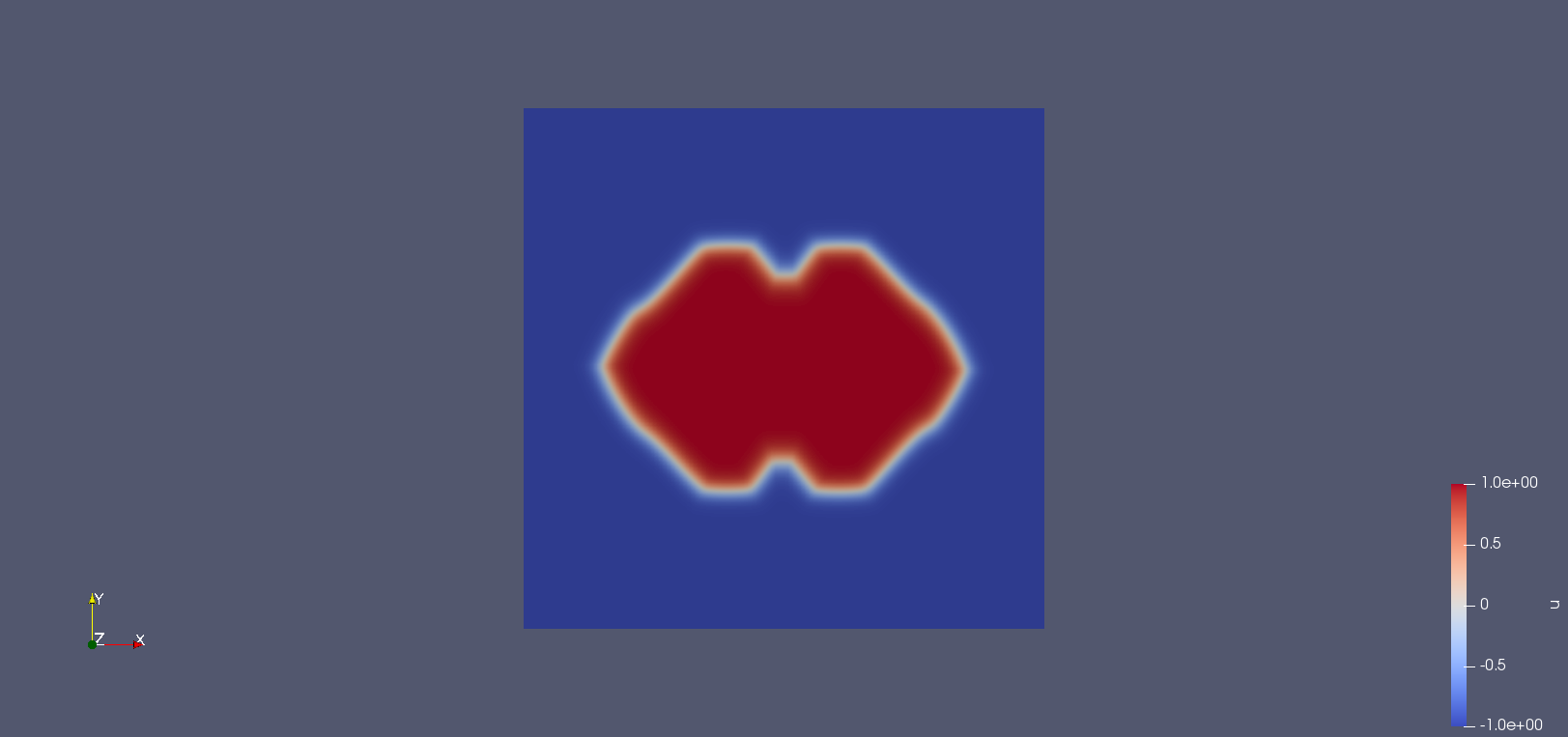}
	\includegraphics[trim={540px 112px 540px 112px}, clip, scale=.08]{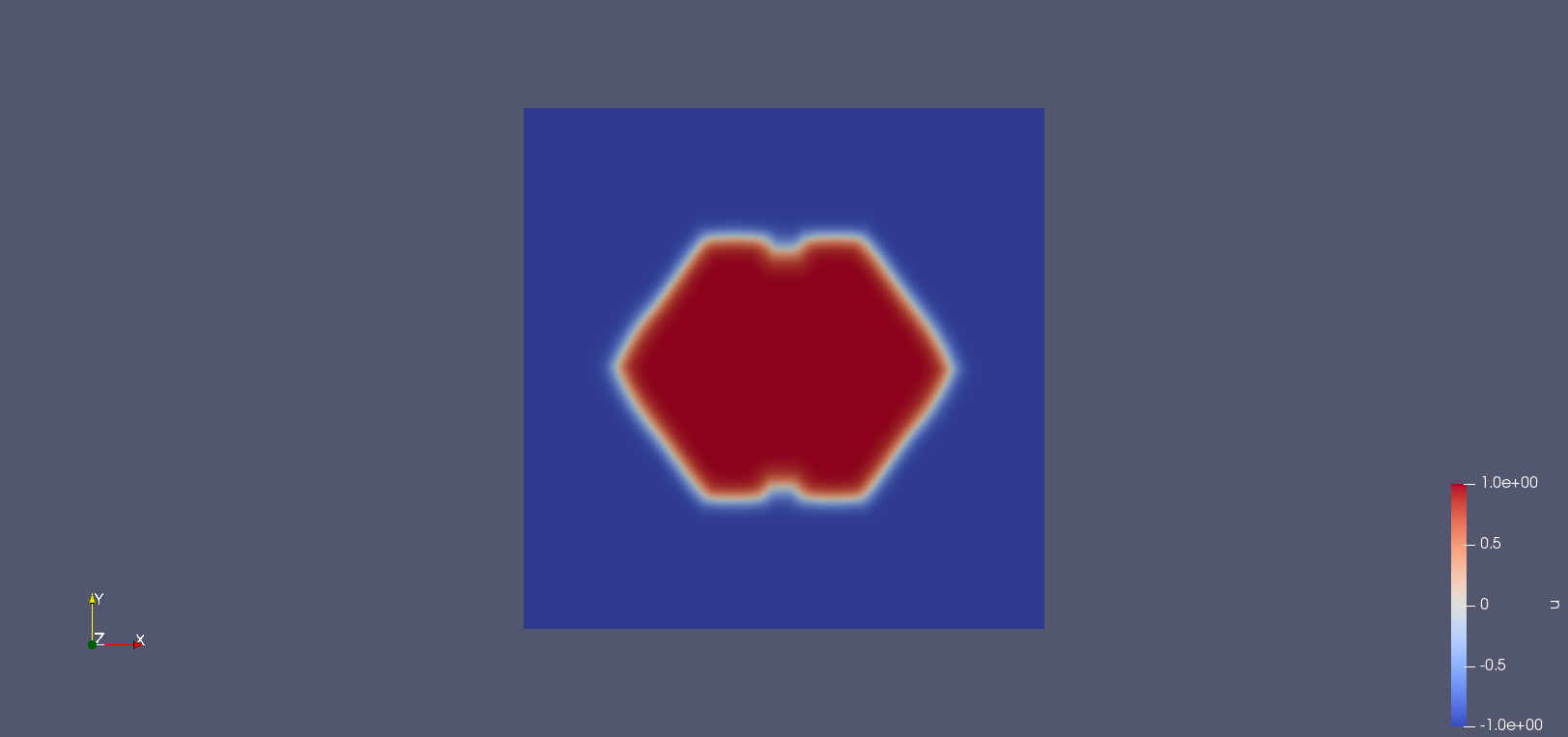}
	\includegraphics[trim={540px 112px 540px 112px}, clip, scale=.08]{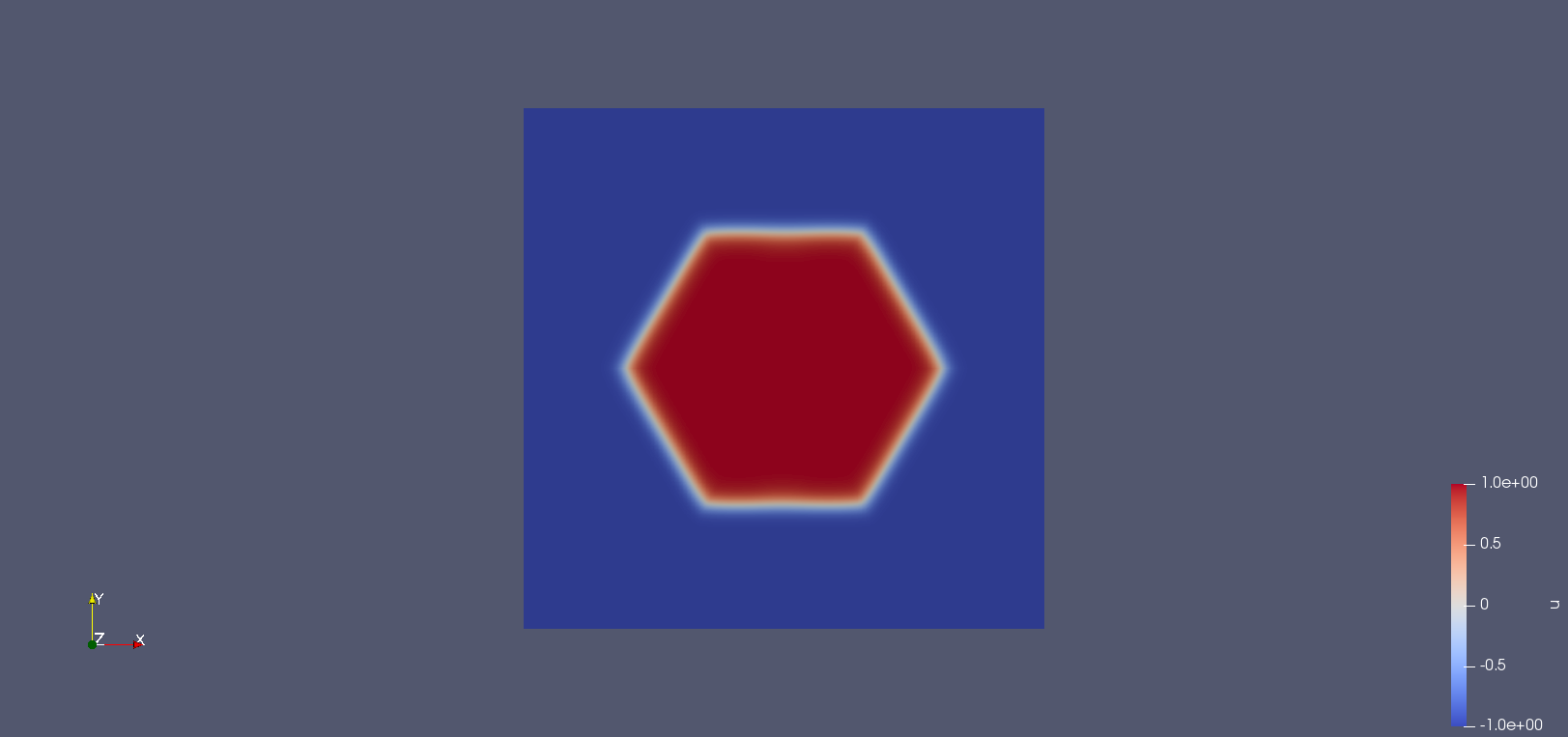}~~~
	\includegraphics[trim={1490px 0px 20px 490px}, clip, scale=.16]{{hexa_merge_state.0000}.png}~\\~\\
	\includegraphics[trim={540px 112px 540px 112px}, clip, scale=.08]{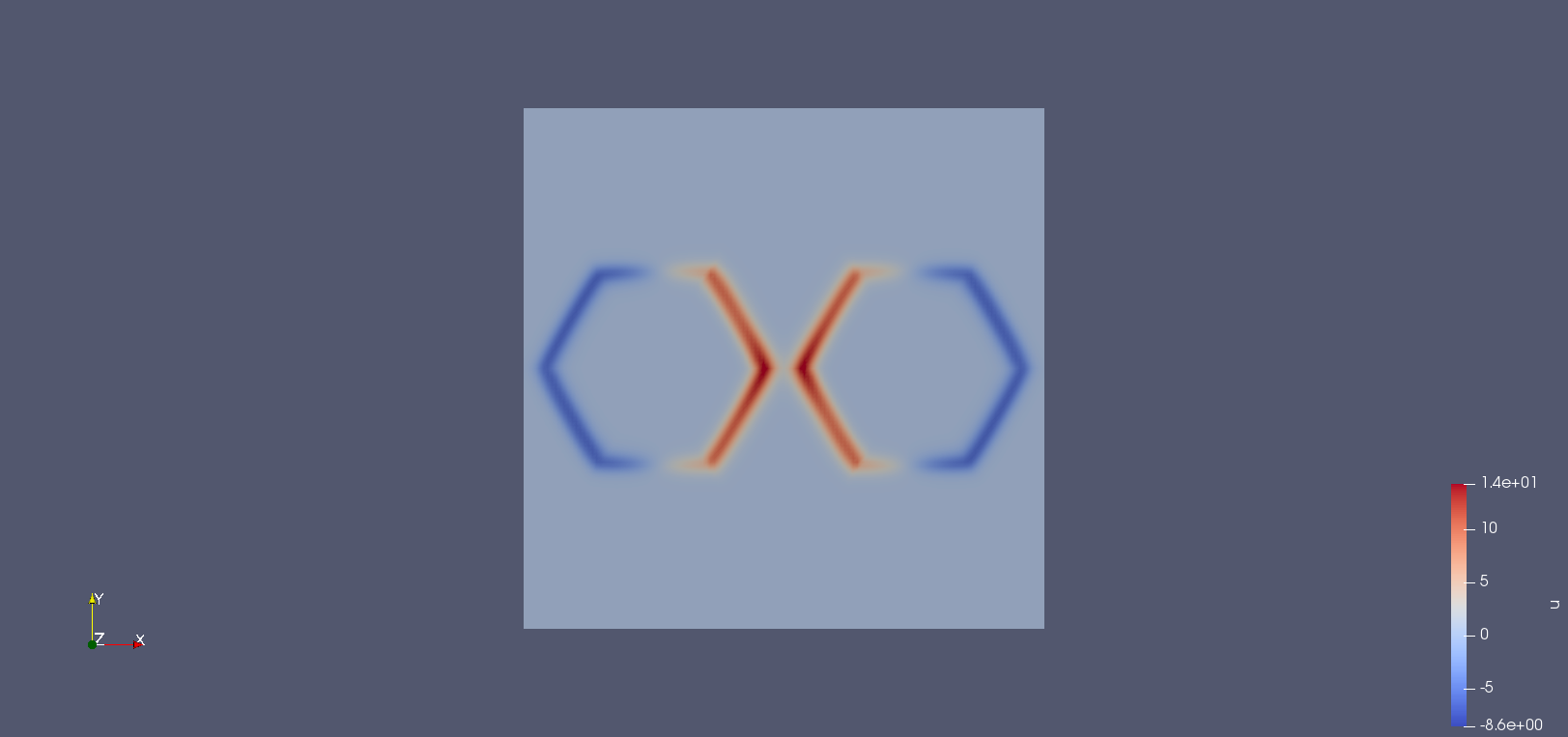}
	\includegraphics[trim={540px 112px 540px 112px}, clip, scale=.08]{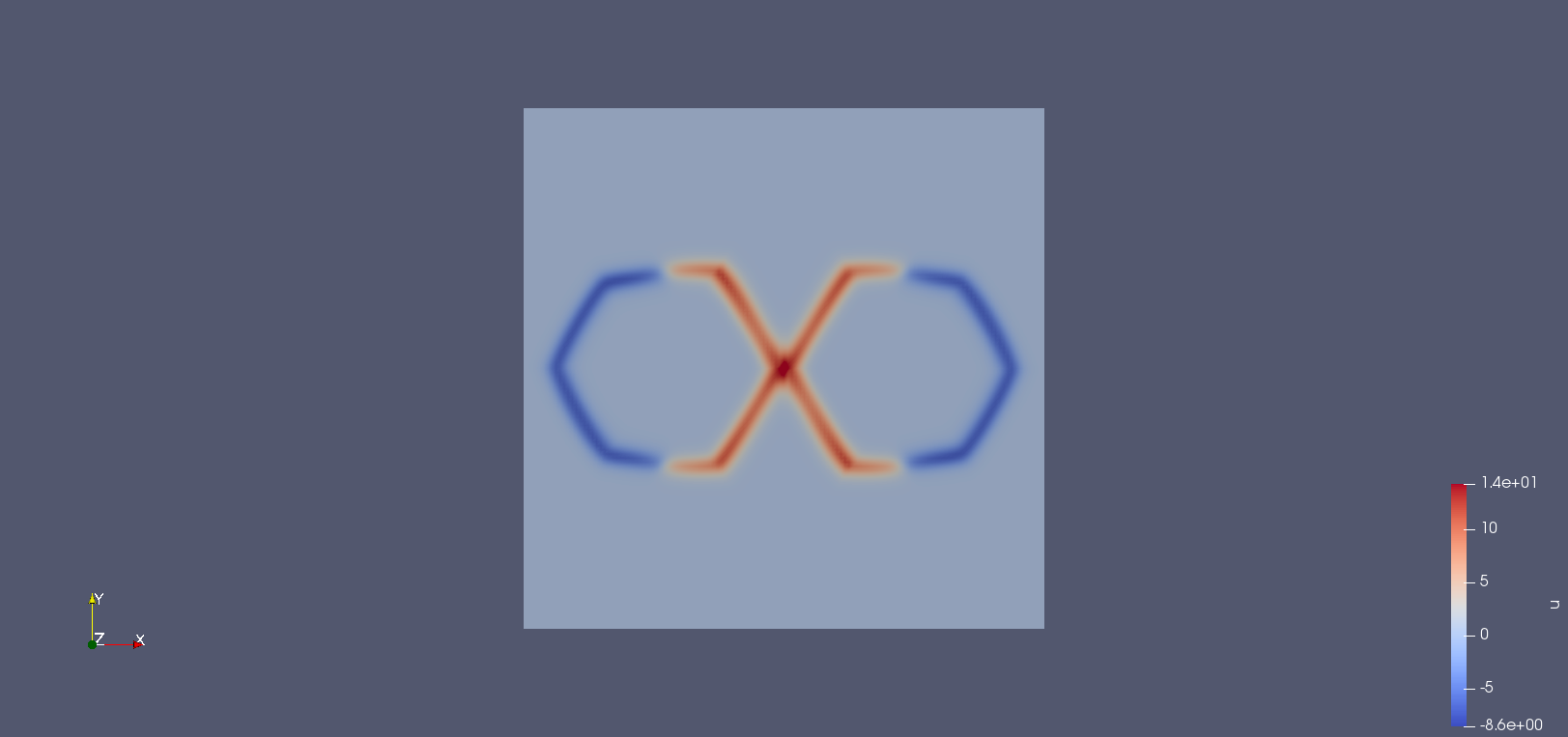}
	\includegraphics[trim={540px 112px 540px 112px}, clip, scale=.08]{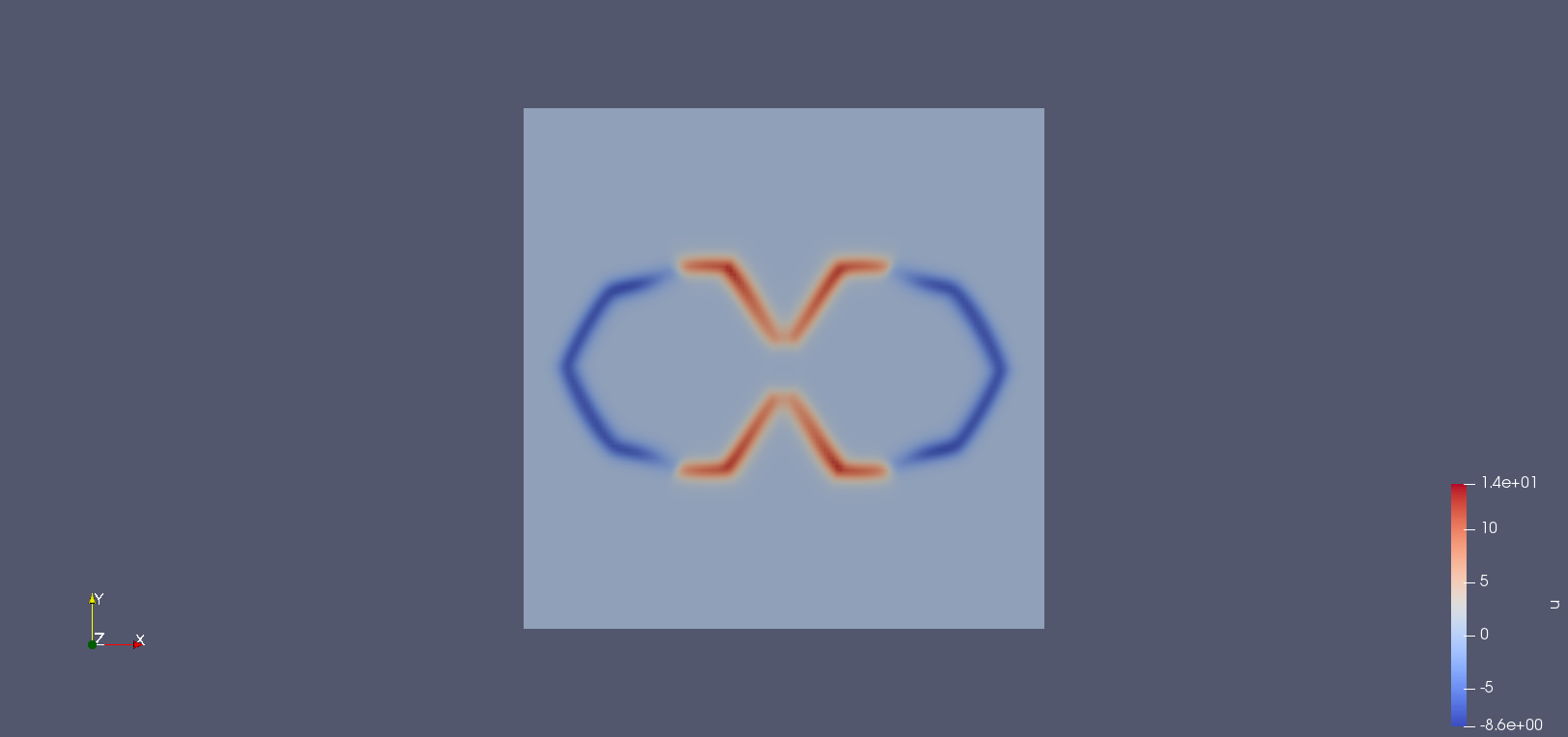}
	\includegraphics[trim={540px 112px 540px 112px}, clip, scale=.08]{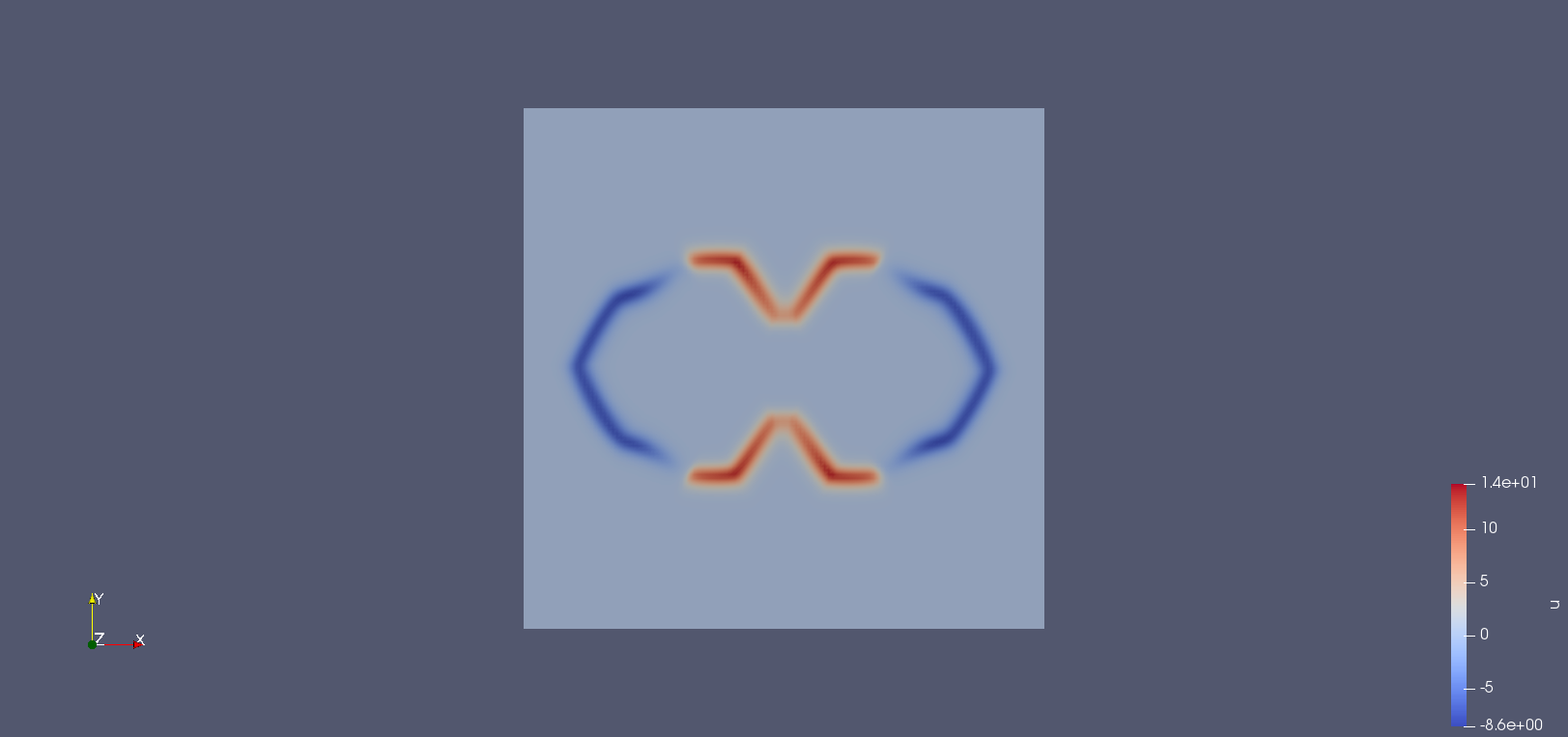}
	\includegraphics[trim={540px 112px 540px 112px}, clip, scale=.08]{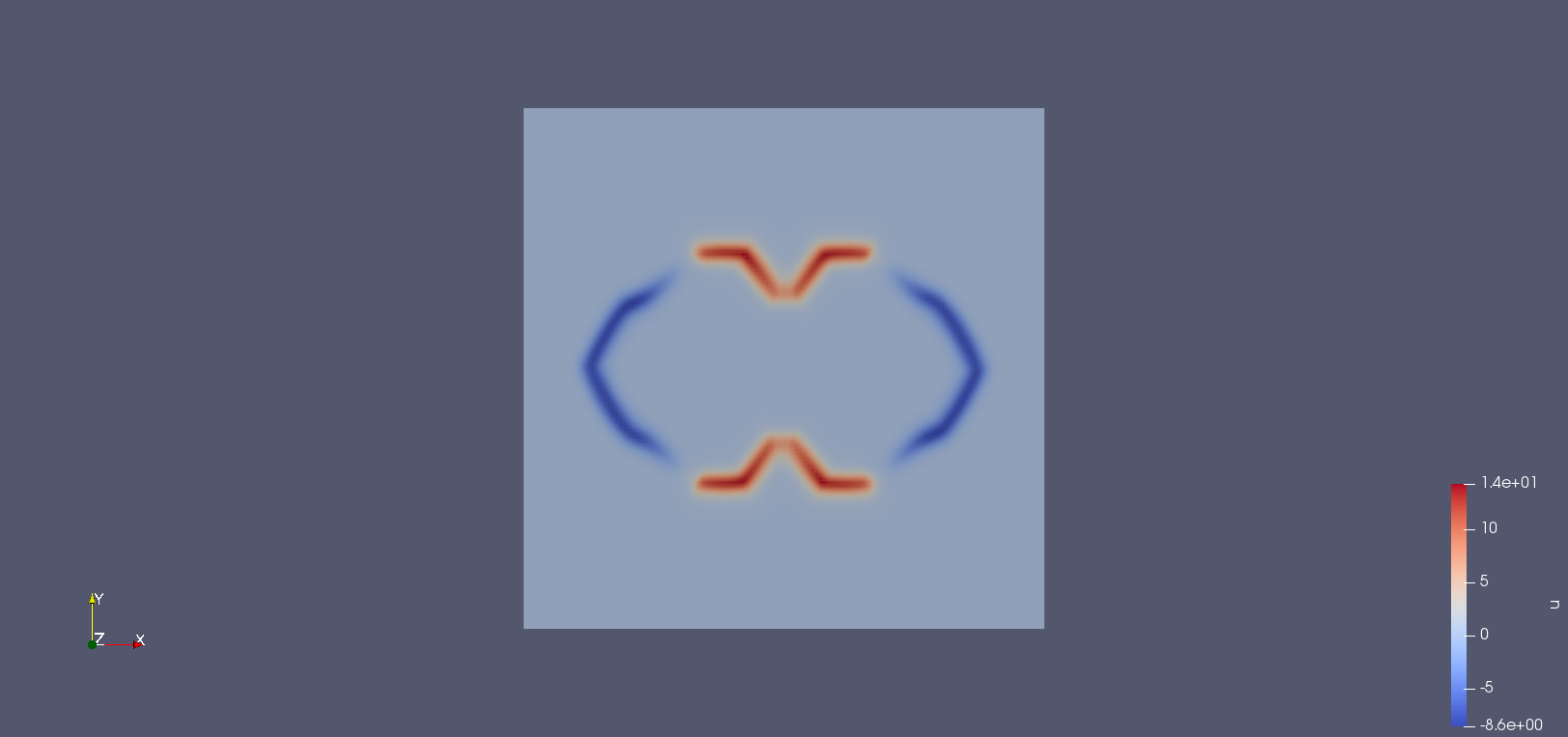}
	\includegraphics[trim={540px 112px 540px 112px}, clip, scale=.08]{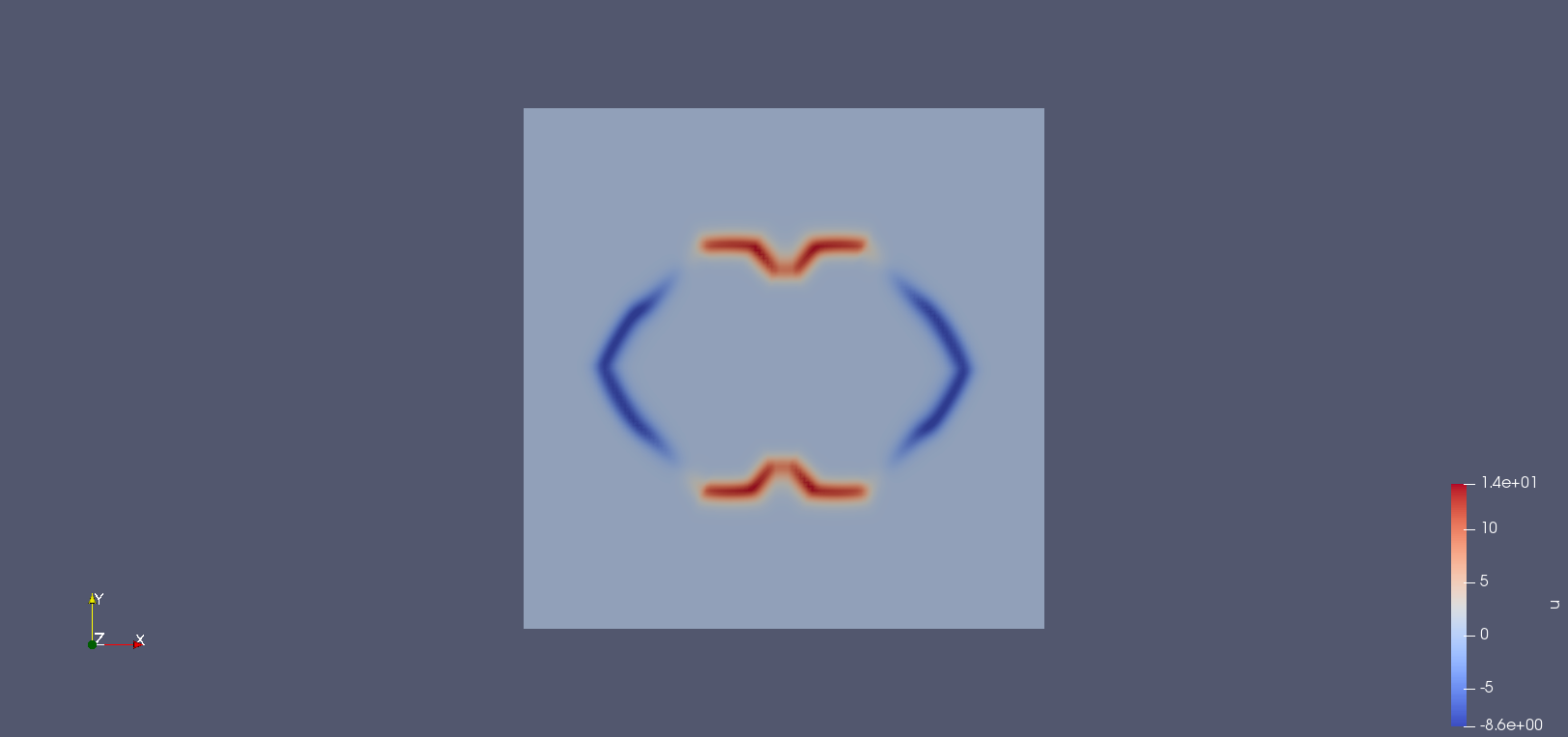}
	\includegraphics[trim={540px 112px 540px 112px}, clip, scale=.08]{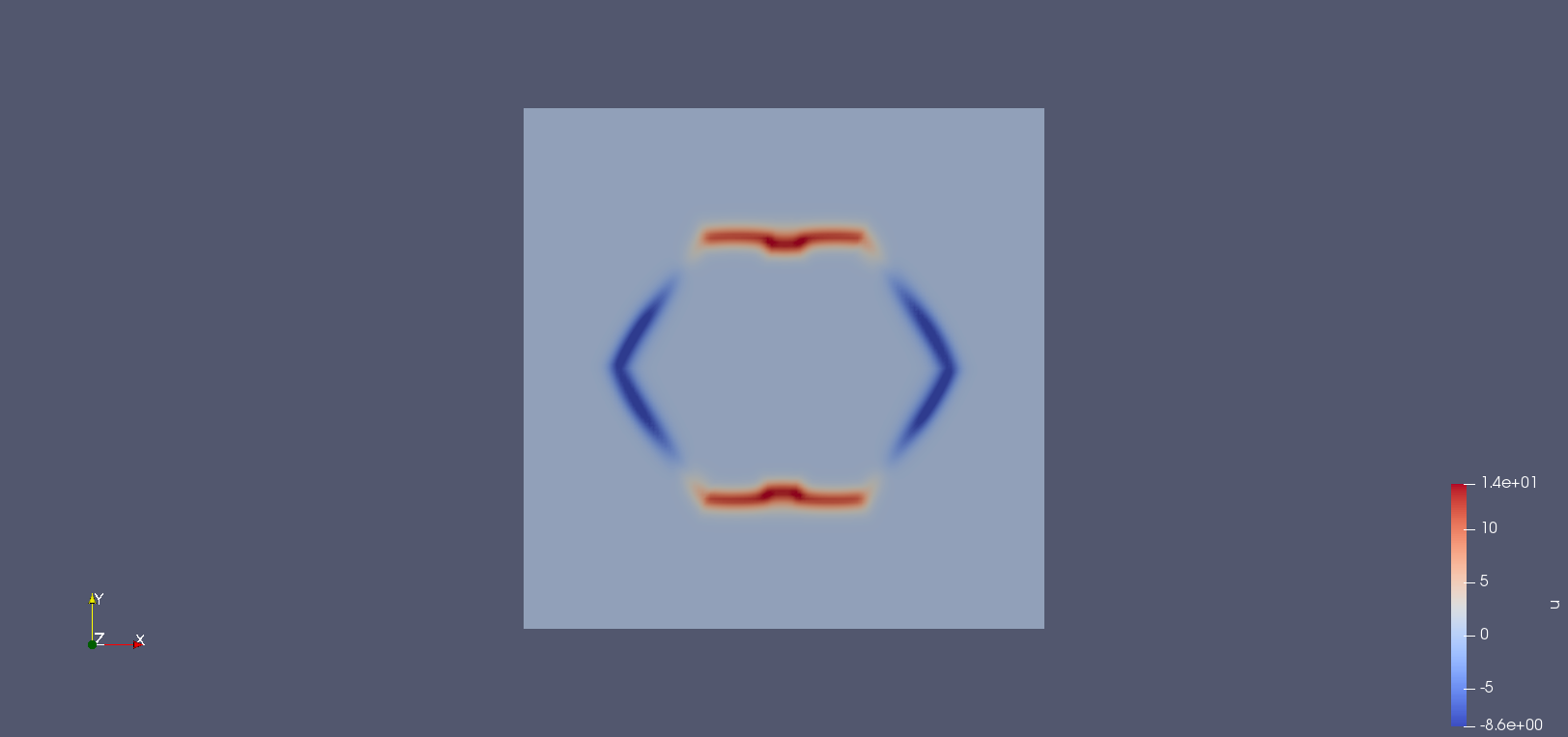}
	\includegraphics[trim={540px 112px 540px 112px}, clip, scale=.08]{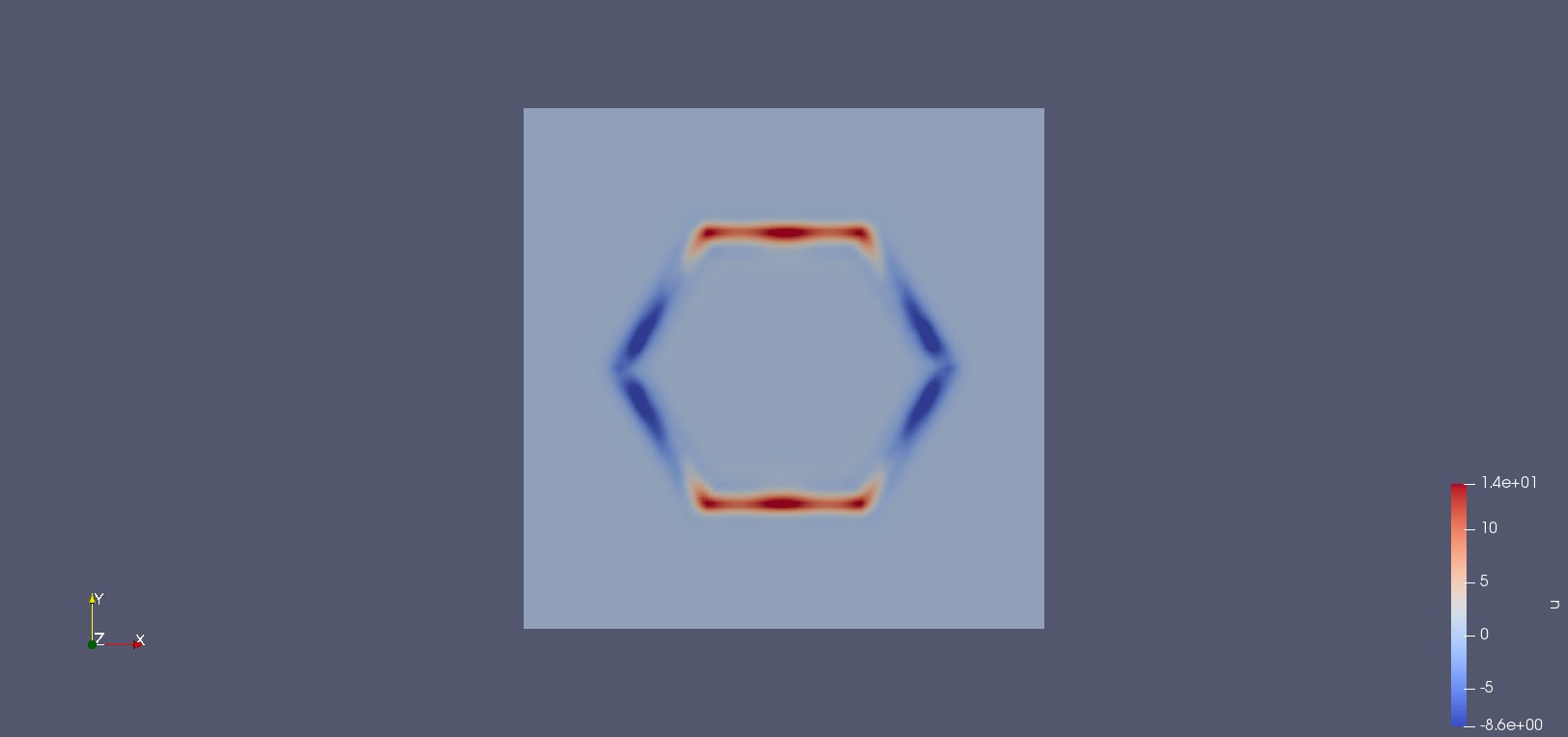}~~~
	\includegraphics[trim={1490px 0px 20px 490px}, clip, scale=.16]{{hexa_merge_control.0000}.png}  
	\caption{Result for `merge hexagon' using the hexagon anisotropy.
	\label{fig:hexa_merge}}
\end{center}    
\end{subfigure}
\caption{State (above) and corresponding control (below) for the solution of merge geometries.
  %at $t = timesnotlookedupsofar$.} % Hinweis: die Zahlen entsprechen ungefähr diesen Zeitschritten, ich hab die Bilder von 0-7 darauf umgerechnet und auf "schöne" Zahlen gerundet
	\label{fig:merge}}
\end{center}    
\end{figure}
\fi

\begin{table}[htbp]
\begin{tabular}{|l|l|l|l|}
		\hline
		& iso & $l_1$  & hexa \\ \hline
		$j(u)$    & 0.103955 & 0.0562286 & 0.0884921  \\ \hline
		$j_1+j_2$ & 0.00409254 + 0.0998625 & 0.000728494 + 0.0555001 & 0.00115402 + 0.0873381  \\ \hline
	\end{tabular}
	\caption{Values of the cost functional for splitting geometries.
	\label{tab:comp_splitting} }
\end{table}

\begin{table}[htbp]
	\begin{tabular}{|l|l|l|l|}
		\hline
		& iso & $l_1$  & hexa \\ \hline
		$j(u)$    & 0.0666414 & 0.0496374 & 0.0588482  \\ \hline
		$j_1+j_2$ & 0.00274715 + 0.0638943 & 0.0010941 + 0.0485433 & 0.000905395 + 0.0579428  \\ \hline
	\end{tabular}
        \caption{Values of the cost functional for merging geometries.
	\label{tab:comp_merge}}
\end{table}

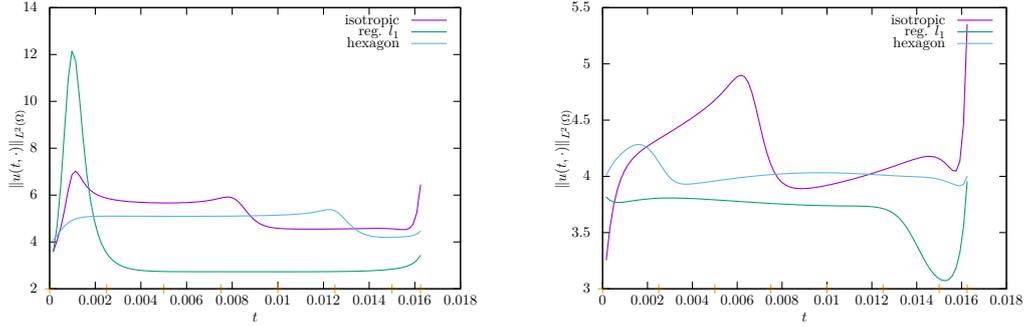
\begin{figure}
  \begin{center}
\begin{subfigure}{.48\textwidth}
  \begin{center}
 	\begin{tikzpicture}[gnuplot, scale=0.5, every node/.style={scale=0.5}]
%% generated with GNUPLOT 5.2p2 (Lua 5.3; terminal rev. 99, script rev. 102)
%% Di 11 Mai 2021 18:11:36 CEST
\path (0.000,0.000) rectangle (12.500,8.750);
\gpcolor{color=gp lt color border}
\gpsetlinetype{gp lt border}
\gpsetdashtype{gp dt solid}
\gpsetlinewidth{1.00}
\draw[gp path] (1.136,0.985)--(1.316,0.985);
\draw[gp path] (11.947,0.985)--(11.767,0.985);
\node[gp node right] at (0.952,0.985) {$2$};
\draw[gp path] (1.136,2.228)--(1.316,2.228);
\draw[gp path] (11.947,2.228)--(11.767,2.228);
\node[gp node right] at (0.952,2.228) {$4$};
\draw[gp path] (1.136,3.470)--(1.316,3.470);
\draw[gp path] (11.947,3.470)--(11.767,3.470);
\node[gp node right] at (0.952,3.470) {$6$};
\draw[gp path] (1.136,4.713)--(1.316,4.713);
\draw[gp path] (11.947,4.713)--(11.767,4.713);
\node[gp node right] at (0.952,4.713) {$8$};
\draw[gp path] (1.136,5.956)--(1.316,5.956);
\draw[gp path] (11.947,5.956)--(11.767,5.956);
\node[gp node right] at (0.952,5.956) {$10$};
\draw[gp path] (1.136,7.198)--(1.316,7.198);
\draw[gp path] (11.947,7.198)--(11.767,7.198);
\node[gp node right] at (0.952,7.198) {$12$};
\draw[gp path] (1.136,8.441)--(1.316,8.441);
\draw[gp path] (11.947,8.441)--(11.767,8.441);
\node[gp node right] at (0.952,8.441) {$14$};
\draw[gp path] (1.136,0.985)--(1.136,1.165);
\draw[gp path] (1.136,8.441)--(1.136,8.261);
\node[gp node center] at (1.136,0.677) {$0$};
\draw[gp path] (2.337,0.985)--(2.337,1.165);
\draw[gp path] (2.337,8.441)--(2.337,8.261);
\node[gp node center] at (2.337,0.677) {$0.002$};
\draw[gp path] (3.538,0.985)--(3.538,1.165);
\draw[gp path] (3.538,8.441)--(3.538,8.261);
\node[gp node center] at (3.538,0.677) {$0.004$};
\draw[gp path] (4.740,0.985)--(4.740,1.165);
\draw[gp path] (4.740,8.441)--(4.740,8.261);
\node[gp node center] at (4.740,0.677) {$0.006$};
\draw[gp path] (5.941,0.985)--(5.941,1.165);
\draw[gp path] (5.941,8.441)--(5.941,8.261);
\node[gp node center] at (5.941,0.677) {$0.008$};
\draw[gp path] (7.142,0.985)--(7.142,1.165);
\draw[gp path] (7.142,8.441)--(7.142,8.261);
\node[gp node center] at (7.142,0.677) {$0.01$};
\draw[gp path] (8.343,0.985)--(8.343,1.165);
\draw[gp path] (8.343,8.441)--(8.343,8.261);
\node[gp node center] at (8.343,0.677) {$0.012$};
\draw[gp path] (9.545,0.985)--(9.545,1.165);
\draw[gp path] (9.545,8.441)--(9.545,8.261);
\node[gp node center] at (9.545,0.677) {$0.014$};
\draw[gp path] (10.746,0.985)--(10.746,1.165);
\draw[gp path] (10.746,8.441)--(10.746,8.261);
\node[gp node center] at (10.746,0.677) {$0.016$};
\draw[gp path] (11.947,0.985)--(11.947,1.165);
\draw[gp path] (11.947,8.441)--(11.947,8.261);
\node[gp node center] at (11.947,0.677) {$0.018$};
\draw[gp path] (1.136,8.441)--(1.136,0.985)--(11.947,0.985)--(11.947,8.441)--cycle;
\node[gp node center,rotate=-270] at (0.276,4.713) {$\|u(t,\cdot)\|_{L^2(\Omega)}$};
\node[gp node center] at (6.541,0.215) {$t$};
\node[gp node right] at (10.479,8.107) {isotropic};
\gpcolor{rgb color={0.580,0.000,0.827}}
\draw[gp path] (10.663,8.107)--(11.579,8.107);
\draw[gp path] (1.234,1.974)--(1.331,2.271)--(1.429,2.613)--(1.526,3.068)--(1.624,3.602)%
  --(1.722,4.008)--(1.819,4.109)--(1.917,3.994)--(2.014,3.833)--(2.112,3.697)--(2.210,3.597)%
  --(2.307,3.524)--(2.405,3.470)--(2.502,3.429)--(2.600,3.398)--(2.698,3.373)--(2.795,3.353)%
  --(2.893,3.337)--(2.990,3.323)--(3.088,3.312)--(3.186,3.302)--(3.283,3.294)--(3.381,3.287)%
  --(3.478,3.281)--(3.576,3.277)--(3.674,3.273)--(3.771,3.269)--(3.869,3.267)--(3.966,3.265)%
  --(4.064,3.264)--(4.162,3.263)--(4.259,3.263)--(4.357,3.263)--(4.454,3.265)--(4.552,3.267)%
  --(4.650,3.270)--(4.747,3.273)--(4.845,3.278)--(4.942,3.284)--(5.040,3.292)--(5.138,3.301)%
  --(5.235,3.313)--(5.333,3.328)--(5.430,3.346)--(5.528,3.367)--(5.626,3.389)--(5.723,3.409)%
  --(5.821,3.419)--(5.918,3.408)--(6.016,3.364)--(6.114,3.283)--(6.211,3.172)--(6.309,3.048)%
  --(6.406,2.930)--(6.504,2.830)--(6.602,2.752)--(6.699,2.695)--(6.797,2.655)--(6.894,2.627)%
  --(6.992,2.608)--(7.090,2.595)--(7.187,2.587)--(7.285,2.580)--(7.382,2.576)--(7.480,2.573)%
  --(7.578,2.571)--(7.675,2.569)--(7.773,2.569)--(7.870,2.568)--(7.968,2.568)--(8.066,2.568)%
  --(8.163,2.568)--(8.261,2.568)--(8.358,2.569)--(8.456,2.570)--(8.554,2.571)--(8.651,2.572)%
  --(8.749,2.573)--(8.846,2.574)--(8.944,2.575)--(9.042,2.577)--(9.139,2.579)--(9.237,2.580)%
  --(9.334,2.582)--(9.432,2.583)--(9.530,2.585)--(9.627,2.586)--(9.725,2.587)--(9.822,2.587)%
  --(9.920,2.585)--(10.018,2.583)--(10.115,2.578)--(10.213,2.572)--(10.310,2.563)--(10.408,2.556)%
  --(10.506,2.557)--(10.603,2.588)--(10.701,2.699)--(10.798,3.005)--(10.896,3.740);
\gpcolor{color=gp lt color border}
\node[gp node right] at (10.479,7.799) {reg. $l_1$};
\gpcolor{rgb color={0.000,0.620,0.451}}
\draw[gp path] (10.663,7.799)--(11.579,7.799);
\draw[gp path] (1.234,1.984)--(1.331,2.645)--(1.429,3.736)--(1.526,5.163)--(1.624,6.548)%
  --(1.722,7.289)--(1.819,7.036)--(1.917,6.087)--(2.014,4.988)--(2.112,4.049)--(2.210,3.335)%
  --(2.307,2.817)--(2.405,2.445)--(2.502,2.178)--(2.600,1.985)--(2.698,1.845)--(2.795,1.742)%
  --(2.893,1.666)--(2.990,1.610)--(3.088,1.568)--(3.186,1.536)--(3.283,1.513)--(3.381,1.495)%
  --(3.478,1.481)--(3.576,1.470)--(3.674,1.463)--(3.771,1.457)--(3.869,1.452)--(3.966,1.449)%
  --(4.064,1.446)--(4.162,1.444)--(4.259,1.442)--(4.357,1.441)--(4.454,1.441)--(4.552,1.440)%
  --(4.650,1.440)--(4.747,1.440)--(4.845,1.439)--(4.942,1.439)--(5.040,1.439)--(5.138,1.439)%
  --(5.235,1.439)--(5.333,1.439)--(5.430,1.439)--(5.528,1.439)--(5.626,1.439)--(5.723,1.439)%
  --(5.821,1.438)--(5.918,1.438)--(6.016,1.438)--(6.114,1.438)--(6.211,1.438)--(6.309,1.438)%
  --(6.406,1.438)--(6.504,1.438)--(6.602,1.438)--(6.699,1.438)--(6.797,1.438)--(6.894,1.438)%
  --(6.992,1.437)--(7.090,1.437)--(7.187,1.437)--(7.285,1.437)--(7.382,1.437)--(7.480,1.437)%
  --(7.578,1.437)--(7.675,1.438)--(7.773,1.438)--(7.870,1.438)--(7.968,1.438)--(8.066,1.439)%
  --(8.163,1.439)--(8.261,1.440)--(8.358,1.440)--(8.456,1.441)--(8.554,1.442)--(8.651,1.443)%
  --(8.749,1.444)--(8.846,1.445)--(8.944,1.447)--(9.042,1.448)--(9.139,1.450)--(9.237,1.453)%
  --(9.334,1.455)--(9.432,1.458)--(9.530,1.462)--(9.627,1.466)--(9.725,1.471)--(9.822,1.477)%
  --(9.920,1.484)--(10.018,1.492)--(10.115,1.502)--(10.213,1.514)--(10.310,1.529)--(10.408,1.549)%
  --(10.506,1.574)--(10.603,1.610)--(10.701,1.661)--(10.798,1.740)--(10.896,1.870);
\gpcolor{color=gp lt color border}
\node[gp node right] at (10.479,7.491) {hexagon};
\gpcolor{rgb color={0.337,0.706,0.914}}
\draw[gp path] (10.663,7.491)--(11.579,7.491);
\draw[gp path] (1.234,2.283)--(1.331,2.404)--(1.429,2.531)--(1.526,2.648)--(1.624,2.737)%
  --(1.722,2.798)--(1.819,2.839)--(1.917,2.866)--(2.014,2.884)--(2.112,2.896)--(2.210,2.903)%
  --(2.307,2.907)--(2.405,2.910)--(2.502,2.911)--(2.600,2.912)--(2.698,2.912)--(2.795,2.912)%
  --(2.893,2.912)--(2.990,2.912)--(3.088,2.912)--(3.186,2.911)--(3.283,2.911)--(3.381,2.911)%
  --(3.478,2.911)--(3.576,2.910)--(3.674,2.910)--(3.771,2.910)--(3.869,2.909)--(3.966,2.909)%
  --(4.064,2.909)--(4.162,2.908)--(4.259,2.908)--(4.357,2.908)--(4.454,2.908)--(4.552,2.908)%
  --(4.650,2.908)--(4.747,2.908)--(4.845,2.908)--(4.942,2.908)--(5.040,2.908)--(5.138,2.908)%
  --(5.235,2.909)--(5.333,2.909)--(5.430,2.910)--(5.528,2.910)--(5.626,2.911)--(5.723,2.912)%
  --(5.821,2.912)--(5.918,2.913)--(6.016,2.914)--(6.114,2.916)--(6.211,2.917)--(6.309,2.918)%
  --(6.406,2.920)--(6.504,2.922)--(6.602,2.924)--(6.699,2.926)--(6.797,2.928)--(6.894,2.931)%
  --(6.992,2.934)--(7.090,2.937)--(7.187,2.940)--(7.285,2.944)--(7.382,2.949)--(7.480,2.954)%
  --(7.578,2.960)--(7.675,2.967)--(7.773,2.975)--(7.870,2.985)--(7.968,2.998)--(8.066,3.013)%
  --(8.163,3.031)--(8.261,3.051)--(8.358,3.072)--(8.456,3.088)--(8.554,3.085)--(8.651,3.052)%
  --(8.749,2.985)--(8.846,2.892)--(8.944,2.786)--(9.042,2.684)--(9.139,2.596)--(9.237,2.524)%
  --(9.334,2.469)--(9.432,2.429)--(9.530,2.399)--(9.627,2.378)--(9.725,2.364)--(9.822,2.355)%
  --(9.920,2.350)--(10.018,2.348)--(10.115,2.348)--(10.213,2.351)--(10.310,2.356)--(10.408,2.362)%
  --(10.506,2.370)--(10.603,2.382)--(10.701,2.401)--(10.798,2.435)--(10.896,2.511);
\gpcolor{rgb color={0.902,0.624,0.000}}
\gpsetpointsize{4.00}
\gppoint{gp mark 1}{(1.136,0.985)}
\gppoint{gp mark 1}{(2.638,0.985)}
\gppoint{gp mark 1}{(4.139,0.985)}
\gppoint{gp mark 1}{(5.641,0.985)}
\gppoint{gp mark 1}{(7.142,0.985)}
\gppoint{gp mark 1}{(8.644,0.985)}
\gppoint{gp mark 1}{(10.145,0.985)}
\gppoint{gp mark 1}{(10.896,0.985)}
\gpcolor{color=gp lt color border}
\draw[gp path] (1.136,8.441)--(1.136,0.985)--(11.947,0.985)--(11.947,8.441)--cycle;
%% coordinates of the plot area
\gpdefrectangularnode{gp plot 1}{\pgfpoint{1.136cm}{0.985cm}}{\pgfpoint{11.947cm}{8.441cm}}
\end{tikzpicture}
%% gnuplot variables
        \caption{Results for `splitting' cf. \cref{fig:iso_splitting,fig:l1_splitting,fig:hexa_splitting}.
        \label{fig:unorm_split}}
  \end{center}
\end{subfigure}
%\hspace{.05\textwidth}
\begin{subfigure}{.48\textwidth}
  \begin{center}
 	\begin{tikzpicture}[gnuplot, scale=0.5, every node/.style={scale=0.5}]
%% generated with GNUPLOT 5.2p2 (Lua 5.3; terminal rev. 99, script rev. 102)
%% Di 11 Mai 2021 18:11:36 CEST
\path (0.000,0.000) rectangle (12.500,8.750);
\gpcolor{color=gp lt color border}
\gpsetlinetype{gp lt border}
\gpsetdashtype{gp dt solid}
\gpsetlinewidth{1.00}
\draw[gp path] (1.320,0.985)--(1.500,0.985);
\draw[gp path] (11.947,0.985)--(11.767,0.985);
\node[gp node right] at (1.136,0.985) {$3$};
\draw[gp path] (1.320,2.476)--(1.500,2.476);
\draw[gp path] (11.947,2.476)--(11.767,2.476);
\node[gp node right] at (1.136,2.476) {$3.5$};
\draw[gp path] (1.320,3.967)--(1.500,3.967);
\draw[gp path] (11.947,3.967)--(11.767,3.967);
\node[gp node right] at (1.136,3.967) {$4$};
\draw[gp path] (1.320,5.459)--(1.500,5.459);
\draw[gp path] (11.947,5.459)--(11.767,5.459);
\node[gp node right] at (1.136,5.459) {$4.5$};
\draw[gp path] (1.320,6.950)--(1.500,6.950);
\draw[gp path] (11.947,6.950)--(11.767,6.950);
\node[gp node right] at (1.136,6.950) {$5$};
\draw[gp path] (1.320,8.441)--(1.500,8.441);
\draw[gp path] (11.947,8.441)--(11.767,8.441);
\node[gp node right] at (1.136,8.441) {$5.5$};
\draw[gp path] (1.320,0.985)--(1.320,1.165);
\draw[gp path] (1.320,8.441)--(1.320,8.261);
\node[gp node center] at (1.320,0.677) {$0$};
\draw[gp path] (2.501,0.985)--(2.501,1.165);
\draw[gp path] (2.501,8.441)--(2.501,8.261);
\node[gp node center] at (2.501,0.677) {$0.002$};
\draw[gp path] (3.682,0.985)--(3.682,1.165);
\draw[gp path] (3.682,8.441)--(3.682,8.261);
\node[gp node center] at (3.682,0.677) {$0.004$};
\draw[gp path] (4.862,0.985)--(4.862,1.165);
\draw[gp path] (4.862,8.441)--(4.862,8.261);
\node[gp node center] at (4.862,0.677) {$0.006$};
\draw[gp path] (6.043,0.985)--(6.043,1.165);
\draw[gp path] (6.043,8.441)--(6.043,8.261);
\node[gp node center] at (6.043,0.677) {$0.008$};
\draw[gp path] (7.224,0.985)--(7.224,1.165);
\draw[gp path] (7.224,8.441)--(7.224,8.261);
\node[gp node center] at (7.224,0.677) {$0.01$};
\draw[gp path] (8.405,0.985)--(8.405,1.165);
\draw[gp path] (8.405,8.441)--(8.405,8.261);
\node[gp node center] at (8.405,0.677) {$0.012$};
\draw[gp path] (9.585,0.985)--(9.585,1.165);
\draw[gp path] (9.585,8.441)--(9.585,8.261);
\node[gp node center] at (9.585,0.677) {$0.014$};
\draw[gp path] (10.766,0.985)--(10.766,1.165);
\draw[gp path] (10.766,8.441)--(10.766,8.261);
\node[gp node center] at (10.766,0.677) {$0.016$};
\draw[gp path] (11.947,0.985)--(11.947,1.165);
\draw[gp path] (11.947,8.441)--(11.947,8.261);
\node[gp node center] at (11.947,0.677) {$0.018$};
\draw[gp path] (1.320,8.441)--(1.320,0.985)--(11.947,0.985)--(11.947,8.441)--cycle;
\node[gp node center,rotate=-270] at (0.276,4.713) {$\|u(t,\cdot)\|_{L^2(\Omega)}$};
\node[gp node center] at (6.633,0.215) {$t$};
\node[gp node right] at (10.479,8.107) {isotropic};
\gpcolor{rgb color={0.580,0.000,0.827}}
\draw[gp path] (10.663,8.107)--(11.579,8.107);
\draw[gp path] (1.416,1.750)--(1.512,2.599)--(1.608,3.188)--(1.704,3.601)--(1.800,3.898)%
  --(1.896,4.115)--(1.992,4.280)--(2.088,4.408)--(2.183,4.512)--(2.279,4.600)--(2.375,4.677)%
  --(2.471,4.747)--(2.567,4.812)--(2.663,4.874)--(2.759,4.934)--(2.855,4.994)--(2.951,5.053)%
  --(3.047,5.113)--(3.143,5.172)--(3.239,5.232)--(3.335,5.293)--(3.431,5.355)--(3.527,5.419)%
  --(3.623,5.483)--(3.718,5.549)--(3.814,5.618)--(3.910,5.688)--(4.006,5.762)--(4.102,5.839)%
  --(4.198,5.922)--(4.294,6.009)--(4.390,6.104)--(4.486,6.206)--(4.582,6.314)--(4.678,6.426)%
  --(4.774,6.532)--(4.870,6.615)--(4.966,6.651)--(5.062,6.604)--(5.158,6.449)--(5.253,6.176)%
  --(5.349,5.807)--(5.445,5.389)--(5.541,4.977)--(5.637,4.614)--(5.733,4.319)--(5.829,4.096)%
  --(5.925,3.936)--(6.021,3.824)--(6.117,3.749)--(6.213,3.699)--(6.309,3.669)--(6.405,3.651)%
  --(6.501,3.643)--(6.597,3.642)--(6.693,3.647)--(6.788,3.656)--(6.884,3.668)--(6.980,3.683)%
  --(7.076,3.700)--(7.172,3.718)--(7.268,3.738)--(7.364,3.759)--(7.460,3.782)--(7.556,3.805)%
  --(7.652,3.829)--(7.748,3.854)--(7.844,3.879)--(7.940,3.905)--(8.036,3.932)--(8.132,3.960)%
  --(8.228,3.988)--(8.323,4.017)--(8.419,4.047)--(8.515,4.078)--(8.611,4.109)--(8.707,4.140)%
  --(8.803,4.173)--(8.899,4.206)--(8.995,4.239)--(9.091,4.273)--(9.187,4.307)--(9.283,4.340)%
  --(9.379,4.374)--(9.475,4.405)--(9.571,4.435)--(9.667,4.461)--(9.763,4.482)--(9.858,4.495)%
  --(9.954,4.497)--(10.050,4.484)--(10.146,4.451)--(10.242,4.396)--(10.338,4.315)--(10.434,4.215)%
  --(10.530,4.121)--(10.626,4.108)--(10.722,4.370)--(10.818,5.353)--(10.914,7.997);
\gpcolor{color=gp lt color border}
\node[gp node right] at (10.479,7.799) {reg. $l_1$};
\gpcolor{rgb color={0.000,0.620,0.451}}
\draw[gp path] (10.663,7.799)--(11.579,7.799);
\draw[gp path] (1.416,3.417)--(1.512,3.332)--(1.608,3.291)--(1.704,3.276)--(1.800,3.277)%
  --(1.896,3.286)--(1.992,3.300)--(2.088,3.315)--(2.183,3.330)--(2.279,3.343)--(2.375,3.354)%
  --(2.471,3.364)--(2.567,3.373)--(2.663,3.381)--(2.759,3.387)--(2.855,3.392)--(2.951,3.394)%
  --(3.047,3.395)--(3.143,3.394)--(3.239,3.392)--(3.335,3.390)--(3.431,3.388)--(3.527,3.385)%
  --(3.623,3.382)--(3.718,3.378)--(3.814,3.373)--(3.910,3.367)--(4.006,3.361)--(4.102,3.354)%
  --(4.198,3.349)--(4.294,3.343)--(4.390,3.337)--(4.486,3.332)--(4.582,3.325)--(4.678,3.319)%
  --(4.774,3.312)--(4.870,3.305)--(4.966,3.298)--(5.062,3.292)--(5.158,3.286)--(5.253,3.280)%
  --(5.349,3.274)--(5.445,3.268)--(5.541,3.262)--(5.637,3.256)--(5.733,3.250)--(5.829,3.244)%
  --(5.925,3.239)--(6.021,3.234)--(6.117,3.230)--(6.213,3.225)--(6.309,3.220)--(6.405,3.216)%
  --(6.501,3.212)--(6.597,3.208)--(6.693,3.204)--(6.788,3.201)--(6.884,3.198)--(6.980,3.196)%
  --(7.076,3.193)--(7.172,3.191)--(7.268,3.189)--(7.364,3.187)--(7.460,3.185)--(7.556,3.184)%
  --(7.652,3.183)--(7.748,3.183)--(7.844,3.182)--(7.940,3.182)--(8.036,3.180)--(8.132,3.178)%
  --(8.228,3.175)--(8.323,3.169)--(8.419,3.162)--(8.515,3.150)--(8.611,3.134)--(8.707,3.110)%
  --(8.803,3.077)--(8.899,3.032)--(8.995,2.971)--(9.091,2.891)--(9.187,2.789)--(9.283,2.664)%
  --(9.379,2.516)--(9.475,2.349)--(9.571,2.168)--(9.667,1.981)--(9.763,1.797)--(9.858,1.628)%
  --(9.954,1.480)--(10.050,1.360)--(10.146,1.272)--(10.242,1.217)--(10.338,1.200)--(10.434,1.233)%
  --(10.530,1.339)--(10.626,1.555)--(10.722,1.945)--(10.818,2.628)--(10.914,3.827);
\gpcolor{color=gp lt color border}
\node[gp node right] at (10.479,7.491) {hexagon};
\gpcolor{rgb color={0.337,0.706,0.914}}
\draw[gp path] (10.663,7.491)--(11.579,7.491);
\draw[gp path] (1.416,3.994)--(1.512,4.182)--(1.608,4.333)--(1.704,4.455)--(1.800,4.556)%
  --(1.896,4.641)--(1.992,4.712)--(2.088,4.766)--(2.183,4.802)--(2.279,4.814)--(2.375,4.796)%
  --(2.471,4.745)--(2.567,4.661)--(2.663,4.546)--(2.759,4.406)--(2.855,4.254)--(2.951,4.108)%
  --(3.047,3.980)--(3.143,3.882)--(3.239,3.816)--(3.335,3.779)--(3.431,3.763)--(3.527,3.760)%
  --(3.623,3.764)--(3.718,3.772)--(3.814,3.783)--(3.910,3.794)--(4.006,3.807)--(4.102,3.820)%
  --(4.198,3.834)--(4.294,3.848)--(4.390,3.861)--(4.486,3.874)--(4.582,3.888)--(4.678,3.900)%
  --(4.774,3.913)--(4.870,3.925)--(4.966,3.936)--(5.062,3.947)--(5.158,3.958)--(5.253,3.968)%
  --(5.349,3.978)--(5.445,3.987)--(5.541,3.996)--(5.637,4.004)--(5.733,4.012)--(5.829,4.019)%
  --(5.925,4.026)--(6.021,4.032)--(6.117,4.037)--(6.213,4.042)--(6.309,4.047)--(6.405,4.051)%
  --(6.501,4.054)--(6.597,4.057)--(6.693,4.059)--(6.788,4.061)--(6.884,4.062)--(6.980,4.063)%
  --(7.076,4.063)--(7.172,4.063)--(7.268,4.062)--(7.364,4.061)--(7.460,4.059)--(7.556,4.056)%
  --(7.652,4.054)--(7.748,4.050)--(7.844,4.047)--(7.940,4.043)--(8.036,4.038)--(8.132,4.033)%
  --(8.228,4.028)--(8.323,4.023)--(8.419,4.017)--(8.515,4.011)--(8.611,4.006)--(8.707,4.000)%
  --(8.803,3.994)--(8.899,3.988)--(8.995,3.983)--(9.091,3.978)--(9.187,3.973)--(9.283,3.968)%
  --(9.379,3.964)--(9.475,3.960)--(9.571,3.956)--(9.667,3.952)--(9.763,3.946)--(9.858,3.939)%
  --(9.954,3.929)--(10.050,3.915)--(10.146,3.896)--(10.242,3.871)--(10.338,3.840)--(10.434,3.804)%
  --(10.530,3.766)--(10.626,3.731)--(10.722,3.715)--(10.818,3.756)--(10.914,3.964);
\gpcolor{rgb color={0.902,0.624,0.000}}
\gpsetpointsize{4.00}
\gppoint{gp mark 1}{(1.320,0.985)}
\gppoint{gp mark 1}{(2.796,0.985)}
\gppoint{gp mark 1}{(4.272,0.985)}
\gppoint{gp mark 1}{(5.748,0.985)}
\gppoint{gp mark 1}{(7.224,0.985)}
\gppoint{gp mark 1}{(8.700,0.985)}
\gppoint{gp mark 1}{(10.176,0.985)}
\gppoint{gp mark 1}{(10.914,0.985)}
\gpcolor{color=gp lt color border}
\draw[gp path] (1.320,8.441)--(1.320,0.985)--(11.947,0.985)--(11.947,8.441)--cycle;
%% coordinates of the plot area
\gpdefrectangularnode{gp plot 1}{\pgfpoint{1.320cm}{0.985cm}}{\pgfpoint{11.947cm}{8.441cm}}
\end{tikzpicture}
%% gnuplot variables
        \caption{Results for `merging' cf. \cref{fig:iso_merge,fig:l1_merge,fig:hexa_merge}.
          \label{fig:unorm_merge}}
  \end{center}
  \end{subfigure}
  \caption{Time evolution $\|u(t,\cdot)\|_{L^2(\Omega)}$
    when topology changes are present.
	\label{fig:u_norm_others}}
      \end{center}
\end{figure}

\section*{Acknowledgements}

The authors gratefully acknowledge the support by the RTG 2339 “Interfaces, Complex Structures, and Singular Limits” of the German Science Foundation (DFG).

% deprecated file
%\nocite{*} %Um alle Einträge des Literaturverzeichnisses auszudrucken, auch wenn sie nicht zitiert werden

%\clearpage
\printbibliography
%\clearpage
%\listoffigures
\clearpage
\end{document}